\documentclass[publmath,draft,numbook]{svjour}

\usepackage{graphicx}
\usepackage{epsf}
\usepackage[active]{srcltx}

\reversetheomheadings     

\usepackage{publmath}                             

\begin{document}
\newcommand{\Z}{{\mathbb Z}}
\newcommand{\R}{{\mathbb R}}
\newcommand{\Q}{{\mathbb Q}}
\newcommand{\C}{{\mathbb C}}
\newcommand{\lms}{\longmapsto}
\newcommand{\lra}{\longrightarrow}
\newcommand{\hra}{\hookrightarrow}
\newcommand{\ra}{\rightarrow}
\newcommand{\sgn}{\rm sgn}
\spnewtheorem{corollary-definition}{Corollary-Definition}[definition]{\bf}{\it}
\title{ MODULI SPACES OF LOCAL SYSTEMS AND HIGHER TEICHM\"ULLER THEORY.}
\author{
\firstname{Vladimir} FOCK\and \firstname{Alexander} GONCHAROV%
}                     
\institute{
V.V.\\
ITEP, \\
B. Cheremushkinskaya 25,\\
 117259 Moscow, Russia,\\
 fock@math.brown.edu 
\and 
A.G.\\Department of Mathematics,\\
Brown University, 
Providence RI 02906, USA,\\ 
sasha@math.brown.edu}
\date{date}
%
\maketitle
\tableofcontents
\begin{abstract}
Let $G$ be a split semisimple algebraic group over $\Q$ with trivial center. 
Let $S$ be a compact oriented surface, with or without boundary. We define {\it positive} 
representations of the fundamental group of $S$ to $G(\R)$, construct explicitly all positive representations, 
and prove that they are faithful, discrete, and positive hyperbolic; the moduli space of positive representations 
is a topologically trivial open domain in the space of all representations.
When $S$ have holes, we defined two moduli spaces closely related to the moduli spaces of $G$-local systems on $S$. 
We show that they carry a lot of interesting structures. In particular we define a distinguished collection of 
coordinate systems, equivariant under the action of the mapping class group of $S$. We prove that their transition 
functions are subtraction free. Thus we have positive structures on these moduli spaces. Therefore we can take 
their points with values in any positive semifield. Their positive real points provide the two higher 
Teichm\"uller spaces related to $G$ and $S$, while the points with values in the tropical semifields provide 
the lamination spaces. We define the motivic avatar of the Weil-Petersson form for one of these spaces. It 
is related to the motivic dilogarithm.
\end{abstract}
\section{Introduction}
\label{intro}

{\bf Summary}. Let $S$ be a compact oriented surface with $\chi(S)<0$, with or without boundary. 
We define {\it positive} representations of $\pi_1(S)$ to a split 
semi-simple real Lie group $G(\R)$ 
with trivial center, construct explicitly 
 all positive representations, and prove that they have the following properties: Every positive 
representation is faithful, its image in $G(\R)$ is discrete, 
and the image  of any non-trivial non-boundary conjugacy class is 
positive hyperbolic;  the moduli space ${\cal L}^+_{G, S}$ of positive representations 
is a topologically trivial domain of dimension 
$-\chi(S) {\rm dim}G$. 
When $G=PGL_2$ we recover the classical 
Teichm\"uller space for $S$. Using total positivity in semi-simple Lie groups \cite{L1}, we 
introduce and study {\it positive 
configurations of flags}, which play a key role in our story. 
The limit sets of positive representations are positive curves 
in the real flag variety related to $G$. 

When $S$ has holes, in addition to the Teihm\"uller 
 space ${\cal L}^+_{G, S}$ we define  a closely related pair 
of dual Teichm\"uller spaces ${\cal X}^+_{G,S}$ and ${\cal A}^+_{G,S}$. 
In order to do this,  we introduce a {\it dual pair of moduli spaces} ${\cal X}_{G,S}$ and ${\cal A}_{G,S}$ 
closely related to the moduli space ${\cal L}_{G,S}$ of $G(\C)$-local systems on $S$. 
We show that, unlike the latter, each of them is rational. Moreover, each of them carries a  
positive atlas,  equivariant with respect to the action of the mapping class group of $S$. 
Thus we can define the 
 sets of their  $\R_{>0}$-points, which are nothing else but the above dual pair of Teichm\"uller spaces. 
We identify each of them with $\R^{-\chi(S) {\rm dim}G}$. 
The sets of their tropical points provide us with a 
dual pair of lamination spaces -- for $G = PGL_2$ we recover Thurston's laminations. 
 The moduli space ${\cal X}_{G,S}$  
is Poisson. The space ${\cal A}_{G,S}$   carries a degenerate $2$-form $\Omega$. It 
comes from  a class in the $K_2$-group of the function field of  ${\cal A}_{G,S}$. 

The Teichm\"uller spaces are related as follows. One has  a canonical projection 
$\pi^+: {\cal X}^+_{G,S} \to {\cal L}^+_{G,S}$ as well as its canonical splitting 
$s:  {\cal L}^+_{G,S}\hra {\cal X}^+_{G,S}$. 
The projection $\pi^+$ is a ramified covering with the structure group $W^n$, where $W$ is the Weyl group of $G$. 
There is a map ${\cal A}^+_{G,S} \to {\cal X}^+_{G,S}$, provided by taking 
the quotient along the null foliation of $\Omega$. We show that higher Teichm\"uller spaces behave nicely 
with respect to  cutting and gluing of $S$. We define a (partial) 
completion of the higher Teichm\"uller spaces, which for $G=PGL_2$ coincides 
with the Weil-Petersson completion of the classical Teichm\"uller space. 

We conjecture that there 
is a duality between the ${\cal X}$- and ${\cal A}$-moduli spaces, which interchanges 
$G$ with its Langlands dual. In particular, it predicts that there exists a canonical basis 
in the space of functions on one of the moduli space, paramatrised by the integral 
tropical points of the other, with a number of remarkable properties. 
We constract such a basis for $G=PGL_2$. The pair $({\cal X}_{G,S}, {\cal A}_{G,S})$, at least for $G=PGL_m$,   
forms an (orbi)-cluster ensemble in the sense of 
\cite{FG2}, and thus can be quantised, as explained in loc. cit.. 
\vskip 3mm

{\bf 1 An algebraic-geometric approach to higher  Teichm{\"u}ller theory}. 
Let $S$ be a compact  oriented surface with $n\geq 0$ holes. 
The  Teichm{\"u}ller space ${\cal T}_S$ is the moduli space of complex
structures on  $S$ modulo diffeomorphisms isotopic to the
identity. We will assume that $S$ is hyperbolic. 
Then the Poincar{\'e} uniformisation theorem 
identifies ${\cal T}_S$ with the space of all  faithful 
representations $\pi_1(S) \to PSL_2(\R)$ with discrete image,  
modulo conjugations by $PSL_2(\R)$. 
The canonical metric on the hyperbolic plane descends to a complete 
hyperbolic metrics of curvature $-1$ on $S$. The space ${\cal T}_S$ is the moduli space 
of such metrics on $S$ modulo isomorphism. 
. 

If $n>0$, the space ${\cal T}_S$ has a boundary; in fact it is a manifold with corners.  
Let us explain why. 
Given a point $p\in {\cal T}_S$,  
a boundary component of $S$ 
is {\it cuspidal} 
if its  neighborhood is isometric to a cusp. The points of the boundary of ${\cal T}_S$ 
parametrise the metrics on $S$ with at least one cuspidal boundary component. 
There is a natural $2^n:1$ 
cover $\pi^+: {\cal T}^+_S \to {\cal T}_S$, ramified at the boundary of ${\cal T}_S$. 
It is described as follows. For each non-cuspidal boundary component of $S$ there is a unique 
geodesic homotopic to it. 
 The space ${\cal T}^+_S$ parametrises pairs ($p\in {\cal T}_S$, 
plus a choice of orientation of
every non-cuspidal boundary component). 
The orientation of $S$ provides orientations of all boundary components and hence orientations 
of the corresponding geodesics. Thus we get a canonical embedding $s: {\cal T}_S \hookrightarrow {\cal T}^+_S$.  Let ${\cal T}^u_S\subset {\cal T}_S$ be the subspace parametrising the metrics 
all of whose boundary components are cuspidal. It is the deepest corner 
of ${\cal T}_S$. Its quotient  
$
{\cal M}_S:= {\cal T}^u_S/\Gamma_S
$ 
by the mapping class group $\Gamma_S:= 
{\rm Diff}(S)/{\rm Diff}_0(S)$ is the classical moduli space ${\cal M}_{g,n}$. 
In particular it has a complex structure. 
\vskip 3mm
The Teichm{\"u}ller space is by no means 
the set of real points of an algebraic variety. Indeed, discreteness 
is not a condition of an algebraic nature. Traditionally the Teichm{\"u}ller theory is 
considered as a part of
analysis and geometry. 
In this paper we  show that  there is
an algebraic-geometric approach to 
the Teichm{\"u}ller theory. Moreover, we show that it admits 
a generalization where $PSL_2$ is
replaced by a split semi-simple algebraic group $G$ over $\Q$ with trivial center.   
Our approach owes a lot to Thurston's ideas in geometry of surfaces. 
So perhaps the emerging theory could be called higher Teichm{\"u}ller-Thurston theory. 
Here are the main features of our story.

\vskip 3mm

Let $G$ be a split reductive algebraic group over $\Q$. 
Let $\widehat S$ be a pair consisting of 
a compact oriented surface $S$ with $n$ holes and a finite 
(possibly empty) collection of marked points on the boundary considered modulo 
isotopy.  
We define two new moduli spaces, denoted ${\cal X}_{G,\widehat S}$ and 
${\cal A}_{G,\widehat S}$. 
{\it We always assume that $G$ has trivial center (resp. is simply-connected) 
for the  ${\cal X}$- (resp. 
${\cal A}$-) moduli space.} 
In particular, when  $\widehat S$ is a disk 
with marked points on the boundary we get  
moduli spaces of cyclically ordered configurations  of points 
in the flag variety ${\cal B}: = G/B$ and twisted cyclic configurations of points of 
the principal affine variety ${\cal A}:= G/U$.

We show that if $S$ does have  holes then these moduli spaces carry 
 interesting additional structures. In particular, for each of them we define  
a distinguished collection of coordinate systems, equivariant under the action 
of the mapping class group of $S$. We  
prove that their transition functions are subtraction-free, 
providing  a  
{\it positive 
atlas} on the corresponding moduli space. 

We say that $X$ is a {\it positive variety} if it admits a positive atlas,
i.e. a distinguished
collection of coordinate systems as above. 
If $X$ is a positive variety, we can take points  of $X$ 
with values in any {\it semifield},
i.e. a set equipped with operations of addition, multiplication and division, 
e.g. $\R_{>0}$. We define {\it higher Teichm{\"u}ller spaces} ${\cal X}^+_{G,
 \widehat S}$ and 
${\cal A}^+_{G,\widehat S}$ as the sets of  $\R_{>0}$-points of the 
corresponding positive spaces 
${\cal X}_{G, \widehat S}$ and ${\cal A}_{G,\widehat S}$. 
Precisely, they  
consist of the real points 
of the corresponding moduli spaces whose coordinates in
 one, and hence in any, of the constructed coordinate systems are positive. 
Our approach for general $G$ uses George Lusztig's theory  
of positivity in semi-simple Lie groups  \cite{L1}-\cite{L2}.  
For $G = SL_m$ we have an elementary and self-contained approach 
developed in Sections 9-10. 
\vskip 3mm
We prove that the representations $\pi_1(S) \to G(\R)$ underlying the points of 
${\cal X}^+_{G, \widehat S}$ are faithful, discrete and positive hyperbolic. 
Using positive configurations of flags, we define positive representations of $\pi_1(S)$ 
for closed surfaces $S$, and show that they have the same properties. 
For $G = PSL_2(\R)$ we recover the Fuchsian representations. If $S$ is
closed, a component in the moduli space of completely reducible $G(\R)$-local system was defined
and studied by Nigel Hitchin \cite{H1} by a 
completely different method. A partial  generalization  to punctured 
 $S$ see in \cite{BAG}. We explain the relationship between our work and the work of Hitchin in Section 1.11. 

\vskip 3mm
The classical Teichm{\"u}ller spaces appear as follows.  
We show that in the absence of  marked points, i.e. when 
$S = \widehat S$, the space 
${\cal X}^+_{PGL_2, S}$ 
is identified with  
the Teichm{\"u}ller space  ${\cal T}^+_S$, 
and  ${\cal A}^+_{SL_2, S}$ is identified with the 
 decorated Teichm{\"u}ller space  ${\cal T}^d_S$
defined by Robert Penner \cite{P1}. 
The restrictions of our canonical coordinates to the 
corresponding Teichm{\"u}ller spaces are well known coordinates there.  
 Our two moduli spaces, their positivity, and especially the 
discussed below motivic data 
seem to be new even in the classical setting.  

\vskip 3mm
It makes sense to consider points of our moduli spaces with values in an arbitrary
semifield. An interesting situation appears when we consider 
exotic 
{\it tropical} semifields ${\Bbb A}^t$. Let us recall their definition. 
 Let ${\Bbb A}$ be one of the three sets $\Z$, $\Q$, $\R$. 
Replacing the operations 
of multiplication, division  
and addition in ${\Bbb A}$ by the operations of addition, subtraction  and 
taking the maximum we obtain the  tropical semifield  ${\Bbb A}^t$.
We show that the set of points of the moduli space ${\cal A}_{SL_2, S}$
 with values in the 
{tropical} semifield $\R^t$ is identified  with Thurston's transversal measured laminations on $S$ \cite{Th}. 
The $\Q^t$- and $\Z^t$-points of ${\cal A}_{SL_2, S}$ 
give us rational and integral laminations. 
Now taking the points of the other moduli spaces with values in tropical
semifields, we get generalizations of
the lamination spaces. 
We show that for any positive variety $X$ the  projectivisation of the space 
$X(\R^t)$   
serves as the Thurston-type boundary for the   space
$X^+:= X(\R_{>0})$.

\vskip 3mm
We suggest that there exists a remarkable duality between the 
${\cal A}$-moduli space for the group $G$ and the ${\cal X}$-moduli space for the 
Langlands dual group $^L G$. Here is one  of the (conjectural) manifestations  
of this duality.  
Let $X$ be a positive variety. A rational function $F$ on $X$ is called {\it a
  good positive Laurent polynomial} if it is a Laurent polynomial with 
positive integral coefficients in every coordinate system of the positive
atlas on $X$. We conjecture that 
the set of $\Z^t$-points of 
the ${\cal A}$/${\cal X}$-moduli space for $G$ should parametrize 
a canonical basis 
in the space of all good positive Laurent polynomials on 
the ${\cal X}$/${\cal A}$-moduli space for the 
Langlands dual group $^L G$, for the positive atlases we defined. 
For $G= SL_2$ or $PGL_2$ we elaborate this theme in Section 12.

\vskip 3mm
If $S$ has no holes, the moduli space ${\cal L}_{G, S}$ has a natural 
symplectic structure. If $n>0$, the symplectic structure is no longer available, and the 
story splits into two different directions. First, the moduli space ${\cal L}_{G, S}$ 
has a natural Poisson structure, studied  in [FR]. There is a natural Poisson structure on 
the moduli space ${\cal X}_{G, \widehat S}$ such that the canonical projection   
$\pi: {\cal X}_{G, \widehat S} \to {\cal L}_{G, \widehat S}$ is Poisson.  
Second, we introduce a degenerate symplectic structure 
on 
the moduli space ${\cal A}_{G, \widehat S}$, 
an analog of the Weil-Petersson form. To define the Weil-Petersson form  
we construct its $K_2$-avatar, that is a class in $K_2$ whose image under the regulator map 
gives the $2$-form. We show that 
it admits  an even finer motivic avatar 
related to the motivic dilogarithm and the  
second motivic Chern class $c_2^{\cal M}$ of the universal $G$-bundle on $BG$. 

Amazingly an explicit construction of a cocycle for the class  $c_2^{\cal M}$ for $G=SL_m$ 
given in \cite{G3} delivers the canonical coordinates on the
${\cal A}_{SL_m, \widehat S}$-space. Namely, in Section 15 we give a simple  general construction 
which translates a class $c_2^{\cal M}$ to a class in $K_2$ of 
${\cal A}_{G,\widehat S}$. Applying it to the explicit cocycle 
defined in \cite{G3}, we get a class in $K_2$ of  
${\cal A}_{SL_m, \widehat S}$, which shows up
  written in the canonical coordinates! This construction was our starting point.  

The Poisson structure and the 
Weil-Petersson form are
 compatible in a natural way: the quotient of the  
${\cal A}$-space  along the null foliation of the 
Weil-Petersson form  is embedded into the ${\cal X}$-space 
as a symplectic leaf of the Poisson structure. 
The representations of $\pi_1(S)$ corresponding to its points are the ones 
with unipotent monodromies around each hole. 
\vskip 3mm
Recall an embedding $SL_2 \hra G$, as a principal $SL_2$--subgroup,   
well-defined up to a conjugation. 
It leads to  embeddings ($G'$ denotes the adjoint group)
$
{\cal A}_{SL_2, S} \hra 
{\cal A}_{G, S}, \quad {\cal X}_{PGL_2, S} \hra {\cal X}_{G', S} 
$ 
and their counterparts for the Teichm{\"u}ller and  lamination  spaces. 
\vskip 3mm
Summarizing, we suggest that for a surface $S$ with boundary 
\vskip 3mm\noindent
{\it The pair 
of moduli spaces 
 $({\cal A}_{G, \widehat S}, {\cal X}_{G', \widehat S})$
is 
the algebraic-geometric avatar of  
higher Teichm{\"u}ller theory}. 
\vskip 3mm

After the key properties of the $({\cal X}, {\cal A})$-moduli spaces were 
discovered and worked out, we started to look for an adequate 
framework describing their
 structures. 
It seems that it 
is provided by the object which we call an {\it (orbi)-cluster
  ensemble} and develop in \cite{FG2}. 
It generalizes the {\it cluster algebras}, discovered and studied in a remarkable
series of papers by S. Fomin and A. Zelevinsky \cite{FZI},
\cite{FZII}, \cite{FZ3}. A cluster ensemble is a pair of 
positive spaces ${\cal X}_{\cal E}$ and 
${\cal A}_{\cal E}$ defined for any integral-valued 
skew-symmetrizable function ${\cal E}$ on  the square of a finite
set. The cluster algebra closely related to the positive space 
 ${\cal A}_{\cal E}$. The orbi-cluster ensemble is a bit more general structure.

For any $G$ there exists an 
 orbi-cluster ensemble $({\cal A}_{{\cal E}(G, \widehat S)}, {\cal X}_{{\cal E}(G, \widehat S)})$
related to the pair of moduli spaces 
$({\cal A}_{G, \widehat S}, {\cal X}_{G, \widehat S})$. For $G = SL_2$ the
corresponding pairs of positive spaces coincide. 
For other $G$'s the 
precise relationship between them is more complicated. 
In Section 10 we explain this for $G = SL_m$.  
A more technical general case is treated in a sequel to this paper. The fact that 
Penner's decorated Teichm\"uller space is related   to a cluster algebra 
was independently
observed in \cite{GSV2}. 

The relationship between the  higher Teichm\"uller theory and
 cluster ensembles is very fruitful in both directions. We show in \cite{FG2} that cluster ensembles 
have essentially all properties of the two moduli spaces described above,
including  the relation with the motivic dilogarithm. 
Our point is that  
\vskip 3mm\noindent
{\it A big part of  
the higher Teichm\"uller theory can be generalized in the framework of 
 cluster ensembles, providing a guide for investigation of the latter.}
\vskip 2mm
 
For example, we suggest in Section 4
loc. cit. that there 
exists a duality 
between the ${\cal X}_{\cal E}$ and ${\cal A}_{\cal E}$ spaces, which in the special 
case of our pair 
of moduli spaces includes the one discussed above. 
This duality should also include the canonical bases of Lusztig and 
the Laurent phenomenon of Fomin and
Zelevinsky   \cite{FZ3} as special cases. 

We show in \cite{FG2} that  
any cluster ensemble can be quantized. 
In a sequel to this paper  we  apply this to quantize the moduli space 
${\cal X}_{G, \widehat S}$. The quantization is given 
by means of the quantum dilogarithm. Remarkably the
quantization is governed by the motivic avatar of the Weil-Petersson form on
the ${\cal A}_{G, \widehat S}$ moduli space.

The rest of the Introduction contains a more detailed discussion of main
definitions and results.

\vskip 3mm 
{\bf 2. The two moduli spaces}. Let 
$ 
S = \overline S - D_1 \cup ... \cup D_n
$, 
where $\overline S$ is a 
compact oriented two-dimensional manifold without boundary of genus $g$, 
and $D_1, ..., D_n$ are non intersecting open
disks. 
The boundary $\partial S$ of $S$ is a disjoint union of circles. 

\begin{definition} \label{4.7.03.1}
A marked surface $\widehat S$ is a pair $(S, \{x_1, ..., x_k\})$, where $S$ 
is a compact oriented surface 
 and  $\{x_1, ..., x_k\}$ is a finite (perhaps empty)  set of distinct 
boundary points, considered up to an  isotopy. 
\end{definition}
We have  $\widehat S = S$  
if there are no boundary points. 
Let us define the {\it punctured boundary} 
$\partial \widehat S$ of 
$\widehat S$ by 
$
\partial \widehat S := \partial S  - \{x_1, ..., x_k\}
$. 
It is a union of circles and (open) intervals. 
Let $N$ be the number of all 
connected components of $\partial \widehat S$.  We say that a 
marked surface $\widehat S$ is  
{\it hyperbolic} if $g>1$, 
or $g=1, N>0$, or $g=0, N \geq 3$. 
We will usually assume that $\widehat S$ is hyperbolic, and $N>0$.

\vskip 3mm
Let $H$ be a group. 
Recall that an $H$--local system is a principle $H$-bundle 
 with a flat connection. There is a well known bijection 
between the 
homomorphisms of the fundamental group 
$\pi_1(X, x)$ to $H$ modulo $H$-conju\-gation,  
  and the isomorphism classes of $H$--local systems on 
 a nice topological space $X$. It associates to a local system its 
monodromy representation. 
\vskip 3mm
Let $G$ be a split reductive algebraic group over $\Q$. 
Then ${\cal L}_{G, S}$ is the moduli space of $G$-local 
systems on $S$. It is an  algebraic stack over $\Q$, and its  
generic part is an algebraic variety over $\Q$. 
We introduce two other moduli spaces, denoted 
 $ {\cal A}_{G, \widehat S}$ and ${\cal X}_{G,\widehat S}$, closely 
related to ${\cal L}_{G, S}$.  
They are also algebraic  stacks over $\Q$.

The flag variety ${\cal B}$ parametrizes Borel subgroups
in $G$. Let $B$ be a Borel subgroup. Then ${\cal B} = G/B$. Further,  
$U := [B,B]$ is
 a  maximal unipotent subgroup in $G$. 

Let ${\cal L}$ be a $G$-local system on $S$. We assume that 
$G$ acts on ${\cal L}$ from the right. We define 
the associated {\it flag bundle} ${\cal L}_{\cal B}$ and {\it principal affine bundle} 
${\cal L}_{\cal A}$ by setting 
\begin{equation} \label{4.7.03.2fx}
{\cal L}_{\cal B}:= {\cal L}\times_G {\cal B} = 
{\cal L}/B; \qquad {\cal L}_{\cal A}:= {\cal L}/U.
 \end{equation}
\begin{definition} \label{4.7.03.2}
Let $G$ be a split reductive group. 
A {\rm framed}  $G$-local system on $\widehat S$ is a pair $({\cal L},
\beta)$, where  ${\cal L}$ is a  $G$-local system 
on $S$, and  $\beta$ a flat section of the restriction of 
${\cal L}_{\cal B}$ to the punctured boundary $\partial \widehat S$. 
The space 
${\cal X}_{G, \widehat S}$ is the moduli space of 
framed $G$--local systems on $\widehat S$. 
\end{definition}

 Assume now that $G$ is simply-connected. The maximal length
 element $w_0$ of the Weyl group of $G$ has a natural lift  
to $G$, denoted $\overline w_0$. 
Set $s_G:= {\overline w}^2_0$. Then one shows that 
 $s_G$ is in the center of $G$ and $s^2_G =e$.
In particular $s_G=e$ if the order of the center of $G$ is odd,
 i.e. $G$ is of type $A_{2k}, E_6, E_8, F_4, G_2$. If 
$s_G \not = e$, e.g. $G = SL_{2k}$, the definition of 
the moduli space ${\cal A}_{G, \widehat S}$ is more subtle.   
So in the Introduction we do not consider this case. 

 \begin{definition} \label{4.7.03.2m} 
Let $G$ be a simply-connected 
split semi-simple algebraic group, and $s_G =e$. 
A {\rm decorated}  $G$-local system on $\widehat S$ is a pair 
$({\cal L}, \alpha)$, where ${\cal L}$ is a $G$-local system 
on $S$, and $\alpha$ a flat section of the restriction of 
${\cal L}_{\cal A}$ to  $\partial \widehat S$. The space 
${\cal A}_{G, \widehat S}$ is the moduli space of 
decorated $G$-local systems on $\widehat S$. 
\end{definition}
 
Let $X$ be a $G$--set. The elements of $G \backslash X^n$ are called 
{\it configurations} of $n$ elements in $X$. The coinvariants of the 
cyclic translation are called  {\it cyclic configurations}
of $n$ elements in $X$. 
\vskip 3mm 
{\bf Example}. Let  $\widehat S$ be a 
disk with $k\geq 3$ boundary points, so ${\cal L}$ is a 
trivial $G$-local system.  
Then 
${\cal X}_{G, \widehat S}$ is identified with the cyclic 
configuration space of $k$ flags in $G$, that is coinvariants of the
cyclic shift on 
$G\backslash {\cal B}^k$. 
Indeed, choosing a trivialization of ${\cal L}$, we identify  its fibers with  $G$. Then  
the section $\beta$ over the punctured boundary of the disk determines a 
cyclically ordered collection of 
flags.  A change of the trivialization amounts 
to a  left diagonal action of  $G$. 
So a cyclic 
configuration of $k$ flags is well defined. It 
determines the isomorphism class of the framed 
$G$-local system on the marked disk. 
Similarly, assuming $s_G =e$, the space 
 ${\cal A}_{G, \widehat S}$ is identified with the 
cyclic configuration space of $k$ points in the principal affine space
 ${\cal A}:= G/U$. 
\vskip 3mm 

Forgetting $\beta$ we get a canonical projection 
$$
\pi: {\cal X}_{G, \widehat S} \lra {\cal L}_{G, S}.
$$
If $\widehat S = S$, it is  a finite map at the generic point of degree  
$|W|^n$, where $W$ is the Weyl group of $G$. 

We show that the moduli spaces ${\cal X}_{G, \widehat S}$ and 
${\cal A}_{G, \widehat S}$ have some 
 interesting additional structures. 
\vskip 3mm
{\bf 3. The Farey set, configurations of flags, and framed local systems on $\widehat S$.} 
A {\it cyclic structure} on a set $C$ is defined so that every point of the set 
provides an order of $C$, and the orders corresponding to different points are related in the obvious way. 
A {\it cyclic set} is a set with a cyclic structure. An example is given by any subset of the 
oriented circle. 

\begin{figure}[ht]
\centerline{\epsfbox{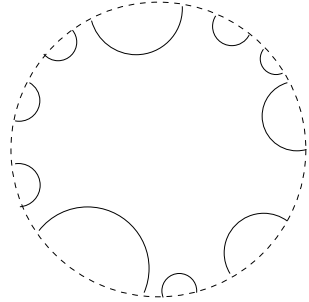}}
\caption{The Farey set: a cyclic $\pi_1(S)$-set ${\cal F}_{\infty}(S)$.}
\label{fg84}
\end{figure}

{\it The Farey set of $S$}. It is a cyclic $\pi_1(S)$-set assigned to a topologicval surface $S$. 
Let us give three slightly different versions of its  definition. 

(i) Shrinking holes on $S$ to punctures, we get a homotopy equivalent surface $S'$.  
The universal cover of $S'$ is an open disc. The punctures on $S'$ 
give rise to a countable subset on the boundary of this disc. This subset 
inherits a cyclic structure from the boundary of the disc. 
The group $\pi_1(S)$ acts on it by the deck transformations. The obtained 
cyclic $\pi_1(S)$-set is called the {\it Farey set} and denoted by ${\cal F}_{\infty}(S)$. 
Choosing a complete hyperbolic structure on $S$, we identify the  Farey set with a subset of 
the boundary of the hyperbolic plane (the latter is called the absolute).

\begin{figure}[ht]
\centerline{\epsfbox{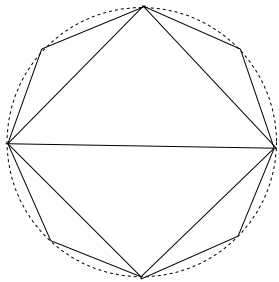}}
\caption{The Farey triangulation of the hyperbolic disc.}
\label{Far}
\end{figure}

(ii) An {\it ideal triangulation}
 $T$ of $S'$ 
is a triangulation of $S'$ with vertices 
at the punctures. Let us lift an ideal triangulation $T$ on $S'$ to the universal cover ${\cal H}$ of $S'$. 
The obtained triangulation $\widetilde T$ of the hyperbolic plane ${\cal H}$ 
is identified with the Farey traingulation, see Figure \ref{Far}. 
The set of vertices of $\widetilde T$ 
inherits a cyclic structure from the absolute of ${\cal H}$. 
The obtained cyclic $\pi_1(S)$-set does not depend on the choice of a triangulation $T$, and is 
canonically isomorphic to 
the {Farey set} ${\cal F}_{\infty}(S)$. 

(iii) Here is a definition of the Farey set via the surface $S$. 
Choose a hyperbolic structure with geodesic boundary on $S$. Its universal 
 cover is identified with the hyperbolic plane with an infinite number of geodesic half 
discs removed, see Figure \ref{fg84}. The set  ${\cal F}_{\infty}(S)$ 
is identified with the set of these removed geodesic half 
discs. The latter  has  an obvious  cyclic $\pi_1(S)$-set structure, see Figure \ref{fg84},  
which does not depend on the choice of a hyperbolic structure on $S$. 

\vskip 3mm 
{\it The Farey set of a marked surface $\widehat S$ with boundary.}  
Choose a hyperbolic structure with geodesic boundary on $S$. Then the set  ${\cal F}_{\infty}(\widehat S)$ 
is the preimage of the punctured boundary $\partial \widehat S$ on the universal cover. 
It has an obvious  cyclic $\pi_1(S)$-set structure,  
which does not depend on the choice of a hyperbolic structure on $S$.

A framed local system on $\widehat S$ gives rise to a map 
\begin{equation} \label{4.1.04.10facoma}
\beta_{\rho}: {\cal F}_{\infty}(\widehat S) \lra {\cal B}(\C) \qquad \mbox{modulo the
  action of $G(\C)$}, 
\end{equation}
which is equivariant with respect to the action of the group 
$\pi_1(S)$: the latter acts on the flag variety via 
the monodromy representation $\rho$ of the local system. 
The following Lemma is straitforward. 

\begin{lemma} \label{fras} There is a natural bijection between the 
framed $G(\C)$-local systems on $\widehat S$ and $\rho$-equivariant maps 
(\ref{4.1.04.10facoma}), where where $\rho: \pi_1(S) \to G(\C)$ is a representation, modulo 
$G(\C)$-conjugation. 
\end{lemma}
A similar interpretation of the 
moduli space ${\cal A}_{G, \widehat S}$ is given in  Section 8.6. 

 Here is a reformulation of this Lemma, which will serve us below. 

\vskip 3mm
{\it Configurations}. Let $C$ be a set, $G$ a group, and $X$ a set. Then we define 
the set 
$$
{\rm Conf}_C(X):= {\rm Map}(C, X)/G
$$
of configurations of points in $X$ parametrised by the set $C$, 
where the action of $G$ on ${\rm Map}(C, X)$ is induced by the one 
on $X$. 
In particular if $C = \{1, ..., n\}$, we denote by ${\rm Conf}_{n}(X)$ the corresponding 
configuration space of $n$ points in $X$.

Let $\pi$ be a group acting on $C$. Then it acts 
 on ${\rm Conf}_C(X)$. So there is  a set of $\pi$-equivariant maps 
$$
{\rm Conf}_{C, \pi}(X) = {\rm Conf}_C(X)^{\pi}.
$$
Let $(\psi, \rho)$ be a pair  consisting of a group homomorphism
 $\rho: \pi \to G$, and a $(\pi, \rho)$-equivariant 
map of sets $\psi: C \to X$, i.e. for any $\gamma \in \pi$ one 
has $\psi(\gamma c) = \rho(\gamma)\psi(c)$. 
The group $G$ acts by conjugation on these pairs. 
\begin{lemma} \label{q1-} Assume that $G$ acts freely on ${\rm Map}(C, X)$. 
Then \newline ${\rm Conf}_{C, \pi}(X)$ is identified with the set of pairs 
$(\psi, \rho)$, where $\psi: C \to X$ is a $(\pi, \rho)$-equivariant 
map of sets,  modulo $G$-conjugations.
\end{lemma} 

{\bf Proof}. Choose a map $\mu: C \to X$ representing a configuration $m$. 
Then $\mu \circ \gamma^{-1}: C \to X$ represents the configuration $\gamma(m)$. 
Since $\gamma(m)=m$, there is an element $g_{\gamma}\in G$ such that $\mu \circ \gamma^{-1} = g \mu$. 
Since the group $G$ acts on ${\rm Map}(C, X)$ freely, $g_\gamma$ is uniquely defined. 
So we get a group homomorphism $\rho: \pi \to G$ given by $\gamma \lms g_\gamma$. 
Changing the representative $\mu$ amounts to conjugation. 
The lemma is proved. 
\vskip 3mm
Using this, we can rephrase Lemma \ref{fras} as the following canonical isomorphism:
\begin{equation} \label{fras10}
{\cal X}_{G, \widehat S}(\C) = {\rm Conf}_{{\cal F}_{\infty}(\widehat S), \pi_1(S)}( {\cal B}(\C)).
\end{equation} 
This way we define points of the space ${\cal X}_{G, \widehat S}(\C)$ without
 even mentioning representations of $\pi_1(S)$: they appear 
thanks to Lemma \ref{q1-}.

\vskip 3mm
{\bf 4. Decomposition Theorem.} 
Let us picture 
 the boundary components of $S$ without  marked points as punctures. 
In the presence of marked points let us choose a single distinguished point 
on each connected component of the punctured boundary $\partial \widehat S$.
 We define an {\it ideal triangulation}
 $T$ of $\widehat S$ 
as a triangulation of $S$ with all vertices 
at the punctures and distinguished points. Each triangle of this triangulation
can be considered as a disc with three marked points. 
These points are located on the edges of the triangle, one per each edge. 
The restriction of an element of ${\cal X}_{G, \widehat S}$ to 
a triangle $t$ of an ideal triangulation provides a framed $G$-local system 
on the marked triangle $\widehat t$. 
It is  an element of ${\cal X}_{G, \widehat t}$. So we get a projection 
$p_t: {\cal X}_{G, \widehat S} \to {\cal X}_{G, \widehat t}$. 
An edge of the triangulation $T$ is called {\it internal} 
if it is not on the boundary of $S$. Given an oriented internal edge 
$e$ of $T$ we  define a rational projection $p_e: {\cal X}_{G, \widehat S} \to H$. 
Let us orient all internal edges of an ideal triangulation.
Then the collection of projections $p_t$ and $p_e$ provide us a rational map 
\begin{equation} \label{9.11.03.10}
\pi_{\bf T} : {\cal X}_{G, \widehat S} \lra \prod_{\mbox{triangles $t$ 
of $T$}}{\cal X}_{G, \widehat t} 
\times H_{G}^{\mbox{\{internal edges of $T$\}}}.
\end{equation} 

\begin{theorem} \label{4.25.04.10} Assume that 
$G$ has trivial center. Then 
for any ideal triangulation $T$ of $\widehat S$, the map (\ref{9.11.03.10}) is a birational isomorphism. 
\end{theorem}

\begin{figure}[ht]
\centerline{\epsfbox{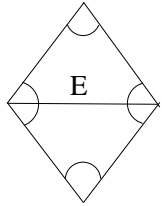}}
\caption{Defining the edge invariant.}
\label{fg86}
\end{figure}

Several different proofs of this theorem are given in Section 6. A similar 
Decomposition Theorem for the space ${\cal A}_{G, \widehat S}$ 
is proved in Section 8, see (\ref{9.28.03.11yx}).
\vskip 3mm 
{\bf Example}. Let $G=PGL_2$. Then for a triangle $\widehat t$ 
the moduli space ${\cal X}_{G, \widehat t}$ is a point, 
and $H = {\Bbb G}_m$. So the Decomposition Theorem in this case states that 
for any ideal triangulation $T$ of $\widehat S$ there is a birational isomorphism
\begin{equation} \label{9.11.03.10asas}
\pi_{\bf T} : {\cal X}_{PGL_2, \widehat S} \lra {\Bbb G}_m^{\mbox{\{internal edges of $T$\}}}.
\end{equation}
It is defined as follows. Take an internal edge $E$ of the tringulation. 
Then there are two triangles of the triangulation sharing this edge, which form a rectangle $R_E$, 
see Figure \ref{fg86}. 
In this case ${\cal L}_{\cal B}$ is a local system of projective lines. A framing provides 
a section of this local system over each corner of the rectangle. Since the rectangle is contractible, 
they give rise to a cyclically ordered configuration of $4$ points on ${\Bbb P}^1$. 
The cross-ratio of these $4$ points, counted from a point corresponding to a vertex of the edge $E$, 
and normalised so that $r(\infty, -1, 0, x) = x$, provides a rational projection  
${\cal X}_{PGL_2, \widehat S}\to {\Bbb G}_m$ corresponding to the edge $E$. 
\vskip 3mm

Let us fix the following data in $G$. 
Let $(B^+, B^-)$ be a pair of opposite 
Borel subgroups. Let $(U^+, U^-)$ be the corresponding pair of maximal unipotent subgroups in $G$. Then  
$B^+\cap B^-$ is identified with the Cartan group $H$. The Cartan 
group  acts on $U^\pm$ by conjugation. 

Recall that given a $G$-set $X$, a configuration of $k$ elements in $X$ is a $G$-orbit in $X^k$. 
Denote by 
${\rm Conf}_k({\cal B})$ the configuration space of $k$-tuples of flags in generic position in $G$, 
meaning that any two of them are opposite to each other. 

If we choose an order of the vertices of the marked triangle $\widehat t$, then there 
is a canonical  birational isomorphism 
${\cal X}_{G', \widehat t} \stackrel{\sim}{\lra} {\rm Conf}_3({\cal B})$. We set $u\cdot B:= uBu^{-1}$. 
Then there is a birational isomorphism 
$$
U^+/H\stackrel{\sim}{\lra} {\rm Conf}_3({\cal B}), \qquad u_+\lms (B^-, B^+, u_+\cdot B^-).
$$
Since the action of $H$ preserves $B^+$ and $B^-$, the formula on the left 
provides a well defined map from $U^+/H$. This plus Decomposition Theorem \ref{4.25.04.10} 
implies the following

\begin{corollary} \label{4.25.04.10s}  
The moduli space ${\cal X}_{G', \widehat S}$ is rational. 
\end{corollary}

Notice that in contrast with this, in general the classical moduli space ${\cal L}_{G, S}$ is by no means rational. In the simplest 
case when $S$ is the sphere with three points removed, description 
of ${\cal L}_{GL_n, S}$ is equivalent to the classical unsolved 
problem of classification of pairs of invertible matrices. 

To explain how the Decomposition Theorem leads to higher Teichm\"uller spaces we recall some 
background on total positivity, and then introduce  a key notion of {positive configurations of flags}. 

\vskip 3mm
{\bf 5. Total positivity, positive configurations of flags and higher Teichm\"uller spaces}. 
A real matrix is called totally positive if all its minors are positive numbers. 
An upper triangular matrix is called totally positive if all its minors which are not identically zero are 
positive. For example 
$$
\left (\begin{matrix}1& a&b\\ 0&1&c\\0&0&1\end{matrix}\right )\quad \mbox{is totally positive if and only if 
$a>0$, $b>0$, $c>0$, $ac-b>0$}.  
$$
The theory of totally positive matrices was developed by  
F.R. Gantmacher and M.G. Krein \cite{GKr}, and by I. Schoenberg \cite{Sch}, in 
1930's. Parametrizations of totally positive unipotent matrices 
were found by A. Whitney \cite{W}. This theory 
was generalized to arbitrary semi-simple real \, Lie \, groups 
by Lusztig \cite{L1}. 
Denote by $U^+(\R_{>0})$ (resp. $H(\R_{>0})$) the subset of all totally positive elements 
in $U^+$ (resp. $H$).

\vskip 3mm
The following definition, followed by Theorem \ref{pcfr},  is one of the crucial steps in our paper.  
\begin{definition} \label{pcf}
A configuration of flags $(B_1, \ldots, B_k)$ is positive if it can be written as 
\begin{equation} \label{pcfs}
\Bigl(B^+, B^-, u_1\cdot B^-, u_1u_2\cdot U^+, \ldots , (u_1 ... u_{k-2})\cdot B^-\Bigr). \qquad 
\mbox{where $u_i \in U^+(\R_{>0})$}
\end{equation}
\end{definition}
One can show that two configurations (\ref{pcfs}) are equal if and only if 
one is obtained from the other by the action of the group $H(\R_{>0})$. 

The definition of totally positive matrices is rather non-invariant: it uses matrix elements. 
Definition \ref{pcf} seems to be even less invariant. However 
it turns out that  
the configuration space of positive flags has remarkable and rather non-trivial properties. 
Let us denote 
by ${\rm Conf}^+_k({\cal B})$ the space of positive configuration of $k$ flags.

\begin{theorem} \label{pcfr}
The cyclic shift $(B_1, \ldots , B_k) \lms (B_2, \ldots , B_k, B_1)$ and reversion 
$(B_1, \ldots , B_k) \lms (B_k, \ldots , B_1)$
preserve the set of positive configurations of flags. 
\end{theorem} 
We suggest that the configuration spaces ${\rm Conf}^+_k({\cal B})$ 
are the geometric objects reflecting the properties of totally positive 
elements in $G(\R)$. 
\begin{figure}[ht]
\centerline{\epsfbox{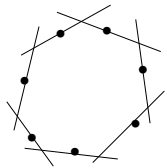}}
\caption{A positive configuration of six flags in ${\Bbb P}^2(\R)$.}
\label{fgo3-15}
\end{figure}
\vskip 3mm
The positive configurations of flags in $PGL_m$ have the following simple geometric description. 
Recall (cf. \cite{Sch}) that 
a curve $C \subset \R{\Bbb P}^{m}$ is {\it convex} if 
any hyperplane intersects it in no more
than $m$ points.   Convex curves 
appeared first in classical problems of analysis (\cite{GKr}). 
\vskip 3mm
\begin{theorem} \label{pcfrsd} A configuration of $k$ real flags $(F_1, \ldots , F_k)$ 
in $\R{\Bbb P}^{m}$ is positive if and only if there exists a (non-unique!) smooth convex curve $C$ 
in $\R{\Bbb P}^{m}$ such that the flag $F_i$ 
is an osculating flag at a point $x_i\in C$, 
and the order of the points $x_1, \ldots, x_k$ is compatible with an orientation of $C$. 
\end{theorem}

\vskip 3mm
{\bf Examples}. The real flag variety for $PGL_2$ is $S^1$. 
A configuration of $n$ points on $S^1$ is positive if and only if 
their order is compatible with an orientation of $S^1$. 
Another example see on Fig. \ref{fgo3-15}.

\vskip 3mm

Generalising Definition \ref{pcf}, and using Theorem \ref{pcfr}, 
we introduce the following key definition. 

\begin{definition} \label{cycuhts} Let $C$ be a cyclic set. A map 
$\beta:C \lra {\cal B}(\R)$ is called positive if for 
any finite cyclic subset $x_1, ..., x_n$ of $S$ the configuration of flags 
$(\beta(x_1), ..., \beta(x_n))$ is positive. 

The space of 
positive configurations of flags parametrsied by a cyclic set $C$ is denoted by  ${\rm Conf}^+_C({\cal B})$. 
\end{definition}
By Theorem \ref{pcfr}, changing a cyclic structure on $C$ to the opposite does not affect 
positivity of a map. 

\noindent \quad \ \ Using positive configurations of real flags we define the higher Teichm\"uller ${\cal X}$-space.  
\begin{definition} \label{q1+++} Let $\widehat S$ be a marked surface with 
boundary, and $G$ a group with trivial center. Then 
the higher Teichm\"uller ${\cal X}$-space is given by 
\begin{equation} \label{q1+++---}
{\cal X}^+_{G, \widehat S}= {\rm Conf}^+_{{\cal F}_{\infty}(\widehat S), \pi_1(S)}({\cal B}).
\end{equation}
\end{definition}
Since $G$ has trivial center, it acts freely on the configurations of $n>2$ flags. 
So by  Lemmas \ref{fras} 
 and \ref{q1-}  a point of ${\cal X}^+_{G,\widehat S}$ is nothing else but a  pair $(\psi, \rho)$,  
where $\rho$ is a representation of $\pi_1(S)$ to $G(\R)$, and $
\psi: {\cal F}_{\infty}(\widehat S) \lra {\cal B}(\R) 
$ is a positive $(\pi_1(S), \rho)$-equivariant 
map. 

There is a similar definition of the higher Teichm\"uller ${\cal A}$-space. 

\vskip 3mm
{\bf Example}. When $\widehat S$ is a disc with $k$ boundary points,  
the  higher Teichm{\"u}ller space ${\cal X}_{G, \widehat S}^+$ is the 
positive configuration space ${\rm Conf}^+_k({\cal B})$. 
\vskip 3mm
Before we proceed any further, let us introduce few simple definitions  relevant 
to  total positivity.

\vskip 3mm
{\bf 6. Positive atlases}. 
Let $H$ be a split algebraic torus over $\Q$, so $H \stackrel{\sim}{=} {\Bbb G}_m^N$. 
 Following [BeK], a rational function $f$ on $H$ is called {\it positive}  if it can be written 
as $f = f_1/f_2$ where  $f_1, f_2$ are linear combinations of characters 
with positive integral coefficients. 
A {\it positive rational map} between two split algebraic tori $H_1, H_2$ 
is a rational map $f: H_1 \to H_2$ such that for any character $\chi$ 
of $H_2$ the composition $\chi \circ f$ is a positive rational function 
on $H_1$. A composition of positive rational functions is positive . 
\vskip 3mm 
 {\bf Example}. The map $x= a+b, y=b$ is positive, but its inverse 
$b=y, a=x-y$ is not. 
\vskip 3mm 
We define a {\it positive divisor} in $H$ as 
the divisor of a positive rational function on $H$. 
Here is our working definition of a positive atlas. 

\begin{definition} \label{7.28.03.1} 
A {\it positive atlas} on an irreducible  scheme/stack $X$ over $\Q$ 
is  a family of  birational isomorphisms
\begin{equation} \label{7.29.03.2}
\psi_{\alpha}: H_{\alpha} \lra X, \quad \alpha \in {\cal C}_X,
\end{equation}
called positive coordinate systems on $X$, 
parametrised by a non-empty set ${\cal C}_X$, 
such that: 

i) each $\psi_{\alpha}$ is regular at the complement to a positive divisor 
in $H_{\alpha}$;

ii)  for any  ${\alpha}, {\beta} \in {\cal C}_X$ the transition map 
$
\psi_{\beta}^{-1} \circ \psi_{\alpha}: H_{\alpha} \lra H_{\beta}
$
is a positive rational map. 
\end{definition}
A {\it positive scheme} is 
a scheme equipped with a positive 
atlas. 
In particular, a positive variety is rational.  
Given a positive scheme $X$, we associate to it an open domain  
$
X({\R_{>0}}) 
\hra X(\R)
$ called  
  the {\it positive part} of $X$. Namely, a point $p \in X(\R)$ lies in 
$X({\R_{>0}})$ if and only if its coordinates in one, and hence any, positive 
coordinate 
system on $X$ are positive.  
Composing  $\psi^{-1}_{\alpha}$ with 
an isomorphism $H_{\alpha} \stackrel{\sim}{\to} {\Bbb G}_m^N$ we get 
  coordinates $\{X^{\alpha}_i\}$ 
on $X$,  providing an isomorphism of manifolds 
\begin{equation} \label{7.29.03.3}
\varphi_{\alpha}: X(\R_{>0}) \stackrel{\sim}{\lra} \R^{{\rm dim} X}, 
\quad p \in X(\R_{>0}) \lms \{\log X^{\alpha}_i(p)\}\in \R^{{\rm dim} X}.
\end{equation}

Let $\Gamma$ be a group of automorphisms of  $X$. We say that a 
positive atlas (\ref{7.29.03.2}) on $X$ is {\it $\Gamma$-equivariant} if 
 $\Gamma$ acts on the set ${\cal C}_X$, and 
for every $\gamma \in \Gamma$ there is an isomorphisms of algebraic tori 
$i_{\gamma}: H_{\alpha} \stackrel{\sim}{\lra} H_{\gamma(\alpha)}$ making the 
following diagram commutative:
\begin{equation} \label{7.20.03.1}
\begin{array}{ccc}
H_{\alpha} & \stackrel{\psi_{\alpha}}{\lra}& X\\
\downarrow i_{\gamma}& &\downarrow \gamma_*\\
H_{\gamma(\alpha)}&\stackrel{\psi_{\gamma(\alpha)}}{\lra}  & X
\end{array}
\end{equation}
Here $\gamma_*$ is the automorphism of $X$ provided by $\gamma$. 
If a positive atlas on $X$ is $\Gamma$-equivariant, then the group 
$\Gamma$ acts on $X({\R_{>0}})$.

\vskip 3mm
{\bf 7. Positive structure of the moduli spaces}. Recall that the mapping class group $\Gamma_{S}$ of $S$ is isomorphic to the group of all outer 
automorphisms of $\pi_1(S)$ preserving the conjugacy classes provided by 
simple loops around 
each of the punctures. It acts naturally on the moduli spaces 
${\cal X}_{G, \widehat S}$, ${\cal A}_{G, \widehat S}$, and ${\cal L}_{G, S}$.

\begin{theorem} \label{7.8.03.1}
Let $G$ be a split semi-simple simply-connected algebraic group, and $G'$ the
corresponding adjoint group. 
Assume that $n>0$. 
Then each of the moduli spaces ${\cal X}_{G', \widehat S}$ and 
${\cal A}_{G,\widehat S}$ 
has a positive $\Gamma_{S}$--equivariant atlas.
\end{theorem}
Here is an outline of our construction of the positive atlas on 
${\cal X}_{G', \widehat S}$.

\vskip 3mm
The work \cite{L1} provides a positive atlas on $U^+$. The set of its $\R_{>0}$-points is 
nothing else but the subset $U^+(\R_{>0})$ discussed in Section 1.5. It induces 
a positive atlas on the quotient $U^+/H$. Thus we get a positive atlas 
on the moduli space ${\cal X}_{G', \widehat t}$. The group $S_3$ permuting the vertices of the triangle $t$ acts on it. We show in Theorem \ref{2.18.03.2a} that 
the group $S_3$ acts by positive birational maps on ${\cal X}_{G', \widehat t}$. Thus 
the $S_3$-orbits of coordinate charts of the 
initial positive atlas give us a new,  $S_3$-equivariant atlas on 
the moduli space ${\cal X}_{G', \widehat t}$. 
Then the  right hand 
side of (\ref{9.11.03.10}) 
provides a positive atlas given by the product of the ones 
on  ${\cal X}_{G', \widehat t}$ and $H$. We prove 
that the transition functions between the atlases provided by different
triangulations are positive. Therefore the 
collection of positive atlases provided by the birational isomorphisms $\{\pi_{T}\}$ for all ideal
triangulations of $\widehat S$ gives rise to 
 a (bigger) positive atlas on  ${\cal X}_{G', \widehat S}$.  
It is $\Gamma_S$-equivariant by the construction. 
 
\vskip 3mm 
For $G= SL_m$ we have a  stronger result. 
A positive atlas is called {\it regular} if the maps  $\psi_{\alpha}$ are  
regular open embeddings. 
\begin{theorem} \label{7.8.03.1q}
The moduli spaces ${\cal X}_{PGL_m, \widehat S}$ and 
${\cal A}_{SL_m,\widehat S}$ 
have regular positive $\Gamma_{S}$--equivariant atlases. 
\end{theorem}
\vskip 3mm 
Theorem \ref{7.8.03.1} is proved in Sections 6 and 9 using the results of Section 5. Theorem \ref{7.8.03.1q} is proved in Sections 9-10
using  simpler constructions of positive atlases   
on the $SL_m$-moduli spaces given in Section 9. 
These constructions can be obtained 
as special cases of the ones given in  Sections 5, 6, and 8, 
but they have some extra nice properties special for $SL_m$, and do not rely 
 on Sections 5, 6, and 8. In particular our $X$- (resp. $A$-)  coordinates on the 
configuration space of triples
of flags (resp. affine flags) 
in an $m$-dimensional vector space  are manifestly invariant under the cyclic 
(resp. twisted cyclic) 
shift. One can not expect this property for a general group $G$, and this
creates some technical issues,  
making formulations and proofs heavier.  In addition to that,  
the construction of positive coordinate systems given in Section 9 
is simpler and can be formulated using the classical language of projective
geometry. 
We continue the discussion of the $SL_m$ case in 
Section 1.15. 

\vskip 3mm 
 {\bf 8. First applications of positivity: higher Teichm{\"u}ller 
spaces and laminations}. Below $G$ stands for a 
split simply-connected semi-simple algebraic group over $\Q$, and $G'$ 
is the corresponding adjoint group. 
\begin{definition} \label{7.29.03.1}
The  
{higher Teichm{\"u}ller spaces} ${\cal X}^+_{G',\widehat S}$ and ${\cal 
A}^+_{G,\widehat S}$ are the real positive parts of the moduli spaces 
${\cal X}_{G',\widehat S}$ and 
${\cal A}_{G,\widehat S}$:
$$
{\cal X}^+_{G',\widehat S}:={\cal X}_{G',\widehat S}(\R_{>0})  \hra {\cal X}_{G',\widehat S}(\R), \qquad 
{\cal A}^+_{G,\widehat S}:={\cal A}_{G,\widehat S}(\R_{>0}) \hra {\cal A}_{G,\widehat S}(\R).
$$ 
\end{definition}

The elements of the space ${\cal X}^+_{G',\widehat S}$ are called 
{\it positive framed $G(\R)$-local systems on $\widehat S$}. 

\begin{theorem} \label{dgdgd} Definition \ref{7.29.03.1} of the 
space ${\cal X}^+_{G', \widehat S}$ 
is equivalent to Definition \ref{q1+++}. 

Thus there is a canonical bijection between 
the positive framed $G(\R)$-local systems $({\cal L}, \beta)$ on $\widehat S$ and 
positive $\pi_1(S)$-equivariant maps 
\begin{equation} \label{4.1.04.10facom}
\Phi_{\cal L, \beta}: {\cal F}_{\infty}(\widehat S) \lra {\cal B}(\R) \qquad \mbox{modulo the
  action of $G(\R)$}. 
\end{equation}
where $\pi_1(S)$ acts on the flag variety  via the monodromy 
representation of the local system ${\cal L}$. 
\end{theorem} 
Theorem \ref{dgdgd} follows from Lemma \ref{fras} plus the 
definition (ii) of the Farey set ${\cal F}_{\infty}(S)$  or its version for $\widehat S$, 
which uses ideal triangulations of $\widehat S$.  
There is a similar interpretation of the space ${\cal A}^+_{G,\widehat S}$. 
 
\vskip 3mm

It follows that the spaces ${\cal X}^+_{G',\widehat S}$ and ${\cal 
A}^+_{G,\widehat S}$ are isomorphic to $\R^N$ for some integer $N$. In particular, if 
$\widehat S = S$, i.e. there are no marked points on the boundary of $S$, we have isomoprphisms
\begin{equation} \label{11.9.05.10}
{\cal X}^+_{G', S} \stackrel{\sim}{\lra} \R^{-\chi(S){\rm dim}G}, \quad 
{\cal A}^+_{G, S} \stackrel{\sim}{\lra} \R^{-\chi(S){\rm dim}G}
\end{equation}

\vskip 3mm
The mapping class 
group  $\Gamma_{S}$ acts on ${\cal X}^+_{G',\widehat S}$ 
and ${\cal A}^+_{G,\widehat S}$.  
Let 
\begin{equation} \label{11.9.05.10sd}
{\cal L}_{G', S}^+:= \pi\Bigl({\cal X}^+_{G', S}\Bigr) \subset {\cal L}_{G', S}(\R)
\end{equation}
We call the points of  ${\cal L}_{G', S}^+$ 
{\it positive local systems} on $S$, and their monodromy representations 
{\it positive representations} of $\pi_1(S)$. 
There is a natural embedding
$$
{\cal L}^+_{G', S} \hra {\cal X}_{G', S}^+
$$
splitting the canonical projection ${\cal X}^+_{G', S} \to {\cal L}_{G', S}^+$.

\vskip 3mm
In Section 11  we show that  our construction,   
specialized to the case  $G= SL_2$  for the
${\cal A}$-space and $G = PGL_2$ for the ${\cal X}$-space, gives the classical 
Teichm{\"u}ller spaces ${\cal T}^d(S)$ and ${\cal T}^+(S)$:
\begin{theorem} \label{3.27.02.1}
Let $S$ be a hyperbolic surface with $n>0$ holes. Then 

a) The 
  Teichm{\"u}ller space ${\cal T}^+_S$ is isomorphic to ${\cal X}^+_{PGL_2,
  S}$.
Further, ${\cal L}_{PGL_2, S}^+$ is the classical 
Teichm{\"u}ller space ${\cal T}_S$. 

b) The decorated Teichm{\"u}ller space ${\cal T}^d_S$ is isomorphic to 
 ${\cal A}^+_{SL_2, S}$. 
\end{theorem}
\ Combining this theorem with the isomorphisms
\ ${\cal X}^+_{PGL_2, S} \stackrel{\sim}{=} \R^{6g-6+3n}$\  and 
${\cal A}^+_{SL_2, S} \stackrel{\sim}{=} \R^{6g-6+3n}$ 
we arrive at the isomorphisms of manifolds 
$$
{\cal T}^+_S   \stackrel{\sim}{=} \R^{6g-6+3n} \quad \mbox{and} \quad {\cal T}^d_S 
  \stackrel{\sim}{=} \R^{6g-6+3n}. 
$$ 
The first  gives an elementary proof of the classical 
Teichm{\"u}ller theorem for surfaces with holes.

We introduce {\it integral, rational and real 
 $G$--lamination spaces} as the  tropical limits of the 
positive  moduli spaces ${\cal X}_{G', \widehat S}$ and ${\cal A}_{G, \widehat S}$, 
i.e. the sets of points of these spaces in the semifields $\Z^t, \Q^t, \R^t$. 
The projectivisations of the real lamination spaces are 
  the Thurston-type boundaries of the corresponding Teichm{\"u}ller spaces.

\vskip 3mm 
{\bf 9.  Positivity, hyperbolicity and discreteness 
of the monodromy representations}.   
Consider the universal $G$-local system  on
the moduli space $S \times {\cal X}_{G, S}$. 
Its fiber over $S \times p$ is the local system corresponding to the
point $p$ of the moduli space ${\cal X}_{G, S}$.
Let ${\Bbb F}_{G, S}$ be the field of rational functions on the moduli
space ${\cal X}_{G, S}$. The
monodromy of  the universal local system around a loop on $S$ is a conjugacy class in
$G({\Bbb F}_{G, S})$.

Observe that given a positive variety $X$, there is a 
subset $\Q_+(X)$ of the field of rational functions of $X$ consisting of the 
functions which are positive in one, and hence in any of the 
positive coordinate systems on $X$. Clearly  $\Q_+(X)$ is a semifield, 
the {\it positive semifield} of the positive variety $X$.

Theorem \ref{7.8.03.1} delivers a positive atlas on the moduli space ${\cal X}_{G, S}$.
Therefore we have the corresponding positive semifield   ${\Bbb F}^+_{G, S}$ 
in the field ${\Bbb F}_{G, S}$. 
The group $G$ has a 
positive atlas provided by the birational isomorphism $G \stackrel{\sim}{\lra} U^- H U^+$ 
(the Gauss decomposition) and the 
positive atlases on the factors.  So we 
have a well defined  subset $G({\Bbb F}^+_{G, S})$. 

\begin{theorem} \label{2.1.04.0} 
The  monodromy of the universal $G$-local system on $S \times {\cal
  X}_{G, S}$ around 
any non-trivial non-boundary loop on $S$ can be conjugated in $G({\Bbb
  F}_{G, S})$  to an element of $G({\Bbb F}^+_{G, S})$. 

The monodromy around a boundary component can be conjugated to an element 
of $B^{\pm}({\Bbb F}^+_{G, S})$, where $B^{\pm}$ means one of the Borel subgroups 
$B^+$ or $B^-$.  
\end{theorem}

Observe that although this theorem is about monodromies of 
local systems on $S$, it can not be even stated without 
introducing the moduli space ${\cal X}_{G, S}$ and a positive atlas on it. 
Here is a corollary 
formulated entirely in classical terms.

Let $H^0$ be the subvariety of the Cartan group $H$ where 
the Weyl group $W$ acts without fixed points. An element of $G(\R)$ is called {\it 
positive hyperbolic} if it is conjugated to an element of $H^0(\R_{>0})$. 
For example, for  $GL_n(\R)$ these are the elements 
with distinct real positive eigenvalues. 

\begin{theorem}\label{2.1.04.0f} 
The monodromy of a positive local system on a surface $S$ with boundary is faithful. 
Moreover its monodromy 
around any non-trivial non-boundary loop is positive hyperbolic. 
\end{theorem}

It follows from Theorem \ref{2.1.04.0}. 
Indeed, Definition \ref{7.29.03.1} and 
Theorem \ref{2.1.04.0} imply that the monodromy of a positive representation 
around a homotopy nontrivial non-boundary loop is conjugated 
to an element of $G(\R_{>0})$. By
 the Gantmacher-Krein theorem  \cite{GKr} for $G= GL_n$, and 
Theorem 5.6 in \cite{L1} in general, any element of 
$G(\R_{>0})$ is positive hyperbolic. To prove faithfulness observe that 
the identity element belongs neither to $G(\R_{>0})$ nor to 
$B^{\pm}(\R_{>0})$.

\vskip 3mm

\begin{theorem} \label{7.8.03.2} Let $G(\R)$ be a split real semi-simple  
Lie group  with trivial center. Then the image of a  positive representation of $\pi_1(S)$ is 
a discrete subgroup in $G(\R)$. 
\end{theorem}

Theorem \ref{7.8.03.2}  is proved in Section 7.1 using Theorem \ref{dgdgd}. Since the  isomorphisms $\psi_{\alpha}$
describing the positive atlas on  ${\cal X}_{G, \widehat S}$ 
are defined explicitly, 
we obtain a  parametrisation of  
a class of discrete subgroups 
in $G(\R)$ which generalizes the Fuchsian subgroups of $PSL_2(\R)$.

\vskip 3mm 
{\bf 10. \!Universal higher Teichm{\"u}ller spaces and positive $G(\R)$-opers}. 
Theorem \ref{dgdgd}  suggests a definition of the {\it universal higher Teichm{\"u}ller spaces} 
${\cal X}_{G}^+$. 

Let us furnish ${\Bbb P}^1(\Q)$ with a cyclic structure 
provided by the natural embedding ${\Bbb P}^1(\Q) \subset {\Bbb P}^1(\R)$ and 
an  orientation of ${\Bbb P}^1(\R)$. 
The cyclic set of the vertices of the Farey triangulation is identified 
with ${\Bbb P}^1(\Q)$. Thus we can identify the cyclic sets ${\cal F}_{\infty}(S)$ and 
${\Bbb P}^1(\Q)$.

\begin{definition} \label{uhts}
The universal higher Teichm{\"u}ller space ${\cal X}_{G}^+$ 
consists of all positive maps
\begin{equation} \label{rho}
\beta: {\Bbb P}^1(\Q) \lra {\cal B}(\R) \qquad \mbox{modulo the
  action of $G(\R)$}.
\end{equation} 
\end{definition}
For $G=PGL_2$ we get the universal Teichm{\"u}ller space 
considered by Bers \cite{Bers} and Penner \cite{P3}. 

For a surface $S$ with boundary, the higher Teichm{\"u}ller spaces 
${\cal X}_{G, S}^+$ can be  embedded to the universal one as follows. 
Let $\Delta$ be a torsion free subgroup of $PSL_2(\Z)$  
and $S_{\Delta}:= {\cal H}/\Delta$, where ${\cal H}$ is the hyperbolic
plane.  
Then $S_{\Delta}$ is equipped with a distinguished ideal
triangulation, given by the image of the Farey triangulation of 
the hyperbolic plane (Figure (\ref{Far})). 
 Conversely, the subgroup $\Delta \subset PSL_2(\Z)$ 
is determined uniquely up to conjugation by an ideal triangulation of $S$. 
The set of vertices of the Farey triangulation is identified with ${\Bbb P}^1(\Q)$. 
Theorem \ref{dgdgd} implies:

\begin{corollary} \label{4.1.04.10}
The space ${\cal X}_{G, S_\Delta}^+$ parametrises pairs 
(a representation $\rho: \Delta \to G(\R)$, a  $\rho$-equivariant map
(\ref{rho})). 
\end{corollary}

The universal higher Teichm{\"u}ller  space ${\cal X}_{G}^+$ is equipped with a 
positive coordinate atlas. 
Observe that the set of vertices of the Farey triangulation 
is identified with ${\Bbb
  P}^1(\Q)$. Using the Farey triangulation  
just the same way as ideal triangulations of
$S$, we obtain the Decomposition Theorem for 
the universal higher Teichm{\"u}ller  space ${\cal X}_{G}^+$, providing 
coordinate systems on ${\cal X}_{G}^+$:

\begin{theorem} \label{1.29.04.3a} There exists a canonical isomorphism 
 \begin{equation} \label{1.29.04.7a} 
\varphi_G: {\cal X}^+_G \stackrel{}{\lra} \prod_{{\mbox{$t$: Farey
  triangles}}}{{\cal X}^+_{G, \widehat t}} ~~\times\prod_{{\mbox{Farey
  diagonals}}}  H(\R_{>0}).
\end{equation} 
\end{theorem} 
In  Section 8.5 we define universal higher Teichm{\"u}ller  spaces 
${\cal A}_{G}^+$  and prove  similar results for them.

\vskip 3mm
{\it Positive continuous  maps} 
\begin{equation} \label{rho1}
S^1 \lra {\cal B}(\R)
\end{equation} 
are objects of independent interest. 
We prove in Section 7.8 
that positive smooth maps (\ref{rho1})  are integral curves  
of a canonical non-integrable distribution on the flag variety, provided by the
simple positive roots, whose dimension 
equals  the rank of $G$. So they can be viewed as the $G(\R)$-opers,   
 in the sense 
of Beilinson and Drinfeld, on the circle. We call them 
{\it positive $G(\R)$-opers}. 

It is well known that 
$PGL_m(\R)$-opers are nothing else 
but smooth curves in $\R{\Bbb P}^{m-1}$: a smooth projective curve 
gives rise to its osculating curve in the flag variety.  
 We show in Section 9.12 
that { positive $PGL_m(\R)$-opers} 
correspond  to smooth convex curves in $\R{\Bbb P}^{m-1}$. 
More generally, positive continuous maps  (\ref{rho1}) give rise to 
$C^{1}$-smooth convex curves. 
Thus positive continuous maps (\ref{rho1}) generalize  
  projective 
convex curves  
to the case 
of an arbitrary reductive group $G$.

\vskip 3mm
{\bf 11.  Higher Teichm\"uller spaces for closed surfaces}. \footnote{Most of the results of Sections 1.11-1.13
 were obtained in 2005; Theorem \ref{H1} and the results of Section 1.18 
were obtained in 2006.} 
Given a surface $S$, with or without boundary, 
there is a countable cyclic $\pi_1(S)$-set ${\cal G}_{\infty}(S)$ defined as follows. 
Choose a hyperbolic structure with geodesic boundary on $S$, and lift all 
geodesics on $S$ to the universal cover. The universal cover 
can be viewed as a part of the hyperbolic plane ${\cal H}$ (Fig. \ref{fg84}). 
The endpoints of the preimages of non-boundary geodesics 
form a subset ${\cal G}_{\infty}'(S)$ of the absolute $\partial{\cal H}$. 
It inherits from the absolute 
a structure of the cyclic $\pi_1(S)$-set, 
which does not depend on the choice of a hyperbolic structure on $S$. Set $$
{\cal G}_{\infty}(S):= 
{\cal F}_{\infty}(S)\cup {\cal G}_{\infty}'(S).
$$ Using the definition (iii) of the 
Farey set (Section 1.3), one easyly sees that 
${\cal G}_{\infty}(S)$ is cyclic $\pi_1(S)$-set extending  ${\cal F}_{\infty}(S)$ and ${\cal G}_{\infty}'(S)$. 
When $S$ is closed, the  Farey set ${\cal F}_{\infty}(S)$ is empty. 
A hyperbolic structure on $S$ provides ${\cal G}_{\infty}(S)$
with a topology induced from $\partial {\cal H}$. It does not depend on the choice of the hyperbolic structure. 

\vskip 3mm
According to Theorem \ref{dgdgd},  
a framed positive local system $({\cal L}, \beta)$ on a surface $S$ with boundary 
provides a $\pi_1(S)$-equivariant positive map (\ref{4.1.04.10facom}).

By Theorem \ref{2.1.04.0}, the monodromy $M_{\gamma}$ 
of a positive local system along a non-boundary
 closed geodesic $\gamma$ 
on $S$ is conjugated to an
element of $G(\R_{>0})$. Therefore there exists a distinguished flag preserved by the monodromy along $\gamma$: 
for $G= PGL_m(\R)$ it is provided by the 
ordering of the eigenspaces of $M_{\gamma}$ in which their  
eigenvalues increase; for the general case see Theorem 8.9 in \cite{L1}. 
Thus, going to the universal cover of $S$, we get a $\pi_1(S)$-equivariant map 
\begin{equation} \label{4.1.04.10f}
\Psi_{{\cal L}, \beta}: {\cal G}_{\infty}(S) \lra {\cal B}(\R)
\end{equation} 
extending the map (\ref{4.1.04.10facom}). 
\begin{theorem} \label{4.1.04.10abc}
Let $S$ be a compact surface with boundary. Then the map (\ref{4.1.04.10f}) 
is a positive map. Moreover, it is a continuous map. 
\end{theorem}

The notion of positive 
configurations of flags and Theorems \ref{dgdgd} and \ref{4.1.04.10abc} suggest the following definition 
of higher Teichm\"uller spaces for compact surfaces $S$ without boundary.

\begin{definition} \label{4.1.04.10a}
Let $S$ be a compact surface with or without boundary. 
A representation $\rho: \pi_1(S) \to G(\R)$ 
is {\em positive}, if there exists a  positive $\rho$-equivariant map (\ref{4.1.04.10f}). 
The moduli space of positive representations is denoted by ${\cal L}^+_{G,S}$. 
\end{definition}

\begin{corollary} \label{4.1.04.10an}
For 
surfaces with boundary  Definition \ref{4.1.04.10a} is equivalent to the one we used before:
 the $\rho$-equivariant maps (\ref{4.1.04.10f}) modulo $G(\R)$-conjugation 
are in bijection with 
points of  ${\cal X}_{G, S}^+$.
\end{corollary}
Indeed, the restriction of the map (\ref{4.1.04.10f}) to ${\cal F}_{\infty}(S)$ provides 
a $\rho$-equivariant map
(\ref{4.1.04.10facom}). Since ${\cal F}_{\infty}(S)$ is dence in ${\cal G}_{\infty}(S)$, 
Theorem \ref{4.1.04.10abc} tells us that 
it extends uniquely from ${\cal F}_{\infty}(S)$ to 
${\cal G}_{\infty}(S)$.

\vskip 3mm
Here is an interpretation of positive representations  as $\pi_1(S)$-equivariant 
configurations of real flags parametrised by the cyclic set ${\cal G}'_{\infty}(S)$ which 
follows from Lemma \ref{q1-} and  
Definition \ref{4.1.04.10a}:
\begin{equation} \label{q1+}
{\cal L}^+_{G,S}= {\rm Conf}^+_{{\cal G}'_{\infty}(S), \pi_1(S)}({\cal B}).
\end{equation}

\vskip 3mm
The following results show  that the moduli space of 
positive representations have the basic features of the classical Teichm\"uller spaces. 
In particular, we extend Theorems \ref{2.1.04.0f}, \ref{7.8.03.2}  and  
the isomorphisms (\ref{11.9.05.10}) to the case of surfaces 
without boundary: 

\begin{theorem} \label{4.1.04.10b}
Let $G$ be a split real semi-simple Lie group with trivial center, and 
$S$ a compact surface without boundary. Then 

(i) A positive representation $\rho: \pi_1(S) \to G(\R)$ 
is faithful, its image is a discrete subgroup in $G(\R)$, and the image of any non-trivial element 
is positive hyperbolic. 

(ii) The moduli space ${\cal L}^+_{G,S}$ 
of positive representations is diffeomorphic to $\R^{-\chi(S) {\rm dim}G}$. 
\end{theorem}
\vskip 3mm

There is an embedding of the classical 
Teichm\"uller space into ${\cal L}^+_{G,S}$ provided by 
a Fuchsian subgroups in 
a principal $PGL_2(\R)$-subgroup of $G(\R)$. 

\vskip 3mm
When $S$ has no boundary, Nigel Hitchin \cite{H1} studied 
topology of the space of
representations $\pi_1(S) \to G(\R)$ modulo conjugation 
such that the adjoint action of $\pi_1(S)$ on the Lie algebra of $G(\R)$ is completely reducible. 
He proved that the component of this space containing 
(via the principal embedding $PGL_2(\R) \hra G(\R)$) the classical 
Teichm\"uller space is diffeomorphic to $\R^{-\chi(S) {\rm dim}G}$. 
It is called the Hitchin component. 
Furthermore, given a complex structure on $S$, let us denote by by ${\Bbb S}$ the corresponding Riemann surface. 
Hitchin proved that there is an isomorphism ($r={\rm rk}(G)$)
\begin{equation} \label{4.22}
\mbox{\rm the Hitchin component}  = H^0({\Bbb S}, \oplus_{i=1}^r\Omega^{d_i}_{\Bbb S})
\end{equation}
where $d_1, ..., d_r$ are the so-called exponents of $G$, e.g. $2, 3, ..., m-1$ when $G = PGL_m$, and 
$\Omega_{\Bbb S}$ is the sheaf of holomorphic differentials on ${\Bbb S}$.  
However the isomorphism (\ref{4.22}) is not invariant under the mapping class group action. 

Hitchin's approach is analytic. It is based on deep analytic results  
of Hitchin \cite{H2},  Corlette \cite{C}, Simpson \cite{S} 
and Donaldson \cite{D}, 
and  makes an essential use of a complex structure on $S$. 
It is very natural but rather non-explicit. 
Our approach 
 is different: it is combinatorial, explicit,
and  does not use a complex structure on $S$. 
The two approaches are completely independent.

\vskip 3mm
In the case when $S$ has no boundary, and $G=PGL_m(\R)$, a class of 
representations of $\pi_1(S)$, called Anosov representations,  
has been defined by Francois Labourie \cite{Lab} using his  Anosov structures, 
the boundary at infinity $\partial_{\infty}\pi_1(S)$ of $\pi_1(S)$
and (hyper)convex curves in $\R{\Bbb P}^{m-1}$. Labourie (with a complement by O. Guichard 
\cite{Gui}) proved that  Anosov representations are 
exactly the ones from Hitchin's component for $PGL_m(\R)$. 
Let us  clarify the connection with his work. 

Recall the following description 
of the  boundary at infinity $\partial_{\infty}\pi_1(S)$ of $\pi_1(S)$. 
Choose a hyperbolic structure with geodesic boundary on $S$. Take the universal cover $D$ of $S$. 
If $S$ has no boundary, it is ${\cal H}$. Otherwise it is obtained by removing from 
${\cal H}$ half discs 
bounded by the preimages of the boundary geodesics on $S$, see Figure \ref{fg84}. 
Now $\partial_{\infty}\pi_1(S)$ is the intersection of the absolute of ${\cal H}$ with the 
closure of $D$. Its cyclic structure is induced by the one of the absolute. The group $\pi_1(S)$ acts by 
 deck transformations. 
For us it is important only that it is a cyclic $\pi_1(S)$-set, which 
does not depend on the choice of hyperbolic structure on $S$.  
So there are  canonical inclusions of cyclic $\pi_1(S)$-sets, where the first two are countable sets: 
$$
{\cal F}_{\infty}(S) \subset {\cal G}_{\infty}(S)\subset \partial_{\infty}\pi_1(S).
$$

\begin{theorem} \label{4.1.04.10ab} For a positive representation the map (\ref{4.1.04.10f}) 
extends uniquely to a continuous map
\begin{equation} \label{4.1.04.10fg}
\overline \Psi_{\rho}: {\partial}_{\infty}\pi_1(S) \lra {\cal B}(\R).
\end{equation} 
\end{theorem} 
The image of this map is the {\it limit set}
 for the corresponding framed positive representation. 
Using 
a homeomorphism $\partial_{\infty}\pi_1(S,p) \to S^1$, we conclude that positive $G(\R)$-opers 
appear as the limit set maps for positive representations for closed surfaces $S$.  

Using Theorems 1.3 and \ref{4.1.04.10ab} one can show that, 
for a surface $S$ without boundary and $G=PGL_m(\R)$, Anosov representations \cite{Lab} 
and our positive representations are the same. 

\begin{theorem} \label{H1} Assume that $S$ has no boundary, and $G$ has trivial center. 
Then the moduli space of positive representations ${\cal L}^+_{G,S}$ coincides with the Hitchin component 
in the representation space of $\pi_1(S)$ to $G(\R)$. 
\end{theorem}

This gives  another proof of the fact that 
for $G= PGL_m(\R)$ 
the positive representations are the same as the Anosov 
representations. 

\vskip 3mm
Here is the scheme of the proof of Theorem \ref{H1}. 
For a general group $G$ the cutting and gluing results, reviewed below, 
imply that positive representations form an open connected subset 
of the Hitchin component. It  contains the classical Teichm\"uller space. 
So to complete the proof we have to show that the space ${\cal L}_{G,S}^+$ is closed 
in the Hitchin component, which is done in Section 7.9.

\vskip 3mm
{\bf 12. Cutting and gluing}. Let $S$ be a surface, with or without boundary, with $\chi(S)<0$. 
The moduli space ${\cal L}^+_{G,S}$ of positive $G(\R)$-local systems on $S$ was defined for surfaces with boundary in (\ref{11.9.05.10sd}), and for  closed $S$
 in Definition \ref{4.1.04.10a}. 

Let $\gamma$ be a non-trivial loop on $S$. 
We assume that $\gamma$ is not homotopic to a boundary component of $S$. 
Denote by $S'$ the surface obtained by cutting $S$ along $\gamma$. It has two boundary components, $\gamma_+$ and 
$\gamma_-$, whose orientations are induced by the one of $S'$.  The surface $S'$  has one or two components, 
each of them of negative Euler characteristic. Denote by ${\cal L}^+_{G,S'}(\gamma_+, \gamma_-)$ 
the subspace of 
 ${\cal L}^+_{G,S'}$  given by the following condition:
\begin{equation} \label{11.13.05.3df}
\parbox{11cm}{\it The monodromies along the loops $\gamma_+$ and 
$\gamma_-$ are positive hyperbolic and mutually inverse.}
\end{equation}

\begin{theorem} \label{11.13.05.1}
Let 
$S$ be a surface with $\chi(S)<0$, and $S'$ is obtained by cutting along a loop $\gamma$, as above. 
Then the restriction from $S$ to $S'$ provides us a principal 
$H(\R_{>0})$-bundle
\begin{equation} \label{11.13.05.3}
{\cal L}^+_{G,S} \lra {\cal L}^+_{G,S'}(\gamma_+, \gamma_-), \qquad {\cal L} \lms {\cal L}|_{S'}.
\end{equation}
The space \ ${\cal L}^+_{G,S}$\  is a connected topologically trivial domain of dimension\  $-\chi(S){\rm dim}G$.
\end{theorem}

The first claim of the theorem  just means that 

\begin{enumerate}

\item
The restriction of a positive local system on $S$ to $S'$ is positive.
\item
The image of the restriction map consists of all positive local systems on $S'$ satisfying the 
constraint (\ref{11.13.05.3df}) on the monodromies around the oriented loops  
$\gamma_+$ and $\gamma_-$. 
\item
Given a positive local system on $S'$ satisfying the constraint (\ref{11.13.05.3df}), 
one can glue it to a positive local system on $S'$, and the group $H(\R_{>0})$ acts simply transitively 
on the set of the gluings. 
\end{enumerate}
 
In particular, for every  non-trivial non-boundary loop $\gamma$ on $S$, there is an action $t_{\gamma}$ 
of the group $H(\R_{>0})$ on ${\cal L}^+_{G,S}$ without fixed points. 
One can show that this action is Hamiltonian for the canonical Poisson structure on ${\cal L}^+_{G,S}$: 
the Hamiltonians for this action are provided by  the monodromy along $\gamma$, understood as a map 
${\cal L}^+_{G,S}\to H(\R_{>0})/W$. 

\vskip 3mm 
{\bf 13. Configuration spaces and laminations}. 
We define 
the lamination spaces for closed surfaces in  Section 6.9. Here is how it is going. 
\vskip 3mm
For an arbitrary countable cyclic set $C$, 
there is a projective limit of positive spaces ${\rm Conf}_{C}({\cal B})$, 
obtained by taking the projective limit 
over finite subsets $C'\subset C$ 
(we underline  $\underline {\cal B}$  to stress that it is an algebraic variety): 
$$
{\rm Conf}_{C}(\underline {\cal B}):= \lim_{\longleftarrow} 
{\rm Conf}_{C'}(\underline {\cal B}).
$$
If a group $\pi$ acts on $C$, it acts on it, so taking the invariants under the action of the group $\pi$ 
we get a subspace ${\rm Conf}_{C, \pi}({\cal B}):= 
{\rm Conf}_{C}({\cal B})^\pi$. So for any semifield $K$ there is a well-defined set 
$$
{\rm Conf}_{C, \pi}(\underline {\cal B})(K) = {\rm Conf}_{C}(\underline {\cal B})(K)^\pi.
$$

So for a closed surface $S$ we define the lamination spaces by taking the tropical points:
\begin{equation} \label{q1++}
\parbox{9cm}{\it The space of $G$-laminations on $S$ with coefficients in ${\Bbb A}$:= 
${\rm Conf}_{{\cal G}_{\infty}(S), \pi_1(S)}(\underline {\cal B})({\Bbb A}^t)$.}
\end{equation} 
The space ${\cal L}_{G,S}$ 
has a similar interpretation, see (\ref{4.1.04.10a}) for closed $S$ and (\ref{q1+++---}) for marked 
surfaces $\widehat S$ with boundary.  So, for any $S$, both  
the higher Teichm\'uller space and laminations spaces can 
be defined as the sets of $\R_{>0}$ and tropical points of certain  configuration spaces.

\vskip 3mm 
{\bf 14. The Weil-Petersson form on the moduli space 
${\cal A}_{G, \widehat S}$
and its $K_2$ and motivic avatars}. We define 
the Weil-Petersson form $\Omega_{G, \widehat S}$ on the moduli space 
${\cal A}_{G, \widehat S}$ by constructing  
a  $K_2$-class on 
${\cal A}_{G, \widehat S}$. We suggest that
$$
 \mbox{\it many interesting symplectic structures 
can be upgraded to their $K_2$-avatars}. 
$$
Let us explain how 
a $K_2$-class provides a symplectic form. 
Recall that thanks to the Matsumoto theorem \cite{Mi} 
the group $K_2(F)$ of  a field $F$ is the quotient of the abelian group 
$\Lambda^2F^*$ by the subgroup generated by the elements  
$(1-x) \wedge x$, $x \in F^* -\{1\}$, called the Steinberg relations. 
If $F = \Q(Y)$ is the field of 
rational functions on a variety $Y$ there is a homomorphism
$$
d\log: K_2(\Q(Y)) \lra \Omega_{\rm log}^2({\rm Spec}(\Q(Y))); 
\quad f_1 \wedge f_2 \lms d\log f_1 \wedge d\log f_2.
$$
Applying it to a class $w\in K_2(\Q(Y))$ we get a rational $2$-form $d\log(w)$ 
with logarithmic singularities. Let us formulate a condition on $w$ which guarantees that the $2$-form $d\log(w)$ is non singular, and hence defines a degenerate symplectic structure of $Y$. 
Let  $Y_1$  be the set of all irreducible divisors in a variety $Y$. 
There is a {\it tame symbol} homomorphism: 
\begin{equation} \label{12.20.02.1es}
{\rm Res}: K_2(\Q(Y)) \lra \bigoplus_{D \in Y_1} \Q(D)^*; \qquad \{f,g\} \lms 
{\rm Rest}_D\Bigl( f^{v_D(g)}/g^{v_D(f)}\Bigr)
\end{equation}
where ${\rm Rest}_D$ is restriction to the generic point of $D$, and 
 $v_D(g)$ is the order of zero of $g$ at the generic point of the divisor $D$. 
Set 
\begin{equation} \label{12.20.02.1}
\begin{array}{rcl}H^2(Y, \Z_{\cal M}(2))&:=& {\rm Ker}({\rm Res}) \hookrightarrow K_2(\Q(Y)),\\
H^2(Y, \Q_{\cal M}(2))&:=& H^2(Y, \Z_{\cal M}(2))\otimes \Q.
\end{array}
\end{equation}
 The  map $d\log$  transforms the tame symbol to the residue map. So 
its restriction to  (\ref{12.20.02.1}) provides a homomorphism
$$
d\log: H^2(Y, \Z_{\cal M}(2)) \lra \Omega^2(Y_{0})
$$ 
where $Y_0$ is the nonsingular part of  $Y$. 
We define a class
\begin{equation} \label{9.5.03.1}
W_{G, \widehat S} \in H^2({\cal A}_{G, \widehat S}, \Q_{\cal M}(2)) ^{\Gamma_{S}}
\end{equation}
and set
$$
\Omega_{G, \widehat S}:= d\log(W_{G, \widehat S}) 
\in \Omega^2({\cal A}_{G, \widehat S})^{\Gamma_{S}}.
$$

We also define the class $W_{G, S}$ for a compact oriented 
surface $S$. One can show that for $G = SL_2$ the form 
$\Omega_{SL_2,  S}$ coincides 
with W. Goldman's symplectic structure \cite{Gol1} on the space 
${\cal L}_{SL_2, S}$. Its restriction to the Teichm{\"u}ller 
component is the Weil-Petersson symplectic structure. 
The $K_2$-avatar of the Weil-Petersson form 
is a new object even in the classical case $G= SL_2$. 
\vskip 3mm
{\it Motivic avatar of the Weil-Petersson form}. 
We show that the $K_2$-class 
 $W_{G,\widehat  S}$ is a manifestation 
of a richer algebraic structure described by using 
 the weight two  motivic 
cohomology understood via the Bloch-Suslin complex, 
and the motivic dilogarithm. Namely, we construct a  
 {\it motivic avatar} 
$$
{\Bbb W}_{G,\widehat  S} \in H^2_{\Gamma_{S}}({\cal A}_{G,\widehat S}, \Q_{\cal M}(2))
$$
of the form $\Omega_{G, \widehat S}$. It lives 
 in the $\Gamma_{S}$--equivariant weight two 
motivic cohomology of ${\cal A}_{G, \widehat S}$.

The $2$-form $\Omega_{G, \widehat S}$ as well as its 
avatars are lifted from similar objects on the moduli space ${\cal U}_{G,
  \widehat S}$ of the unipotent framed $G$-local systems on $\widehat S$. 
The corresponding $2$-form on the moduli space ${\cal U}_{G,
  \widehat S}$ is non degenerate, and thus provides a symplectic structure on
this space. 

Our construction uses  as a building block  the second motivic 
Chern class $c_2^{\cal M}$ of the universal $G$--bundle over 
the classifying space $BG_{\bullet}$. 
When $G= SL_m$, 
there is an explicit cocycle for the 
class $c_2^{\cal M}$ constructed 
in \cite{G3}. Applying to it our general procedure 
we arrive at  an explicit cocycle for the class ${\Bbb W}_{SL_m, \widehat S}$. 
 Magically this construction delivers  canonical coordinates 
on the space  ${\cal A}_{SL_m, \widehat S}$. 
Their properties are outlined in the next subsection.

\begin{figure}[ht]
\centerline{\epsfbox{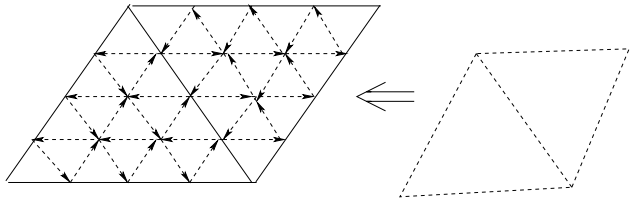}}
\caption{The $4$-triangulation arising from  a triangulation of a surface.}
\label{fg32}
\end{figure}

\vskip 3mm 
{\bf 15. \!The (orbi)-cluster ensemble structure in the\! $PGL_m/SL_m$ case}. 
Let $T$ be an  ideal  triangulation of $\widehat S$. 
Take the triangle 
$$
x+y+z=m, \quad x,y,z \geq 0
$$
 and consider its triangulation 
given by the lines $x=p$, $y=p$, $z=p$ where 
$0 \leq p \leq m$ is an integer. 
The {\it $m$--triangulation} of 
a triangle is a triangulation isotopic to this one. 
Let us $m$--triangulate each triangle of the triangulation 
$T$,  so each side gets $m+1$ vertices.   
We get 
a subtriangulation, called the 
{\it $m$--triangulation} of  $T$.   
An edge of the $m$--triangulation of $T$ 
is called an {\it internal edge} if it does not lie 
on a side of the original  triangulation $T$. 
The orientation of $S$  provides orientations of 
the internal edges. Indeed, 
take a triangle $t$ of the triangulation $T$. 
The orientation of $S$ provides an orientation of its boundary. 
An internal edge $e$ sitting inside of $t$ is parallel to a certain side of $t$, so the orientation of this side induces an orientation of $e$. 
This is illustrated on Figure \ref{fg32}. 
Here on the right there are  two adjacent triangles of the triangulation 
$T$. On the left we show the $4$--triangulation of these triangles. 
Their internal edges are oriented accordantly to 
the cyclic structure coming from the clockwise orientation of the plain.

 Given an  ideal  triangulation $T$ of $\widehat S$, we define 
 a canonical coordinate system $\{\Delta_i\}$
 on the space ${\cal A}_{SL_m, \widehat S}$, parameterized by the set 
\begin{equation} \label{11.21.02.1}
{\rm I}_m^{T} := \parbox{8cm}{$\{$vertices of the $m$--triangulation of 
$T$$\}$ - $\{$vertices at the  punctures 
of $S$$\}$}.
\end{equation}
Further, we define a canonical coordinate system $\{X_j\}$ 
on  ${\cal X}_{PGL_m, \widehat S}$ 
parameterised by the set
$$
{\rm J}_m^{T}:= {\rm I}_m^{T} - \mbox{$\{$ the 
vertices at the boundary 
of $S$ $\}$}.
$$ 
\begin{figure}[ht]
\centerline{\epsfbox{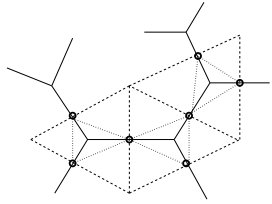}}
\caption{The set ${\rm J}_2^{T}$ is naturally identified with 
the set of all internal edges of $T$. }
\label{fg42}
\end{figure}
\begin{definition} \label{11.25.02.1}
A skew symmetric 
$\Z$--valued function 
$\varepsilon_{p q}$ on the set of vertices of the $m$-triangulation of $T$ 
is given by
\begin{equation} \label{11.25.02.1we}
\varepsilon_{p q}:= \#\{\mbox{oriented  edges 
from $p$ to $q$}\} - 
\#\{\mbox{oriented  edges from  $q$ to $p$}\}. 
\end{equation}
\end{definition} 
The Weil-Petersson form and the Poisson  structure 
  in these coordinates are given by 
\begin{equation} \label{5.17.04.11}
\Omega_{SL_m, \widehat  S} = \sum_{i_1, i_2} \varepsilon_{i_1 i_2}d\log \Delta_{i_1} \wedge d\log \Delta_{i_2}, \qquad 
\{X_{j_1}, X_{j_2}\}_{PGL_m,\widehat  S} = \sum_{j_1 j_2} \varepsilon_{j_1, j_2}
X_{j_1} X_{j_2}.
\end{equation}
The $K_2$-avatar of the Weil-Petersson form is given by 
$$
{ W}_{SL_m,\widehat S} = \sum_{i_1, i_2} \varepsilon_{i_1 i_2}\{\Delta_{i_1}, \Delta_{i_2}\}. 
$$

\begin{figure}[ht]
\centerline{\epsfbox{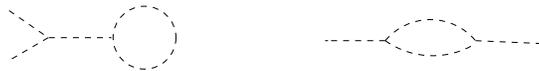}}
\caption{Special graphs, called a virus and an eye.}
\label{fgo-120}
\end{figure}

 We call an ideal triangulation $T$ of $\widehat S$ is {\it special} 
if its dual graph  contains, as a part, one of the graphs shown on Figure \ref{fgo-120}. 
A marked surface $\widehat S$ is {\it
  special} if it admits a special triangulation, and {\it regular} otherwise. 

\begin{lemma} \label{5.16.04.1}
A marked hyperbolic surfaces $\widehat S$ is special if and only if it
 contains at least two holes, and in addition one of the holes has no marked
 points on the boundary. 
\end{lemma}

The sets 
${\rm I}_m^{T}$ and ${\rm J}_m^{T}$ can be defined for any, 
not necessarily finite,
triangulation $T$. An especially interesting example 
is provided by the Farey triangulation. It is obviously regular. It has no boundary, 
so the ${\rm I}$ and ${\rm J}$ sets coincide, and are denoted 
by ${\rm I}_m$.  Since  
the function $\varepsilon_{pq}$ for the set ${\rm I}_m$ 
has finite support for every given $p$, 
all structures of the cluster ensemble but the $W$-class make
sense. 

\begin{theorem} \label{10.29.03.100}
a) The universal 
pair of positive spaces \newline $({\cal X}_{PGL_m}, {\cal A}_{SL_m})$ 
is a part of the cluster ensemble for the function $\varepsilon_{pq}$ 
on the set ${\rm I}_m$. 

b) If  $\widehat S$ 
is regular, the pair of positive  spaces 
$({\cal X}_{PGL_m, \widehat S}, {\cal A}_{SL_m, \widehat S})$ is a part of the cluster 
ensemble for the function $\varepsilon_{pq}$ 
from Definition \ref{11.25.02.1}.
\end{theorem}

We assign to every triangulation $T$ of $\widehat S$ a coordinate system on the ${\cal X}$- and ${\cal A}$-space   
For every $T$ the function 
$\varepsilon_{pq}$ from Definition \ref{11.25.02.1} describes 
the Poisson structure on ${\cal X}_{PGL_m, \widehat S}$ and the $W$-class 
on ${\cal X}_{SL_m, \widehat S}$ in  
the corresponding coordinate system 
via formulas (\ref{5.17.04.11}). 
However the transition functions between the coordinate systems assigned to 
special triangulations 
are slightly different than the ones prescribed by the cluster ensemble data. 
We show that the pair of positive spaces $({\cal X}_{PGL_m, \widehat S}, {\cal A}_{SL_m,
\widehat S})$ is related to a more general  
 orbi-cluster ensemble structure.

\begin{theorem} \label{10.29.03.100a}
If $\widehat S$ is special, the pair of positive spaces 
\newline $({\cal X}_{PGL_m, \widehat S}, {\cal A}_{SL_m,
  \widehat S})$ is a part of the orbi-cluster 
ensemble for the function $\varepsilon_{pq}$ 
from Definition \ref{11.25.02.1}. 
\end{theorem}
The precise meaning of this theorem is discussed in Section 10. 
In a sequel to this paper we  will prove a similar result 
for an arbitrary $G$, quantize the higher Teichm{\"u}ller spaces 
and as a 
result construct an infinite dimensional unitary 
projective representation of the group  $\Gamma_{S}$. 

\vskip 3mm

{\bf 16. The cluster mapping class group}. 
As we show in \cite{FG2}, every 
cluster ensemble comes with a group of symmetries, called the {\it cluster mapping group}. 
Here is how it is defined (a bit more general definition see in loc. cit.). 
Let ${\cal X}$ be a rational Poisson variety equipped with a rational 
coordinate system $\{X_i\}$ parametrised by a set $I$, $i\in I$, so that the 
Poisson structure is given by 
$$
\{X_i, X_j\}=\varepsilon_{ij} X_i X_j\quad 
\mbox{where $\varepsilon_{ij} \in \Z$}.
$$   
Then for every $k \in I$ there 
is a rational  automorphism $\mu_k$ of  ${\cal X}$, called a {\it mutation in the direction $k$},
 given by 
\begin{equation} \label{5.11.03.1xa}
\mu_k^*X_{i} = \left\{\begin{array}{lll} X_k^{-1}& \mbox{ if } & i=k \\
    X_i(1+X_k^{-\sgn \varepsilon_{ik}})^{-\varepsilon_{ik}} & \mbox{ if } 
&  i\neq k. 
\end{array} \right.
\end{equation}
The 
Poisson structure in the mutated coordinates $X_i':= \mu_k^*(X_i)$ is again quadratic, and given by 
$\{X'_i, X'_j\}=\varepsilon'_{ij} X'_i X'_j$, where $\varepsilon'_{ij} \in \Z$. The Poisson tensor 
$\varepsilon'_{ij}$ is calculated from the one $\varepsilon_{ij}$ by an explicit formula. 
The mutated coordinate system is parametrized by the same set $I$, so 
we can apply the above procedure again and again.  Further, let a 
$\sigma: I \to I$ be a bijection. It gives rise to an automorphism $\sigma: {\cal X} \to {\cal X}$ 
called a {\it symmetry}, so that  $X_i':= \sigma^*X_i := X_{\sigma(i)}$, 
 with a new Poisson tensor $\varepsilon'_{ij}:= 
\varepsilon_{\sigma(i)\sigma(j)}$. 
We define an {\it ${\cal X}$-cluster transformation} as a composition of mutations and symmetries. 
It is a birational Poisson automorphism $\alpha$ of ${\cal X}$  such that 
 $\{\alpha^*X_i, \alpha^*X_j\} = \varepsilon^\alpha_{ij}X_iX_j$. 
The cluster transformations $\alpha$  preserving the $\varepsilon_{ij}$-tensor, i.e.  
$\varepsilon^\alpha_{ij} = \varepsilon_{ij}$, form a group. The {\it cluster mapping class group} is 
the image of this group in the group of birational automorphisms of ${\cal X}$. 

The moduli space ${\cal X}_{G,S}$ admits such a coordinate system 
(with slight modifications when $G$ is not simply-laced). For $G=PGL_m$ the corresponding data 
was described in Section 1.15. 
We prove that in our situation the cluster mapping class group, denoted $\Gamma_{G, \widehat S}$, 
contains the classical mapping class group $\Gamma_{S}$ of $S$. As a result we show that 
the classical mapping class group acts by cluster transformation of our moduli spaces, and this it 
is given in a very explicit form.  If $G=PGL_2$, 
the two mapping class groups coincide. It would be extremely interesting to determine the group 
$\Gamma_{G, \widehat S}$ in general.

\vskip 3mm
{\bf Examples}.   1. 
Let $G$ be a group of type $G_2$, 
and let $\widehat S$ be a disc with three marked points at the boundary. 
Then the classical mapping class group is $\Z/3\Z$. 
The cluster mapping class group is an infinite quotient of the braid group of type $G_2$, 
which conjecturally coincides with the latter (\cite{FG5}). 

2. Let $S$ be an annulus with $2k$ marked points on one component of the boundary, $k \geq 1$. 
The classical mapping class group is $\Z/2k\Z$. 
Let $G$ be an arbitrary split semi-simple group over $\Q$. 
 The braid group ${\Bbb B}_G$ of type $G$ acts by cluster transformations of the corresponding 
moduli space. The center acts trivially. 
It is likely that the cluster mapping class group in this case is the quotient of 
the braid group ${\Bbb B}_G$ by its center. (\cite{FG6}). 
\vskip 3mm
Here is a concrete problem. 
Recall (see Section 1.9) that, given a loop $\gamma$ on $S$, there is an action 
of the group $H(\R_{>0})$ on the Teichm\"uller space ${\cal L}_{G, S}^+$. 
It always contains a subgroup $\Z$ generated by the Dehn twist along $\gamma$. 
If $G = PGL_2$, then $H(\R_{>0})$ is isomorphic to $\R$, so the quotient by the subgroup generated by 
the Dehn twist is compact. It is no longer compact in all other cases.  
Is there  an  abelian subgroup of $\Gamma_{G, S}$, which is of rank ${\rm dim}H$, 
lies in $H(\R_{>0})$,  and is cocompact there?

\vskip 3mm
It follows from \cite{FG2} 
that the cluster mapping class group $\Gamma_{G, \widehat S}$ acts by birational automorphisms 
of the non-commutative 
$q$-deformation of the moduli space 
${\cal X}_{G, \widehat S}$. Moreover, this leads to a construction of an 
 infinite-dimensional unitary representation of 
$\Gamma_{G, \widehat S}$, and we conjecture that it provides an example of 
infinite-dimensional modular functor.  This and other indications strongly suggest that the 
 ``true'' moduli space related to the pair $(G, \widehat S)$ 
should be the quotient ${\cal X}^+_{G, \widehat S}/\Gamma_{G, \widehat S}$, or, better, just 
the pair $({\cal X}^+_{G, \widehat S}, \Gamma_{G, \widehat S})$.

\vskip 3mm 
{\bf 17. Duality conjectures}. We suggest that there exists  a remarkable duality interchanging the 
${\cal A}$- and ${\cal X}$-moduli spaces on $\widehat S$, 
which changes the group $G$ to its Langlands dual $^L G$. So the dual pairs are
$$
{\cal A}_{G, \widehat S}\quad \mbox{and} \quad  {\cal X}_{^L G, \widehat S}. 
$$
Observe that the Langlands dual to an adjoint simple group 
is simply-con\-nected, and vice versa.  Here is one of the manifestations of this duality. 

Let ${\cal X}$ be a positive space. 
Let ${\Bbb L}_+[{\cal X}]$ be the subset of all rational functions $F$ on ${\cal X}$ 
such that for every positive coordinate system $\psi_{\alpha}$ the function $F$ 
written in the corresponding coordinates
is a Laurent polynomial with positive integral coefficients. In other words 
$\psi_{\alpha}^*F$ is a regular positive function on the torus 
$H_{\alpha}$. Then    ${\Bbb L}_+[{\cal X}]$ is evidently a  semiring. 
Let ${\bf E}({\cal X})$ be its set of {\it extremal elements}, i.e. elements 
$F \in  {\Bbb L}_+[{\cal X}]$ which can not be decomposed into a sum 
$F = F_1 + F_2$ of two non-zero elements  $F_i \in {\Bbb L}_+[{\cal X}]$. 

Observe that for any semifield  $K$ the set ${\cal X}(K)$ is determined by a single 
coordinate system of the positive atlas on ${\cal X}$, which provides an isomorphism 
${\cal X}(K) \stackrel{\sim}{\lra} K^{\rm dim{\cal X}}$. On the other hand 
the semiring ${\Bbb L}_+[{\cal X}]$ shrinks 
when we get more coordinate systems in our positive atlas. In particular it may be empty. 

Let $L$ be a set. Denote by $\Z_+\{L\}$ the abelian 
semigroup generated by $L$. Its elements are finite expressions  
$\sum_i n_i \{l_i\}$ where $l_i \geq 0$, where 
$\{l_i\}$ is the generator corresponding to $l_i \in L$. 

The group $W^n$ acts naturally on both moduli spaces: this is obvious for the ${\cal X}$-space, but 
rather surprising for the ${\cal A}$-space -- see Section 12.6 where the latter action is discussed.

\begin{conjecture} \label{10.10.03.10q}
{\it Let $G$ be a connected, simply-connected, split semi-simple algebraic group. Let $S$ be a surface with boundary. 
Then there exist  canonical isomorphisms of sets, 
which are equivariant with respect to the action of the (cluster) mapping class group, and 
intertwine the  two actions of the group $W^n$:

\begin{equation} \label{10/9/03/3w}
{\cal A}_{G,   S}(\Z^t) \stackrel{\sim}{\lra} {\bf E}({\cal X}_{^L G,   S}), 
\qquad {\cal X}_{^L G,   S}(\Z^t) \stackrel{\sim}{\lra} {\bf E}({\cal A}_{G,   S}), 
\end{equation}

Moreover they extend to isomorphisms of the semirings
\begin{equation} \label{10/9/03/3ww}
\Z_+\{{\cal A}_{G,   S}(\Z^t)\} \stackrel{\sim}{\lra} {\Bbb L}_+({\cal X}_{^L G,   S}), 
\qquad \Z_+\{{\cal X}_{^L G,   S}(\Z^t)\} \stackrel{\sim}{\lra} {\Bbb L}_+({\cal A}_{G,   S}). 
\end{equation}}
\end{conjecture} 

\vskip 3mm
A considerable part of Conjecture \ref{10.10.03.10q} for $G= SL_2$ is proved in Section 12: we define 
the above canonical maps, and prove most of their properties. However  
we can not show that the functions assigned to laminations are extremal elements.

As explained in Section 4 of \cite{FG2}, there are several other versions 
of the duality conjectures. In one of them   
$\Z^t$ is replaced by the real tropical field $\R^t$, while ${\Bbb L}_+$ is 
replaced by functions on $\R_{>0}$ points of the dual space. 
So in this form we are looking for a pairing between the $\R^t$-points 
of one of the spaces with the  $\R_{>0}$-points of the other. 

There are quantum versions of these conjectures where ${\Bbb L}_+({\cal X}_{G,   S})$ 
is replaced by the corresponding semiring on the non-commutative 
q-deformation  ${\cal X}^q_{G,   S}$, or by the non-commutative $\ast$-algebra of functions on the Teichmuller space ${\cal X}_{G,   S}(\R_{>0})$. In the first case it is paired with 
${\cal A}_{G,   S}(\Z^t)$, while in the second with ${\cal A}_{G,   S}(\R^t)$. 

The Duality Conjecture  for surfaces without boundary is discussed in Section 13.4. 

Summarizing, the spirit of these conjectures is this: we are looking for a description of 
one of the spaces, perhaps made non-commutative (in several different ways), via tropicalization 
of the dual space. The tropical points are points of the Thurston type 
boundary. It would be very interesting to 
have a duality between the two moduli spaces, which degenerates to our duality when one 
of the spaces is replaced by its tropicalization.  This reminds us of the mirror duality.  

\vskip 3mm 
{\bf 18. Completions of Teichm\"uller spaces and canonical bases.} In Section 13 
we define a (partial for $G \not = PSL_2$) completion of the higher Teichm\"uller space. It is 
equipped with an action of the mapping class group of $S$. Its strata are parametrised by 
{\it simple laminations}, that is the isotopy classes of collections of simple disjoint non-isotopic loops on $S$, 
non-isotopic to any of the boundary components. 
The stratum corresponding to a simple lamination $l$ is given by an appropriate 
 Teichm\"uller space 
for the surface $S-l$. In the classical case 
when $G=PGL_2(\R)$ 
we get a new construction of the Weil-Petersson completion of 
the Teichm\"uller space which goes back to L. 
Bers \cite{Bers2} - see \cite{Wo} and references therein. 
Its quotient under the mapping class group is the 
Knudsen-Deligne-Mumford moduli space $\overline {\cal M}_{g,n}$. 
For example when $S$ is a punctured torus, the 
Teichm\"uller space is identified with the hyperbolic disc, and 
the completion is obtained by adding a countable set cusps to its boundary, identified with  
${\Bbb P}^1(\Q)$. 

We conjecture the existence of the canonical map for surfaces without boundary. 
Its restriction to the boundary component corresponding to a simple lamination $l$ on $S$  should be given by 
the canonical map from Conjecture \ref{10.10.03.10q} for the surface $S-l$. 
This requirement should determine uniquely the canonical map for surfaces without boundary.

\vskip 3mm 
 {\bf 19. Coda}. We suggest that  
investigation of the pair of moduli spaces 
$( {\cal A}_{G,\widehat S}, {\cal X}_{G',\widehat S})$ can be viewed as a 
marriage of representation theory and the theory of surfaces. 
Indeed,  
when the surface $\widehat S$ is simple, i.e. is a disc 
with marked points on the boundary or an annulus, we recover 
many aspects of the representation theory for the group $G$,   
quite often  
from a new point of view: positivity 
in $G(\R)$, Steinberg varieties, invariants in the tensor products of finite-dimensional
representations,  quantum groups, canonical bases, 
$W$-algebras, etc. 
On the other hand, when $G = SL_2$ but $\widehat S$ is general we 
get the theory of Riemann surfaces. 

\noindent \qquad Moreover, we expect that the\! theory of quasifuchsian and Kleinian groups 
also admits a generalization, where $PGL_2(\C)$ is replaced by an arbitrary 
semi-simple complex Lie group -- see Section 7.12.  The higher Teichm\"uller spaces should 
describe the corresponding deformation spaces just the same way as 
in the Bers double uniformization theorem  
and its generalization to Kleinian groups.

\begin{figure}[ht]
\centerline{\epsfbox{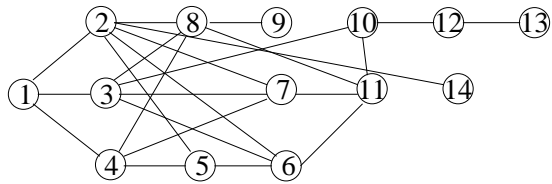}}
\caption{Leitfaden}
\label{fgoleit}
\end{figure}
 
\vskip 3mm 
\noindent \qquad {\bf Acknowledgments}. 
N. Reshetikhin informed us about his 
work in prog\-ress with R. Kashaev,  containing a result similar 
 to our decomposition (\ref{9.28.03.11yx}). 
After the first version of this work was posted in the ArXiv, we
learned about a very interesting preprint of F. Labourie \cite{Lab} 
who, for surfaces $S$ without holes, 
defined, using a dynamical systems approach, 
 a class of representations $\pi_1(S) \to PSL_m(\R)$, and proved that they are discrete, hyperbolic, and 
coincide with the ones defined by Hitchin.

 This work was developed  when 
V.F. visited Brown University. He would like to thank Brown University 
for hospitality and support. He wishes to thank grants 
 RFBR 01-01-00549, 0302-17554 and NSh 1999.2003.2. 
A.G. was 
supported by the  NSF grants DMS-0099390 and DMS-0400449. 
The paper was essentially written when A.G. enjoyed the hospitality of IHES (Bures sur Yvette)  
and MPI (Bonn)
during the Summer and Fall of 2003. 
He is grateful to NSF, IHES and MPI for the support. 
We are grateful to J. Bernstein,  M. Gromov, 
D. Kazhdan, M. Kontsevich  
for many stimulating discussions, to F. Labourie for 
the explanation of his approach in \cite{Lab}, 
and to M. King and D. Dumas for proofreading a part of the paper.  

We are extremely grateful to the referee,  who read very carefully the paper, corrected several errors,  
and made a lot of  remarks and suggestions, incorporated into the text, to improve the exposition. 

\section{The moduli spaces ${\cal A}_{G, \widehat S}$ and 
${\cal X}_{G, \widehat S}$}
\label{S3}


Let us replace the holes on $S$ by punctures, getting a surface $S'$. The local systems on $S$ and $S'$ are the same. Let us equip the surface $S'$
with a complex structure. Then we get a complex 
algebraic curve ${\Bbb S}$. 
The moduli space of local systems on ${\Bbb S}$
has two different algebraic structures: the monodromic 
and the De Rham structures. Nevertheless the corresponding complex 
analytic spaces 
are isomorphic. The monodromic structure 
depends only on the fundamental group $\pi_1({\Bbb S}(\C))$. Therefore it is determined by the topological surface $S'$ underlying the curve ${\Bbb S}$. 
The corresponding moduli stack is denoted ${\cal L}_{G, {S}}$. 
To define the De Rham structure recall 
the Riemann-Hilbert correspondence 
([Del])
between the isomorphism classes of 
$G$--local systems on $S' = {\Bbb S}(\C)$ and algebraic $G$--connections 
with regular singularities on the algebraic curve ${\Bbb S}$. 
The moduli space 
${\cal L}_{G, {\Bbb S}}$ of regular $G$-connections on 
 ${\Bbb S}$ is an algebraic stack over $\Q$.  
Its algebraic structure depends on the algebraic structure 
of the curve ${\Bbb S}$. It is the De Rham algebraic 
structure.  
The Riemann-Hilbert correspondence provides an isomorphism 
of the complex 
analytic spaces ${\cal L}_{G, S}(\C)$ and ${\cal L}_{ G, {\Bbb S}}(\C)$.

\vskip 3mm

{\bf 1. The moduli spaces ${\cal X}_{G, \widehat S}$}. 
The flag variety  ${\cal B}$ 
parametrises all Borel subgroups in $G$. Choosing 
a Borel subgroup $B$ we identify it with $G/B$. 
Let ${\cal L}$ be a $G$-local system on $S$. Recall 
the associated {\it flag bundle} $
{\cal L}_{\cal B}:= 
{\cal L} \times_G {\cal B}$.  We call the  elements of 
a fiber of ${\cal L}_{\cal B}$ 
over  $x$ {\it flags over $x$ in ${\cal L}$}. 
A reduction 
of a $G$--local system ${\cal L}$ on a circle to $B$ is just the same as 
a choice of a flag in a fiber  of ${\cal L}$ invariant under the monodromy.

\vskip 3mm
{\bf Examples}. a) A  flag in a 
vector space $V$ is given by a filtration
\begin{equation} \label{11.4.02.2}
0 = V_0 \subset V_1 
\subset V_2 \subset ... \subset V_m =V, \quad {\rm dim}V_i =i.
\end{equation} 
The flag variety for $G= GL_m$ parametrizes all flags in $V$.

b) Let $G = GL_m$. Recall the bijective correspondence between the 
$GL_m$--local systems and local systems of 
$m$-dimensional  vector spaces. Namely, let 
 $V$ be a standard $m$-dimensional 
representation of $GL_m$.  A $GL_m$--local system ${\cal L}$  
provides a 
 local system of vector spaces ${\cal L}':= {\cal L} \times_{\rm GL_m} V$. 
Conversely, given an $m$-dimensional 
local system  ${\cal L}'$ we recover 
the corresponding $GL_m$--local system by taking 
all bases in fibers of ${\cal L}'$. 
The flag bundle ${\cal L}_{\cal B}$ consists of all flags in the fibers of the  
 local system  ${\cal L}'$. 
If $G = SL_m$,   ${\cal L}'$ is a local system with 
an invariant volume form. 
\vskip 3mm

\begin{definition} \label{4.7.03.2d}
Let $G$ be an arbitrary split reductive group. 
A {\it framed  $G$-local system on $\widehat S$} is a pair 
$({\cal L}, \beta)$ where ${\cal L}$ is a $G$-local system 
on $S$ and $\beta$ is a flat section of the restriction of 
${\cal L}_{\cal B}$ to the punctured boundary $\partial \widehat S$. 
${\cal X}_{G, \widehat S}$ is the moduli space of 
framed $G$--local systems on $\widehat S$. 
\end{definition}

One can rephrase this definition as follows:

i)  If $C_i$ is a boundary component 
without marked points, we choose a flag $F_x$  
over a certain point $x \in C_i$   
invariant under the monodromy around $C_i$. Equivalently, we choose 
a reduction to the subgroup $B$ of the restriction 
of the $G$-local system ${\cal L}$ 
to the boundary component $C_i$. 
The parallel transport of the flag $F_x$ provides the restriction of the section 
$\beta$ to  $C_i$, and vice versa.

ii) If there is a non empty 
set $\{x_1, ..., x_p\}$ of boundary points on the component 
$C_i$, the complement 
$C_i - \{x_1, ..., x_p\}$ is a union of $p$ arcs. Let us pick 
a point $y_j$ inside of 
each of these arcs, and choose a flag ${F}_{j}$ over  $y_j$. 
These flags may not be invariant under the monodromy.
The flags ${F}_{j}$ determine the restriction of the  flat section 
$\beta$ to the corresponding arc, and vice versa. 

\vskip 3mm
{\it The Cartan group of $G$}. A split reductive algebraic group $G$ determines 
the Cartan group $H$. Namely, let $B$ be a Borel subgroup in $G$. 
Then $H := B/U$ where $U$ is the unipotent radical of $B$. If $B'$ is 
another Borel subgroup then $B$ and $B'$ are conjugate in $G$, and this conjugation induces a canonical isomorphism 
$B/U \lra B'/U'$. So the quotient $ B/U$ is canonically isomorphic 
to the Cartan subgroup $H$ in $G$. Let 
\begin{equation} \label{9.26.03.1}
i_B: H \to B/U 
\end{equation}
be the canonical isomorphism. The Cartan group $H$ 
is not considered 
as a subgroup of $G$.   

If $S = \widehat S$, i.e. there are no marked points, the canonical projection $B \to H= B/U$ provides a natural projection
\begin{equation} \label{9.26.03.2}
\kappa: {\cal X}_{G, S} \lra H^{\{\mbox{punctures of $S'$}\}}.
\end{equation}

\vskip 3mm
{\it The moduli space ${\cal X}_{G, \widehat S}$ as a quotient stack}. 
Let $F$ be an arbitrary field. Choose a point $x \in S$. Let $g$ be the genus
of the 
surface $\overline S$. 
Let us choose generators $\alpha_1, ..., \alpha_g, \beta_1, ..., \beta_g, 
\gamma_1, ..., \gamma_n \in \pi_1(S, x)$ such that 
\begin{equation} \label{10.25.03.10}
\prod_{j=1}^g[\alpha_j, \beta_j]\prod_{i=1}^n\gamma_i = 1
\end{equation} 
and where $\gamma_i$ belongs to the conjugacy class in $\pi_1(S, x)$ determined
by the $i$-th hole in $S$. It is well known that such a system of generators
exists, and that (\ref{10.25.03.10}) is the only relation between 
them in $\pi_1(S, x)$. 

\begin{definition} \label{10.25.03.12}
A {\rm framed representation} $\pi_1(S, x)\to G(F)$ is the following data:
\begin{equation} \label{10.25.03.11}
\begin{array}{c}\Bigl(X_{\alpha_1}, ..., X_{\alpha_g}, X_{\beta_1}, ..., X_{\beta_g}, X_{\gamma_1},
..., X_{\gamma_n}, \{B_j\}\Bigr),\\
X_{?}\in G(F), \quad j \in
\pi_0(\partial \widehat S), \quad B_j \in {\cal B}(F)
\end{array}
\end{equation}
where $X_{?}$ satisfy (\ref{10.25.03.10}), and 
if $j$ corresponds to a connected component $C_j$ of the boundary of $S$ 
without 
marked points on it, then $X_{\gamma_j} \in B_j(F)$. 
\end{definition}
Obviously the 
elements of the set (\ref{10.25.03.11}) are the  $F$-points of an affine
algebraic variety denoted by $\widetilde {\cal X}_{G, \widehat S}$. The
algebraic group 
$G$ acts on this variety, and, by definition,  the stack ${\cal X}_{G, \widehat
  S}$ 
 is the quotient stack 
$\widetilde {\cal X}_{G, \widehat S}/G$.  A version of definition of the stack 
${\cal X}_{G, S}$ see in Section 12.5. 

\vskip 3mm
{\bf Examples}. 1. Let $\widehat D_n$ be a disc with $n$ marked points on the boundary. 
Then $\widetilde {\cal X}_{G, \widehat D_n} = {\cal B}^n$, and 
${\cal X}_{G, \widehat D_n} = G\backslash {\cal B}^n$ is the space of configurations of $n$
flags. 

2. Let $S$ be an annulus. Then $\widetilde {\cal X}_{G, S} = 
\{(g, B_1, B_2)\}$ 
where $B_i$ are Borel subgroups in $G$ and $g \in B_1, g\in B_2$. 
The quotient ${\cal X}_{G, S}:= \widetilde {\cal X}_{G, S}/G$ is the Steinberg
variety, 
see \cite{CG}.

\vskip 3mm
{\bf 2. Pinnings}. 
Let us introduce some notations 
following  \cite{L1}. 
Let $G$ be a split reductive simply-connected connected algebraic group over $\Q$. 
Let $H$ be a split maximal torus of $G$ and $B^+, B^-$ a pair 
of opposed Borel subgroups containing $H$, with unipotent radicals 
$U^+, U^-$. Let $U_i^+$ ($i \in I$) be the simple root subgroups of $U^+$ and 
let $U_i^-$ be the corresponding  root subgroups of $U^-$. 
Here $I$ is a finite set indexing the simple roots. 
Denote by $\chi'_i:H \lra {\Bbb G}_m$ the simple root corresponding to $U_i^+$. 
Let $\chi_i:{\Bbb G}_m \lra H$ be the simple coroot corresponding to $\chi'_i$. 
Let ${\Bbb G}_a:= {\rm Spec}\Q[t]$ be the additive one dimensional algebraic group. 
We assume that  for each $i \in I$ we are given isomorphisms $x_i: {\Bbb G}_a \lra U_i^+$ 
and $y_i: {\Bbb G}_a \lra U_i^-$ such that the maps
$$
\left (\begin{matrix}1& a\\ 0&1\end{matrix}\right ) \lms x_i(a), \quad 
\left (\begin{matrix}1& 0\\ b&1\end{matrix}\right )\lms y_i(b), \quad 
\left (\begin{matrix}t& 0\\ 0&t^{-1}\end{matrix}\right )  \lms  \chi_i(t) 
$$
 provide a homomorphism 
$\varphi_i: SL_2 \lra G$. 
The datum $(H, B^+, B^-, x_i, y_i; i \in I)$ is called a {\it pinning} 
for $G$. 
The group $G$ acts by conjugation on pinnings, making it into 
a principal homogeneous space for the group 
$G/{\rm Center}(G)$. 
In particular any two pinnings for $G$ are conjugate in $G$.

There is a unique involutive antiautomorphism $\Psi: G \lra G$ such that 
for all $i \in I$ and $t \in H$ one has 
$\Psi(x_i(a)) =  y_i(a)$, $\Psi(y_i(a)) =  x_i(a)$ and $\Psi(t) =  t$. 
\vskip 3mm
{\bf Example}. If $G = SL_m$ then $B^+, B^-$ are the subgroups of upper and 
lower triangular matrices, so $H$ is the subgroup of diagonal matrices. 
$\Psi$ is the transposition. 
If $G = PSL(V)$ where $V$ is a vector space 
then a choice of pinning is the same thing as  a choice of a projective 
basis in $V$.  
\vskip 3mm

Let $X^*(T)$ and $X_*(T)$ be the groups of characters $T \to {\Bbb G}_m$ 
and cocharacters ${\Bbb G}_m\to T $
for a split torus $T$. There is a pairing $$
<\cdot, \cdot>: X^*(T) \times X_*(T) \lra \Z, \quad 
\chi'\circ \chi : {\Bbb G}_m \lra {\Bbb G}_m, t \lms t^{<\chi', \chi>}.
$$

Let $A = (a_{ij})$ 
be the Cartan matrix given by $a_{ij}:= <\chi'_j,  \chi_i>$.   
Recall the weight lattice $P\subset X^*(H)$. It is 
the subgroup of all $\gamma \in X^*(H)$ such that 
$<\gamma, \chi_i> \in \Z$ for all $i \in [1, ..., r]$. 
Let $\{\omega_1, ..., \omega_r\}$ be the ${\Z}$--basis of fundamental weights 
given by $<\omega_j, \chi_i> = \delta_{ij}$.  
\vskip 3mm
{\bf 3. The element $s_G$}. A choice of pinning for $G$ provides an inclusion of the Cartan 
group $H$ into $G$. Abusing notation we will denote the 
obtained  subgroup by $H$ and call it the Cartan subgroup corresponding 
to the pair $(B^-, B^+)$ of opposite Borel subgroups.  Then the element
$$
\overline s_i:= y_i(1) x_i(-1) y_i(1) = \varphi_i
\left( \begin{matrix}0& -1\\ 1&0\end{matrix}\right)
$$
lifts to ${\rm Norm}_G(H)$ the generator $s_i:= s_{\chi_i} \in W$. 
It is well known that the elements $\overline s_i$, $i \in I$, 
satisfy the braid relations. Therefore we can associate to each $w \in W$ 
its standard representative $\overline w \in {\rm Norm}_G(H)$ in such a way 
that for any reduced decomposition $w = s_{i_1} ... s_{i_k}$ one has 
$
\overline w = \overline s_{i_1} ... \overline s_{i_k}
$. Let $w_0$ be the maximal length element of $W$. 
We set $s_G:= \overline w_0^2$.    
Then evidently $s_G^2 = e$. 
\vskip 3mm
{\bf Example}. If $G=  SL_m$ then $s_{SL_m} = (-1)^{m-1} e$.  
\vskip 3mm
Let $R$ be the system of roots for the group $G$, and  $R = R_+ \cup R_-$ 
its decomposition into a union of the sets of  positive and negative
roots. Then for any element $w \in W$ its length $l(w)$ equals to  
$|R_- \cap w(R_+)|$. For any root $\alpha = \sum n_i \alpha_i$, where $\{\alpha_i\}$ is the set of positive simple roots,  we set 
$\chi_{\alpha}(t):= \prod \chi_{i}(t)^{n_i}$. 
\begin{lemma} \label{10.05.03.q1}
 Let $w \in W$. Then 
$$
\overline w \overline {w^{-1}} = \prod_{\alpha \in R_- \cap w(R_+)} \chi_{\alpha}(-1).
$$
\end{lemma}

{\bf Proof}. By induction on the length $l(w)$. If $l(w)=1$  this is
clear. Assume 
$w = w_i w'$ where $l(w') = l(w) -1$. Then by the induction assumption 
$$
\overline w \overline {w^{-1}} = \overline w_i \overline {w'} \overline
{{w'}^{-1}}
\overline {w^{-1}_i} = \left[w_i\Bigl( \prod_{\alpha \in R_- \cap w'(R_+)}
\chi_{\alpha}(-1)\Bigr) 
w_i^{-1}\right]\overline w_i\overline {w_i^{-1}} = $$ $$ = \prod_{\alpha \in R_- \cap w'(R_+)}
\chi_{w_i(\alpha)}(-1)\chi_{\alpha_i}(-1).
$$
Since $w_i(R_- \cap w'(R_+)) \cup -\alpha_i = R_- \cap w(R_+)$, the lemma is proved. 
\vskip 3mm

\begin{corollary} \label{10.05.03.1}
$s_G$ is a central element of $G$ of order $2$.  
\end{corollary}

{\bf Proof}. Let $2\rho$ be the sum of all positive roots. Then 
(see Proposition 29, Chapter VI, paragraph 1 in \cite{Bo})  we have 
$\rho = \omega_1 + ... + \omega_r$. Denote by
$\chi'_{2\rho}$ 
the character of $H$ corresponding to $2\rho$. Then for any coroot $\chi$ 
one has $<\chi'_{2\rho}, \chi> \in 2\Z$. Since $w^{-1}_0 = w_0$, this
combined with Lemma \ref{10.05.03.q1} implies that $s_G$ commutes 
with the elements $x_i(t)$, and thus is a central element. The corollary is
proved. 

There is a well defined up to conjugation embedding $p: SL_2
\hookrightarrow G$,  called the principal $SL_2$ subgroup of $G$. 
Its Lie algebra is described as follows. 
If we choose a pinning in $G$, 
the image of the element $\left( \begin{matrix}0& 1\\ 0&0\end{matrix}\right )$ (resp. $\left( \begin{matrix}0& 0\\ 1&0\end{matrix}\right )$) of $sl_2$ under the differential of the map $p$ 
is $\sum_{i \in I}e_{i}$ (resp. $\sum_{i \in I}f_{i}$), where $e_i$ (resp. $f_i$) is the 
generator of ${\rm Lie}(U)$ (resp. ${\rm Lie}(U^-)$) 
corresponding to a simple root parametrised by $i \in I$. 
\begin{lemma} \label{12.7.03.1}
$s_G= p(s_{SL_2})$  for the principal $SL_2$ subgroup of $G$.
\end{lemma}

{\bf Proof}. An easy exercise. 
\vskip 3mm
{\bf 4. The moduli spaces ${\cal A}_{G, \widehat S}$}. 
Let  us choose a maximal unipotent subgroup $U$  
in $G$.  The normalizer of $U$ in $G$ is 
a Borel subgroup $B$.  Conversely, $U$ is the commutant of  $B$. 
The Cartan group $H$ acts from the right 
on the principal affine variety ${\cal A}:= G/U$: if $A = gU \in {\cal A}$ 
we set 
\begin{equation} \label{9.26.03.3}
A\cdot h =  gU \cdot h:=  gi_B(h) U
\end{equation}
Since $H$ is commutative it,  of course, can be considered as a left action. 
 This action is free. 
   The quotient ${\cal A}/H$ is identified with 
the flag variety ${\cal B}= G/B$. 
The canonical projection ${\cal A}  \to {\cal B}$ is a principal $H$-bundle.   
 More generally, let $X$ be a   
right principal homogeneous space for $G$. 
Then $ X/U$ is an ${\rm Aut}X$--homogeneous space. 
 Recall the {\it principal affine bundle} 
${\cal L}_{\cal A}:=  {\cal L}/U$ 
associated to a $G$--local system ${\cal L}$ on $S$.

Let $T'_S$ be 
the punctured tangent bundle to the surface $S$, 
i.e. the tangent bundle with the zero section
removed. Its fundamental group $\pi_1(T'S, x)$ is a central extension
of $\pi_1(S, y)$ by $\Z$, where $x \in T_yS $:
\begin{equation} \label{chzhq}
0 \lra \Z \lra \pi_1(T'S, x) \lra \pi_1(S, y) \lra 0.
\end{equation}
Let $T'_yS$ be the punctured tangent space at a 
point $y$. The inclusion $T'_yS \hra T'S$ induces an isomorphism of 
$\pi_1(T'_yS) \stackrel{\sim}{=} \Z$ with the central subgroup $\Z$ in
(\ref{chzhq}). We denote by $\sigma_S$ a generator of this central subgroup.
It is well defined up to sign.

\begin{definition} \label{10.14.03.1}
A twisted $G$-local system on $S$ is a local system on $T'S$ with  
the monodromy $s_G$ around  $\sigma_S$. 
\end{definition}
Observe that since $s_G$ 
is of order two, this does not depend on the choice of the generator $\sigma_S$. 

Let $\overline \pi_1(S, y)$ be the quotient of $\pi_1(S, y)$ by the central
subgroup $2\Z \subset \Z$, so it is a central extension
$$
0 \lra \Z/2\Z \lra \overline \pi_1(T'S, x) \lra \pi_1(S, y) \lra 0. 
$$
We denote by $\overline \sigma_S$ the order two central element in the
subgroup $\Z/2\Z$. 
The twisted local systems on $S$ are in bijective correspondence with 
the representations $\rho: \overline \pi_1(T'S, x) \to G$, considered modulo
conjugation, such that $\rho(\overline \sigma_S) = s_G$. 
\vskip 3mm
{\bf Remark}. Let us choose a complex structure on $S$. Then the uniformisation provides an
embedding $i: \pi_1(S, y) \hra PSL_2(\R)$. It is easy to see that the group 
$\overline \pi_1(T'S, x)$ is isomorphic to the preimage of the subgroup 
$i(\pi_1(S, y))$ in $SL_2(\R)$. In particular $\overline \sigma_S$ 
corresponds to the element $-e \in SL_2(\R)$. 
\vskip 3mm
The group $\overline \pi_1(T'S, x)$ is isomorphic to the direct product 
$\Z/2\Z \times \pi_1(S, y)$, although this isomorphism is by no means
canonical. The set of all isomorphisms 
$$
\overline \pi_1(T'S, x) \lra \Z/2\Z \times \pi_1(S, y)
$$ 
which send the central subgroup $\Z/2\Z $ on the left to the left factor 
on the right, and reduce to  the canonical isomorphism modulo 
the $\Z/2\Z$-subgroups, 
is a principal homogeneous space over ${\rm Hom}(\pi_1(S, y), \Z/2\Z)$. 
Therefore choosing such an isomorphism we can identify the space  of 
twisted $G$-local systems
on $S$ with the space of $G$-local systems on $S$. Observe that a choice of
such an isomorphism is the same thing as a choice of a spin structure 
on $S$.

Let ${\bf C}_i$ be a little annulus containing a boundary component $C_i$ of 
$S$ as
a component of
its boundary. Observe that $\pi_1(T'{\bf C}_i)$ is canonically identified 
with $\Z \times \Z$: the first factor is $\pi_1(T'_x{\bf C}_i)$, and the
second is generated by the tangent vectors to $C_i$. 

Let $x_1, ..., x_p$  be the 
 marked points on the component $C_i$. Let us present the annulus ${\bf C}_i$ as a product 
$C_i \times [0,1]$, and let ${\bf C}'_i:= ({ C}_i - \{x_1, ..., x_p\})\times
[0,1]$. 
Then ${\bf C}'_i$ has $p$ connected components.

\begin{definition} \label{5.7.02.10xs} 
Let $G$ be a split simply-connected reductive group.  
Let ${\cal L}$ be a $G$-local system on $T'S$ representing 
a twisted local system on $S$.  
A {\em decoration} on ${\cal L}$ 
is a choice of locally constant section $\alpha$ of the restriction of the
principal affine bundle 
${\cal L}_{\cal A}$ 
to $\cup_i {\bf C}'_i$. 

A {\rm decorated}  twisted $G$--local system on a marked surface 
$\widehat S$ is  a pair $({\cal L}, \alpha)$, where $\cal L$ is 
 a twisted $G$-local system on $S$, and $\alpha$ is a decoration on ${\cal L}$. 
The space 
${\cal A}_{G, \widehat S}$ is the moduli space of  
 decorated  twisted $G$--local systems on $\widehat S$. 
\end{definition}
When $s_G =e$ this definition reduces to the one 
given in the Definition \ref{4.7.03.2m}. 
Just like in the ${\cal X}$-case,  the moduli space ${\cal A}_{G, \widehat S}$ 
can be understood as the quotient moduli stack. We leave the obvious details
to the reader. 

Let ${\rm Conf}_n({\cal A}):= G \backslash {\cal A}^n$ be the configuration
space of $n$ affine flags in $G$. 
\begin{definition} \label{4.21.03.157} The twisted cyclic shift map 
$S: {\rm Conf}_n({\cal A}) \lms 
{\rm Conf}_n({\cal A})$ is given by
$$
S: (A_1, A_2,..., A_{n}) \lms  (A_2, ..., A_n, s_G A_1 ).
$$
\end{definition}
Since $s_G$ is a central element one has 
$A \cdot s_G = s_G A$. Therefore $S^n = {\rm Id}$ because 
$$
S^n (A_1, A_2,..., A_n) = ( s_G A_1, s_G A_2, ...,  s_G A_n) \sim  (A_1, A_2,..., A_n).
$$
The coinvariants of the twisted cyclic map acting on the space ${\rm Conf}_n({\cal
  A})$ are called {\it twisted cyclic configurations} of $n$ affine flags 
in $G$.  Let $\widetilde {\rm Conf}_n({\cal A})$ be the moduli space  
of twisted cyclic configurations of $n$  affine 
flags in $G$. 

\begin{lemma} \label{10.20.03.1} 
Let $\widehat D_k$ be 
 a disc $D$ with $k$ marked points at the boundary. Then ${\cal A}_{G, \widehat
 D_k}$ 
is identified with the twisted cyclic configuration space 
$\widetilde {\rm Conf}_k({\cal A})$. 
\end{lemma}

{\bf Proof}. 
The punctured tangent bundle to a disc is homotopy equivalent to a
circle. So there is a unique up to isomorphism 
 twisted $G$-local system on $\widehat D_k$. Consider a path in $T'D$ 
obtained by choosing a non zero tangent vector field
on  a simple loop near the boundary of $D$. The monodromy around this
path is $s_G$. The marked points $x_1, ..., x_k$ on the boundary $C = \partial
D$ divide $C$ into a union of the arcs $C^*_1, ..., C^*_k$. 
We assume they follow the boundary clockwise.  
Let us choose a  vector $v_i$ tangent to $D$ at a certain point of the arc
$C^*_i$. We assume that the vectors $v_i$ look inside of the disc $D$. For
each vector $v_i$ there is a canonical homotopy class of paths from $v_i$ to
$v_1$ in $T'D$ obtained by counterclockwise transport of $v_i$ along $C$
towards $v_1$. Restricting the decoration section $\alpha$ to $v_1, ..., v_k$
and then translating the obtained affine flags in the fibers 
over $v_2, ..., v_k$ to the
fiber over $v_1$ we get a configuration of affine flags $(A_1, ..., A_k)$ over
$v_1$.  
If we choose $v_2$  instead of $v_1$ as the initial vector, this
configuration will be replaced by its twisted cyclic shift 
$(A_2, ..., A_k, s_G A_1)$. So we identified ${\cal A}_{G, \widehat
 D_k}$ with $\widetilde {\rm Conf}^*_k({\cal A})$. 
The lemma is proved. 
 \vskip 3mm
Recall the adjoint group  $G'$ corresponding to the group $G$. 
Since the canonical map $G \to G'$ kills the center, and $s_G$ is in the
center, a twisted $G$-local system ${\cal L}$ on $T'S$ provides a $G'$-local
system 
${\cal L}'$ on
$S$. 
The canonical projection $p: {\cal L}_{\cal A} \lra {\cal L}'_{\cal B}$ 
provides a map 
\begin{equation} \label{7.18.03.1}
p: {\cal A}_{G, \widehat S} \lra {\cal X}_{G', \widehat S}, \qquad 
({\cal L}, \alpha) \lms ({\cal L}', \beta), \quad \beta:= p(\alpha).
\end{equation}

\begin{lemma} \label{5.6.02.11} 
Let $k$ be the number of boundary points. Then \newline
$
{\rm dim} {\cal X}_{G', \widehat S} + k \cdot {\rm dim}H = {\rm dim} 
{\cal A}_{G, \widehat S}
$. 

 If $k=0$, one has 
${\rm dim} {\cal X}_{G', S} = {\rm dim} {\cal A}_{G, S} = 
{\rm dim} {\cal L}_{G, S} 
= 
(2g+n-2){\rm dim}G$. 
\end{lemma}

{\bf Proof}. The map (\ref{9.26.03.2}) is surjective over the generic point. 
Consider a boundary component $C_i$ without boundary points. 
The condition that a $G$--local system has unipotent 
monodromy around $C_i$ decreases the dimension by  ${\rm rank}(G)$, 
while adding a decoration at the puncture $s$ we increase the dimension by
 ${\rm dim} H = {\rm rank}(G)$. On the other hand 
over an arc the space of $\alpha$'s is $G/U$, while the space of $\beta$'s is 
$G/B$. The lemma is proved.

\vskip 3mm

{\bf 5.  Remarks, complements  and  examples}. 1. 
If $\widehat S = S$, i.e. there is no boundary points, it may be  
 convenient to assume that $S = \overline S - \{s_1, ..., s_n\}$ 
is a surface with $n$ punctures. 
For each punctures  there is a 
well defined isotopy class of  
a little oriented circle making a simple counterclockwise loop  
around $s$. Let ${\cal L}_s$ be the isomorphism class of the 
restriction of a $G$--local system ${\cal L}$ to it. 
 A {\rm framed}   $G$--local system on $S$ is then 
a $G$--local system ${\cal L}$ 
equipped with the following additional structure: 
for every puncture $s$ we choose a reduction of the 
local system ${\cal L}_s$ to the Borel subgroup $B$.  

 One can always avoid talking about holes $D_i$ in $\overline S$.
Indeed, let $S$ be a surface with $n$ punctures. 
Denote by  ${\rm Sp}_{s_i}{\cal L}$  the specialization of the local system ${\cal L}$ 
to the punctured tangent  space $T_{s_i}\overline S - \{0\}$ at $s_i$.
It is a local system on the punctured tangent space.  
The marked points are replaced 
by a collection of distinct rays in the punctured tangent  bundle 
$T_{s_i}\overline S - \{0\}$. 
A framed structure on ${\cal L}$ is given by the 
following data at the punctures $s_i$: 
restrict the local system ${\rm Sp}_{s_i}{\cal L}$ to the complement 
of the chosen rays, and take  its flat  section $\beta$. 
This language is convenient when we work with an algebraic structure 
on $S$. 
\vskip 3mm
2. Let us look more closely at the projection  (\ref{7.18.03.1}) 
in the case when $\widehat S = S$. 
The monodromy of a  $G$--local system  on a circle is 
described by a conjugacy class 
in $G$. The monodromy provides a bijective correspondence between the 
 conjugacy classes in $G$ and the isomorphism classes of  
 $G$--local systems on a circle. 
A  $G$--local system ${\cal L}$ on $ S$ is {\it  unipotent}  
 if for every boundary component $C_i$ the monodromy 
of the restriction of  ${\cal L}$ to ${C_i}$ is unipotent. 
 
\begin{definition} \label{11.4.02.1} 
${\cal U}_{G, S}$ is the moduli space of 
the   {\it framed 
unipotent} $G$--local systems on $S$. 
\end{definition}

Observe that a  $B$-local system  ${\cal L}_s$ on a circle has a 
unipotent monodromy if and only if it can be reduced to a $U$-local system. 
A decoration on ${\cal L}_s$ is  the same  as  
a reduction of the structure group to the subgroup $U$. 
So a decoration on a $B$--local system ${\cal L}_s$ 
on a circle exists if and only if ${\cal L}_s$ is unipotent. 
Thus a decorated $G$--local system on $S$ 
is unipotent. 
The right action of $H$ on ${\cal L}_s$ provides an action of $H$ on the 
set of all  decorations, making it into a principal homogeneous space over $H$.

Therefore forgetting the decoration we get a canonical  projection 
\begin{equation} \label{5.8.02.1} 
p': {\cal A}_{G, S} \lra {\cal U}_{G,S}.
\end{equation}
It provides ${\cal A}_{G, S}$ with a 
structure of a principal $H_G^{\mbox{\{punctures\}}}$-bundle over $ {\cal U}_{G,
  S}$, where $H_G$ is the Cartan group of $G$.  
Observe that the map $\kappa$ is the
projection onto the  group $H_{G'}^{\mbox{\{punctures\}}}$. 
There is a commutative diagram, where  
$i$ is the natural closed embedding:  
\begin{equation} \label{4.15.02.1f1}
\begin{array}{ccccccc}
{\cal A}_{G, S} 
&&& & &&\\
&&&& &&\\
p' \downarrow &&\searrow p & & &&\\
&&&& &&\\
 {\cal U}_{G', S} && \stackrel{i}{\lra}&&{\cal X}_{G', S} &
 \stackrel{\kappa}{\lra} H_{G'}^{\mbox{\{punctures\}}}.&
\end{array}
\end{equation}
  One has 
$
{\cal U}_{G', S} = p({\cal A}_{G, S}) = \kappa^{-1}(e)
$. 
\vskip 3mm

3. {\it Affine flags}. The points of ${\cal L}_{\cal B}$ parametrize 
the flags in the fibers of ${\cal L}$. In many cases, e.g. when the center of $G$ is trivial, there is a similar interpretation of the principal affine bundle 
 based on the notion of an affine flag. 

Let $A$ be a direct sum of  one dimensional $H$-modules corresponding to all 
simple roots. Let ${\cal U}$ be the Lie algebra of $U$.
The $H$-module 
${\cal U}/[{\cal U}, {\cal U}]$ is isomorphic to $A$. 
  \begin{definition} \label{5.6.02.5} Suppose that the center of $G$ is trivial. 
An {\it affine flag} is a pair 
\begin{equation} \label{5.6.02.1}
\{B, i: A \stackrel{\sim}{\lra} {\cal U}/[{\cal U}, {\cal U}]\}
\end{equation}  
 where $B$ is the Borel subgroup containing $U$, and 
$i$ is an isomorphism of $H$-modules.  
\end{definition}

The group $G$ acts by conjugation on affine flags.
 This action is transitive. Since the center of $G$ is trivial, 
   the stabilizer of 
the affine flag (\ref{5.6.02.1}) is $U$. So  
the principal affine bundle parametrises the affine flags in the fibers of ${\cal L}$. 
The same conclusion can be made for the classical groups with non trivial center. 
Here are some examples.  

\vskip 3mm
{\it The classical affine flags}. 
i). $G = GL_m$. 
Recall that 
an affine flag in a vector space $V$ is given by a flag (\ref{11.4.02.2}) in $V$ 
plus a choice of non zero vectors $v_i$ at each quotient $V_i/V_{i-1}$ 
for all $i=1, ..., m$. Equivalently, 
an affine flag in $V$ can be defined by a sequence of non zero 
volume elements $\eta_i \in {\rm det}V_i$:
$$
\eta_1:= v_1, \quad \eta_2:= v_1 \wedge v_2,\quad  ..., \quad \eta_m:= 
v_{1} \wedge ... \wedge v_{m}.
$$ 
Alternatively,  consider the dual volume forms  
$\omega_i \in {\rm det}V_i^*$, $<\eta_i, \omega_i>=1$. 
We may parametrise the affine flags by the sequence of volume forms 
$(\omega_1, ..., \omega_m)$. 
The variety of all affine flags is identified with $GL_m/U$ where 
$U$ is the subgroup of  upper triangular matrices in $GL_m$ 
with units at the diagonal. 

\vskip 3mm
ii). {\it The $SL_m$ case}.  
An affine flag in $V$ determines a volume form 
 in $V$. 
By an affine flag in a vector space with a given volume form 
$\omega$ we understand an affine flag whose 
volume form is the form $\omega$. 
The variety all affine flags with a given volume form is 
identified with $SL_m/U$. 
\vskip 3mm
{\bf Example}.  An  affine flag in a two 
dimensional vector space $V$ 
is a pair
 $\{$a non zero vector $v_1$ and a non zero $2$--form $\omega$ in $V$$\}$. 
The variety of all affine flags in a two dimensional 
symplectic vector space is identified with ${\Bbb A}^2 -\{0\}$. 
\vskip 3mm





\section{Surfaces, graphs, mapping class groups and groupoids}
\label{surf}

This Section contains some background material, suitably adopted or
generalized. We start 
from the ribbon graphs and
their generalizations needed to describe two
dimensional surfaces with marked points on the boundary. We use it throughout
the paper. Then we introduce 
a convenient combinatorial model of the classifying orbi-space  for 
mapping class groups which will be used in Section 15.  
In the end we recall some basic facts about the classical Teichm{\"u}ller spaces 
which will be used in Section 12 and some other parts of the paper. 
\vskip 3mm 
{\bf 1. Ribbon graphs and hyperbolic surfaces}. 
Recall that a 
{\it flag} of a graph $\Gamma$ is   
a pair $(v, E)$ where $v$  is a vertex of $\Gamma$,  $E$ is an edge of $\Gamma$,  
and $v$ is an an endpoint of $E$.  
A {\it ribbon graph} is a finite   
graph such that  for each vertex a cyclic order of the set of all 
flags sharing this vertex is chosen. 
An example of a ribbon graph is provided by a graph drawn on an oriented surface $S$. 
The cyclic order of the flags sharing the same vertex is given by the 
counterclockwise ordering of the flags. 

We say that an oriented path on a ribbon graph {\it turns left} 
at a vertex if we come to the vertex along the edge which is next 
(with respect to the cyclic order) to the edge  
used to get out of the vertex. 
A {\it face path} on a ribbon graph is a closed oriented path 
turning left at each vertex, which ends as soon as the path starts to repeat
itself. 
To visualize a face path imagine a ribbon graph glued of actual thin
ribbons.  Then a face path is a boundary component of such a ribbon graph. 
For example the punctured torus has one face path. 
A flag determines a face path. 
This face path can be visualized by a strip bounded by the face path. 
Conversely a face path plus a vertex (resp. edge) 
on this path determine uniquely the flag containing 
this vertex (resp. edge) situated on the face path. 

Given a ribbon graph $\Gamma$ we can get an oriented compact surface $S_{\Gamma, c}$ 
with a ribbon graph  isomorphic to $\Gamma$ lying on the surface. Namely,  
we take a disc for each of the face paths and glue its boundary to the graph along the face 
path. If we use punctured discs in the above construction we get 
an oriented surface $S_{\Gamma}$ which can be obviously retracted onto $\Gamma$. 
Clearly $S_{\Gamma}$ is obtained by making $n$ punctures on  $S_{\Gamma, c}$, where 
$n= F(\Gamma)$.
Any surface with negative Euler characteristic can be obtained as $S_{\Gamma}$ 
for certain ribbon graph $\Gamma$. Indeed, consider a triangulation 
of the surface with the vertices at the punctures, and take the dual graph. 



\vskip 3mm 
\noindent \qquad {\bf 2. Ideal triangulations of marked hyperbolic surfaces 
and marked trivalent graph}. 
Let $\widehat S$ be a marked hyperbolic surface. If its  boundary component $C_i$ 
has no boundary points on it, let us 
picture the corresponding hole on $S$ as a puncture. 

A {\it marked trivalent graphs
 on $S$ of type $\widehat S$} is a graph $\Gamma$  on $S$ with the 
following properties:

i) The vertices of $\Gamma$ are of valence three or one.

ii) The set 
of univalent vertices coincides with the set of  marked points 
on the boundary $\partial S$.

iii) The surface $S$ can be shrinked onto $\Gamma$. 

We define  an {\it ideal  
triangulation} of  $\widehat S$ as a triangulation of $S$ with the 
vertices either at the punctures or at the boundary arcs, such that 
each arc carries exactly one vertex of the triangulation, and 
each puncture serves as a vertex.  
The edges of an ideal triangulation either lie inside of  
$S$, or belong to the boundary of $S$. 

The duality between triangulations and trivalent 
graphs provides a bijection between the 
isotopy classes of  ideal triangulations of $\widehat S$ 
and isotopy classes of 
marked trivalent graphs on $\widehat S$. 
The ends of the dual graph are dual to the boundary edges 
of the ideal triangulation.  

A {\it marked ribbon graph} 
is a graph with trivalent and univalent vertices 
equipped with a cyclic order of the edges sharing each trivalent vertex. 
Given a marked ribbon graph one can construct the corresponding 
marked surface.

Consider a marked graph $(\overline \Gamma, v)$ which has 
 just one $4$-valent vertex in addition to
 trivalent and univalent vertices. Then there are two marked 
trivalent graphs $\Gamma$ and $\Gamma'$ 
related to $(\overline \Gamma)$, as shown on the picture. 
We say that $\Gamma'$ is obtained from $\Gamma$ 
by a flip at the edge $E$ shrinked to $v$. A similar operation on 
triangulations consists of changing an edge in a $4$-gon of the triangulation. 
It is well known that any two isotopy classes of marked trivalent graphs 
on a surface are related by a sequence of flips. Indeed, here is an argument, which assumes for simplicity that 
there are no marked points. Then, presenting the surface as a $4g$-gon with glued sides, we reduce the claim 
to triangulations of the $4g$-gon punctured in $n-1$ points, where it is obvious.

\begin{figure}[ht]
\centerline{\epsfbox{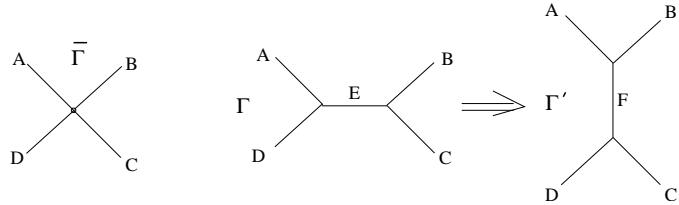}}
\caption{Flip at the edge $E$.}
\label{fg11}
\end{figure}
\vskip 3mm 
{\bf 3. The mapping class group of $S$}. 
The mapping class  group $\Gamma_S$ is the quotient of the group 
of all orientation preserving diffeomorphisms of $S$ modulo its connected 
component  consisting of the diffeomorphisms isotopic to the identity. 
It has the following well known  group theoretic description. 
Every hole $D_i$ on $S$ provides  
a conjugacy class in the 
fundamental group $\pi_1(S,x)$.  
The mapping class group $\Gamma_S$ consists of all 
outer automorphisms of $\pi_1(S,x)$  
preserving  each of these 
 conjugacy classes.
The group $\Gamma_S$ acts 
by automorphisms of the ${\cal A}$ and ${\cal X}$  moduli spaces. 
We will use a combinatorial description of 
 the mapping class group as the automorphism group of 
objects of the {\it modular groupoid} of $S$ defined below. 

\vskip 3mm 
{\bf 4. The graph complex $G_{\bullet}(S)$ of  
a hyperbolic surface $S$}.  It is a version of the complex used by 
Boardman and Kontsevich \cite{Ko}. To define it, recall  the $\Z/2\Z$--torsor ${\rm Or}_{\Gamma}$ 
of orientations of a graph $\Gamma$, 
given by invertible elements in the determinant of the 
free abelian group generated by the edges of $\Gamma$. 
Let $G_{k}(S)$ be the abelian group of $\Z$--valued functions 
$\varphi(\Gamma, \varepsilon_{\Gamma})$ with finite support on 
the set of pairs $(\Gamma, \varepsilon_{\Gamma})$, 
where $\Gamma$ is an isotopy class of 
 an embedded graph on  $S$ such that $S$ can be retracted to it, 
$\varepsilon_{\Gamma} \in {\rm Or}_{\Gamma}$, one has 
$\varphi(\Gamma, -\varepsilon_{\Gamma}) = -\varphi(\Gamma, 
\varepsilon_{\Gamma})$,  and  
$k$ is given by 
\begin{equation} \label{12.14.02.5}
k -3 = \sum_{v \in V(\Gamma)} ({\rm val}(v) -3).  
\end{equation}
Here ${\rm val}(v)$ is the valency of a vertex $v$ of $\Gamma$. 
Observe that if all except one vertex of $\Gamma$ are of valency three, then 
$k$ is the valency of this one vertex. 
The differential $\partial: G_{\bullet}(S) \lra G_{\bullet-1}(S)$ 
is provided by contraction of the edges of $\Gamma$:
\begin{equation} \label{12.15.02.1}
\partial \varphi(\Gamma, \varepsilon_{\Gamma}):= \sum_{E', \Gamma'}
\varphi(\Gamma', \varepsilon_{\Gamma'}).
\end{equation}
Here the sum is over all pairs $(E', \Gamma')$, where $E'$ is an edge of 
$\Gamma'$,  such that 
 the ribbon graph $\Gamma'/E'$ obtained 
by shrinking of $E'$ is isomorphic $\Gamma$, and 
$\varepsilon_{\Gamma'} = E' \wedge \varepsilon_{\Gamma}$. 
Then $\partial^2 =0$. 
We place the trivalent graphs in the degree 
$3$, and consider  $G_{\bullet}(S)$ as a homological complex.

This normalisation will be useful in Section 14. 
In this Section we will use the shifted complex $G_{\bullet}(S)[3]$, 
where the trivalent graphs are in degree $0$. 
The shifted graph complex $G_{\bullet}(S)[3]$ is the chain complex of 
a finite dimensional 
contractible polyhedral complex ${\Bbb G}(S)$, called the {\it modular complex} of
$S$, 
constructed in the next
subsection. 
\vskip 3mm 
{\bf 5. The modular complex of $S$}. 
Recall the $n$--dimensional Stasheff polytope $K_n$. 
It is  a convex polytope
whose combinatorial structure is described as follows. 
Consider plane trees with internal vertices of valencies $\geq 3$, and with $n+3$ ends lying on an oriented circle. 
The ends are marked by the set $\{1, ..., n+3\}$,  so that its 
 natural order is compatible with the circle orientation.  
The codimension $p$ faces of $K_n$ are 
parametrized by the isotopy classes of such plane trees with  
$p+1$ internal vertices. 
Let  $F_T$ be the face corresponding to a tree $T$. 
Then $F_T$ belongs to the boundary of $F_{\overline T}$ 
if and only if $ \overline T$ is obtained from $T$ by shrinking 
some internal edges. 
The Stasheff polytope $K_n$  can be realized as a cell in 
the configuration space of 
$n+3$ cyclically ordered points on a circle, see \cite{GM} for a survey. 
\vskip 3mm 
{\bf Example}. $K_0$ is a point, $K_1$ is a line segment, 
and $K_2$ is a pentagon. 
\vskip 3mm 

Let us define a finite dimensional 
contractible polyhedral complex ${\Bbb G}(S)$, called 
the modular complex of $S$. The faces of the modular  
complex ${\Bbb G}(S)$  
are  parametrized by the  isotopy classes  of 
graphs  $\Gamma$ 
on $S$, homotopy equivalent to $S$. 
The polyhedron $P_{\Gamma}$ corresponding to $\Gamma$ is 
a product of the Stasheff polytopes: 
\begin{equation} \label{12.21.02.10}
P_{\Gamma} : = \prod_{v \in V(\Gamma)}K_{{\rm val}(v)-3}.
\end{equation}
Here $V(\Gamma)$ is the set of vertices of $\Gamma$, and 
${\rm val}(v)$ is the valency of a vertex $v$. 
The face $P_{\Gamma}$ is at the boundary of the face $P_{\overline \Gamma}$ 
if and only if $ \overline \Gamma$ is obtained from $\Gamma$ by 
shrinking of some edges. 
Precisely, the inclusion $ P_{\Gamma} \hookrightarrow P_{\overline 
\Gamma}$
is provided by the same data for the Stasheff polytopes appearing 
in (\ref{12.21.02.10}). 
The mapping class group ${\Gamma}_S$ acts on ${\Bbb G}(S)$ preserving 
the polyhedral structure.

\vskip 3mm 
{\bf Examples}. 1. The vertices of ${\Bbb G}(S)$ correspond to the 
isotopy classes of the 
trivalent graphs on  $S$.  The edges -- to the graphs with one 
$4$--valent vertex.  
The two dimensional faces are either pentagons corresponding to 
graphs with a unique $5$--valent vertex, or squares corresponding to 
graphs with two $4$--valent vertices. 
\vskip 3mm 

{2}. Let $S$ be a torus punctured at a single point. The 
mapping class group of $S$ is isomorphic to $SL_2(\Z)$. 
Then ${\Bbb G}(S)$ is the tree dual to the classical modular 
triangulation of the hyperbolic plane:

\begin{figure}[ht]
\centerline{\epsfbox{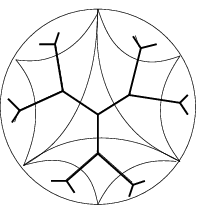}}
\caption{The tree dual to the modular triangulation.}
\label{fg54}
\end{figure}

\begin{proposition} \label{12.02.15.30s}
 The shifted graph complex ${ G}_{\bullet}(S)[3]$ is the chain complex of 
the modular complex ${\Bbb G}(S)$. 
\end{proposition}

{\bf Proof}. The faces of the modular complex are glued according to the 
boundary map (\ref{12.15.02.1}). 
An orientation  $\varepsilon_{\Gamma}$ of $\Gamma$ 
 provides 
 an orientation of the face $P(\Gamma)$. The lemma now 
follows from the very definition. 
\vskip 3mm

\begin{proposition} \label{12.02.15.30}
 The polyhedral complex ${\Bbb G}(S)$ is contractible.
The stabilisers of the mapping class group acting 
on ${\Bbb G}(S)$ are finite. 
\end{proposition}

The proof is deduced from Strebel's theory of quadratic differentials 
on Riemann surfaces [Str], see also \cite{Ha}. 
 One can show that there exists 
a ${\Gamma_S}$-equivariant contraction of the 
Teichm{\"u}ller space to the modular complex of $S$. 
So the modular complex is a "combinatorial version" 
of the Teichm{\"u}ller space. 

  \vskip 3mm 
{\bf 6. Modular groupoid and  combinatorial 
description of the mapping class group}. 
The quotient ${\Bbb G}(S)/\Gamma_S$ 
is a polyhedral orbifold. 
By proposition \ref{12.02.15.30} its orbifold fundamental group 
is isomorphic to $\Gamma_S$. 
The 
{\it modular groupoid of $S$} is a 
 combinatorial version of the Poincare groupoid of this  orbifold. 
Precisely, the objects of the modular groupoid of $S$ 
are $\Gamma_S$--orbits of the vertices of the modular complex ${\Bbb G}(S)$, 
i.e. isomorphism classes of   the trivalent ribbon graphs of type $S$.  
The morphisms are represented by paths on the $1$-skeleton composed with 
the  graph automorphisms, i.e. the elements of $\Gamma_S$ stabilizing  the vertices. 
The two paths represent the same morphism if an only if they are homotopic 
in the $2$--skeleton ${\Bbb G}_{\leq 2}(S)/\Gamma_S$.

Here is a description of the modular groupoid by generators and relations. 
Recall the flip determined by an oriented edge of the modular complex, see 
Figure  \ref{fg11}. 
Every morphism in the modular groupoid is a composition of the elementary ones, given by 
flips and graph symmetries. There are two types 
of  relations between the elementary 
morphisms provided by the two types of the two dimensional faces of ${\Bbb G}(S)$:

i) ({\it Square}). The flips performed at two disjoint edges commute.

ii) ({\it Pentagon}). The composition of the five flips 
corresponding to the (oriented) boundary of a pentagonal face of ${\Bbb G}(S)$ is zero. 

They together with the graph symmetries generate all the 
relations between the morphisms in the modular groupoid. 
\vskip 3mm 
{\bf 7. Background on the Teichm{\"u}ller spaces}.  
The classical Teichm{\"u}ller space ${\cal T}_S$ for a hyperbolic surface $S$ 
can be defined as the moduli space of any of the following four different objects 
related to  $S$:

1) Complex structures on $S$.

2) Faithful representations $\pi_1(S) \to PSL_2(\R)$ 
with discrete image, modulo $PSL_2(\R)$-conjugation.

3) Complete Riemannian metrics on ${\rm int}(S)$ of curvature $-1$, called  
{\it complete hyperbolic metrics}.

4) Hyperbolic metrics on $S$  with geodesic
boundary $\partial S$ or a cusp. 

Indeed, the uniformisation theorem implies 
$1) \iff 2)$.  Taking the quotient of the hyperbolic plane ${\cal H}$ by
the action of the image of $\pi_1(S)$ in $PSL_2(\R)$ we
go from 2) to 3). A Riemannian metric on $S$ provides a complex structure on
$S$, and hence the arrow $3) => 1)$. Taking the universal cover of 
$S$ we see the equivalence $2) \iff 3)$. 
For each hole on a surface $S$ with complete hyperbolic structure 
there is a unique geodesic, called the {\it boundary geodesic},  homotopic to the boundary of
the hole, 
except the case when the hole degenerates to a puncture, i.e. the
corresponding end of $S$ is a cusp. Cutting out the ends of the surface
$S$ along the boundary  geodesics we get a surface with geodesic
boundary. This leads to the equivalence $3) \iff 4)$. 

The universal cover of a surface $S$ equipped with a hyperbolic metric with
geodesic boundary can be realized as the hyperbolic plane minus 
an infinite collection of disjoint discs bounding geodesics. 
These geodesics project onto the boundary geodesics on $S$. 
Precisely, let $S'$ be a surface with
 complete hyperbolic metric corresponding to $S$, i.e. 
$S = S' - T_1\cup  ... \cup  T_n$ where $T_i$ is
 the tube cut out from $S'$ by the boundary geodesic corresponding to the
 $i$-th hole. The preimage of each  boundary geodesic 
is a disjoint collection of geodesics on ${\cal H}$, 
with the  monodromy group $\Gamma$ acting transitively on it.  
The tube  $T_i$ lifts to a union of 
discs bounding these geodesics. Cutting out these discs we get a
universal cover ${\cal H}_S$ of $S$ sitting inside of ${\cal H}$. 
The boundary of ${\cal H}_S$ is 
the limit set of $\Gamma$. If the lengths of the boundary geodesics are non zero, it is 
a Cantor set.

If the number $n$ of holes is positive, 
the Teichm{\"u}ller space ${\cal T}_S$ is a manifold with corners. 
Choosing orientations of all boundary
geodesics we get a version of the Teichm{\"u}ller space denoted ${\cal T}^+_S$. 
Forgetting orientations of the geodesics we get a $2^n:1$ covering 
${\cal T}^+_S \to {\cal T}_S$. On the other hand the orientation of $S$
provides canonical orientations of the geodesics, and hence an embedding 
${\cal T}_S \hra {\cal T}^+_S$.

\begin{figure}[ht]
\centerline{\epsfbox{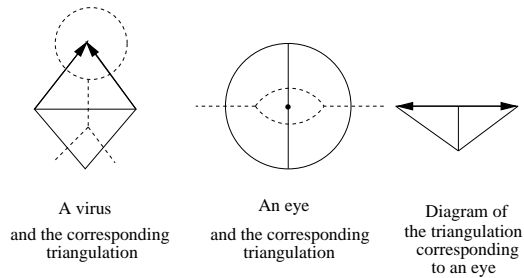}}
\caption{A virus and an eye, with the corresponding triangulations.}
\label{fgo-100}
\end{figure}

\vskip 3mm 
{\bf 8. Special graphs and flips}. Recall that a trivalent graph on $S$ is
{\it special} 
if it contains, as a part, one of the two graphs shown on Figure
\ref{fgo-100} by dotted lines. 
We call them 
an eye and a virus. The 
 triangulations dual to the eye and virus are shown on Figure \ref{fgo-100} by
 solid lines. 
The sides shown by arrows are identified. The triangulations 
dual to the graphs containing a virus  are
 characterised by the property that they 
contain a triangle with two sides identified. 
To obtain a diagram of the 
triangulation corresponding to an eye we cut one of the edges of the
triangulation,
 the top one on the middle picture, 
  getting the right diagram on  Figure \ref{fgo-100}.

\begin{figure}[ht]
\centerline{\epsfbox{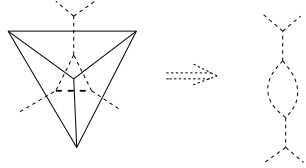}}
\caption{An eye is obtained by a flip at the bottom edge on the left.}
\label{fgo-101}
\end{figure}

\begin{figure}[ht]
\centerline{\epsfbox{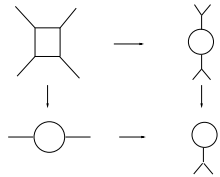}}
\caption{A flip from an eye to a virus is a folding of a composition
  of two flips on regular graphs.}
\label{fgoeye}
\end{figure}

Figures \ref{fgo-101} and \ref{fgoeye} 
show that the special graphs can be obtained by flips from the  usual graphs. 
Namely, on Figure \ref{fgo-101} we obtain an eye by making a flip at the
bottom horizontal edge of the graph on the left. The bottom of 
Figure \ref{fgoeye} shows that a flip at the leg of a virus produces an
eye, and vice versa. A flip at the loop of a virus is not defined.

Figure \ref{fgoeye} shows that a flip from an eye to a virus, lifted to a $2:1$ cover 
of the surface ramified at the puncture at the center of the eye, turns out
to be a composition of flips at the edges $E_1$ and $E_2$ of a regular graph. 
The edges $E_1$ and $E_2$ are the preimage of the edge $E$. The non-trivial
automorphism of the cover is given by the central symmetry. Below we always
treat 
a special triangulation of $S$ by going to a cover of $S$
where the triangulation becomes regular.


\begin{figure}[ht]
\centerline{\epsfbox{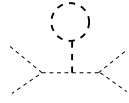}}
\caption{Adding a virus to a graph.}
\label{fgo-121}
\end{figure}

{\bf Proof of Lemma \ref{5.16.04.1}}. 
Let $S$ be a surface with a trivalent graph $\Gamma$ on $S$, which is homotopy
equivalent to $S$. Let us make a puncture of $S$, getting a new surface
$S'$. Then adding a virus to $\Gamma$ as shown on Figure \ref{fgo-121}, we get a
trivalent graph $\Gamma'$ homotopy equivalent to $S'$. So $S'$ is special. 
Hence any marked  
surface with at least two holes, such that the boundary of one of the holes is
without
marked points, is special. 
Let us prove the converse claim. A virus separates the surface into two
domains, each of them containing a whole. Moreover the one of them surrounded
by the head of the virus has no marked points on its boundary. It works the same for an eye. The lemma is
proved. 
\vskip 3mm 

\section{Elements of positive algebraic geometry}
\label{posgeom}


{\bf 1. Semifields and positive varieties}. 
Following  \cite{BFZ96}, we define  a {\it semifield}  as a   
 set $P$  endowed with the  operations of
 addition and multiplication, which have the following properties:

i) Addition in $P$ is commutative and associative.

ii) Multiplication makes $P$ an abelian group.

iii) distributivity: $(a+b)c = ac + bc$ for $a,b,c \in P$. 

Here are   important examples of  semifields:

a) $P = \Q_{>0} = \{x \in \Q| x>0\}$,

b) $P = \R_{>0} = \{x \in \R| x>0\}$,

c) Let $K$ be a field and $P\subset  K$  a  semifield inside of $K$, 
with the addition and multiplication in $P$ induced from $K$, 
like in the examples a) and b). Then there is a new 
semifield
$$
K((\varepsilon))_{P}:= \{f=a_s\varepsilon^s + 
a_{s+1}\varepsilon^{s +1} + ... | a_i \in P\} \subset K((\varepsilon))
$$
 consisting of all Laurent power series with the leading coefficient from $P$. 
For example starting from  $\Q_{>0}$ we get a 
semifield denoted $\Q((\varepsilon))_{>0}$. Similarly one defines 
$\R((\varepsilon))_{>0}$ starting from $\R_{>0}$.

d) The tropical semifields:  
$P = \Z$, $P = \Q$ or $P = \R$ with the multiplication $\otimes $ 
and addition $\oplus$ given by 
\begin{equation} \label{10.24.03.1}
a \otimes b := a+b, \quad a \oplus b := {\rm max}(a,b).
\end{equation}

Observe that in the examples a),  b) and c) the semifield $P$ is naturally 
embedded to a field, while in the example d) this is impossible since 
$a  \oplus a =a$. 

\begin{lemma} \label{8.27.03.13} Assigning to $f \in 
K((\varepsilon))_{P}$ the minus degree of its leading coefficient we get a homomorphism 
of  semifields
$$
-{\rm deg}: K((\varepsilon))_{P} \lra \Z^t.
$$
\end{lemma}

 {\bf Proof}. Clear. 
\vskip 3mm 
There is a unique way to 
assign to a split torus $H$ an abelian group $H(P)$ so that 
${\Bbb G}_m(P) = P$ and $H \to H(P)$ is a functor from the category 
of split algebraic tori to the one of abelian groups. One has $H(\Z^t) =
X_*(H)$ and 
$H(P) = X_*(H)\otimes_{\Z}P$, where $P$ is considered here as an abelian
group.  

Let ${\rm Pos}$ be the category whose objects are split algebraic tori
and morphisms are positive rational maps. 
A positive rational map $\varphi: H_1 \to H_2$ gives rise to a map of sets 
$\varphi_*: H_1(P) \to H_2(P)$, providing us a functor  
${\rm Pos} \lra {\rm Sets}$. 
In particular if $\varphi$  is invertible 
in the category ${\rm Pos}$ then $\varphi_*$ is an isomorphism. Therefore given a positive variety $(X, \{\psi_{\alpha}\})$ we get a collection of abelian groups $H_{\alpha}(P)$ 
and set isomorphisms $\psi_{\alpha, \beta}: H_{\alpha}(P) \to H_{\beta}(P)$ between them. 

Recall the definition of a positive structure on a scheme/stack 
given in Section 1.3. 
\begin{definition} \label{8.28.03.10}
\!\!Let $P$ be a\!  semifield and $(X, \{\psi_{\alpha}\})$ a\! positive scheme. 
 Then its $P$--part $X(P)$ is a set whose elements are collections  
$\{p_{\alpha} \in H_{\alpha}(P)\}$ such that  
$\psi_{\alpha, \beta}(p_{\alpha}) = 
p_{\beta}$ for any 
$\alpha, \beta \in {\cal C}_X$. 
\end{definition}
Decomposing $H_{\alpha}$ into a product of ${\Bbb G}_m$'s 
we get an 
isomorphism
$$
\widetilde \varphi_{\alpha}: P^{{\rm dim}H_{\alpha}} 
\stackrel{\sim}{\lra} X(P).
$$

The $P$-part $X(P)$ of $X$ depends only on the positive 
birational class of   $X$, and compatible positive 
structures give rise to the same $P$-part. 

If  $P=K_{>0}$ is a  semifield in a field $K$ 
we can realize the $K_{>0}$--part of $X$ inside of 
the set of $K$-points of $X$.  
Indeed, 
thanks to the condition i) of definition \ref{7.28.03.1} 
the map $\psi_{\alpha}$ is regular at $H_{\alpha}(K_{>0})$. Therefore  
there is a subset  
\begin{equation} \label{4.03.23.8}
X^+_{K_{>0}}:= \psi_{\alpha}(H_{\alpha}(K_{>0})) \hookrightarrow X(K).
\end{equation}
Thanks to the condition ii) of definition \ref{7.28.03.1} it 
does not depend on $\alpha$. 

We say that a positive atlas $(\{\psi'_{\alpha}\}, {\cal C}')$ 
 on $X$ {\it dominates} 
another  one $(\{\psi_{\alpha}\}, {\cal C})$ if there is an inclusion 
$i: {\cal C} \hra {\cal C}'$ such that for any $\alpha \in {\cal C}$ 
the coordinate systems $\psi_{\alpha}: H_{\alpha} \to X$ 
and $\psi'_{i(\alpha)}: H'_{i(\alpha)} \to X$  are identified. 
So the positive atlases on $X$ are partially ordered. 
If two positive atlases ${\cal C}'$ and ${\cal C}''$ dominate 
the same one ${\cal C}$, then there is a new positive 
atlas ${\cal C}'\cup {\cal C}''$ dominating all of them. 
Therefore for a given positive 
atlas on $X$ there is the {\it maximal positive 
atlas} on $X$  dominating it. 

Let $(X, \{\psi_{\alpha}\}, {\cal C})$ be a positive scheme. 
Then  $\psi_{\alpha}(H_{\alpha})$ 
is a Zariski open subset of $X$. 
Their union over all $\alpha \in {\cal C}$ 
is a regular Zariski open subset of $ X$. 
Applying this construction to the maximal positive atlas 
on $X$ dominating the original one we get 
a regular Zariski open  subset 
$
X^{\rm reg} = X^{\rm reg}_{\cal C} \hra X
$ 
 called the 
{\it regular part} of a positive variety $X$.

The $P$-part of $X$ can be  defined in a more general set up of positive spaces. 
\vskip 3mm 
{\bf 2. A category of positive spaces}. 
Recall that a groupoid is a category where all morphisms are isomorphisms. 
Below we will assume that the set of
isomorphisms between every two objects of a groupoid  is non-empty, i.e. all groupoids are 
connected.  
Let $ {\cal G}$ be a connected groupoid and $\alpha$ its object. 
The {\it fundamental group} $\Gamma_{ {\cal G}, \alpha}$ 
of the groupoid $ {\cal G}$ based at $\alpha$ is the
group of automorphisms of the object $\alpha$. The fundamental groups
based at different objects are isomorphic, but the
isomorphism is well defined only up to an inner automorphism.

 \begin{definition} \label{10.31.03.1} \!A {\em positive space} ${\cal X}$ is a pair consisting of a groupoid $ {\cal G}_{\cal X}$ and a functor 
\begin{equation} \label{11.6.03.99}
\psi_{\cal X}:  {\cal G}_{\cal X}\lra {\rm Pos}.
\end{equation}
The groupoid ${\cal G}_{\cal X}$ is called the coordinate groupoid of the positive
space. 
\end{definition}
Thus for every object $\alpha$ of the coordinate groupoid $ {\cal G}_{\cal X}$ 
there is an algebraic torus $H_{\alpha}$, 
and for every morphism $f: \alpha \lra \beta$ in the groupoid 
there is a positive birational isomorphism 
$\psi_{f}: H_{\alpha} \lra  H_{\beta}$.

\noindent \qquad Let ${\cal X}$ and ${\cal Y}$ be two positive spaces defined via the coordinate groupoids 
$ {\cal G}_{\cal X}$ and $ {\cal G}_{\cal Y}$ and functors $\psi_{\cal X}$ and 
 $\psi_{\cal Y}$. A morphism ${\cal X} \to {\cal Y}$ 
is given by a functor $\mu_{{\cal X}, {\cal Y}}:  {\cal G}_{\cal X} \to 
 {\cal G}_{\cal Y}$ 
and a natural transformation of functors 
$m: \psi_{\cal X} \to \psi_{\cal Y} \circ \mu_{{\cal X}, {\cal Y}}$. 
A morphism is called a {\it strict morphism}  if 
for every object $\alpha \in  {\cal G}_{\cal X}$ the map $m_{\alpha}: \psi_{\cal X}(\alpha) \to 
\psi_{\cal Y} \circ \mu_{{\cal X}, {\cal Y}}(\alpha)$ is a homomorphism of algebraic tori.

Given a semifield $P$ and a positive space ${\cal X}$ defined via 
(\ref{11.6.03.99}) we get a functor
$$
\psi_{\cal X}^P:  {\cal G}_{\cal X}\lra {\rm Sets}, \qquad \alpha \lms H_{\alpha}(P), 
\quad H_{\alpha}:= \psi_{\cal X}(\alpha).
$$
It describes the set 
${\cal X}(P)$ of 
$P$-points of ${\cal X}$. The fundamental group of the groupoid
$ {\cal G}_{\cal X}$ acts on ${\cal X}(P)$.
\vskip 3mm 
{\bf Example 1}. A positive variety $X$  provides a functor 
(\ref{11.6.03.99}) as follows.  The fundamental group of the coordinate
groupoid $ {\cal G}_{\cal X}$ is trivial, so it is just  a
set. 
Namely, it is the set parametrising 
positive coordinate charts of the atlas defining the positive structure of $X$. 
The functor $\psi_{\cal X}$ is given by  $\psi_{\cal X}(\alpha):= H_{\alpha}$  and 
 $\psi_{\cal X}(\alpha\to \beta):= \psi_{\alpha, \beta}$. 

\vskip 3mm 
{\bf Example 2}. Let $X$ be a $\Gamma$-equivariant positive scheme. 
It   provides a positive space ${\cal X}$ given by a functor 
(\ref{11.6.03.99}) such that the fundamental group of the groupoid  
of $ {\cal G}_{\cal X}$ is isomorphic to $\Gamma$.   
\vskip 3mm

{\bf 3. The 
semiring of good positive Laurent polynomials for a positive 
space ${\cal X}$}. Let $T_{i, \alpha}$ be the coordinates on the tori $H_{\alpha}$. 
Then $\Z[H_{\alpha}]$ is the ring of Laurent polynomials in 
$T_{i, \alpha}$. 
Let us say that a Laurent polynomial $F$ in the variables $T_{i, \alpha}$ 
is {\it good} 
it for every morphism $f: \beta \to \alpha$ in ${\cal G}_X$ the rational
function $f^*F$ is a Laurent polynomial in the coordinates
$T_{i, \beta}$. Denote by   ${\Bbb L}({\cal X})$ the ring 
of all good Laurent polynomials for 
${\cal X}$. We say that $F$ is a  {\it good positive Laurent
  polynomial} if for every morphism $f: \beta \to \alpha$ in
${\cal G}_X$ the function $f^*F$ is a Laurent polynomial  in the coordinates
$T_{i, \beta}$ with positive integer
 coefficients.  Denote by ${\Bbb L}_+({\cal X})$ the semiring of good 
positive Laurent polynomials for  
${\cal X}$. 

Let $\Q_+({\cal X})$ be the semifield of all positive
rational functions for ${\cal X}$. Set  
$
\widetilde {\Bbb L}_+({\cal X}) = \Q_+({\cal X}) \cap {\Bbb L}({\cal X})
$. It contains ${\Bbb L}_+({\cal X})$, but may be bigger, as the following
example shows: $x^2-x+1 = (x^3+1)/(x+1)$.

The semiring ${\Bbb L}_+({\cal X})$ determines the set of {\it extremal elements} ${\bf E}({\cal X})$, that is the elements which can not be decomposed into a sum of two non zero 
good 
positive Laurent polynomials. Any element of 
${\Bbb L}_+({\cal X})$ is a linear combination of the extremal elements with 
non negative integer coefficients. 

\vskip 3mm 
{\bf Remark}. If we 
replace a positive atlas on a scheme $X$  by one dominating it, the set of 
points of $X$ with values in a semifield $P$ does not change, while   
the (semi)ring of good (positive) Laurent polynomials can get smaller. 
\vskip 3mm 

{\bf 4. The schemes ${\cal X}^{*}$ and ${\cal X}^{**}$}. 
Let ${\cal X} = (\psi_{\cal X}, {\cal G}_{\cal X})$ be a positive space. 
Let $\Gamma:= \Gamma_{{\cal X}, \alpha}$ be the fundamental group 
of the groupoid ${\cal G}_{\cal X}$. Below we define 
two $\Gamma$-equivariant positive schemes, ${\cal X}^{**} \hra {\cal X}^{*}$,
corresponding to ${\cal X}$ in the sense that ${\cal X}$ is the positive space 
assigned to the positive schemes ${\cal X}^*$ and ${\cal X}^{**}$. 
We start with the scheme  ${\cal X}^{*}$. 

\begin{lemma} \label{12.6.03.1}
Let ${\cal X}$ be a positive space. Then there exists a unique 
positive regular $\Gamma$-equivariant scheme ${\cal X}^{*}$ corresponding to
${\cal X}$ which is ''minimal'' in the following sense. 
For any positive regular $\Gamma$-equivariant scheme ${\cal Y}$ corresponding to
${\cal X}$ there exists a unique morphism ${\cal X}^{*} \hra {\cal Y}$. 
\end{lemma}

{\bf Proof}. Let us assume first that the fundamental group of 
the coordinate groupoid of ${\cal X}$ is trivial. Then the scheme   
${\cal X}^{*}$ is obtained by gluing the coordinate tori 
$H_{\alpha}$ according to  the birational isomorphisms $\psi_{\alpha, \beta}$. 
For any field $K$ the set of its $K$-points 
is given by collections of elements $p_{\alpha} \in H_{\alpha}(K)$ such that 
$\psi_{\alpha, \beta}p_{\alpha} = p_{\beta}$. The general case is obtained by
considering the universal covering groupoid for ${\cal G}_{\cal X}$.
The lemma is proved. 
\vskip 3mm

The ring ${\Bbb L}[{{\cal X}}]$ is identified with the ring 
of global regular functions 
on  ${\cal X}^{*}$. 

Let $\{\psi^{*}_{\alpha}\}$ be the positive atlas of the scheme
${\cal X}^{*}$. 
The scheme ${\cal X}^{**}$ is defined as the subscheme of ${\cal X}^{*}$ 
given by the intersection of all coordinate charts of this atlas: 
$$
{\cal X}^{**}:= \cap_{\alpha \in {\cal
    C}_{\cal X}}{\rm Im}\psi^{*}_{\alpha}(H_{\alpha}) \subset {\cal X}^{*}
$$
It is a $\Gamma$-equivariant positive scheme: the positive atlas for ${\cal X}^{*}$
induces the one for ${\cal X}^{**}$.  
For any $\alpha \in {\cal C}_{\cal X}$ 
the scheme ${\cal X}^{**}$ is 
the complement to a possibly infinite collection of positive divisors in 
the torus $H_{\alpha}$. 

If the fundamental group of the groupoid ${\cal G}_{\cal X}$ is trivial, 
 ${\cal X}^{**}$ is characterised by the following universality property. 
It is the maximal scheme such that for every $\alpha$ there exists 
a regular open embedding $i_{\alpha}: {\cal X}^{**} \hra H_{\alpha}$ and
$\psi_{\alpha, \beta} i_{\alpha} = i_{\beta}$. 
Reversing the arrows we arrive at the universality property for  the scheme 
${\cal X}^{*}$.

Let $(X, \{\psi_{\alpha}\}, {\cal C})$ be a positive regular scheme/stack.\newline 
Then  $X^{*} = \cup_{\alpha \in {\cal
    C}}{\rm Im}\psi_{\alpha}(H_{\alpha})$ and one has
$
X^{**} \subset X^{*}  \subset X
$.

\begin{lemma} Let $(X, \{\psi_{\alpha}\}, {\cal C})$ be a positive 
regular scheme/stack, such that  $X$ is noetherian. Then $X^{*}$ can be  covered 
by a finite collection of 
the open subsets $\psi_{\alpha}(H_{\alpha})$.
\end{lemma} 

{\bf Proof}. Since $X$ is noetherian, any open subset of $X$ is quasi-compact, 
i.e. any cover by open 
subsets admits a finite subcover. The subset  $X^{*}$ is open, and covered by 
the open subsets $\psi_{\alpha}(H_{\alpha})$. Therefore there exists 
a finite subcover. The lemma is proved. 
\vskip 3mm
Therefore for a positive regular noetherian scheme $X$ 
the ring ${\Bbb L}(X)$ as well as the semiring ${\Bbb L}_+(X)$ are 
 determined by a finite number of conditions. 

\vskip 3mm

{\bf 5. $\Z$--, $\Q$-- and $\R$--tropicalizations
 and the Thurston boundary of a positive space}. 
Let ${\Bbb A}$ be either $\Z$, or $\Q$, or $\R$. Denote by ${\Bbb A}^t$ 
the corresponding tropical  semifield. Then a positive 
rational map $\varphi: H_1 \to H_2$ gives rise to a {\it piecewise-linear} map  
$\varphi_*: H_1({\Bbb A}^t) \to H_2({\Bbb A}^t)$.

Let ${\cal X}$ be a positive space. 
Then we have the manifold ${\cal X}(\R_{>0})$, and the ${\Bbb A}$-{\it tropicalization} 
${\cal X}({\Bbb A}^t)$  
of ${\cal X}$. They are related as follows. First of all, 
one has 
$
{\cal X}({\Z^t}) \hra {\cal X} ({\Q^t}) \hra {\cal X} ({\R^t}) 
$. 
The multiplicative group ${\Bbb A}^*_{>0}$ of positive elements in 
${\Bbb A}$ acts by automorphisms of the tropical semifield ${\Bbb A}^t$. 
Let us denote by $\{0\}$ the zero element in the abelian group 
$H({\Bbb A}^t)$. Let ${\Bbb A}$ be either $\Q$ or $\R$. Set 
$$
{\Bbb P}H({\Bbb A}^t) := \Bigl(H({\Bbb A}^t)- \{\{0\}\}\Bigr)/{\Bbb A}^*_{>0}.
$$
Let $\varphi: H \to H'$  be 
a positive rational map. Then  
the piecewise-linear map $\varphi_*: H({\Bbb A}^t) \to H'({\Bbb A}^t)$ 
commutes with the multiplication by ${\Bbb A}^*_{>0}$, and  
$\varphi_*\{0\} = \{0\}$. 
 Thus we get a morphism
$$
\overline \varphi_*: {\Bbb P}H({\Bbb A}^t) \stackrel{}{\lra}
{\Bbb P}H'({\Bbb A}^t).
$$

\begin{definition} \label{8.28.03.10q}
Let ${\cal X}$ be a positive space. Let ${\Bbb A}$ be either $\Q$ or $\R$. 
 The projectivisation ${\Bbb P}{\cal X} ({\Bbb A}^t)$ of the ${\Bbb A}$-tropicalisation 
of ${\cal X}$ is 
$$
{\Bbb P}{\cal X}({\Bbb A}^t) := \{\{p_{\alpha} \in {\Bbb P}H_{\alpha}({\Bbb A}^t)\} \quad | \quad \mbox{ 
$\overline \psi_{\alpha, \beta}(p_{\alpha}) = 
p_{\beta}$ for any 
$\alpha, \beta \in {\cal G}_X$}\}.
$$ 
\end{definition}

Observe that ${\Bbb P}H_{\alpha}({\R}^t)$ is a sphere, and the maps 
$\psi_{\alpha, \beta}$ are
piecewise-linear maps of spheres. The set 
${\Bbb P}{\cal X}({\Q}^t)$ is an everywhere  dense subset of 
${\Bbb P}{\cal X}({\R}^t)$. 

The degree map from Lemma \ref{8.27.03.13} provides a  map 
$$
-{\rm deg}: {\cal X}(\R((\varepsilon))_{>0}) \lra {\cal X}(\Z^t).
$$
An $\R((\varepsilon))_{>0}$-point $x(\varepsilon)$ of ${\cal X}$, 
considered modulo action of the 
group of positive reparametrizations of the formal $\varepsilon$-line, 
i.e. modulo substitutions $\varepsilon = \phi(\varepsilon)$, where 
$\phi(\varepsilon) \in \R((\varepsilon))_{>0}$, gives rise to 
a well defined element of ${\Bbb P}{\cal X}({\Q^t})$, 
provided  that ${\rm deg}~x(\varepsilon) \not = \{0\}$. 
Thus ${\Bbb P}{\cal X}({\R^t})$ 
serves as the Thurston type boundary of ${\cal X}({\R_{>0}})$.


\section{Positive configurations of flags}
\label{posflag}


The first two subsections contain background material 
about positivity in algebraic semi-simple groups, due to Lusztig, and
Fomin and Zelevinsky. The rest of the Section is devoted to 
a definition and investigation of  positive structures 
on configurations of three and four generic flags in $G$.  

\vskip 3mm
We use throughout 
the paper the following conventions. The {\it moduli spaces}, that is stacks, 
 of configurations of flags are denoted by 
${\rm Conf}({\cal B})$. The corresponding {\it varieties} of 
flags in generic position are denoted by ${\rm Conf}^*({\cal B})$. 
If we care only birational structure of the configuration space, we use the simpler notation ${\rm Conf}$. 
We use ${\rm Conf}^*({\cal B})$ only when we need it as a variety. 
Finally, ${\rm Conf}^+({\cal B})$ denotes the set of $K$-positive configurations, 
which is a subset of the $K$-points of 
${\rm Conf}^*({\cal B})$. 
Similar conventions are used for the configurations of affine flags.

\vskip 3mm
{\bf 1. A positive regular structure on  $U$ and $G$}. 
We recall some basic results of G. Lusztig \cite{L1}, \cite{L2}. 
The only difference is that we present them 
as results about existence and properties of  
regular positive structures on   $U$ and  $G$.  
The proofs in loc.cit.  give exactly what we  need. 
However, as explained in Section 5.8, one needs to be careful when 
relating this point of view with the one used by Lusztig. 

Let $W$ be the Weyl group of $G$. 
Let ${s}_i \in W$ be the simple reflection corresponding to  $x_i, y_i$. 
Let $w_0$ be the maximal length element in the Weyl group $W$, and  
$w_0= s_{i_1} ... s_{i_n}$ its reduced expression into a product 
of simple reflections. It is encoded by the sequence ${\bf i}:= \{i_1, ..., i_n\}$ of 
simple reflections in the decomposition of $w_0$.

\begin{proposition} \label{7.22.03.1}
a)Each reduced expression of $w_0$, encoded by ${\bf i}$,   provides an open regular 
embedding 
\begin{equation} \label{7.22.03.2}
\psi_{\bf i}: {\Bbb G}_m^n \hra U, \qquad (a_1, ..., a_n) 
\lms x_{i_1}(a_1) ... x_{i_n}(a_n).
\end{equation}

b) The collection of embeddings $\psi_{\bf i}$, when ${\bf i} $ run through the set of all reduced 
expression of $w_0$, provides a positive structure on $U$. 

c) Replacing $x_{i_p}$ by $y_{i_p}$  in (\ref{7.22.03.2}) 
we get an open regular 
embedding  $\psi_{\bf \overline i}: {\Bbb G}_m^n \hra U^-$, providing 
a positive regular structure 
on $U^-$. The map $\Psi: U^- \to U$ is a positive regular isomorphism. 
\end{proposition}

{\bf Proof}. Proposition 2.7 in \cite{L1} and the proof of its part a) 
prove the proposition.  
\vskip 3mm

{\bf 2. Positive regular structures on double Bruhat cells {\rm \cite{FZ99}}}.  
 The double Bruhat cell corresponding to an element $(u,v) \in W \times W$ 
is defined by 
\begin{equation} \label{9.19.03.20}
G^{u, v} := BuB \cap B^-vB^- .
\end{equation} 
A {\it double reduced word} for the element $(u, v) \in W \times W$ is a 
reduced word 
for an element $(u, v)$ of the Coxeter group $W \times W$. 
We will use the indices $\overline 1, \overline 2, ..., \overline r$ 
for the simple reflections 
in the first copy of $W$, and $1, 2, ...,  r$ for the second copy. 
Then a double reduced word ${\bf j}$ for $(u, v)$ is simply a shuffle of 
a reduced word for $u$ written in the alphabet $[\overline 1, 
\overline 2, ..., \overline r]$, 
and a reduced word for $v$ written in the alphabet $[1, 2, ...,  r]$. 
Recall that $l(u)$ is the 
length of an element $u \in W$. We set $x_{\overline {i}}(t):= 
y_{{i}}(t)$. The following result is  
Theorem 1.2 in \cite{FZ99}. 
\begin{proposition} \label{8.28.03.1}
a) \!Given a\! double reduced word ${\bf j} = (j_1, ..., j_{l(u) + l(v)})$ 
for $(u, v)$ 
and $0 \leq k \leq l(u) + l(v)$ 
there is a regular open embedding
$$
x_{{\bf j}, k}: H \times {\Bbb G}_m^{l(u) + l(v)} \hra G^{u, v} 
$$
$$
(h; t_1, ..., t_{l(u) + l(v)}) \lms x_{i_1}(t_1) ...x_{i_{k}}(t_{i_{k}}) 
h x_{i_{k+1}}(t_{i_{k+1}})  ... x_{i_{l(u) + l(v)}}(t_{l(u) + l(v)}). 
$$

b) The collection of all  regular open embeddings $x_{{\bf j}, k}$ 
provide a regular positive structure on the double Bruhat cell 
$G^{u, v}$. 
\end{proposition}
Indeed, if $k=0$ this is exactly what stated in loc. cit. Moving 
$h$ we do not change the positive structure since the characters of $H$ 
are positive 
by definition. It follows from Proposition \ref{8.28.03.1} 
that ${\rm dim} G^{u,v} = r + l(u) + l(v)$. 
\vskip 3mm
{\it The generalized minors and double Bruhat cells {\rm \cite{FZ99}}}. 
Let ${\bf i} = (i_1, ..., i_n)$ be a reduced decomposition of $w_0$. Set 
$$
\overline {\overline s}_i:= \varphi_i
\left (\begin{matrix}0& 1\\ -1&0\end{matrix}\right ), \qquad \overline {\overline w}_0 := 
\overline {\overline s}_{i_1} ... \overline {\overline s}_{i_n}. 
$$
Then $\overline {\overline w}_0$ does not depend on choice of a reduced decomposition. 
Using both reduced decompositions $(i_1, ..., i_n)$ and $(i_n, ..., i_1)$ one 
immediately sees that 
$$
\overline {\overline w}_0 {\overline w}_0 = e, \quad \overline {\overline w}_0  = s_G 
{\overline w}_0, 
\quad \Psi(\overline {\overline w}_0) 
 = {\overline w}_0, \quad \Psi({\overline w}_0) 
 = \overline {\overline w}_0.
$$
Let $G_0 := U^-HU$. Any element $x \in G_0$ admits unique decomposition 
$x = [x]_- [x]_0 [x]_+$ with $[x]_- \in U_-$, $[x]_0 \in H$, $[x]_+ \in U_+$.  
Let $i \in [1, ..., r]$. Then according to \cite{FZ99},
the  generalized minor
$\Delta_{u\omega_i, v\omega_i}$ is the regular function on $G$ 
whose restriction to the open subset $\overline u G_0 \overline v^{-1}$ is given by 
\begin{equation} \label{9.19.03.1yt}
\Delta_{u\omega_i, v\omega_i}(x) = 
\omega_i([ {\overline u} ^{-1} x \overline v]_0).
\end{equation}
It was shown in loc. cit. that it depends only on the weights 
$u\omega_i, v\omega_j$. For example if $G = SL_{r+1}$ then for 
permutations $u,v \in S_{r+1}$ the function 
$\Delta_{u\omega_i, v\omega_i}$ is given by the $i \times i$ minor 
whose rows (columns) are labeled by the elements of the set $u([1, ..., i])$ 
(resp., $v([1, ..., i])$). 
Let ${\bf i}$ be a reduced decomposition of $w_0$ and 
$k \in [1, l(u) + l(v)]$. We set $\varepsilon(i_l) = -1$ if $l$ is from the alphabet $\overline 1, \ldots, \overline r$, and $+1$ otherwise. Following (\cite{BFZ03}, 
Section 2.3), we set 
\begin{equation} \label{9.19.03.1}
u_{\leq k} = u_{\leq k}({\bf i}) = \prod_{l = 1, ..., k, \quad \varepsilon(i_l) 
= -1}s_{|i_l|}
\end{equation}
\begin{equation} \label{9.19.03.2}
v_{>k} = v_{>k}({\bf i}) = \prod_{l = l(u) + l(v), ..., k+1, \quad \varepsilon(i_l) = +1}s_{|i_l|}
\end{equation}
where the index $l$ is increasing (decreasing) in (\ref{9.19.03.1}) (resp., 
(\ref{9.19.03.2})). 
For $k \in -[1, ..., r]$ we set 
$u_{\leq k} = e$ and $v_{>k} = v^{-1}$. Finally, for 
$k \in -[1,r] \cup [1, l(u) + l(v)]$, set
\begin{equation} \label{bfz03}
\Delta(k; {\bf i}):= \Delta_{u_{\leq k}\omega_{|i_k|}, v_{>k}\omega_{|i_k|}}.
\end{equation}
Let us set
$$
F({\bf i}) = \{\Delta(k; {\bf i}): k \in -[1,r] \cup [1, l(u) + l(v)]\}.
$$

Let 
$$
G^{u,v}({\bf i}):= \{g \in G^{u,v}: \Delta(g) \not = 0 \mbox{ for all } 
\Delta \in F({\bf i})\}.
$$
The lemma below is a part of Lemma 2.12 in \cite{BFZ03}, 
and follows directly from \cite{FZ99}. 

\begin{lemma} \label{9.19.03.3}
The map 
$
G^{u,v} \to {\Bbb A}^{r + l(u) + l(v)}$ given by $g \mapsto 
(\Delta(g))_{\Delta \in F({\bf i})}
$ restricts to an isomorphism
\begin{equation} \label{9.19.03.4}
G^{u,v}({\bf i}) \stackrel{\sim}{\lra} {\Bbb G}_m^{r + l(u) + l(v)}.
\end{equation}
\end{lemma}

Let 
\begin{equation} \label{9.19.03.41q}
\Delta_{{\bf i}}: {\Bbb G}_m^{r + l(u) + l(v)} {\hra} G^{u,v}
\end{equation}
be the regular open embedding provided by Lemma \ref{9.19.03.3}. 
The theorem below is a reformulation of the results in \cite{FZ99}. 
\begin{theorem} \label{9.19.03.41}
The collection of 
embeddings $\{\Delta_{{\bf i}}\}$, when ${\bf i}$ runs \newline through the set of all 
reduced decompositions of $(u,v)$,  provides a positive regular structure 
on the double Bruhat cell $G^{u,v}$. It is compatible with the 
positive regular structure 
defined in Proposition \ref{8.28.03.1}.
\end{theorem}

\vskip 3mm
{\bf 3.  A positive regular structure on the intersection of two 
opposite open Schubert cells in ${\cal B}$}. 
Recall the two opposite Borel subgroups $B^-, B^+$. They determine 
the two open Schubert  cells ${\cal B}_-$ and ${\cal B}_+$
 in the flag variety ${\cal B}$, parametrising
the Borel subgroups of $G$ in generic position to $B^-$ and $B^+$. One has 
\begin{equation} \label{7.22.03.15}
{\cal B}_-:= \{ u_- B^+ u_-^{-1}\},\quad u_- \in U^-; \qquad {\cal B}_+:= \{u_+ B^- u_+^{-1} \},\quad u_+ \in U^+.
\end{equation}
 The formula  (\ref{7.22.03.15}) provides  isomorphisms
$
\alpha^-: U^- \stackrel{\sim}{\lra} {\cal B}_-, \quad 
\alpha^+: U^+ \stackrel{\sim}{\lra} {\cal B}_+
$. 
 Let us set 
$$
{\cal B}_{*}:= {\cal B}_-  \cap {\cal B}_+, \quad 
U^+_*:= U^+ \cap B^-w_0B^-, \quad U_*^-:= U^- \cap B w_0 B.
$$
\begin{lemma} \label{8.27.03.10}
a) There are canonical isomorphisms
$$
\alpha_*^+: U_*^+ \stackrel{\sim}{\lra} {\cal B}_*, \quad u_+ \lms u_+B^-u_+^{-1}, 
\qquad 
\alpha_*^-: U_*^- \stackrel{\sim}{\lra} {\cal B}_*, \quad u_- \lms u_-B^+u_-^{-1}.
$$

b) $\Psi: U^-_* \lra U^+_*$ is an isomorphism. 
\end{lemma}

{\bf Proof}. a) The map $\alpha_*^+$ is evidently injective. 
Let us check the surjectivity. Given a pair $(u_+, u_-)$ as in (\ref{7.22.03.15}), 
such that 
$u_+B^-u_+^{-1} = u_-B^+u_-^{-1}$, let us show that $u_-w_0B^- \cap U^+$ is 
non empty, and hence is a point. Indeed, conjugating $B^-$ by $u_-w_0B^-$ we get 
$u_-B^+u_-^{-1}$. By the assumption conjugation of $B^-$ by $u_+$ gives the same result. 
Thus $u_+ \in u_- w_0B^-$. The part b) is obvious. The lemma is proved. 
\vskip 3mm

Lemma \ref{8.27.03.10} implies that  there is an isomorphism
$$
\phi: U_*^- \stackrel{\sim}{\lra}  U_*^+, \qquad  \phi:= (\alpha_*^+)^{-1}\alpha_*^-.
$$

\begin{proposition} \label{8.27.03.1}
a) For every reduced decomposition  ${\bf i}$ of $w_0$  
there are regular open embeddings
\begin{equation} \label{8.23.03.02}
\alpha^+_* \circ \psi_{\bf i}: ({\Bbb G}_m)^n \hra {\cal B}_+; \qquad 
\alpha^-_* \circ \psi_{\bf \overline i}: ({\Bbb G}_m)^n \hra {\cal B}_-.
 \end{equation}

b) The maps $\{\alpha^+_* \circ \psi_{\bf i}\}$
  provide a positive regular structure on ${\cal B}_*$. Similarly 
the maps $\{\alpha^-_* \circ \psi_{\bf \overline i}\}$ 
provide another positive regular structure on ${\cal B}_*$. 

c) The map $\phi$ is positive. Therefore 
the two positive structures from the part b) are compatible, and thus provide 
a positive regular structure on ${\cal B}_*$ dominating them. 
\end{proposition}

{\bf Proof}. a) According to  Theorem 1.2 in \cite{FZ99} there is 
a regular open embedding 
$$x_{\bf i}:H \times {\Bbb G}_m^n \to 
G^{e, w_0}:= B^+ \cap B^- w_0 B^-, \qquad 
x_{\bf i}(h; t_1, ..., t_n):= h 
x_{i_1}(t_1) \cdot ... \cdot x_{i_n}(t_n).
$$
Putting $h=e$ we land in $U^+ \cap B^-w_0B^-$. 

b) Follows from Proposition \ref{7.22.03.1}.

c) In s. 3.2 of \cite{L2} there is 
an algorithm for computation of a map $\phi': U^- \lra U^+$, such that 
\begin{equation} \label{8.27.03.111}
uB^+u^{-1} = \phi'(u) B^- \phi'(u)^{-1}, \quad u \in U^-.
\end{equation}
Such a map is unique. It follows from loc. cit. that it is a positive rational map
for the positive structures on $U^+$ and $U^-$ given by Proposition \ref{7.22.03.1}. 
Its restriction to $U^-_*$ is our $\phi$.  The Proposition is proved. 
\vskip 3mm

{\bf 4. A positive regular variety ${\cal V}$}.  
Let ${\cal V}$ be the quotient of the affine variety 
$U_*$ by the action of $H$:
\begin{equation} \label{7.27.03.1}
{\cal V}:= U_*/H:= {\rm Spec}\Q[U_*]^H.
\end{equation}
Let us define a coordinate system on ${\cal V}$ 
corresponding to a reduced decomposition  ${\bf i}$ of $w_0$. 
Let $\alpha_1, ..., \alpha_r$ be the set of all 
 simple roots.  Each $i_k \in {\bf i}$ is labeled by a simple root,  
 providing a decomposition 
$$
{\bf i} = {\bf i}(\alpha_1) \cup ... \cup {\bf i}(\alpha_r). 
$$
Precisely,  $i_j \in {\bf i}(\alpha_p)$ if and only if $i_{j} = p\in I$. 
Let $l_k:= |{\bf i}(\alpha_k)|$. So there is a positive regular open embedding 
$$
\psi_{\bf i}: {\Bbb G}_m^{{\bf i}(\alpha_1)} \times ... 
\times {\Bbb G}_m^{{\bf i}(\alpha_r)} \hra U_*.
$$
The group  $H$  acts on the factor 
${\Bbb G}_m^{{\bf i}(\alpha_k)}$ through its character $\chi'_k$. 
Since the center of $G$ is trivial, there is an isomorphism
$
(\chi'_1, ..., \chi'_r): H \stackrel{\sim}{\lra}
 {\Bbb G}_m \times ... \times {\Bbb G}_m
$. 
 So $H$ acts freely on the product. The map 
$\psi_{\bf i}$ transforms this action 
to the action of $H$ on $U_*$ by conjugation. 
 Thus we get a  regular open embedding
$$
\overline \psi_{\bf i}: {\Bbb G}_m^{{\bf i}(\alpha_1)}/{\Bbb G}_m 
\times ... \times {\Bbb G}_m^{{\bf i}(\alpha_r)}/{\Bbb G}_m  \hra U_*/H
$$
where the group ${\Bbb G}_m$ acts on ${\Bbb G}_m^{{\bf i}(\alpha_k)}$ diagonally. 
Let $a_{\alpha_k, 1}, ..., a_{\alpha_k, l_k}$ be the natural coordinates 
on ${\Bbb G}_m^{{\bf i}(\alpha_k)}$. Then the natural coordinates 
on its quotient by the diagonal action of ${\Bbb G}_m$ are 
\begin{equation} \label{7.27.03.5}
t_{\alpha_k, j}:= \frac{a_{\alpha_k, j+1}}{a_{\alpha_k, j}}, \quad 1 \leq j \leq l_k-1.
\end{equation} 

\begin{lemma} \label{7.27.03.12} The collection of regular open embeddings 
$\{\overline \psi_{\bf i}\}$, where ${\bf i}$ run through the set of 
all reduced expressions of $w_0$,  
provides a  regular positive structure on  ${\cal V}$. 
\end{lemma}

{\bf Proof}. Follows immediately from the very definition 
 and 
Proposition \ref{7.22.03.1}. 
\vskip 3mm

Similarly, replacing $U$ by $U^-$ we get a positive variety ${\cal V}^-$. 
It is handy to denote ${\cal V}$ by ${\cal V}^+$. The 
map $\Psi$ evidently commutes with the action of  $H$. 
It follows from the construction of the map $\phi$ given in s. 3.2 of [L2] that 
 $\phi$ also commutes with the action of  $H$. 
Thus they provide isomorphisms 
\begin{equation} \label{7.27.03.15}
 \phi: {\cal V}^- \lra {\cal V}^+, \quad 
\Psi: {\cal V}^- \lra {\cal V}^+.
\end{equation} 
Abusing the notation we denote them the same way as the original 
maps. 
\begin{lemma} \label{7.27.03.13} The maps $ \phi$ and 
$ \Psi$ in (\ref{7.27.03.15}) 
are isomorphisms of positive regular varieties.
\end{lemma}

{\bf Proof}. For $ \Psi$ this is clear. For $\phi$ this follows from 
Proposition \ref{8.27.03.1}. The lemma is proved. 
\vskip 3mm

{\bf 5. A positive structure on configurations of triples of flags}. 
We say that a pair of flags is in {\it generic position} if they 
are conjugated to the pair $(B^+, B^-)$, or, equivalently, 
project to the open $G$-orbit in ${\cal B}\times {\cal B}$. 
Recall that configurations of $n$--tuples 
of flags are by definition the $G$--orbits 
on ${\cal B}^n$. Denote by $\sim$ the equivalence relation provided by  $G$--orbits. 
We say that a configuration of flags is {\it generic} if 
every two flags of the configuration are in generic position. 
The moduli space ${\rm Conf}^*_n({\cal B})$ of generic configurations 
of $n>2$ flags is a variety over $\Q$. 
Indeed, a generic configuration can be presented as
$$
(B^-, B^+, u_1B^+u_1^{-1}, ..., u_kB^+u_k^{-1}), \quad u_k \in U^-_*.
$$
Thus they are parametrised by a Zariski open part of 
 $(U_*^- \times ... \times U_*^-)/H$, where there are $k$ factors in the product. 
The variety ${\rm Conf}^*_n({\cal B})$ is a Zariski open part of the stack ${\rm Conf}_n({\cal B})$.
So the two objects are  birationally isomorpic.

Let us define a positive structure on the configuration space 
${\rm Conf}_3({\cal B})$. 
Let ${\rm pr}: U_*^{\pm} \to {\cal V}^{\pm}$ be the canonical projection. 
Consider a regular open embedding
\begin{equation} \label{7.28.03.10}
{\cal V}^- \hra {\rm Conf}_3({\cal B}), \  v \lms (B^-, B^+, vB^+v^{-1}):= 
(B^-, B^+, uB^+u^{-1}), \  {\rm pr}(u) = v.
\end{equation}
 Since $H$ stabilizes the pair $(B^-, B^+)$, this map 
does not depend on a choice of $u \in U_*^-$ lifting  $v\in {\cal V}^-$.  
\begin{definition} \label{2.3.03.2}
A positive structure on ${\rm Conf}_3({\cal B})$ 
is defined by the map (\ref{7.28.03.10}) and Lemma \ref{7.27.03.12}. 
\end{definition}
According to Lemma  \ref{7.27.03.13} the positive structure is compatible with the one given by 
\begin{equation} \label{7.28.03.10sa}
{\cal V}^+ \hra {\rm Conf}_3({\cal B}),\  v \lms (B^-, B^+, vB^-v^{-1}):= 
(B^-, B^+, uB^-u^{-1}), \quad {\rm pr}(u) = v.
\end{equation}
Yet another natural positive regular structure on ${\rm Conf}_3({\cal B})$ 
is given by a regular open embedding
\begin{equation} \label{7.28.03.101}
{\cal V}^+ \hra {\rm Conf}_3({\cal B}), \qquad v \lms (B^-, v^{-1}B^-v, B^+).
\end{equation}
Recall that the map $\Psi$ acts on the flag variety ${\cal B}$ as well as on the configurations of flags.

\begin{theorem}\label{2.18.03.2a} Let  $G$ be a split reductive group  
with trivial center. Consider a positive atlas on ${\rm Conf}_3({\cal B})$ 
provided by (\ref{7.28.03.10}). Then  one has:

a) The cyclic shift $(B_1, B_2 , B_3) \lms (B_2, B_3, B_1)$  
is a positive isomorphism on ${\rm Conf}_3({\cal B})$. Precisely:
\begin{equation} \label{8.11.03.1}
(B^-, B^+, uB^-u^{-1}) \sim 
(B^+, cB^-c^{-1}, B^-), \quad u \in {U}_*^+, \quad c:= (\phi\Psi)^2(u). 
\end{equation}

b) $\Psi: (B_1, B_2 , B_3) \lms 
(\Psi B_1, \Psi B_2, \Psi B_3)$ provides a positive  isomorphism  on ${\rm Conf}_3({\cal B})$.

c) The two positive regular structures on ${\rm Conf}_3({\cal B})$ given by 
(\ref{7.28.03.10}) and (\ref{7.28.03.101}) are compatible. 

d) The reversing map $(B_1, B_2, B_3) \lms (B_3, B_2, B_1)$ is a
positive isomorphism on ${\rm Conf}_3({\cal B})$.

\noindent 
So the positive structure  on ${\rm Conf}_3({\cal B})$ is invariant 
under the action of the symmetric group $S_3$. 
\end{theorem}

{\bf Proof}. a) 
Consider a configuration 
\begin{equation} \label{2.03..3.5}
(B^-, B^+, u_+B^-u_+^{-1}), \qquad {\rm pr}(u_+) \in {\cal V}_+.
\end{equation}
Let us  compute the composition  $\Psi \circ \Psi$ applied to 
this configuration. Conjugating  (\ref{2.03..3.5}) by $u_+^{-1}$ we get    
$(u_+^{-1}B^-u_+, B^+, B^-)$. It is, of course, equivalent to the initial 
configuration (\ref{2.03..3.5}). 
Then applying $\Psi$ we get 
\begin{equation} \label{2.03..3.10}
(u_-B^+ u_-^{-1}, B^-, B^+), \qquad u_- = \Psi(u_+) \in U^-_*, \quad u_+ \in U^+_*.
\end{equation}
According to (\ref{8.27.03.111}) it  can be written as 
$$
(v_+B^- v_+^{-1}, B^-, B^+)\quad \mbox{for some $v_+ \in U^+_{*}$}, \quad v_+ = \phi\Psi(u_+).
$$
Conjugating by $v_+^{-1}$ we 
get $(B^-, v_+^{-1}B^-v_+, B^+)$. Applying $\Psi$ again we obtain
$$
(B^+, v_-B^+v_-^{-1}, B^-), \qquad v_-:= \Psi(v_+) \in U^-_{*}.
$$
We can write it as 
$$
(B^+, c_+B^-c_+^{-1}, B^-), \qquad c_+:= (\phi\Psi)^2(u_+) \in U^+_{*}.
$$
This configuration is the cyclic shift of a configuration (\ref{2.03..3.5}). On the other hand, 
since $\Psi^2={\rm Id}$,  it coincides 
with (\ref{2.03..3.5}). It remains to use positivity of the maps
$\phi$ (Proposition \ref{8.27.03.1}c)) and $\Psi$. The statement a) is proved.

b) By (\ref{8.27.03.111}) and as $\phi$ is positive (or by using (\ref{7.28.03.10sa})), the configuration (\ref{2.03..3.10}) is a 
cyclic shift of the 
configuration (\ref{2.03..3.5}). 
On the other hand, as the proof of a) shows,  it is 
equivalent to $\Psi(B^-, B^+, u_+B^-u_+^{-1})$. The cyclic shift is positive by a). 
The statement b) is proved. 

c) Conjugating (\ref{2.03..3.5}) by $u_+^{-1}$ we get    
$(u_+^{-1}B^-u_+, B^+, B^-)$. Applying the cyclic shift twice and using a) we get c).

d) By the parts c) and a) the following equivalences provide positive
maps: 
$$
(B^-, B^+, u_+B^-u_+^{-1}) \sim (B^-, v^{-1}_+B^-v_+, B^+) \sim (v^{-1}_+B^-v_+, B^+, B^-). 
$$
Thus it remains to check that map 
$
(B^+, B^-, v^{-1}_+B^-v_+) \mapsto (B^+, B^-, v_+B^-v^{-1}_+) 
$ 
induced by the inversion $v \lms v^{-1}$
is a positive map on ${\rm Conf}_3({\cal B})$. 
Although the inversion map is not positive on $U$, 
it induces a positive map on the quotient $U_*/H$. Indeed, 
inversion changes the order of factors in a reduced decomposition,
which amounts to changing the reduced decomposition, and
replaces each factor $x_i(t)$ by $x_i(-t)$. However this does not
affect the ratios 
$t_i/t_j$ used to define the positive atlas on $U_*/H$, and changing a reduced
decomposition   we get a positive map. The part d) is proved. 
 
The last assertion of the theorem is a corollary of a) and d). 
The theorem is proved.
\vskip 3mm
{\bf 6. A positive structure on ${\rm Conf}_4({\cal B})$}. Consider a 
triangulation of a convex $4$-gon 
with the flags $B_i$ assigned to its vertices, as on the picture. 
Choose in addition an orientation of the internal edge, from $B_1$
to $B_3$. 
Such a data, denoted by $\overline T$, 
 is determined by an ordered configuration of $4$ flags $(B_1, B_2,
B_3, B_4)$. 
 
\begin{figure}[ht]
\centerline{\epsfbox{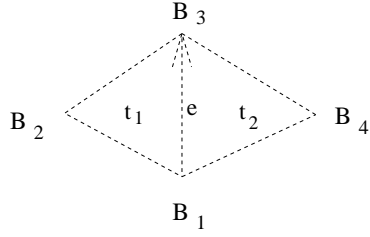}}
\caption{A triangulated quadrilateral with flags assigned to vertices.}
\label{fgo2}
\end{figure}

We will employ the following notation for the group action by conjugation 
on the flag variety: 
$$
u\cdot B:= uBu^{-1}.
$$
\begin{lemma} \label{9.9.03.1}
The map 
\begin{equation} \label{2.03..3.4}
(u_+, u_-)\lms 
(B^-,  u_+^{-1}\cdot B^-,  B^+,  u_-\cdot B^+).
\end{equation}
provides a regular open embedding 
\begin{equation} \label{2.03..3.4q}
i_{\overline T}: {\rm Conf}^*_4({\cal B}) \hra ({U}_*^+ \times {U}_*^-)/H. 
\end{equation}
\end{lemma}

{\bf Proof}. Let us identify  $(B_1,
B_3)$ with the standard pair $(B^-, B^+)$ and choose one of the  pinnings 
related to  $(B^-, B^+)$. 
Then the element 
\begin{equation} \label{2.03..3.4qq}
(B_1, B_2, B_3, B_4) \in 
{\rm Conf}_4({\cal B})
\end{equation}
 can be written  as the right hand side of 
(\ref{2.03..3.4}). The action of the group $H$ by conjugation  on 
the pinnings with the given pair 
$(B_1, B_3)$ is a simple
transitive action.  So the element (\ref{2.03..3.4qq}) 
is described uniquely by an element of $(U_*^- \times
U_*^+)/H$.  The lemma is proved. 
\vskip 3mm

Just like in Section 5.4, there is a positive regular structure on $({U}_*^+ \times {U}_*^-)/H$. Precisely, take a  pair of reduced words 
$({\bf i_1}, {\bf i_2})$  for $w_0$. Using ${\bf i_1}$ to define 
a positive coordinate system on ${U}_*^+$ and ${\bf i_2}$ for ${U}_*^-$, we 
  get a regular open embedding
$$
\overline \psi_{\bf (i_1, i_2)}: {\Bbb G}_m^{{\bf i_1}(\alpha_1) + 
{\bf i_2}(\alpha_1)}/{\Bbb G}_m \times ... 
\times {\Bbb G}_m^{{\bf i_1}(\alpha_r) + 
{\bf i_2}(\alpha_r)}/{\Bbb G}_m \hra (U_*^- \times
U_*^+)/H.
$$
Similar to (\ref{7.27.03.5}), we can 
introduce homogeneous coordinates on each of the factors 
${\Bbb G}_m^{{\bf i_1}(\alpha_k) + 
{\bf i_2}(\alpha_k)}/{\Bbb G}_m$. Introducing the homogeneous coordinates one must take into account the fact that ${\Bbb G}_m$ acts on 
the factor ${\Bbb G}_m^{{\bf i_1}(\alpha_k)}$ through the character 
$\chi_k'$, while on 
the factor ${\Bbb G}_m^{{\bf i_2}(\alpha_k)}$ via its inverse.

\begin{theorem} \label{2.18.03.2wq}
a) The rational map $$
i^{-1}_{\overline T}: ({U}_*^+ \times {U}_*^-)/H \lra {\rm Conf}_4({\cal B})
$$ transforms 
the positive structure on $({U}_*^+ \times {U}_*^-)/H$ to  a positive structure 
on the moduli space ${\rm Conf}_4({\cal B})$. 

b) The cyclic shift on ${\rm Conf}_4({\cal B})$ is a positive rational
map. 

c) The reversion $(B_1, B_2, B_3, B_4)\lms (B_4, B_3, B_2, B_1) $ is a positive rational
map on ${\rm Conf}_4({\cal B})$. 

\noindent
So the positive structure on ${\rm Conf}_4({\cal B})$ has a 
natural dihedral symmetry. 
\end{theorem}

{\bf Proof}. a) The map $i_T$ is not an isomorphism.  
 Its image consists of the elements $(u_+, u_-)$ 
such that the corresponding pair of flags $$(u_+^{-1}B^-u_+, u_-B^+u_-^{-1})$$
is in generic position. So to  prove a) it remains to 
show that the complement to the image of $i_T$  is contained in a positive divisor. 

Let $\overline T'$ be another triangulation of the $4$-gon on Figure \ref{fgo2}, 
obtained from $\overline T$ by changing the oriented edge $B_1 \to B_3$ to the one 
$B_2 \to B_4$. It provides  
another  regular open embedding 
\begin{equation} \label{2.03..3.4qas}
i_{\overline T'}: {\rm Conf}^*_4({\cal B}) \hra ({U}_*^+ \times {U}_*^-)/H. 
\end{equation}
It is induced by $(u_+, u_-)\lms (u_-\cdot B^+, B^-, u_+^{-1}\cdot B^-, B^+)$. 
We need the following proposition. 

\begin{proposition} \label{9.11.03.1} The rational map 
$$
i_{\overline T}\circ i_{\overline T'}^{-1}: ({U}_*^+ \times {U}_*^-)/H  \hra ({U}_*^+ \times {U}_*^-)/H 
$$
as well as its inverse are positive rational maps. 
\end{proposition}

{\bf Proof}. Applying Proposition \ref{8.28.03.1} to 
 ${\bf j_1}:= ({\bf i}, {\bf \overline i})$, $k=n$  and 
${\bf j_2}:= ({\bf \overline i}, {\bf i})$, $k=n$ we get   
\begin{lemma} \label{8.28.03.2}
The map 
$$
\eta: U^+_* \times H \times U^-_* \lra U^-_* \times H \times U^+_*, \  
\eta: (u_+, h, u_-) \lms (a_-, t, a_+),
$$ 
where $u_-hu_+ = a_+ t a_-$ is a positive rational map.
\end{lemma} 
Let $(u_+, u_-)$ be as in (\ref{2.03..3.4}). Then  thanks to Lemma 
\ref{8.28.03.2} there exists a positive rational map 
$$
(u_+, u_-) \lra (a_-, t, a_+), \quad a_{\pm} \in U_*^{\pm}, t \in H\quad 
\mbox{such that $u_+u_- =  a_-a_+t $}.
$$
  Set 
$g := a_-^{-1}u_+ = a_+tu_-^{-1}$. Conjugating (\ref{2.03..3.4}) 
by $g$ we get 
$$
(a_+\cdot B^-,  B^-,  a_-^{-1}\cdot B^+,  B^+) = 
(\phi(a_+)\cdot B^+,  B^-,  \phi(a_-)^{-1}\cdot B^-,  B^+).
$$
Since $\phi$ is a positive isomorphism, this proves that the rational map $i_{\overline T}\circ i_{\overline T'}^{-1}: (u_+, u_-) \lms (\phi(a_-), \phi(a_+))$ is positive. 
Applying this argument three times we show that 
the rational map induced by the cube of the cyclic shift is positive. 
This is equivalent to positivity of the inverse of $i_{\overline T}\circ i_{\overline T'}^{-1}$. 
Proposition \ref{9.11.03.1} is proved. 
\vskip 3mm
We will finish the proof of Theorem \ref{2.18.03.2wq} in the end of the next subsection. 

\vskip 3mm
{\bf 7.  Another description of ${\rm Conf}_4({\cal B})$}. Recall the  
triangulation $T$ of a convex $4$-gon 
with the flags $B_i$ assigned to its vertices as on the picture 
and the internal edge  oriented from $B_1$
to $B_3$. 
\begin{proposition} \label{2.18.03.4} There exists a canonical birational isomorphism
\begin{equation} \label{2.18.03.4123}
\varphi_T: {\rm Conf}_4({\cal B})  \stackrel{\sim}{\lra} {\rm Conf}_3({\cal B})
\times H \times {\rm Conf}_3({\cal B}) \stackrel{\sim}{\lra} {\cal V}^+ \times H \times {\cal V}^-. 
\end{equation} 
\end{proposition}

{\bf Proof}. 
Let $(B_1, B_2, B_3, B_4) \in {\rm Conf}_4({\cal B})$. 
Then there are projections 
$$
p_1: {\rm Conf}_4({\cal B}) \lra {\rm Conf}_3({\cal B}), \quad (B_1, B_2, B_3, B_4) 
\lms (B_1, B_2, B_3); 
$$
$$
p_2: {\rm Conf}_4({\cal B}) \lra {\rm Conf}_3({\cal B}), \quad (B_1, B_2, B_3, B_4) 
\lms (B_1, B_3, B_4)
$$
corresponding to the two  triangles of the triangulation $T$. 

We are going to  define a rational projection 
\begin{equation} \label{9.29.03.10}
p_{\bf e}: {\rm Conf}_4({\cal B})  \lra H
\end{equation}
corresponding to the internal oriented edge $\bf e$ of the triangulation. 
Unlike the constructions in Section 5.4, it will not depend on the
choice of a reduced decomposition ${\bf i}$ of $w_0$. 

Recall the map (\ref{2.03..3.4q}). Recall that we assume that $G$ has trivial center. 
Let ${\cal R}_{U}$ be the maximal abelian quotient  of the unipotent radical 
$U$ of a Borel subgroup $B$. 
The torus $H_U:= B/U$ acts on ${\cal R}_{U}$. 
Let ${\cal R}_{U}^*$ be the unique dense $H_U$-orbit  in ${\cal R}_{U}$.
It is a principal homogeneous $H_U$--space. 
Let ${\pi}: U \to {\cal R}_{U}$ be the canonical projection. 
It provides a rational map 
$$
p_{\bf e}: ({U}_*^+ \times {U}_*^-)/H \lra 
({\cal R}^*_{U^+} \times {\cal R}^*_{U^-})/H = H, \quad (v_+, u_-) \lms 
{\pi}(\Psi(v_+))/{\pi}(u_-).
$$
Here the right hand side has the following meaning. 
For generic $(v_+, u_-)$ both  $\Psi(v_+)$ and $u_-$ lie 
in the $H$-torsor $ {\cal R}^*_{U^-}$, so their ratio is a well defined element
of $H$.  The Proposition is proved
\vskip 3mm

 \begin{proposition} \label{2.18.03.123} The birational isomorphism 
(\ref{2.18.03.4123}) provides a positive atlas on 
the moduli space ${\rm Conf}_4({\cal B})$. It is 
equivalent to the one defined in  
Section 5.6. 
\end{proposition}

{\bf Proof}. It is an immediate corollary of the following lemma. 

\begin{lemma} \label{9.1.03.1a} 
Let $(s_1, ..., s_m, t_1, ..., t_n)$ be the natural  coordinates 
on ${\Bbb G}_m^{n+m}$.  Set
$$
p_i = \frac{s_{i+1}}{s_i}, \quad 
i=1, ..., m-1, \quad q_j = \frac{t_{j+1}}{t_j},  \quad
j=1, ..., n-1,
$$
$$ X= \frac{t_{1}}{s_m},\quad
 Y = \frac{s_1+ ... + s_m}{t_1
  + ... + t_n}.
$$ 
Then the two coordinate systems  $(\{p_i\}, \{q_j\}, X)$ and
$(\{p_i\}, \{q_j\}, Y)$ 
on the quotient of ${\Bbb G}_m^{n+m}$
by the diagonal action of ${\Bbb G}_m$  are related by positive
transformations.
\end{lemma}

{\bf Proof}. Follows immediately from the  formula
$$
\frac{1 + p_1 + p_1p_2 + ... + p_1 ... p_{m-1}}
{1 + q_1 + q_1q_2 + ... + q_1 ... q_{n-1}} = p_1 ... p_{m-1}XY.
$$
Proposition \ref{2.18.03.123} is proved. 

\vskip 3mm
Let us return to the proof of Theorem \ref{2.18.03.2wq}. 
The map $i_{\overline T'}\circ i_{\overline T}^{-1}$ is certainly non regular 
at the subvariety 
of configurations of flags  $(B_1, B_2, B_3, B_4)$ 
determined by the condition that all pairs except $(B_2, B_4)$ 
are in generic position. But according to 
Proposition \ref{9.11.03.1} 
this subvariety is contained in a positive divisor. 
The part a) of the theorem is proved. 
The part b) follows immediately from Proposition \ref{9.11.03.1}.
To prove the part d), we apply the part b) of Theorem \ref{2.18.03.2wq}, 
then the results of the previous subsection, and then Theorem \ref{2.18.03.2a}.
Theorem \ref{2.18.03.2wq} is proved. 
\vskip 3mm
{\bf 8. A remark}. 
According to Theorem \ref{2.18.03.2wq}, positive 
configurations of four flags have a dihedral symmetry. The following lemma shows that they do not have $S_4$-symmetry. 
Recall the $H$-invariant $p_{\bf e}(B_1, B_2, B_3, B_4)$ of a 
configuration of four flags defined by (\ref{9.29.03.10}). 

\begin{lemma} \label{3.3.04.1} 
If $p_{\bf e}(B_1, B_2, B_3, B_4) \in H(\R_{>0})$, then \newline
$p_{\bf e}(B_1, B_3, B_2, B_4) \in H(\R_{<0})$.
\end{lemma}

\begin{figure}[ht]
\centerline{\epsfbox{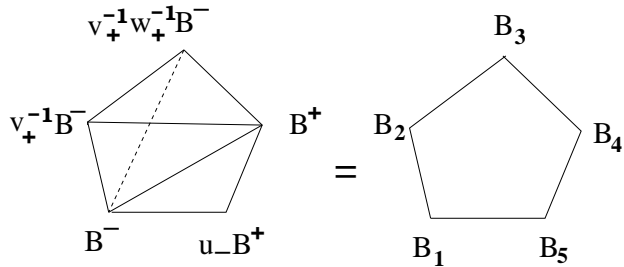}}
\caption{A triangulated pentagon with flags assigned to vertices.}
\label{fgo-39}
\end{figure}

{\bf Proof}. One has 
$$
(B_1, B_2, B_3, B_4) \sim (B^-, v_-^{-1}B^+v_-, B^+, u_+B^-u_+^{-1}) \sim 
$$
$$
\sim \Bigl(B^-, B^+, v_-B^+v^{-1}_-, v_-u_+B^-(v_-u_+)^{-1}\Bigr) 
$$
where $v_- \in U^-(\R_{>0})$ and $u_+ \in U^+(\R_{>0})$.
 Writing $v_-u_+= u'_+hv'_-$, where $u'_+,h,v'_-$ 
are also positive, we get 
$$
(B_1, B_3, B_2, B_4) \sim \Bigl(B^-, v_-B^+v_-^{-1}, B^+, u'_+B^-(u'_+)^{-1}\Bigr). 
$$
Recall the canonical projection $\pi: U^- \to {\cal R}_{U^-}:= U^-/[U^-, U^-]$. 
Elements of the abelian group ${\cal R}_{U^-}$ 
are represented uniquely as $\prod y_{\alpha}(t_{\alpha})$, where $\alpha$ runs 
through all simple positive roots. If $v_- \in U^-(\R_{>0})$ then these
 $t_{\alpha}$-coordinates 
of the element 
$\pi(v_-)$ are in $\R_{>0}$, while for 
$\pi(v^{-1}_-)$ they are in $\R_{<0}$. Now the lemma follows 
immediately from the definition of the $H$-invariant. The lemma is proved. 
\vskip 3mm
{\bf 9. A positive atlas on ${\rm Conf}_5({\cal B})$ and its 
dihedral symmetry}. 
Consider a convex 
 pentagon with vertices labeled by a configuration
$(B_1, ..., B_5)$. Let $T$ be its triangulation. The 
two diagonals of the triangulation share a vertex. We will 
assume the diagonals are oriented out of the vertex. 
 Then there is a birational isomorphism 
\begin{equation} \label{1.11.04.5cd}
\varphi_{T}: {\rm Conf}_5({\cal B})  \stackrel{\sim}{\lra} 
{\rm Conf}_3({\cal B})^3
\times H^2
\end{equation} 
where the factors ${\rm Conf}_3({\cal B})$ 
correspond to triangles with ordered vertices of the
triangulation $T$, and 
the factors $H$ are assigned to the oriented diagonals of the triangulation. 
Precisely, assuming that the diagonals share the vertex $B_1$, 
the map is given by
$$
\begin{array}{rcl}
(B_1, ..., B_5) &\lms& (B_1, B_2, B_3)\times (B_1, B_3, B_4) \times(B_1, B_4, B_5)\times \\
& &\times p_{{\bf e_{13}}}(B_1, B_2, B_3, B_4) \times p_{{\bf e_{14}}}(B_1, B_3, B_4, B_5).
\end{array}
$$
The decomposition (\ref{1.11.04.5cd}) provides a positive atlas 
of the space ${\rm Conf}_5({\cal B})$. We are going to show that 
it is invariant under the natural action of the dihedral group $D_5$.
So we have to show that it is invariant under the  
flips and reversing the order 
of the flags. The claim that a dihedral symmetry map followed by 
projection to a factor ${\rm Conf}_3({\cal B})$ is positive immediately reduced 
to dihedral symmetry of positive atlases on the spaces 
${\rm Conf}_n({\cal B})$ for $n=3,4$. So it remains only to check 
similar claims for projections onto the $H$-factors. 
More specifically, the only statement not covered 
by the dihedral symmetry of positive atlases on the spaces 
${\rm Conf}_n({\cal B})$ for $n=3,4$ is this: a flip 
on the quadrilateral $(B_1, B_2, B_3, B_4)$ of the triangulation on 
Figure \ref{fgo-39}, followed by  projection  
to the $H$-invariant assigned to the diagonal  $B_1B_4$, 
is a positive map.

It is convenient to use the following parametrisation of the space \newline
${\rm Conf}_5({\cal B})$. 
Consider the map
\begin{equation} \label{4.25.04.2}
\frac{U_*^-\times U_*^+\times U_*^+}{H} \lra 
{\rm Conf}_5({\cal B});
\end{equation}
$$
 (u_-,v_+, w_+) \lms (B^-, v_+^{-1}\cdot B^-, v_+^{-1}w_+^{-1}\cdot B^-, B^+,u_-\cdot B^+).
$$
This map is evidently a birational isomorphism. The positive atlases on $U_*^{\pm}$
provide in the standard way 
a positive atlas on the quotient on the left, and hence on 
${\rm Conf}_5({\cal B})$. 
Observe that for the triangulation on Figure \ref{fgo-39}, the vertex shared by the two diagonals is $B_4$,
 while for 
the one used in (\ref{1.11.04.5cd}) 
it is $B_1$. 

\begin{lemma} \label{4.25.04.1} 
a) The positive atlas on 
${\rm Conf}_5({\cal B})$ given by the birational isomorphism 
(\ref{4.25.04.2}) 
is equivalent to the 
standard one corresponding to the triangulation shown on Figure 
\ref{fgo-39}. 

b) The flip at the punctured diagonal on Figure 
\ref{fgo-39} is a positive automorphism of ${\rm Conf}_5({\cal B})$. 
Thus the cyclic shift is a positive automorphism of ${\rm Conf}_5({\cal B})$. 

c) The reversion $(B_1, ..., B_5) \lms (B_5, ..., B_1)$ is a 
positive automorphism of ${\rm Conf}_5({\cal B})$. 
\end{lemma}

{\bf Proof}.   a) It is obvious that, taking into account the 
$S_3$-symmetry of the positive atlas on  ${\rm Conf}_3({\cal B})$, 
the configurations of triples of flags 
assigned to the three triangles of the triangulation 
shown on Figure \ref{fgo-39} are parametri\-sed by 
$u_-, v_+, w_+$, considered modulo the $H$-conjugation. So it remains to 
investigate the $H$-invariants assigned to the edges. 

The two quadrilaterals of the triangulation 
provide two configurations of $4$ flags. The first is 
$$
(B_1, B_2, B_3, B_4) \sim (v_+\cdot B^-, B^-, w_+^{-1}\cdot B^-, B^+) = 
(\phi(v_+)\cdot B^+, B^-, w_+^{-1}\cdot B^-, B^+). 
$$
Applying the cyclic shift we 
recover the standard positive atlas on the configurations 
$(B_1, B_2, B_3, B_4)$. Further, the configuration corresponding to the second
quadrilateral is 
$$
(B_1, B_2, B_4, B_5) = (B^-, v_+^{-1}\cdot B^-, B^+, u_-\cdot B^+)
$$
is already in the standard form. 

b) After the flip at the diagonal $B_1B_3$ we get
$$
(B_1, B_3, B_4, B_5) \sim (B^-, (w_+v_+)^{-1}\cdot B^-, B^+, u_-\cdot B^+).
$$
c) Follows immediately from the properties of ${\rm Conf}_n({\cal B})$ 
for $n=3,4$.
The lemma is proved. 
\vskip 3mm
{\bf Remark}. The configuration $(B_1, B_3, B_4, B_5)$ is 
parametrised by the pair $(w_+v_+, u_-)$ considered modulo the $H$-conjugation. 
So the group product in $U^+$ appears naturally here. 
\vskip 3mm
{\bf 10. Positivity for the principal $PGL_2$-embedding}. Recall the principal $PGL_2$-embedding 
$i: PGL_2 \hra G$, well-defined up to a conjugation. It can be normalized so that 
the generator of $U_{>0}(PGL_2)$ goes to the element 
$\sum_{\alpha}X_\alpha$, where ${\rm exp}(X_{\alpha})= x_{\alpha}(1)$ and the sum is over all positive simple roots. 

\begin{proposition} \label{cb1}
The principal embedding respects the total positivity: one has 
$i(PGL_2(\R_{>0})) \subset G(\R_{>0})$ and  a similar statement for the totally positive parts of the 
maximal unipotent subgroups $U^+$. 
\end{proposition}

{\bf Proof}. 
Let $\lambda$ be a highest weight which is positive on every 
simple coroot. Denote by $V_{\lambda}$  
the corresponding irreducible representation of $G$. Let ${\Bbb B}_{\lambda}$ 
be the canonical basis in  $V_{\lambda}$. 
Let  $\eta_{\lambda}$ be a lowest weight vecor in $V_{\lambda}$.  
According to Proposition 5.4 of \cite{L2}, for a simply-laced $G$, 
$u \in U^+_{>0}$ if and only if $u (\eta_{\lambda})$ is a linear combination with $>0$ 
coefficients of the elements in the canonical basis ${\Bbb B}_{\lambda}$. 
According to Lusztig, the generators $X_{\alpha}$ of ${\cal G}$ act in the canonical basis 
by matrices with $\geq 0$ coefficients. Thus the same is true for ${\rm exp}(\sum_{\alpha}X_\alpha)$. 
In fact these coefficients are $>0$: indeed, take a monomial $c\cdot X_{\alpha_1} \ldots 
X_{\alpha_k}$, $c>0$, 
in the expansion of the above exponential corresponding to a weight 
$\eta_{\lambda} - (\alpha_1 + \ldots +\alpha_k)$ in $V_\lambda$ such that 
$\eta_{\lambda} - (\alpha_1 + \ldots +\alpha_p)$ is a weight for every $p<k$. 
The non simply-laced case is treated by folding. 
The proposition is proved. 

\begin{corollary} \label{cb2}
There is a canonical map ${\rm Conf}^+_k({\Bbb P}^1) \hra {\rm Conf}^+_k({\cal B})$ provided by 
a principal $PGL_2$-embedding. 
\end{corollary}

\vskip 3mm
{\bf 11.  Positive configurations of $K$-flags: two points of view}. 
Given a semifield $P$ there are  
the $P$-parts $U_{P}^{\pm}$, $H_{P}$, $G_{P}$ 
of the corresponding positive varieties. 
They are monoids. 
One has 
$
G_{P} = U_{P}^+H_{P} U_{P}^- = U_{P}^-H_{P}U_{P}^+
$. 

Let  $K_{>0} \subset K$ be a 
semifield in a field $K$.  Then there are two ways 
to define $K_{>0}$-positive configurations of flags. The first, preferable one,
 is this: 
we furnish the space ${\rm Conf}_n({\cal B})$ with a positive atlas and 
take its $K_{>0}$-points. For $n= 3,4,5$ the positive 
atlas has been defined above, and the general case is deduced from this
 in Section 6. 
This method works for an arbitrary semifield. 
 Below we elaborate another, more naive approach, where the fact that $K_{>0}$ 
is a semifield in a field $K$ is important.  We show that
 the two approaches lead to two {\em a priori} 
different notions, which nevertheless coincide 
if $K_{>0} = \R_{>0}$. 

Recall the configuration space ${\rm Conf}_n({\cal B}(K))$ of $n$ tuples of flags over $K$.  
\begin{definition}\label{2.18.03.20q}
a) A  configuration
$(B_1, B_2 , B_3) \in {\rm Conf}_3({\cal B}(K))$ 
is $K_{>0}$-positive  
 if 
$$
(B_1, B_2 , B_3) \sim (B^-, B^+, u\cdot B^+) \quad \mbox{
for some $u \in U^-_{K_{>0}}$}.
$$ 

b) A  configuration
$(B_1, B_2 , B_3, B_4) \in {\rm Conf}_4({\cal B}(K))$ 
is $K_{>0}$-positive  
 if 
$$
(B_1, B_2 , B_3, B_4) \sim (B^-, v^{-1}\cdot B^-, B^+, u\cdot B^+)\ \mbox{ for some $v \in U^+_{K_{>0}}$ and $u \in U^-_{K_{>0}}$}.
$$  
\end{definition}
  The set of $K_{>0}$-positive 
configurations of flags for $n=3,4$ is denoted by ${\rm Conf}^+_n({\cal
  B}(K))$. Repeating the  proofs of Theorems \ref{2.18.03.2a}
and 
\ref{2.18.03.2wq} we get the following:

\begin{corollary}\label{2.18.03.2} Let $K_{>0}$ be a 
semifield in a field $K$. Let  $G$ be as in Theorem
\ref{2.18.03.2a}. Then the sets ${\rm Conf}^+_n({\cal B}(K))$ 
for $n=3,4$ have a  natural 
dihedral symmetry. 
\end{corollary}

\begin{definition}\label{2.18.03.20qa}
Choose a triangulation $T$ of a convex  $n$-gon whose vertices 
are labeled by flags $(B_1, ..., B_n)$. 
A configuration $(B_1, ... , B_n)\in \newline {\rm Conf}_n({\cal B}(K))$ is
$K_{>0}$-positive 
if the four flags at the vertices of each quadrilateral of the triangulation $T$
form a $K_{>0}$-positive configuration of flags. 
\end{definition}

Corollary \ref{2.18.03.2} shows that only dihedral  order of 
vertices of each quadrilateral  matters. 
Any two triangulations of a convex  $n$-gon are related by a
sequence of flips. Using Corollary \ref{2.18.03.2} and the proof of 
Lemma \ref{4.25.04.1} we see that the 
notion of $K_{>0}$-positivity does not change under a flip. So we get 
the following proposition, telling us that Definition \ref{2.18.03.20qa} does
not depend on 
choice of a triangulation:
 
\begin{proposition}\label{10.28.03.1}
If a configuration $(B_1, ... , B_n)\in {\rm Conf}_n({\cal B}(K))$ is
$K_{>0}$-positive with respect to one triangulation, it is positive with
respect to any triangulation of the $n$-gon.  
\end{proposition}

Let us denote by ${\rm Conf}_n^+({\cal B}(K))$ the set of 
all $K_{>0}$-positive configurations of $K$-flags. 
\begin{corollary}\label{10.28.03.2}
 The set  ${\rm Conf}_n^+({\cal B}(K))$ is invariant under the natural 
 dihedral group action.
\end{corollary}

{\bf Proof}. Although the dihedral group action  changes the triangulation used to define 
$K_{>0}$-positive configurations, by Proposition \ref{10.28.03.1} it
does not affect $K_{>0}$-positivity. The corollary is proved. 
\vskip 3mm

For an
arbitrary semifield $K_{>0}$ one may have 
\begin{equation} \label{10.27.03.10}
{\rm Conf}^+_n({\cal B}(K)) \not = {\rm Conf}_n({\cal B})(K_{>0}).
\end{equation}
\vskip 3mm
{\bf Example}. The configuration space of three distinct points on $P^1$ is 
defined as follows: 
$$
{\rm Conf}^*_3(P^1):= (P^1
\times P^1 \times P^1 - \mbox{diagonals})/SL_2.
$$
It is a single point. However for a field $F$
the set of $SL_2(F)$-orbits  
on $P^1(F)^3 - \mbox{diagonals}$ is 
$F^*/(F^*)^2$. Indeed, under the action of $SL_2(F)$ the first two points can
be made $0$ and $\infty$, and then the subgroup ${\rm diag}(t, t^{-1})$ 
of diagonal matrices in
$SL_2(F)$ stabilizing these two points acts on $z \in F^*$ by $z \lms t^2 z$. 
For example when $F = \R$ there
are $2$ orbits. Indeed, three ordered distinct points on $P^1(\R)$ provide an
orientation of $P^1(\R)$. Moreover ${\rm Conf}^+_n({P^1}(K)) = 
{\rm Conf}_n({P^1})(K_{>0})$ if and only if $(K_{>0})^2 = K_{>0}$, which may not
be the case. 

Here is a general set up. Assume that we have a map of $F$-varieties $f:X \to
Z$, and an $F$-group $G$ acts on $X$ from the right. Assume that $G$ acts on
fibers of $f$, and this action is transitive on the fibers (over $\overline
F$). Let $z \in Z(F)$. Set $Y:= f^{-1}(z)$. It is a right homogeneous space
for $G$. Let us assume that there exists a point $y \in Y(F)$. Then, setting
$H:= {\rm Stab}_G(y)$, we have $Y = H \backslash G$. So by 
(\cite{Se}, Cor 1 of Prop. 36
in I-5.4) there is canonical isomorphism 
$$
Y(F)/G(F) = {\rm Ker}\Bigl(H^1(F, H) \lra H^1(F, G) \Bigr).
$$ 
The question whether there exists an $F$-point in $Y$ is non trivial if 
$Y$ is not a principal homogeneous space for $G$, and, as was 
pointed out to us M. Borovoi, is treated 
for $p$-adic and number fields in the Subsection 7.7 of \cite{Bor}.  
We will not need this in our story.  

So if $Z:= X/G$ is the quotient
  algebraic variety, there are two different 
sets: $Z(F)$ and
  $X(F)/G(F)$. The canonical map $X(F)/G(F) \lra Z(F)$ can be 
  neither 
surjective nor injective. 

However for positive configurations of real flags 
the situation is as nice as it could be:

\begin{lemma} \label{10.25.03.30}
For $n=3,4,5$, one has $${\rm Conf}^+_n({\cal B}(\R)) = {\rm Conf}_n({\cal
  B})(\R_{>0}),$$ and there is a canonical inclusion 
$$
{\rm Conf}_n({\cal B})(\R_{>0}) \hra {\rm Conf}_n({\cal B}(\R)). 
$$ 
\end{lemma}

{\bf Proof}. It boils down to the following two facts:  \newline
${\cal V}_*(\R_{>0}) = U_*(\R_{>0})/H(\R_{>0})$ and similarly 
$((U^-_* \times U^+_*)/H)(\R_{>0}) = U^-_*(\R_{>0}) \times
U^+_*(\R_{>0})/H(\R_{>0})$. They follow from $(\R_+^*)^m = \R^*_+$. The lemma
is proved. \vskip 3mm

A positive atlas on the space ${\rm Conf}_n({\cal
  B})$ is defined in Section 6. Then it is straightforward 
to check that Lemma 
\ref{10.25.03.30} remains valid for all $n>2$. 
So 
the real positive configurations of flags in the sense of Definition
\ref{2.18.03.20qa} always coincide with the 
$\R_{>0}$-points of the positive space ${\rm Conf}_n({\cal B})$, and are 
embedded into the configuration space  of $n$ real flags in 
generic position modulo the  $G(\R)$-action. 

\section{A  positive structure on the moduli space ${\cal X}_{G, \widehat S}$ }
\label{posAX}


In this Section we give three different proofs, in Sections 6.1, 6.3 and 
 6.5, of the 
decomposition of the moduli space ${\cal X}_{G, \widehat S}$ according to a
triangulation. This provides a positive atlas
 on ${\cal X}_{G, \widehat S}$.

\vskip 3mm
{\bf 0. Basic properties of the positive atlas on ${\rm Conf}_n({\cal B})$ for $n=3,4,5$}. 
Let us summarize the properties of the positive atlases on the spaces 
${\rm Conf}_n({\cal B})$, where $n=3,4,5$,  which are used
below. The proofs are given in Theorem \ref{2.18.03.2a}, Theorem 
\ref{2.18.03.2wq}, and Lemma \ref{4.25.04.1}.

1) The space ${\rm Conf}_3({\cal B})$ has a positive atlas invariant under
   the action of the group $S_3$.

2) There is a birational isomorphism 
\begin{equation} \label{1.11.04.5}
\varphi: {\rm Conf}_4({\cal B})  \stackrel{\sim}{\lra} {\rm Conf}_3({\cal B})
\times H \times {\rm Conf}_3({\cal B}), 
\end{equation} 
where projections to the left and right factors ${\rm Conf}_3({\cal B})$ 
are given by the maps $$
(B_1, B_2, B_3, B_4) \lms (B_1, B_2, B_3) \quad \mbox{and}\quad 
(B_1, B_2, B_3, B_4) \lms (B_1, B_3, B_4).$$ 

3) The birational isomorphism (\ref{1.11.04.5})
 provides a positive atlas on the space ${\rm Conf}_4({\cal B})$. This
 positive atlas  is invariant under
   the cyclic shift and reversion maps given by 
$$
(B_1, B_2, B_3, B_4) \lra (B_4, B_1, B_2, B_3), \qquad (B_1, B_2, B_3, B_4) \lra 
(B_4, B_3, B_2, B_1).
$$

If we picture a configuration $(B_1, B_2, B_3, B_4)$ 
at the vertices of a $4$-gon as
on Figure \ref{fgo2}, then the cyclic invariance boils down to the 
invariance of the
positive atlas on ${\rm Conf}_4({\cal B})$ under a flip, and
hence the invariance under the 
change of an orientation of the diagonal $e$. The invariance 
under the reversion map shows the invariance of the positive atlas 
under the change of an orientation of the plane where the $4$-gon is located. 

4) Consider the convex 
 pentagon whose vertices are labeled by a configuration
$(B_1, ..., B_5)$. Let $T_1$ be its triangulation 
by the diagonals $B_1B_3$
 and $B_1B_4$. 
 Then there is a birational isomorphism 
\begin{equation} \label{1.11.04.5c}
\varphi_{T_1}: {\rm Conf}_5({\cal B})  \stackrel{\sim}{\lra} 
\Bigl({\rm Conf}_3({\cal B})\Bigr)^3
\times H^2 
\end{equation} 
where the factors ${\rm Conf}_3({\cal B})$ correspond to the triangles of the
triangulation $T_1$, and 
the factors $H$ are assigned to the diagonals $B_1B_3$
 and $B_1B_4$ of the triangulation. 

The positive atlas 
of the space ${\rm Conf}_5({\cal B})$ provided by (\ref{1.11.04.5c})
is invariant under the natural action of the dihedral group $D_5$. 

\vskip 3mm
 {\bf 1. A positive atlas on ${\cal X}_{G, \widehat S}$}. Let us 
spell in detail the definition of 
 the map $\pi_{\bf T}$ outlined in Section  1.3. 
 Let $\widehat S$ be a  marked hyperbolic surface.
Let us choose a pair  ${\bf T}$ given by 
\begin{equation} \label{1.11.04.1}
\mathbf T = \left\{\begin{array}{l}\mbox{an ideal triangulation} T \mbox{ of } \widehat S,\\
 \mbox{a choice of orientation of all internal edges of } T
\end{array}
\right\}.
\end{equation} 
Let us assume first that the triangulation is regular in the sense 
of Section 3.8. 
Denote by ${\rm tr}(T)$ the set of the triangles of $T$, and by 
 ${ ed_i}({\bf T})$ 
be the set of the internal edges of ${\bf T}$. 
Denote by ${\cal X}_{G, \widehat t}$ the moduli space of the framed $G$--local
systems on a triangle $t$, considered as a disc with three marked points on
the boundary given by certain three points inside of the edges, one point per
each edge. 
Let $({\cal L}, \beta) \in {\cal X}_{G, \widehat S}$. 
We can view a triangle $t$ of the ideal triangulation 
as a disc with three marked points at the boundary, as we just explained. 
Restricting 
the framed local system $({\cal L}, \beta) \in {\cal X}_{G, \widehat S}$ 
to such a  triangle $t$ 
we get an element of ${\cal X}_{G, \widehat t}$. 
Indeed, if a vertex $s$ of the triangle $t$ belongs to a boundary arc on 
$\widehat S$ then there is a flag at this point. 
If $s$ is a puncture then we take the monodromy invariant flag 
at the intersection of the triangle  with a little disc around the puncture.
 So we get a projection $p_t: 
{\cal X}_{G, \widehat S} \lra {\cal X}_{G, \widehat t}$. 
 Further, each  internal  edge 
${e}$ of the triangulation $T$ determines 
a pair of triangles $t_1, t_2$ sharing ${ e}$. Their union $t_1 \cup t_2$
is a $4$-gon with  four marked points on the boundary, one per each
side. Thus restricting 
$({\cal L}, \beta)$ to $t_1 \cup t_2$ we get an element of ${\cal X}_{G, 4}$. Its
$H$--invariant corresponding to the oriented 
edge $\bf e$ is given by (\ref{9.29.03.10}). It   provides a rational projection
$p_{\bf e}: {\cal X}_{G, \widehat S} \lra H$. The collection of 
projections $\{p_t\}, \{p_{\bf e}\}$ provides a 
 rational map 
\begin{equation} \label{9.1.03.1e}
\pi_{\bf T}: {\cal X}_{G, \widehat S} \lra 
\prod_{t \in {\rm tr}(T)}{\cal X}_{G, \widehat t} 
\times \prod_{{\bf e} \in { ed_i}({\bf T})}H.
\end{equation}
If the triangulation $T$ is special, the above construction works 
unless we consider the edge given by the loop of the virus. 
The construction of the $H$-invariant for this edge 
is obtained by going to the $2:1$ cover of 
$S$ described in Section 3.8, and then using the above construction 
of the $H$-invariant. The same trick, described in detail in the Section 10.7, 
 is used in the proof of the theorems in this Section 
in the part related to this edge. 

\begin{theorem} \label{9.1.03.1n}  Let $G$ be a
 split semi-simple algebraic group with trivial center and $\widehat S$ a  marked hyperbolic
 surface. Then the map $\pi_{\bf T}$ is a birational isomorphism. 
The collection of rational maps $\{\pi^{-1}_{\bf T}\}$, 
where ${\bf T}$ run through all  datas  (\ref{1.11.04.1}),  
provides a $\Gamma_S$--equivariant positive atlas on ${\cal X}_{G,
  \widehat S}$. 
\end{theorem}

{\bf Proof}. Let us prove first the theorem in the special case when 
$\widehat S = \widehat D_n$ is a disc with $n$ marked points on the boundary. We proceed by
the induction on $n$. For $n=3$ this is given by the property 1) 
above. To make the inductive step from $n$ to $n+1$  consider a convex 
$(n+1)$-gon $P_{n+1}$ with the vertices $b_1, ..., b_{n+1}$ taken in
the order compatible with the orientation of the disc. We assign a
configuration of flags 
 $(B_1, ..., B_{n+1})$ to the vertices of $P_{n+1}$, so that $B_i$ is assigned
 to the vertex $b_i$. Let us triangulate
 this $(n+1)$-gon so that $(b_1, b_2, b_3)$ is one of the triangles of the
 triangulations. Let $P_n$ be the $n$-gon obtained by cutting out the triangle
 $(b_1, b_2, b_3)$. By the induction assumption a generic data  
on the right hand side of (\ref{9.1.03.1e}) determines a configuration $(B_1,
B_3, B_4, ..., B_{n+1}) \in {\rm Conf}_{n}({\cal B})$ assigned to the vertices
of the $n$-gon $P_n$. Adding the triangle 
$(b_1, b_2, b_3)$ we add to this data the element $\beta_{1 2 3} \in {\rm
  Conf}_{3}({\cal B})$ corresponding to this triangle, and an $H$-invariant
$h_{13}\in H$ corresponding to the edge $(b_1, b_3)$. 
 Let $(b_1, b_3, b_k)$ be the triangle of the triangulation
 adjacent to the edge $(b_1, b_3)$, and  $\beta_{1 3 k} \in {\rm
  Conf}_{3}({\cal B})$ the invariant assigned to this triangle. (Since we work rationally, 
we may ignore the difference between ${\rm Conf}_{3}({\cal B})$ and ${\cal X}_{G, 3}$). 
According to the property 2) there is a unique configuration $(B'_1, B'_2, B'_3,
B'_k)$ of four flags 
with the invariant  $(\beta_{1 2 3}, h_{13}, \beta_{1 3 k})$. 
Observe that the configurations $(B'_1, B'_3, B'_k)$
and 
$(B_1, B_3, B_k)$ are the
same.

\begin{figure}[ht]
\centerline{\epsfbox{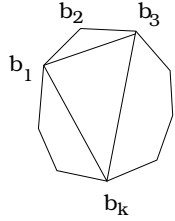}}
\caption{The quadrilateral corresponding to the triangle $b_1b_2b_3$.}
\label{fg60}
\end{figure}

Since $G$ has trivial center, it acts without stable points
on the space of generic triples of flags in $G$.  Therefore 
there exists a unique configuration $(B_1, B_2,  ..., B_{n+1})$ 
such that $(B'_1, B'_2, B'_3, B'_k) = (B_1, B_2, B_3, B_k)$
corresponding to the data provided by the triangulation of the $(n+1)$-gon. 
So, given a  triangulation $T$ of the $(n+1)$-gon,  the 
 map (\ref{9.1.03.1e})  is a birational isomorphism for the disc
$\widehat D_{n+1}$. Thus it provides a positive atlas on ${\rm
  Conf}_{n+1}({\cal B})$. It follows from 3) that the positive atlas 
corresponding to the triangulation $T'$ obtained from $T$ by a flip 
is compatible with this
positive atlas. 
The role of the mapping class group in this case is played by the dihedral group 
generated by the cyclic shift and reversion maps. 
This proves the theorem for the disc. 

Now let us consider the general case. It is sufficient to show that, given an
arbitrary field $F$, a generic point $x$ of the right hand side of (\ref{9.1.03.1e})
provides an $F$-point 
$y$ of the right hand side such that $\pi(x) =y$. Let us define 
a representation $\rho_x:\pi_1(S) \lra G(F)$ corresponding to $x$. 
Let $T$ be an ideal triangulation of $S$. Lifting it to 
a universal cover $\widetilde{S}$ of
$S$ we obtain a triangulation $\widetilde T$. The triangulated 
surface $(\widetilde{S}, \widetilde T)$  
is isomorphic to  the hyperbolic disc ${\cal H}$ equipped with the 
Farey triangulation. Let $(a,b,c)$ be a triangle of $\widetilde T$. 
Let $\gamma \in \pi_1(S)$. Take a
finite polygon $P_{\gamma}$ of the triangulation $\widetilde T$ 
containing the triangles $(a,b,c)$ and $(a',b',c') := \gamma (a, b, c)$.  
The polygon $P_{\gamma}$ inherits a triangulation. So by 
the first part of the proof, the data $y$,
restricted to the triangulated polygon, gives rise to a  unique 
element $y(P_{\gamma}) \in {\rm Conf}_n({\cal B}(F))$ corresponding to it. 
Let $(B_a, B_b, B_c)$ be the corresponding configuration of flags attached to 
the vertices $(a,b,c)$ of the polygon. Then there exists 
an element $g_{\gamma} \in G(F)$ such that 
$$
(B_{a'}, B_{b'}, B_{c'}) = g_{\gamma}(B_a, B_b, B_c).
$$
Such an element is unique since $G$ has trivial center. 
It is easy to see that $\gamma \to g_{\gamma}$ is a representation of 
$\pi_1(S, p)$, and it does not depend on the choice of the polygon. 
Here $p$ is a point inside of the triangle on $S$ obtained by projection of the triangle $(a,b,c)$. 
We leave construction of the framing as an easy exercise. 
The theorem is proved. 
\vskip 3mm
In the next subsection we work out an explicit construction of the framed 
local system corresponding to the generic data in (\ref{9.1.03.1e}).

\vskip 3mm

{\bf 2. The six canonical pinnings determined by a 
generic triple of flags in $G$}. Let us choose two reduced words ${\bf i}^+$
and ${\bf i}^-$ of the element $w_0$. Set
$$
U^-_*({\bf  i}^-):= \psi_{\overline {{\bf i}^-}}({\Bbb G}_m^n) \subset U^-_*, 
\qquad U^+_*({\bf i}^+):= \psi_{\bf i^+}({\Bbb G}_m^n) \subset U^+_*.
$$
The group $H$ acts freely on each of these spaces. So the $H$-orbit of any
element  $u
\in U^-_*({\bf  i^-})$  intersects the subvariety
\begin{equation} \label{U^-_*}
\overline U^-_*({\bf i^-}):= \psi_{\overline {\bf i^-}}\Bigl(e \times 
{\Bbb G}_m^{{\bf i^-}(\alpha_1)-1} \times e \times 
{\Bbb G}_m^{{\bf i^-}(\alpha_2)-1} \times ... \times e \times 
{\Bbb G}_m^{{\bf i^-}(\alpha_r)-1}\Bigr) \subset U^-_*({\bf i^-})
\end{equation} 
in a single point, called the normalized representative $\overline u$. 
Similarly we define the normalized representative $\overline v$ of the
$H$-orbit 
of an element 
$v \subset U^+_*({\bf i^+}):= \psi_{\bf i^+}({\Bbb G}_m^n)$ as the intersection
point of the orbit with the subvariety
\begin{equation} \label{U^+_*}
\overline U^+_*({\bf i^+}):= \psi_{\bf i^+}
\Bigl({\Bbb G}_m^{{\bf i^+}(\alpha_1)-1} \times e
 \times {\Bbb G}_m^{{\bf i^+}(\alpha_2)-1} \times e
 \times ... \times 
{\Bbb G}_m^{{\bf i^+}(\alpha_r)-1} \times e \Bigr) \subset U^+_*({\bf i^+}).
\end{equation}

\begin{figure}[ht]
\centerline{\epsfbox{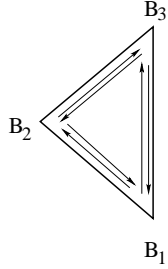}}
\caption{Arrows encode the six pinnings assigned to a configuration of three flags.}
\label{fg-2}
\end{figure}

 Let us  say  that a triple $(B_1, B_2, B_3)$ of flags in $G$ is a sufficiently generic, if one has 
$$
(B_1, B_2, B_3) = g_-\cdot (B^-, B^+, \overline u_- B^+ \overline u_-^{-1}), \quad
u_- \in U^-_*({\bf   i^-})
$$
$$
(B_1, B_2, B_3) = g_+ \cdot (B^-,  \overline v_+^{-1} B^- \overline v_+, B^+), \quad
v_+ \in U^+_*({\bf   i^+}) 
$$
for the uniquely defined elements $g_-, g_+ \in G$. The sufficiently generic triples form a non empty 
Zariski open subset in the configuration space of triples of flags in $G$.  
Therefore such a triple 
$(B_1, B_2, B_3)$ of flags in $G$ 
determines the following two pinnings for $G$:
$$
p_{(B_1, B_2, B_3)}^+(B_1, B_2):= g_-\cdot (B^-, B^+), \qquad p_{(B_1, B_2, B_3)}^-(B_1, B_3):= g_+ \cdot 
(B^-, B^+).
$$ 
Using the cyclic shifts we define the rest of the pinnings. So we get 
the six pinnings for $G$, where the indices $i$ are modulo $3$:
\begin{equation} \label{9.28.03.2}
p_{(B_1, B_2, B_3)}^+(B_i, B_{i+1}), \quad  p_{(B_1, B_2, B_3)}^-(B_i, B_{i-1}). 
\end{equation}
These pinnings are determined by the cyclic configuration 
of the three flags $(B_1, B_2, B_3)$ and the choice of two reduced words 
${\bf   i^-}, {\bf   i^+}$ for $w_0$.  Observe that $p^{\pm}$ corresponds 
to $g^{\mp}$.  
The $+$ and $-$  in (\ref{9.28.03.2}) indicate whether the
corresponding edge is oriented accordantly to the cyclic order of $(B_1, B_2,
B_3)$ or not. It is handy to picture the six pinning related to a triple $(B_1, B_2, B_3)
\in {\rm Conf}_3({\cal B})$ by arrows as on Figure \ref{fg-2}. 
Say that  {\it $U$ is a positive Zariski open subset of a positive variety $X$}
if $X-U$ is contained in a positive divisor of $X$.  
The discussion above proves the following lemma.

\begin{lemma} \label{9.28.03.1}
There exists a positive Zariski open subset 
${\rm Conf}^o_3({\cal B}) \subset {\rm Conf}_3({\cal B})$ such that for any
configuration $(B_1, B_2, B_3) \in {\rm Conf}^o_3({\cal B})$ the above construction
determines the six pinnings (\ref{9.28.03.2}) for $G$. 
Changing a standard pinning for $G$ amounts to a shift of each of 
our the six pinnings by a common element in $G$.
\end{lemma}
We leave a precise description of ${\rm Conf}^o_3({\cal B})$ as the
intersection of several explicitly given positive subvarieties 
of ${\rm Conf}_3({\cal B})$ as a straightforward exercise.

\vskip 3mm
 {\it  The edge invariant}. 
Given a generic configuration  $(B_1, B_2, B_3,  B_4)$ of flags
 in $G$,
such that both $(B_1, B_2, B_3)$ and $(B_1, B_3, B_4)$ are in 
${\rm Conf}^o_3({\cal B})$, there exists a unique element 
$$
h_{({\bf   i^-}, {\bf   i^+})}(B_1, B_2, B_3,  B_4) \in H
$$
 such that 
\begin{equation} \label{9.28.03.4}
p_{(B_1, B_3, B_4)}^+(B_1, B_3)  = 
h_{({\bf   i^-}, {\bf   i^+})}(B_1, B_2, B_3,  B_4) \cdot
p_{(B_1, B_2, B_3)}^-(B_1, B_3).
\end{equation}
The element $h_{({\bf   i^-}, {\bf   i^+})}(B_1, B_2, B_3,  B_4) \in H$ in 
(\ref{9.28.03.4}) 
is called the {\it edge
 invariant}  of the 
configuration of flags  $(B_1, B_2, B_3,  B_4)$ corresponding to the oriented edge
$(B_1, B_3)$. 

\begin{figure}[ht]
\centerline{\epsfbox{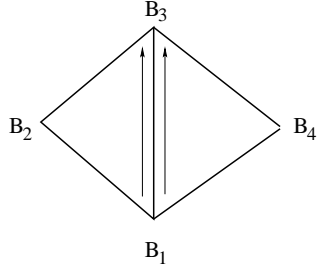}}
\caption{The edge invariant is obtained by comparison of the two pinnings shown by
  arrows.}
\label{fg-7}
\end{figure}
 
\vskip 3mm
{\bf 3. The positive structure of ${\cal X}_{G, \widehat S}$}. 
To describe the rational inverse to the map  $\pi_{\bf T}$, see (\ref{9.11.03.10}), 
we need 
a finer data assigned to a triangulation $T$ of $\widehat S$.

\begin{definition} \label{9.28.03.7} Denote by ${\Bbb T}$ the following data
  on $\widehat S$: 

i) An ideal  triangulation 
$T$ of a marked hyperbolic surface $\widehat S$.

ii) For each internal edge of the triangulation 
$T$ a choice of orientation of this edge. 

iii) For each triangle $t$ of the triangulation 
$T$ a choice of a vertex of this triangle. 

iv) A pair of reduced decompositions $({\bf   i^-}, {\bf   i^+})$ 
of the maximal length element $w_0
\in W$. 
\end{definition}
We can  visualize a choice described in iii) as  a choice of an angle
for each of the triangles $t$.

\begin{theorem} \label{9.1.03.1}  Let $G$ be a
 split semi-simple algebraic group with trivial center. 
 Let $\widehat S$ be a  marked hyperbolic
 surface and ${\Bbb T}$ a data assigned to an ideal triangulation $T$ of
 $\widehat S$, as in  Definition \ref{9.28.03.7}. Then 

a) There exists a regular open embedding
\begin{equation} \label{9.28.03.11}
\psi_{{\Bbb T}}: {{\rm Conf}^o_3({\cal B})}^{\{\mbox{triangles of $T$}\}} \times 
H^{\{\mbox{edges of $T$}\}} \hra {\cal X}_{G, \widehat S}.
\end{equation}

b) The collection of regular open embeddings $\{\psi_{{\Bbb T}}\}$, when 
${\Bbb T}$ runs through all possible datas from Definition \ref{9.28.03.7} 
assigned to $\widehat S$, provides
a $\Gamma_S$--equivariant 
positive structure on ${\cal X}_{G, \widehat S}$. 

c) The positive structure defined in Theorem \ref{9.1.03.1n} is 
compatible with the one defined in b). 
\end{theorem}

{\bf Proof}. a) Recall the six pinnings assigned according to Lemma
\ref{9.28.03.1} to an element of $C \in {\rm Conf}^o_3({\cal B})$. 
Their dependence on the choice of a standard pinning $(B^-, B^+)$ boils down  
to a shift by a common element in $G$. 
Further, let us suppose that we have two elements $C_1, C_2 \in {\rm Conf}^o_3({\cal B})$
and an element $h \in H$. Then there exists a unique configuration 
$(B_1, B_2, B_3, B_4) \in {\rm Conf}_4({\cal B})$ such that 
\begin{equation} \label{8.28.03.10a}
C_1 = (B_1, B_2, B_3), \quad C_2 = (B_1, B_3, B_4), \quad 
h_{({\bf   i^-}, {\bf   i^+})}(B_1, B_2, B_3,  B_4)=h.
\end{equation}
Moreover if we choose one of the six pinnings corresponding to the  configuration 
$(B_1, B_2, B_3)$, then there is a unique way to normalize the six pinnings 
corresponding to the  configuration $(B_1, B_3, B_4)$ so that 
formula (\ref{9.28.03.4})  holds with 
$h_{({\bf   i^-}, {\bf   i^+})}(B_1, B_2, B_3,  B_4)=h$. 
Therefore given (\ref{8.28.03.10a}) and one of the $12$ pinnings shown by
arrows on the picture we determine uniquely the rest of the pinnings. 
 \begin{figure}[ht]
\centerline{\epsfbox{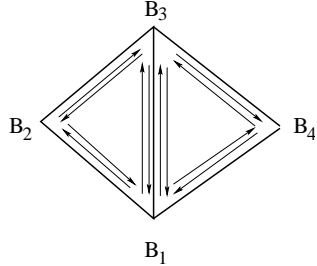}}
\caption{The twelve pinnings assigned to a quadrilateral.}
\label{fg-4}
\end{figure}
Let $\Gamma$ be a graph on a marked hyperbolic surface $\widehat S$ dual to an
ideal triangulation $T$ of $\widehat S$. 
Consider a closed path ${\bf p}$ on the graph. 
We are going to attach an element $g_{\bf p} \in G$ to this path which will give rise
to a homomorphism $\pi_1(S) \to G$. 
Let us consider the triangles of the
triangulation corresponding to the vertices of the path ${\bf p}$, as on the
picture. Let us choose one of these triangles and denote it $t_1$. 
\begin{figure}
\centerline{\epsfbox{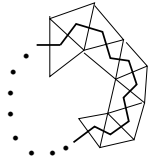}}
\caption{Constructing the local system.}
\label{fg-3}
\end{figure}
The data ${\Bbb T}$ allows to assign to each triangle $t$ a
 configuration of flags $(B_1, B_2, B_3)$ sitting at the vertices 
of  $t$, so that  $B_1$ is assigned to the distinguished vertex (provided by 
the  condition 
iii) in the definition of  ${\Bbb T}$), 
and the orientation of the triangle $t$ given by the order of $B_i$'s 
is the same as the 
one induced by the orientation of the surface  $S$ (clockwise on our pictures). 
It follows from the remarks made about pinnings above that given one of the
six 
canonical 
pinnings associated to the flags sitting at the vertices of the 
triangle $t_1$, we determine uniquely the $6$-tuples of pinnings 
corresponding to the triangles following $t_1$ along the path ${\bf p}$, 
till we return  to the
triangle $t_1$. Thus for a closed path ${\bf p}$ there are two 
$6$-tuples of canonical pinnings  assigned to the triangle $t_1$: the original one, and the
 one obtained after traveling along the loop ${\bf p}$. These two 
$6$-tuples of  pinnings differ by a shift in $G$. 
So we pick up an element $g_{\bf p} \in G$ transforming the original $6$-tuple
of pinnings to the final one. Clearly the map ${\bf p} \lms g_{\bf p}$ is
 multiplicative: ${\bf p}_2{\bf p}_1 \lms g_{\bf p_2}g_{\bf p_1}$. 
Moreover ${\bf p}^{-1} \lms g^{-1}_{\bf p}$. Since  ${\bf p}$ is a path on the
 graph, its homotopy class is determined uniquely by the path itself. 
So we get a group homomorphism  
$\pi_1(S, x) \to G$ where $x$ is a base  point on $\Gamma$ inside of the
 triangle $t_1$. 

Let us define a framed structure on this representation. 
Let $s$ be a puncture on $S$. Consider all triangles of the 
triangulation sharing the vertex $s$. 
\begin{figure}[htb]
\centerline{\epsfbox{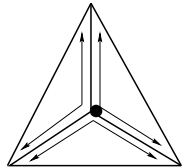}}
\caption{The link of a puncture $s$.}
\label{fg-5}
\end{figure}
Let us choose a path ${\bf p_s}$ from the base 
 point $x$  to a point inside of a certain  triangle $t_s$ sharing the vertex $s$. 
Then given a pinning $P_x$ in the triangle $t_x$  containing $x$ we get the six
 pinning assigned to the triangle $T_s$. Let 
 $P_s$ be one of them corresponding to an arrow going from the vertex $s$. Let 
$(B, B')$ be the pair of Borel subgroups underlying this pinning, so
 that $B$ is assigned to the vertex $s$. Going around  $s$ 
we define pinnings corresponding to the other triangles sharing $s$. 
Observe that all pinnings corresponding to the arrows starting at $s$ 
have the same ''initial''  Borel subgroup,  $B$. Therefore the
 monodromy of the local system corresponding to our representation 
around a small loop ${\bf l_s}$ going around $s$ 
belongs to $B$. So if  $P_s = gP_x$, then the element  
${\bf p_s}^{-1}{\bf l_s}{\bf p_s}$ leaves the Borel subgroup 
$g^{-1}Bg$ invariant. This provides a framed 
structure of our representation at the puncture $s$.

We have defined the map $\psi_{\Bbb T}$. Let us show that it is injective. 
Let us define a left inverse map
$$
\pi_{\Bbb T}: {\cal X}_{G, \widehat S} \lra 
\prod_{t \in {\rm tr}(T)}{\cal X}_{G, \widehat t} 
\times \prod_{{\bf e} \in {\bf ed_i}(\Bbb T)}H.
$$
Here ${\bf ed_i}(\Bbb T)$ is the set of internal edges of the triangulation $T$
oriented as prescribed by the data $\Bbb T$. So the second product is over the
set 
of all internal edges of $T$. 
The $\prod_{t \in {\rm tr}(T)}{\cal X}_{G, \widehat t}$--component of the map
$\pi_{\Bbb T}$ 
is the same as for the map  $\pi_{T}$. 
The $H$-component corresponding to
an oriented internal edge $\bf e$ is the edge invariant from 
 (\ref{9.28.03.4}) assigned to $\bf e$.  
    Notice that it is 
different from the one used in $\pi_{T}$. 
Then by the very construction we have $\pi_{\Bbb T} \circ \psi_{\Bbb T} =
{\rm Id}$. 
Thus the map $\psi_{\Bbb T}$ is injective. 
\begin{lemma} \label{5.18.03.1}
Let $f:X \lra Y$ and $g:Y \lra X$ be maps of irreducible algebraic varieties of the same dimension. Then $g \circ f = {\rm Id}_X$ implies that 
both $f$ and $g$ are birational isomorphisms. 
\end{lemma}

{\bf Proof}. The conditions of the lemma imply that the maps $f,g$ are dominant, 
so there are maps of fields $f^*: \Q(Y) \to \Q(X)$, $g^*: \Q(X) \to \Q(Y)$. 
Since $f^*g^* = {\rm Id}$, the lemma follows. 
\vskip 3mm
 
\begin{lemma} \label{9.30.03.1}  Let $\widehat S$ be a marked hyperbolic surface. 
Denote by $t(T)$ and $e_i(T)$ the number of triangles and internal edges of an ideal
 triangulation $T$ of  $\widehat S$. 
Then one has
\begin{equation}\label{3.28.02.4}
{\rm dim}{\cal X}_{G, \widehat S} = e_i(T) 
{\rm dim}H + t(T) ({\rm dim}U - {\rm dim}H).
\end{equation}
\end{lemma}

This Lemma, of course, follows from the proof given in Section 6.1. However
we will give an independent proof for completeness of the arguments. 
 
{\bf Proof}. Assume for a moment that the number $k$ of marked boundary points
on $\widehat S$ is equal to zero, i.e. $\widehat S = S$ and $e(T) = e_i(T)$ is
the number of all edges. Then
$
\chi(S) = -e(T) + t(T), \quad 2e(T)  = 3t(T)
$. 
Therefore 
$$
e(T) = -3\chi(S), \quad t(T) = -2\chi(S).
$$
 It is well known that ${\rm dim}{\cal L}_{G, S} = -\chi(S){\rm
  dim} G$.
The canonical map ${\cal X}_{G, S}\to {\cal L}_{G, S}$ is a finite map at the
  generic point. This follows from the well known fact that the monodromies of
  a generic $G$-local system on $S$ around the punctures are regular conjugacy
  classes in $G$. This gives ${\rm dim}{\cal X}_{G, S} = -\chi(S){\rm dim}
  G$.  
The equality (\ref{3.28.02.4}) in the case $\widehat S = S$ 
follows immediately from these remarks. To check (\ref{3.28.02.4}) in  general  
observe that adding
a marked point we increase the dimension by ${\rm dim} (G/B)  = {\rm dim} U$. On the other hand 
adding a marked point we add by one the number of internal edges and triangles of a
triangulation. So it remains to check the statement in the case when $S$ is a
 disc or cylinder to have the base for the induction. 
The lemma is proved. 
\vskip 3mm
It follows from Lemma \ref{9.30.03.1} 
both spaces in (\ref{9.28.03.11}) are of the same dimension. 
Thus (Lemma \ref{5.18.03.1}) $\psi_{\Bbb T}$ is a birational isomorphism. 
The part a) of the theorem is proved.

b) It is sufficient to show that 
if $T$ and $T'$ are two triangulations of $\widehat S$ related by a flip, and ${\Bbb T}$ and ${\Bbb T}'$ are two datas refining
these triangulations, then the positive structures defined by $\psi_{\Bbb T}$ 
and $\psi_{\Bbb T'}$ are compatible.  
This follows from Theorem
\ref{2.18.03.2wq} and Lemma \ref{4.25.04.1}. Observe that we need to consider 
configurations of $5$-tuples of flags in order to make sure that the positive structures related by a flip in an edge $E$ 
are compatible: 
indeed, consider the $4$-gon containing $E$ as a diagonal. Then the $H$-coordinates assigned to an  
external edge $F$ of this $4$-gon depend on the flags in the vertices of the $4$-gon containing $F$ as a diagonal. 
One of them is assigned to a vertex outside of the original $4$-gon.  The part b) of the theorem is proved.

c) It follows from Proposition \ref{2.18.03.123}. 
The theorem is proved. 
\vskip 3mm

{\bf 4. The higher Teichm{\"u}ller ${\cal X}$-space 
and the ${\cal X}$-lamination space 
related to $G$}. 
\begin{definition} \label{4.02.03.3} 
Let $\widehat S$ be a marked hyperbolic surface and $G$ a 
split reductive algebraic group with trivial center.  

a) The Teichm{\"u}ller space  ${\cal X}^+_{G, \widehat S}$ is 
the $\R_{>0}$-positive part ${\cal X}_{G, \widehat S}(\R_{>0})$ of the  
positive space ${\cal X}_{G, \widehat S}$. 

b) Let ${\Bbb A}^t$ be one of the tropical semifields 
$\Z^t, \Q^t, \R^t$.  The set 
of ${\cal X}^{{\Bbb A}}$-laminations on $\widehat S$ corresponding to  $G$ 
is 
the set  ${\cal X}_{G, \widehat S}({\Bbb A}^t)$ 
of points of the positive space ${\cal X}_{G, \widehat S}$
with values in the tropical semifield ${\Bbb A}^t$.
\end{definition}

According to Section 4.5 projectivisation 
of the set of real ${\cal X}$-laminations 
can be considered as the Thurston type boundary of the space 
${\cal X}^+_{G, \widehat S}$.

Applying the definition of the regular part of a positive variety (Section
 4.1) 
 in our case we get the regular  
part ${\cal X}^{\rm reg}_{G,\widehat S} \hra {\cal X}_{G,\widehat S}$. 
The  higher Teichm{\"u}ller space is a connected component in 
${\cal X}^{\rm reg}_{G,\widehat S}(\R)$.

Recall from Section 2.1 the set $\widetilde 
{\cal X}_{G, \widehat S}(\R)/G(\R)$ parametrising framed 
representations $\pi_1(S, x) \to G(\R)$ modulo $G(\R)$-conjugation. 
If $\widehat D_3$ is a disc with three marked points, then ordering these
points we get an isomorphism 
$\widetilde 
{\cal X}_{G, \widehat D_3}(\R) = {\cal B}(\R)^3$. 
Let $T$ be an ideal triangulation of $\widehat S$. Then each triangle $t$
of the triangulation provides a restriction map 
$\widetilde {\cal X}_{G, \widehat S}(\R)/G(\R) \to \widetilde 
{\cal X}_{G, \widehat t}(\R)/G(\R)$. We employ quadrilaterals 
of the triangulation in a similar
way.

An element of  
$\widetilde {\cal X}_{G, \widehat S}(\R)/G(\R)$ is called {\it positive} 
if its
restriction to every quadrilateral of a given triangulation $T$ 
belongs to ${\rm Conf}_4^+({\cal B}(\R))$. 

\begin{corollary} \label{4.02.03.1} 
The set of positive elements in $\widetilde {\cal X}_{G, \widehat
  S}(\R)/G(\R)$
does not depend on the choice of triangulation $T$. 
It is isomorphic to ${\cal X}^+_{G, \widehat S}$. So 
there is an embedding 
$
{\cal X}^+_{G, \widehat S} \hra \widetilde 
{\cal X}_{G, \widehat S}(\R)/G(\R)
$.  
\end{corollary}

{\bf Proof}.    
Any two ideal triangulations of a hyperbolic marked 
surface are related by a sequence of flips. 
Corollary \ref{2.18.03.2} plus Lemma \ref{4.25.04.1} imply that a 
flip does not change 
the set of positive configurations. The second claim 
follows from Lemma \ref{10.25.03.30} and the proof of  Theorem 
\ref{9.1.03.1}.  
\vskip 3mm

{\bf 5. A constructive proof of Theorem \ref{9.1.03.1n}}.  Here is a 
convenient way to spell the definition of 
a framed local system given by the map $\psi_{\Bbb T}$, which we will use to study 
its  monodromy. 
To simplify the exposition let us assume that $\widehat S = S$. To describe 
the answer we use the following picture. Starting from an ideal triangulation $T$ of 
$S$, we 
 construct a graph $\Gamma'_T$ embedded into the surface by drawing small edges 
transversal to each side of the triangles and inside each triangle 
connect the ends of edges pairwise by three more edges, as shown on 
the picture. The edges of the little oriented triangles are called 
  $t$-edges. The  edges dual to the edges of the
triangulation $T$ are called  $e$-edges.

\begin{figure}[ht]
\centerline{\epsfbox{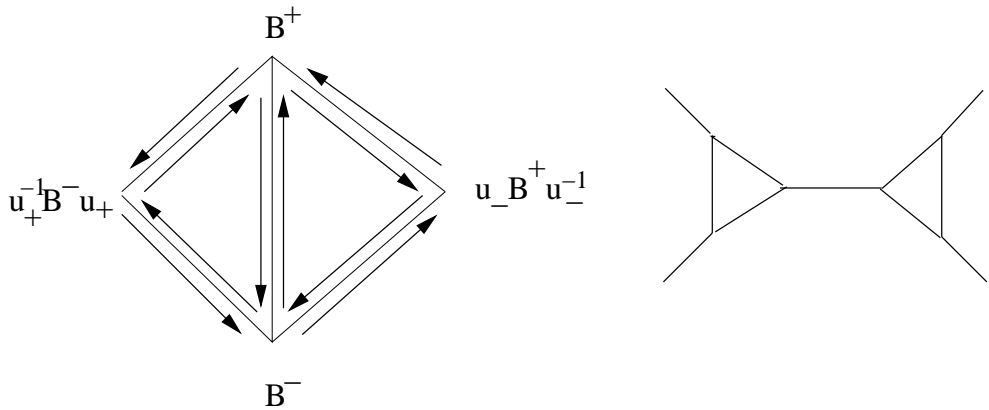}}
\caption{The graph $\Gamma_T'$.}
\label{fgo-10}
\end{figure}

Given a left principal homogeneous $G$-space $P$, a point $p_0 \in P$, and a
configuration $(p_1, p_2)$ of two pints in $P$ (that is, a pair of points in
$P$ considered up to the left diagonal action of $G$), there is an element
$g_{p_0}(p_1, p_2) \in G$ defined as follows:
\begin{equation} \label{3.17.04.1}
g_{p_0}(p_1, p_2) := g_1^{-1}g_2 \in G, \qquad \mbox{where } p_1 = g_1p_0, g_2
= g_2p_0.
\end{equation}
We will apply this construction to the principal homogeneous space of pinnings for $G$. 
Recall that 
we choose a standard pinning $P_0 = (B^-, B^+)$ for $G$. 
 Every edge of the graph $\Gamma'_T$ provides a configuration of two pinnings
in $G$, assigned naturally to the vertices of this edge. Indeed, the vertices
of the graph on the right of Figure \ref{fgo-10} match the pinnings assigned to the 
triangulation on the left. We use 
only the pinnings compatible with the orientation of $S$.  
Therefore we can assign to every oriented edge ${\bf a}$ of the graph
$\Gamma'_T$ 
an element $M({\bf a}) \in G$. Namely, it is the element assigned by (\ref{3.17.04.1}) to the 
configuration of two pinnings attached to the oriented 
edge ${\bf a}$, and the standard
pinning $P_0$. It follows from the very definitions that these elements enjoy
the following properties:

1) $M({\bf a}) M(\overline {\bf a}) = {\rm Id}$, where $\overline {\bf a}$
   denotes the edge as ${\bf a}$ taken with the opposite orientation. 

2) $M({\bf t_1}) M( {\bf t_2})M( {\bf t_3}) = {\rm Id}$, where the oriented
edges 
${\bf t_1},  {\bf t_2},  {\bf t_3}$ make a small  oriented triangle on the
graph $\Gamma'_T$. 

Thanks to 1)  the elements $\{M({\bf a})\}$ define a $G$-local system on the graph
$\Gamma'_T$. Thanks to 2) it is descended to a $G$-local system on the graph
$\Gamma$. Since $\Gamma$ is homotopy equivalent to $S$, we got a $G$-local
system on $S$. 
 
\vskip 3mm
{\bf 6. Proof of Theorem \ref{2.1.04.0}}. Given a loop on $S$, let us shrink it to a loop on
${\Gamma}'_T$. We may assume that this loop contains no 
consecutive $t$-edges: indeed, a composition of two $t$-edges is a
$t$-edge. Thus we may choose an initial vertex on the loop so that the
edges have the pattern $et ... etet$. 
Let $t$ and $e$ be consecutive $t$- and $e$-edges of the graph $\Gamma_T'$. 
Then there is a unique way to orient them so that   
$e$ goes after  $t$,  
i.e. the oriented edge ${\bf e}$ starts at the end of the one  ${\bf t}$. 
We say that a $t$-edge is oriented clockwise if it is 
compatible with the clockwise orientation of the corresponding triangle. 
 Therefore the monodromy is computed as a
product of elements of type  $M({\bf e})M({\bf t})$, where 
${\bf e}$ follows after ${\bf t}$.

 This construction  
defines the universal framed $G$-local system on $\widehat S$. 
Recall the semifield ${\Bbb F}^+_{G, S}$ of positive 
rational functions on the moduli space 
${\cal X}_{G, S}$ for the positive atlas defined above.

\begin{proposition} \label{3.15.04.1} 
Let $t$ and $e$ be two consecutive $t$- and $e$-edges of the graph $\Gamma_T'$. 
Let us orient them so that ${\bf e}$ goes after  ${\bf t}$. 
Then 
$$
M({\bf e})M({\bf t}) \in \quad \left\{ \begin{array}{ll} B^{+}({\Bbb F}^+_{G, S}) & 
\mbox{if ${\bf t}$ is oriented clockwise }
\\B^{-}({\Bbb F}^+_{G, S}) & \mbox{if ${\bf t}$ is oriented counterclockwise.}
\end{array}\right.
$$
\end{proposition}

\begin{figure}[htb]
\centerline{\epsfbox{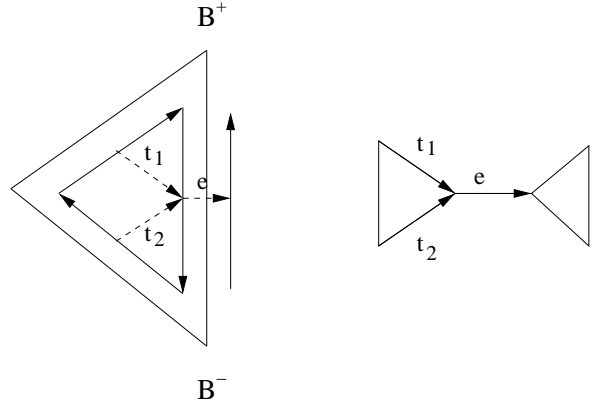}}
\caption{The oriented edges ${\bf t_1}$ and ${\bf t_2}$ on the graph $\Gamma_T'$.}
\label{fgo3-11}
\end{figure}

{\bf Proof}. We have to consider separately the cases when ${\bf t}$ is 
oriented clockwise, or counterclockwise, i.e.  ${\bf t} = {\bf t_1}$ or 
${\bf t} = {\bf t_2}$ where ${\bf t_1}$ and ${\bf t_2}$ are as on Figure 
\ref{fgo3-11}.

We denote by $[B_1, B_2, B_3]$ triples of flags, to distinguish from 
configurations of flags, and by  $g\cdot [B_1, B_2, B_3]$ the result of the action 
of $g \in G$ on the triple $[B_1, B_2, B_3]$.

{\it Case 1}: ${\bf t} = {\bf t}_1$. 
Let $u_+ = h_+(\overline u_+)$ and $u_- = h_-(\overline u_-)$, where 
$h_{\pm} \in H$ acts on $U^{\pm}$ by conjugation. 
We refer to Figure \ref{fgo-10} for the definition of $u_+$ and $u_-$. Then we have 
$$
\left[u_+^{-1}B^-u_+, B^+, B^-\right] = u_+^{-1}\cdot \left[B^-, B^+,
  u_+B^-u_+^{-1}\right]
= u_+^{-1}h_+\cdot \left[B^-, B^+, \overline u_+B^-\overline u_+^{-1}\right]
$$
$$
= u_+^{-1}h_+\cdot \left[B^-, 
B^+, \phi(\overline u_+)B^+ \phi(\overline u_+)^{-1}\right]. 
$$
On the other hand
\begin{equation} \label{3.15.04.2} 
\left[B^-, B^+, u_-B^+u_-^{-1}\right] = h_-\cdot \left[B^-, B^+, \overline 
u_-B^+\overline u_-^{-1}\right].
\end{equation}
Let $h(\phi(\overline u_+))(\overline {\phi(\overline u_+)}) = {\phi(\overline
  u_+)}$. Observe that $h(\phi(\overline u_+)) \in H$. 
Then 
$$
M({\bf e})M({\bf t}_1) = h(\phi(\overline u_+))^{-1}h_+^{-1}u_+h_- = 
h(\phi(\overline u_+))^{-1}\overline u_+h_+^{-1}h_- \in B^+({\Bbb F}^+_{G,S}).
$$
Indeed,  $h_+^{-1}h_-$ is the $H$-invariant, and $h(\phi(\overline
u_+))^{-1}$ is a positive rational function by positivity of $\phi$.

{\it Case 2}: ${\bf t} = {\bf t}_2$. We have 
$$
\left[B^-, u_+^{-1}B^-u_+, B^+\right] = \left[B^-, \phi(u_+)^{-1}B^+\phi(u_+),
  B^+\right] =
$$ 
$$
=\phi(u_+)^{-1}\left[B^-, B^+, \phi(u_+)B^+\phi(u_+)^{-1}\right] =
$$
$$ 
=\phi(u_+)^{-1}h_+\left[B^-, B^+, \phi(\overline u_+)B^+\phi(\overline u_+)^{-1}\right] =
$$ 
$$
=\phi(u_+)^{-1}h_+ h(\phi(\overline u_+))\left[B^-, B^+, 
\overline {\phi(\overline u_+)}B^+{\overline {\phi(\overline u_+)}}^{-1}\right]. 
$$
Comparing with  (\ref{3.15.04.2}), 
 and using the fact that \ $u_+ = h_+\overline u_+ h_+^{-1}$ \ implies 
$\phi(u_+) = h_+\phi(\overline u_+) h_+^{-1}$,  we get 
$$
M({\bf e})M({\bf t}_2) = h(\phi(\overline u_+))^{-1} h_+^{-1} \phi(u_+)h_+ h_+^{-1}h_- =
$$
$$ 
=h(\phi(\overline u_+))^{-1} \phi(\overline u_+)h_+^{-1}h_- \in B^-({\Bbb F}^+_{G,S}).
$$
Indeed,  $h_+^{-1}h_-$ is the $H$-invariant, and $\phi(\overline u_+)$ 
is positive.  The Proposition is proved. 
 
\vskip 3mm
Let us return to the proof of the theorem. 
A loop on $S$ contains elements of just one kind, $B^+$ or $B^-$, 
 if and only if it is homotopic to the loop
around a boundary component. Therefore the monodromy around any non-boundary
loop is obtained as a product of elements both from 
$B^-({\Bbb F}^+_{G,S})$  and $B^+({\Bbb F}^+_{G,S})$. 
Thus it lies in $G({\Bbb F}^+_{G,S})$. The theorem is proved.
\vskip 3mm

{\bf Example}. Let $G = PGL_2$, $u_+ = \left ( \begin{matrix}1& y\\ 0&1\end{matrix}\right )$, $u_- = \left (\begin{matrix}1& 0\\ x&1\end{matrix}\right )$. Then 
$$
(B^-, u_+^{-1}B^-u_+, B^+, u_-B^-u_-^{-1}) = (\infty, -y^{-1}, 0, x) = 
(\infty, -1, 0, xy). 
$$ 
The edge invariant equals  $xy$, and coincides with the cross-ratio 
of the corresponding $4$ points on ${\Bbb P}^1$. 
\vskip 3mm
{\bf 7. The universal  higher Teichm{\"u}ller ${\cal X}$-space related to $G$}. 
\begin{definition} \label{cstr1} A set $C$ is {\em cyclically ordered} if  any of its
  finite
  subset is cyclically ordered, and any inclusion of finite subsets
  preserves the cyclic order. 
\end{definition}
To define a cyclic order on a  set $C$ it is sufficient 
to specify  cyclic orders for all  triples and quadruples of
elements of $C$ so that the following compatibility condition holds: 
forgetting any element of any cyclically ordered quadruple we
get 
a triple whose induced cyclic order is the same as the prescribed one.  
Reversing a cyclic order on $C$ we get a new one. 

According to Theorems \ref{2.18.03.2a} and  \ref{2.18.03.2wq} 
the set of positive $n$-tuples of flags in ${\cal B}(\R)$ is
invariant under the operation of 
reversion of the cyclic order.

\begin{definition} \label{cstr} i) 
A {\it dihedral structure} on $C$ is a cyclic order considered up to reversion. 

ii) A map $\beta$ from a set $C$ with given dihedral structure  to the flag variety
  ${\cal B}(\R)$ is {\em positive} if it maps every cyclically ordered quadruple 
  $(a,b,c,d)$ to a positive 
  quadruple of flags $(\beta(a), \beta(b), \beta(c), \beta(d))$.
\end{definition}
\vskip 3mm
{\bf Remarks}. 1. Positive maps are  injective. We say that $C$ is a {\it positive subset} of the flag variety ${\cal B}(\R)$, if there exists a cyclically ordered 
set $C'$ and a positive map $\beta:C' \to {\cal B}(\R)$ whose image is $C$.

2. There exist quadruples of flags in
$\R {\Bbb P}^2$ such that 
forgetting any of them we get a positive triple, while the
quadruple itself is not positive. For example take four points $(x_1, x_2, x_3, x_4)$
going clockwise on a circle, and consider the tangent
flags  
to the points $(x_1, x_3, x_2, x_4)$ (see [FG3], Lemma 2.3 and 2.4). 
\vskip 3mm

Consider the dihedral 
structure on ${\Bbb P}^1(\Q)$ provided by the 
embedding ${\Bbb P}^1(\Q) \hra {\Bbb P}^1(\R)$.

\begin{definition} \label{1.29.04.3} The universal  Teichm{\"u}ller space ${\cal X}^+_G$ 
     is the quotient of the space of positive maps
\begin{equation} \label{1.29.04.1}
\beta: {\Bbb P}^1(\Q) \lra {\cal B}(\R)
\end{equation} 
by the natural action of the group $G(\R)$ on it. 
\end{definition} 

Let $\varphi: C \to {\cal B}(\R)$ 
be a positive map from a cyclically ordered set 
$C$. 
The group $G(\R)$ acts on the set of maps from $C$ to ${\cal B}(\R)$. 
If a map $\varphi: C \to {\cal B}(\R)$ is positive then for any $g \in G(\R)$ the map 
$g\varphi$ is also positive. The $G(\R)$-orbits on the set of positive maps are called 
{\it positive configurations of flags in ${\cal B}(\R)$}. 

So ${\cal X}^+_G$ is the space of positive configurations of real
flags in $G(\R)$ parametrised by ${\Bbb P}^1(\Q)$. 

\vskip 3mm
{\bf Canonical decomposition of ${\cal
  X}^+_G$}. To parametrise ${\cal
  X}^+_G$, 
let us choose three distinct points on ${\Bbb P}^1(\Q)$,
  called $0, 1, \infty$. Recall the Farey triangulation, understood 
as a triangulation of the hyperbolic disc with a
  distinguished oriented edge (= flag).  
Then we have  canonical identifications 
\begin{equation} \label{1.29.04.6} 
{\Bbb P}^1(\Q) = \Q \cup \infty = \{\mbox{vertices of the Farey
  triangulation}\}.
\end{equation} 
The distinguished oriented edge  
connects $0$ and $\infty$. 
So a point of ${\cal X}_G^+$ gives rise to  
a positive map 
$$
\{\mbox{vertices of the Farey triangulation}\} \lra {\cal B}(\R)
$$
considered up to the action of $G(\R)$. Assigning to
every Farey triangle the corresponding configuration 
of the three flags  at the
vertices, 
and to every Farey geodesic, oriented in a certain way,  the  $H$-invariant of the associated
quadruple of flags, 
we get a canonical map
\begin{equation} \label{1.29.04.7} 
\varphi_G: {\cal X}^+_G \stackrel{}{\lra} {{\cal X}^+_{G, 3}}^{\mbox{\{Farey
  triangles\}}} \times H(\R_{>0})^{\{\mbox{Farey
  diagonals}\}}.
\end{equation} 

\begin{theorem} \label{1.29.04.3q} The map (\ref{1.29.04.7}) is an isomorphism. 
\end{theorem} 

{\bf Proof}. It is completely similar to the one of Theorem
\ref{9.1.03.1n}, and thus is omitted. 
\vskip 3mm
{\bf Example}. When $G=PSL_2$ our universal Teichm\"uller space is essentially
the one considered by L. Bers \cite{Bers} and R. Penner [P3]. 
\vskip 3mm
Recall the {\it Thompson group} ${\Bbb T}$ of all piece-wise linear automorphisms of
${\Bbb P}^1(\Q)$. By definition,  for  every
element $g \in {\Bbb T}$ there is a decomposition of ${\Bbb P}^1(\Q)$
into a union of a finite number of segments, which may overlap only at the
ends,  so that restriction of $g$ to each segment is
given by an element of $PGL_2(\Q)$ (depending on the segment).  It acts
on  ${\cal X}^+_G$ in an obvious way. Each
element of the Thompson group can be presented as a composition of
flips. 

Let  $\Delta \subset PSL_2(\Z)$ be a torsion-free finite index subgroup. 
Denote by ${\cal H}$ the hyperbolic disc equipped with the Farey
triangulation. Set
$S_{\Delta}:= {\cal H}/\Delta$. The surface $S_{\Delta}$
inherits an ideal triangulation $T_{\Delta}$. The higher Teichm\"uller space 
${\cal X}^+_{G, S_{\Delta}}$ is embedded into the universal one 
 ${\cal X}_{G}^+$ 
as the subspace of $\Delta$-invariants. Moreover, let
\begin{equation} \label{1.29.04.7f} 
\varphi_{G, S_{\Delta}}: {\cal X}^+_{G, S_{\Delta}} \stackrel{\sim}{\lra} 
{{\cal X}^+_{G, 3}}^{\mbox{\{triangles of $T_{\Delta}$\}}} \times 
H(\R_{>0})^{\{\mbox{diagonals of $T_{\Delta}$}\}}
\end{equation} 
be the decomposition of the higher Teichm\"uller space according to a
triangulation $T_{\Delta}$ of $S$, where we assume 
an orientation of the edges of $T_{\Delta}$ is
chosen. This isomorphism is obtained by taking the $\R_{>0}$-points 
of the spaces in (\ref{9.1.03.1e}). Then the isomorphism (\ref{1.29.04.7f}) 
can be obtained by taking the ${\Delta}$-invariants of
the isomorphism  (\ref{1.29.04.7}).


\vskip 3mm
{\bf 8. $K$-positive equivariant configurations of flags}. 
Let $C$ be a finite cyclic set 
and $|C|>2$. Let  $G$ be a split semi-simple group over $\Q$.  
Recall the positive structure on the configuration space 
${\rm Conf}_C({\cal B})$ of maps $C \to {\cal B}$ modulo $G$-conjugation.  
An 
inclusion of finite cyclic sets $C_1 \subset C_2$ gives rise to a map of positive spaces, the restriction map: 
\begin{equation} \label{CCCa}
{\rm Res}_{C_2/C_1}: {\rm Conf}_{C_2}(\underline{\cal B}) \lra 
{\rm Conf}_{C_1}(\underline{\cal B}), 
\end{equation} 
Therefore there is a functor ${\rm Conf}$ 
from the category of finite 
cyclic sets, with morphisms given by inclusions 
of sets 
preserving the cyclic structure to the category of positive spaces. 
It takes a cyclic set $C$ to the positive space ${\rm Conf}_C({\cal B})$.

One can extend the functor ${\rm Conf}$ to a functor from 
the category ${\rm Cyc}$ of possibly infinite (countable) cyclic sets to the category whose objects 
are projective limits of positive spaces. We will use only the fact that, for any semifield $K$, 
there is a functor from the category ${\rm Cyc}$ to the category of sets, covariant with respect to $K$. 
It takes a cyclic set $C$ to a set 
${\rm Conf}_C({\cal B})(K)$;  
if $C_1 \subset C_2$ is an inclusion of cyclic sets, there is the restriction  map of sets  
$$
{\rm Res}_{C_2/C_1}: {\rm Conf}_{C_2}({\cal B})(K) \lra 
{\rm Conf}_{C_1}({\cal B})(K).
$$ 

\begin{definition} \label{CCC1} Let $\pi$ be a group acting by automorphisms of a cyclic set 
$C$. It acts therefore on the set ${\rm Conf}_C({\cal B})(K)$. The set of 
invariants of this action  is denoted by ${\rm Conf}_{C, \pi}({\cal B})(K)$. It is  called 
the set of $\pi$-equivariant $K$-positive configurations of flags parametrised by $C$. 
\end{definition}

So an element 
${\cal C}\in {\rm Conf}_C({\cal B})(K)$ is $\pi$-equivariant if and only if for any subset 
$C' \subset C$ the elements ${\rm Res}_{C'}{\cal C}$ and ${\rm Res}_{\pi(C')}{\cal C}$ coincide 
in ${\rm Conf}_{C'}({\cal B})(K)$. 

\vskip 3mm
{\bf Remark 1}. The set ${\rm Conf}_C({\cal B})(K)$ is not necessarily the set of $K$-points 
of a positive space in the sense of the definition from Section 4. It is  the set of $K$-points 
of a {\it generalised positive space}, but we will not pursue its definition here, 
working with its $K$-points only.  

\vskip 3mm
{\bf Remark 2}. Let $\Gamma$ be a group acting by automorphisms on a cyclic set $C$ and 
on a group $\pi$. Then, for any semifield $K$, the group 
 $\Gamma$ acts on the set ${\rm Conf}_{C, \pi}({\cal B})(K)$. 
\vskip 3mm

The following Lemma is an immediate consequence of Lemma \ref{q1-}. 
\begin{lemma} \label{CCC4} Suppose that $G$ has trivial center. 
Let us assume that $K$ is a semifield in a field $F$, e.g. $K=\R_{>0}, F=\R$. 
Then an element ${\cal C}\in {\rm Conf}_{C, \pi}({\cal B})(K)$ gives rise to 
a homomorphism
$$
\rho_{\cal C}: \pi \lra G(F)\qquad \mbox{well defined modulo $G(F)$-conjugation}. 
$$
\end{lemma}

\vskip 3mm 
{\bf 9. An application: higher Teichm\"uller spaces and laminations for closed surfaces}. 
Recall (Sections  1.3, 1.10) the three different 
cyclic $\pi_1(S)$-sets assigned to a surface $S$: 
\begin{equation} \label{11.06.05.7} 
{\cal F}_{\infty}(S) \subset {\cal G}_{\infty}(S) \subset \partial_{\infty}\pi_1(S)
\end{equation}
For a surface without boundary the first one 
 is empty. So for any surface $S$, open or closed, and any semifield $K$,
 there is a set ${\rm Conf}_{{\cal G}_{\infty}(S), \pi_1(S)}( {\cal B})(K)$ with a natural action of the mapping class group on it. 

The higher Teichm\"uller space ${\cal X}_{G,S}^+$ was defined for open 
surfaces in Definition \ref{7.29.03.1}, and for closed surfaces in Definition \ref{4.1.04.10a}. 
Here is an ``algebraic-geometric''  definition, which 
treats 
simulteneously higher Teichm\"uller spaces and lamination spaces for 
open and closed surfaces.

\begin{definition} \label{CCC2} Let $S$ be a surface, open or closed. 
Let $G$ be a split semi-simple algebraic group over $\Q$ with trivial center. 

(i) The higher Teichm\"uller space is the space
 ${\rm Conf}_{{\cal G}_{\infty}(S), \pi_1(S)}({\cal B})({\R}_{>0})$. 

(ii) Let ${\Bbb A}^t$ be one of the three tropical semifields $\Z^t$, $\Q^t$, $\R^t$. 
The set ${\rm L}_{G,S}({\Bbb A})$ of $G$-laminations on $S$ with coefficients in the ring ${\Bbb A}$ is the set 
${\rm Conf}_{{\cal G}_{\infty}(S), \pi_1(S)}({\cal B})({\Bbb A}^t)$. 
\end{definition}

We will show in Section 7 
 that the higher Teichm\"uller space in this definition is 
the same as the one defined in Definitions \ref{7.29.03.1} and \ref{4.1.04.10a}. 
Further, the set ${\rm L}_{PGL_2, S}({\Bbb A})$ coincides with the 
space of Thurston's integral (${\Bbb A}=\Z$), rational (${\Bbb A}=\Q$) or measured (${\Bbb A}=\R$) laminations on $S$ 
\vskip 3mm
{\bf Remark}. For a closed surface $S$, ${\rm Conf}_{{\cal G}_{\infty}(S), \pi_1(S)}( {\cal B})$ is a generalised positive space 
which is not a positive space in the sense of the definition from Section 4. 

\section{Higher Teichm\"uller spaces and their basic properties}
\label{higher}


In this Section we extend the definition of higher Teichm\"uller spaces 
to closed surfaces, and establish basic properties of higher Teichm\"uller spaces 
for arbitrary surfaces, with or without boundary. 
In the very end we propose a conjectural generalisation of the theory of quasifuchsian and 
Kleinian groups 
where $PSL_2(\C)$ is replaced by an arbitrary semi-simple complex Lie group.  

\vskip 3mm
{\bf 1. Discreteness of positive representations: Proof of Theorem \ref{7.8.03.2}}. 
Take a positive framed $G(\R)$-local system on a
hyperbolic surface $S$ with boundary. Let $T$ be an ideal triangulation of $S$. 
Consider a universal cover $\widetilde{S}$ and 
the corresponding triangulation $\widetilde T$. The triangulated 
surface $(\widetilde{S}, \widetilde T)$  
is isomorphic to  the hyperbolic disc ${\cal H}$ with the 
Farey triangulation. 
Choose a point 
$\widetilde{p}$ on $\widetilde{S}$ and a trivialization of the local system 
at $\widetilde p$. Let $V(\widetilde{T})$ be the set of vertices of 
$\widetilde{T}$.  Let us construct a map
\begin{equation} \label{11.9.06.1}
\beta: V(\widetilde{T}) = {\cal F}_{\infty}(S) \rightarrow {\cal B}(\R).
\end{equation}
 Let $\widetilde{v}$ be a vertex of 
$\widetilde{T}$ and $\widetilde{\alpha}$ be a path from $\widetilde{p}$ to 
$\widetilde{v}$. Let $p,v$ and $\alpha$ be the corresponding data on 
$S$. Let $\beta$ be a small counterclockwise loop on $S$ 
surrounding the point $v$ and intersecting with $\alpha$ at one point $r$. 
Consider the flag at the fiber of the local system  over $r$ 
provided by the framing. The part of the path $\alpha$ between $p$ and  $r$ maps it
 to the fiber over $p$,  thus giving the desired map $\beta$.

The map $\beta$ has the following properties.
\vskip 3mm
1. Changing the basepoint $\widetilde{p}$ or trivialization at $\widetilde p$ 
results in a change $\beta 
\rightarrow g\cdot \beta$ for a certain $g \in  G(\R)$. 

2. The map $\beta$ is $\pi_1(S)$-equivariant. 

3. $\beta$ is a positive map. 
\vskip 3mm
The first two properties follow directly from the construction. The third 
one follows from the definition of a positive framed local system if 
the triple / quadruple of vertices of $\widetilde{T}$ constitute vertices of a 
triangle or quadrilateral. If they do not, 
there exists a finite polygon of the triangulation $\widetilde T$ 
containing all these vertices. One can 
change the triangulation of this polygon by a finite number of flips 
so that  the original triple/quadruple  of points become vertices of one 
triangle/quadrilateral of the new triangulation. According to 
Theorem \ref{2.18.03.2} if a local system is positive 
with respect to   one 
triangulation, it is positive with respect to   any other one, hence the property.  

Let us fix a pair $(B^-, B^+)$ of opposite Borel 
 subgroups in $G(\R)$. Consider the following two subsets 
in the real flag variety:
$$
{\cal B}^+ := \{u_-\cdot B^+\}, \qquad {\cal B}^- := \{v^{-1}_-\cdot B^+\}, \quad u_-, v_- \in U^-(\R_{>0}).
$$
They are open in ${\cal B}(\R)$. Let $\overline {\cal B}^+$ 
be the closure of ${\cal B}^+$ in  ${\cal B}(\R)$. 

\begin{lemma} \label{3.8.04.1}
The sets    $\overline {\cal B}^+ $ and ${\cal B}^- $ are disjoint.
\end{lemma}

{\bf Proof}. One has $u_-B^+u_-^{-1} = v^{-1}_-B^+v_-$ 
if and only if  $u_- = v^{-1}_-$, i.e. $u_-v_- =1$.
Let us show that this is impossible. 
Recall the projection $\pi: U^- \to U^-/[U^-, U^-]$. 
Then $\pi(u_-) = \prod_{i}x_{i}(t_i)$ and $\pi(v_-)= \prod_{i}x_{i}(s_i)$, 
where $t_i>0, s_i \geq 0$, and 
the product is over all simple positive roots. 
Thus $0 = \pi(u_-v_-) = \prod_{i}x_{i}(t_i+s_i)$. But $t_i+s_i>0$. 
This contradiction proves the lemma. 
\vskip 3mm

Now fix a triangle of $\widetilde{T}$  and let $a,b,c$ be its vertices. 
The deck transformation 
along a loop $\gamma$ sends the triangle $(abc)$ to a triangle  $(a'b'c')$ 
without common interior points with $(abc)$. 
Without loss of generality we can assume that 
the vertices of the triangle $(a'b'c')$ are at the arc between $a$ and $c$ which  
does not  contain $b$. 
However the point $b'$ may coincide with $a$ or $c$. 
Let us assume first this does not happen. Then, since $\beta$ is positive, 
we have a positive quadruple of flags 
$(\beta(a), \beta(b), \beta(c), \beta(b'))$, and  $\beta(b') = \rho(\gamma) \beta(b)$.
 \begin{figure}[ht]
\centerline{\epsfbox{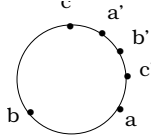}}
\caption{A deck transformation moves the triangle $abc$ to the one $a'b'c'$.}
\label{fg80}
\end{figure}
We may assume that  
$\beta(a) = B^-, \beta(c) = B^+$. Then $\beta(b) \in {\cal B}^-$ and 
$\beta(b') \in {\cal B}^+$. Since ${\cal B}^-$ is open, 
there exists an open neighborhood $O$ of the unit in $G(\R)$ such that 
$g\beta(b) \in {\cal B}^-$ for any $g \in O$. Since 
${\cal B}^+$ is disjoint with ${\cal B}^-$, $\rho(\gamma) \not \in O$. 
In the case when $b'$ coincides with $a$ or $c$, we may 
assume without loss of generality that it coincides with $c$.  Then the argument 
is the same except that in the very end we use the fact that 
${\cal B}^-$ does not contain $B^+$, which is a special case of the lemma. 
The theorem  is 
proved. 
\vskip 3mm
{\bf Proofs of Theorems \ref{4.1.04.10} and \ref{dgdgd}}. Follow from the above
constructions.

\vskip 3mm
{\bf 2. Two results on Higher Teichm\"uller spaces.} According to Theorem \ref{dgdgd},  
a framed positive local system $({\cal L}, \beta)$ on a surface $S$ with boundary 
is described uniquely by a $\pi_1(S)$-equivariant positive map 
\begin{equation} \label{11.06.05.8} 
\Phi_{\cal L, \beta}: {\cal F}_{\infty}(S) \lra {\cal B}(\R) \quad \mbox{modulo $G(\R)$-conjugation}, 
\end{equation}
where $\pi_1(S)$ acts on the flag variety  via the monodromy 
representation of the local system ${\cal L}$. 

\begin{theorem} \label{11.06.05.18} The map (\ref{11.06.05.8})  has a unique 
extension 
to a positive $\pi_1(S)$-equivariant map 
\begin{equation} \label{11.12.05.1} 
\Psi_{\cal L, \beta}: 
{\cal G}_{\infty}(S) \lra {\cal B}(\R).
\end{equation} 
\end{theorem}

We will prove this theorem in Sections 7.3-7.5.

\begin{corollary} \label{CCC3} 
Let $S$ be a surface, with or without boundary. 
Then there is a canonical isomorphism, equivariant with respect to the action  of the mapping class group 
of $S$: 
$
{\cal X}_{G,S}^+ \stackrel{\sim}{\lra} {\rm Conf}_{{\cal G}_{\infty}(S), \pi_1(S)}( {\cal B})(\R_{>0})
$. 
\end{corollary}

{\bf Proof}. Let us present an element of 
${\rm Conf}_{{\cal G}_{\infty}(S),\pi_1(S)}( {\cal B})(\R_{>0})$ 
by a positive map $\Psi: 
{\cal G}_{\infty}(S)\to {\cal B}(\R)$. Since $\R_{>0}\subset \R$,  by Lemma \ref{CCC4} 
it gives rise to a representation 
$\rho: {\cal G}_{\infty}(S)\to G(\R)$, and $\Psi$ is $\rho$-equivariant. 
This proves the claim for closed $S$. For open $S$, one needs in addition Theorem \ref{11.06.05.18}. 
The corollary is proved. 

\begin{corollary} \label{CCC30} 
Let $S$ be a surface, with or without boundary. 
Then the principal $PGL_2$-embedding provides an embedding 
$
{\cal X}_{PGL_2,S}^+ \hra {\cal X}_{G,S}^+
$. 
\end{corollary}

{\bf Proof}. By Corollary \ref{cb2} there is a canonical embedding of the configuration spaces 
${\rm Conf}^+_k({\Bbb P}^1) \hra {\rm Conf}^+_k({\cal B})$. It remains to use Corollary \ref{CCC3}. 
The corollary is proved. 

\vskip 3mm
Recall Definition \ref{4.1.04.10a} of positive local systems for surfaces $S$ with or without boundary. 
It is clear from Theorems \ref{dgdgd} and \ref{11.06.05.18} that,  
in the case when  $S$ has holes, ${\cal L}_{G,S}^+$ is the image of ${\cal X}_{G,S}^+$ 
under the canonical projection ${\cal X}_{G,S}(\R) \lra {\cal L}_{G,S}(\R)$.  
For closed $S$ one has ${\cal L}_{G,S}^+ = {\cal X}^+_{G,S}$.

\begin{theorem} \label{11.20.05.1}
The map (\ref{11.12.05.1}) extends uniquely to a continuous positive map 
\begin{equation} \label{11.12.05.2} 
\overline \Psi_{\cal L, \beta}: 
{\partial}_{\infty}\pi_1(S) \lra {\cal B}(\R).
\end{equation}
\end{theorem}
This map is continuous by Theorem \ref{5.9.03.100} below. 
We will not use it.

Let us prepare the ground for the proof of Theorem \ref{11.06.05.18}. 
\vskip 3mm

{\bf 3. Positive configurations of flags revisited}. 
\begin{theorem} \label{11.04.05.4} 
For any $m+n \geq 1$, the map 
\begin{equation} \label{11.04.05.6}
\begin{array}{c}
\underbrace{U^+ \times \ldots \times U^+}_{\mbox{$m+n$ \rm copies}}/H \lra {\rm Conf}_{m+n+2}({\cal B}),
\qquad (v_1, \ldots , v_m, u_1, \ldots , u_n) \lms \\

\Bigl(B^+, (v_1 ... v_m)^{-1}\cdot B^-, \ldots , (v_{m-1}v_m)^{-1}\cdot B^-,  v^{-1}_m\cdot B^-,\\
         B^-, u_1\cdot B^-, u_1u_2\cdot B^-, \ldots  , u_1 ... u_n\cdot B^-\Bigr),
\end{array}
\end{equation}
where the quotient is by the diagonal action of $H$, 
provides a positive atlas on ${\rm Conf}_{m+n+2}({\cal B})$. 
\end{theorem}

{\bf Proof}. Let us first prove the claim of the theorem in the special case $m=0$. 
When $n=1$, it follows from the very definition of the positive atlas on ${\rm Conf}_{3}({\cal B})$. 
In the case $n=2$ we  define a positive atlas on ${\rm Conf}_{4}({\cal B})$ by using a map 
$$
(U^+ \times U^+)/H \lra {\rm Conf}_{4}({\cal B}), \quad (v,u) \lra (B^+, v^{-1}\cdot B^-, B^-, u\cdot B^-).
$$
Acting by $v$, we transform the last configuration to the one $(B^+, B^-, v\cdot B^-, vu\cdot B^-)$, 
proving the claim. Now let $n>2$. 
Let $(B_1, ..., B_{n+2})$ be a configuration of flags, and $B_1 = B^+, B_2 = B^-$.  
Let us consider an $(n+2)$-gon whose vertices are decorated by the flags $B_i$, 
so that the cyclic order of the vertices coincides with the one induced by the order of the flags.  
Consider the triangulation of the $(n+2)$-gon by the diagonals 
from the vertex decorated by $B^+$, see Figure \ref{fgo-10}. 
Then the claim of the theorem for $n=1$ and $n=2$ implies that the map 
$$
\underbrace{U^+ \times \ldots  \times U^+}_{\mbox{$n$ copies}}/H \lra {\rm Conf}_{n+2}({\cal B});
$$
$$
(u_1, \ldots , u_{n}) \lms (B^+, B^-, u_1\cdot B^-, u_1u_2\cdot B^-, \ldots , u_1 ... u_n\cdot B^-)
$$
is positive.  
In the general case we observe that multiplying the configuration (\ref{11.04.05.6}) by 
the element $v_1...v_m \in U^+$ we transform it to an equivalent configuration
$$
\Bigl(B^+,      B^-, v_1 \cdot B^-, \ldots , v_1 ... v_m\cdot B^-,  v_1 ... v_mu_1\cdot B^-, \ldots  , v_1 ... v_m u_1 ... u_n\cdot B^-\Bigr).
$$

The theorem follows. 
\vskip 3mm

\begin{corollary} \label{10.30.05.3}
 A configuration of flags, finite or infinite, is positive if and only 
if it can be written as 
\begin{equation} \label{10.30.05.4}
\Bigl(B^+, B^-, u_1\cdot B^-, u_1u_2\cdot B^-, \ldots 
 , u_1 ... u_n\cdot B^-, \ldots\Bigr), \quad u_i \in U^+(\R_+).
\end{equation}
\end{corollary} 
This shows that Definition \ref{pcf} is equivalent to the one we used in Section 5. 
\vskip 3mm
\begin{lemma} 
The $(n+2)$-tuple of flags (\ref{10.30.05.4}) representing a positive configuration 
is well defined up to the action of the positive tori $H(\R_{>0})$. 
\end{lemma}

{\bf Proof}. Suppose we have an $(n+2)$-tuple of flags 
corresponding to $(u_1', \ldots , u_n')$, where  $u'_i \in U(\R_{>0})$, representing 
the same configuration (\ref{10.30.05.4}). Then there exists an element of $H(\R)$ 
conjugating one to the other. The projection of  $u$ to the maximal abelian quotient 
$U^+/[U^+, U^+]$ has natural positive coordinates  $u_{\alpha}$ parametrised by the simple roots $\alpha$. 
The element $t \in H(\R)$ acts on $u_{\alpha}$ by multiplication on
$\chi_{\alpha}(t)$. Thus $\chi_{\alpha}(t)>0$ for every simple root $\alpha$. Since $G$ has trivial 
center, this implies that $t \in H(\R_{>0})$. The lemma is proved.

\vskip 3mm
{\bf 4. Semi-continuity of positive subsets of real flag varieties}. 

\begin{definition} \label{6.26.04.1} Let $C$ be a cyclically ordered set. 
  A map $\varphi: C \to {\cal B}(\R)$ is {\em semi-continuous from the left} if 
for any sequence 
of distinct points $c_0, ..., c_n, ...$ in $C$, whose order is compatible with the cyclic
 order of $C$, the limit $\lim_{n \to \infty}\varphi (c_n)$ exists. 

A positive configuration in ${\cal B}(\R)$ is {\em semi-continuous from the left} 
if a  map $\varphi: C \to {\cal B}(\R)$ representing it is semi-continuous from the left. 
\end{definition}

Changing the cyclic order to the opposite one we define maps 
{\em semi-continuous from the right}. 

\begin{theorem}\label{5.9.03.100} A positive map $\varphi: C \to {\cal B}(\R)$ 
is semi-continuous from the left and from the right. 

If $\lim_{n\to \infty}\varphi(c_n)\not = \varphi(c_i)$ for $i \geq 0$, then 
adding the limit flag we get a positive configuration. 
\end{theorem}

Precisely, in the second statement we define a new cyclic set 
$C':= (c_{\infty}, c_0, c_1, \ldots)$, and extend the map $\varphi$ to a map $\varphi'$ on $C'$ by setting 
$\varphi'(c_{\infty}):= \lim_{n\to \infty}\varphi(c_n)$.

\vskip 3mm
{\bf Proof}. Set $B_i := \varphi(c_i)$. By Corollary \ref{10.30.05.3} 
the configuration of flags \newline $(B_2, B_3, B_4,  ..., B_0, B_1)$ can be written as 
\begin{equation} 
(B^-, v_3\cdot B^-, v_4\cdot B^-, ..., v_n \cdot B^-, ..., v_0\cdot B^-, B^+);\ 
v_k = u_3u_4 ... u_k,\  u_i, v_0 \in U^+(\R_{>0}).
\end{equation} 
Moreover we may assume that for every $n \geq 3 $ one has 
$v_0 = v_n v'_n$ with $v'_n \in U^+(\R_{>0})$.

By Proposition 3.2 in \cite{L1}, 
for any finite dimensional representation $\rho$ of $G$ in a vector space $V$, 
any element $u$ of $U(\R_{>0})$ acts in  the canonical basis in $V$ by a 
matrix with non-negative elements. 
The values of each  of the matrix elements in the sequence 
$\rho(v_i)$, $i=3, 4, ...$,  provide a non-decreasing sequence of real numbers 
bounded by the value of the corresponding matrix element of $\rho(v_0)$. 
Thus there exists the limit $\lim_{i \to \infty}\rho(v_i)$. 
Further, if $u \in U(\R_{>0})$ and $v \in U(\R_{\geq 0})$, then $uv \in U(\R_{>0})$. 
The Theorem is proved.
 \vskip 3mm
 

{\bf 5. Proof of Theorems \ref{11.06.05.18} and \ref{11.20.05.1}.} 
Recall the following result of Lusztig (\cite{L1}, Theorem 5.6), 
generalizing the Gantmacher-Krein theorem 
in the case $G = GL_n(\R)$. 

\begin{theorem} \label{10.06.05.2}
 Let $g \in G(\R_{>0})$. Then there exists a unique $\R$-split maximal 
torus of $G$ containing $g$. In particular, $g$ is regular and semi-simple. 
\end{theorem}

\begin{corollary} \label{10.06.05.3}
 Let $g \in G(\R_{>0})$. Then there is a unique real flag $B^{\rm at}$ (the attracting flag) 
such that 
for any flag $B$ in a Zariski open subset of ${\cal B}(\C)$ containing $B^{\rm at}$ 
the sequence $\{g^nB\}$, $n \to \infty$,  converges to the flag $B^{\rm at}$. 
\end{corollary}

{\bf Proof}. It follows from Theorem \ref{10.06.05.2} that $g$ is conjugate 
to the unique element $h\in H(\R_{>0})$ such that $\chi_{\alpha}(h) >1$ 
for every (simple) positive root $\alpha$, where $\chi_{\alpha}$ is the root 
character corresponding to $\alpha$. The flags containing $h$ are  $wB^+w^{-1}$, $w \in W$. 
The flag $B^+$ is the unique flag with the desired property. 
The Lemma is proved. 

\vskip 3mm
{\bf Proof of Theorem \ref{11.06.05.18}}. 
Let us choose an orientation of a closed ge\-o\-de\-sic $\gamma$. Let 
$M_\gamma$ be the monodromy of the local system along the oriented loop $\gamma$. 
Take a lift $\widetilde \gamma$ of  $\gamma$ to the universal cover. 
Let $\mu_\gamma$ be the deck transformation along the oriented path $\widetilde \gamma$. 
Let $p_+$ be the  attracting point for $\mu_\gamma$. It is an end of $\widetilde \gamma$. Then take  
two points $x_+$ and $x_-$ in ${\cal F}_{\infty}(S)$, 
located on the different sides of $\widetilde \gamma$. Applying 
transformations $\mu^n_{\gamma}$ we get two sequences of points,  $\{x^n_+\}$ 
and $\{x^n_-\}$,  converging to  $p_+$ from two different sides. 
The configuration of flags 
$
\Phi_{\cal L, \beta}(x^m_+, x^1_-, \ldots , x^n_-, \ldots )
$ 
 is positive. 
Thus by Theorem \ref{5.9.03.100} there exists $\lim_{n\to \infty}\Phi_{\cal L, \beta}(x^n_-)$. It  
must coincide with the unique attracting flag  $B^{\rm at}(M_{\gamma})$ 
for the monodromy operator $M_{\gamma}$. Similarly for 
the sequence $\{x^n_+\}$. Thus  
$$
\lim_{n\to \infty}\Phi_{\cal L, \beta}(x^n_-) = \lim_{n\to \infty}\Phi_{\cal L, \beta}(x^n_+) 
= B^{\rm at}(M_{\gamma}).
$$
We set $\Psi_{\cal L, \beta}(p_+)$ to be equal to this flag.  
Observe that $B^{\rm at}(M_{\gamma})$ does not coincide with $\Phi_{\cal L, \beta}(x^n_\pm)$.
Indeed, the latter ones are not stable by  $M_{\gamma}$, since $x^n_\pm$ is not stable by $\mu_\gamma$. 
The positivity of $\Psi_{\cal L, \beta}(x^m_+, x^1_-, \ldots , x^n_-, \ldots, p_+)$ follows now 
from Theorem \ref{5.9.03.100}.
The Theorem is proved. 

\vskip 3mm
{\bf Proof of Theorem \ref{11.20.05.1}.} Let us choose a hyperbolic metric on $S$, and 
a geodesic $\alpha$ on ${\cal H}$ ending at a point 
$p_+ \in \partial_\infty\pi_1(S) \subset \partial {\cal H}$. We 
assume that $p_+ \not \in {\cal G}_\infty\pi_1(S)$. 
Choose a fundamental domain ${\cal D}$ for $\pi_1(S)$ acting by the 
deck transformations on ${\cal H}$.  
Choose elements  $g_i \in \pi_1(S)$ such that 
$g_i{\cal D}$ intersects $\alpha$, 
and, as $i \to \infty$, approaches $p_+$. Denote by $\lambda_i$ the distance between the domains 
${\cal D}$ and $g_i{\cal D}$. According to Section 5, the images of non-boundary 
elements of $\pi_1(S)$ under a positive representation $\rho$ are conjugate to 
elements of $G(\R_{>0})$. In particular, this is so 
for $g_i$, so for any real $\lambda$ 
there exists $\rho(g_i)^{\lambda}$. 

\begin{lemma} \label{11.20.05.2} Let $\rho$ be a positive representation of $\pi_1(S)$ to $G(\R)$. 
Then the set $\{\rho(g_i)^{1/\lambda_i}\}$ is bounded. 
Thus there exists a subsequence $\{g_j\}$ for which there exists the limit
 $g_{\alpha}:= \lim_{j\to \infty}\rho(g_j)^{1/\lambda_j}$. 
\end{lemma}

{\bf Proof}. 
The group $\pi_1(S)$ equipped with the word metric is quasiisometric to 
${\cal H}$ (cf. \cite{Bon}). Since it is finitely generated, the Lemma follows. 

\vskip 3mm 

The element $\rho(g_{\alpha})$ is conjugate to an element of $H^0(\R_{>0})/W$. 
Thus there is a unique attracting flag $B^{\rm at}(\rho(g_{\alpha}))$ for it. 
Further, let $x^+$ and $x^-$ be 
the endpoints of two geodesics lifting closed geodesics on $S$. We assume that $x^+$ and $x^-$ are 
 located on different sides of $\alpha$. Then  
$x_j^\pm:= g_jx^\pm$ are also endpoints of liftings of closed geodesics, one has 
$\lim_{j \to \infty}g_jx^\pm = p_+$, and 
\begin{equation} \label{11.20.05.3q} 
\lim_{j\to \infty}\Phi_{\cal L, \beta}(x_j^+) = B^{\rm at}(\rho(g_{\alpha})); \qquad 
\lim_{j\to \infty}\Phi_{\cal L, \beta}(x_j^-) = B^{\rm at}(\rho(g_{\alpha}))
\end{equation}
We define 
$\overline \Phi_{\cal L, \beta}(p_+)$ to be the attracting flag of 
$\rho(g_{\alpha})$. Thanks to (\ref{11.20.05.3q}), 
it is the limit of flags assigned by $\Psi_{{\cal L}, \beta}$  
to the points of ${\cal G}_\infty\pi_1(S)$ converging to $p_+$. 
In particular it is independent of the choice of $\alpha$. 
The theorem is proved. 

\vskip 3mm 
{\bf 6. The cutting and gluing maps on the level of Teichm\"uller spaces.}  
Let $(S, \mu)$ be a surface with a 
hyperbolic metric $\mu$ and a geodesic boundary. Cutting $S$ along a non-boundary 
geodesic $\gamma$ we obtain a surface $S'$ with the induced metric. It has a geodesic boundary. 
Thus we have a map of Teichm\"uller spaces ${\cal T}_S \to {\cal T}_{S'}$, 
which we call the {\it cutting along $\gamma$} map. On the other hand, if $S'$ is a 
surface consisting of one or two components, with geodesic (and non-cuspidal) boundary, 
and the lengths of two boundary components 
$\gamma_1$ and $\gamma_2$ coincide and are different from zero, (i.e. the corresponding ends 
 are not cusps), then  we can glue 
$S'$ by gluing these curves, getting a surface $S$. There is a one-parameter  
family of hyperbolic structures 
on $S$ which, being restricted to $S'$, give the original hyperbolic 
structure there. It is a principal homogeneous space over the additive group $\R$. Thus, if 
${\cal T}_{S'}(\gamma_1, \gamma_2)$ is the subset of the space of hyperbolic structures on $S$ 
such that the lengths of the boundary geodesics $\gamma_1$ and $\gamma_2$ 
 are different from zero and coincide, 
then there is a principal $\R$-fibration
$$
{\cal T}_{S} \lra {\cal T}_{S'}(\gamma_1, \gamma_2).
$$
In particular this allows us to describe ${\cal T}_{S}$ 
if we know ${\cal T}_{S'}$. 

Below we generalize this 
to higher Teichm\"uller spaces. 

\vskip 3mm

Let $\gamma$ be the homotopy class of 
a simple loop on a surface $S$, 
which is neither trivial nor  
homotopic to a boundary component. Cutting $S$ along $\gamma$ we get a new surface $S'$. 
It is either connected or has two connected components. The surface $S'$  has two new boundary components, 
denoted $\gamma_+$ and $\gamma_-$.  Their orientations are 
induced by the orientation of $S'$. 

\vskip 3mm
{\bf The induced framing}. 
Let $({\cal L}, \beta)$ be a positive framed $G(\R)$-local system on $S$. Let 
${\cal L}'$ be the restriction of ${\cal L}$ to $S'$. There is a canonical framing 
$\beta'$ on ${\cal L}'$ extending $\beta$. Namely, 
the framing at the boundary components of $S'$ inherited from $S$ is given by $\beta$. 
According to Theorem \ref{2.1.04.0}, the monodromy $M_{\gamma}$ along the loop $\gamma$ is conjugate to an element of $G(\R_{>0})$. Thus we can apply Corollary \ref{10.06.05.3}, and take the attracting flag 
for the monodromy $\mu_{\gamma}$ as the framing. Changing the orientation 
of $\gamma$ we arrive at a different flag. Thus, using the canonical orientations of the boundary 
loops 
 $\gamma_+$ and $\gamma_-$, we get 
canonical framings assigned to the boundary components $\gamma_+$ and $\gamma_-$. 
We call $({\cal L}', \beta')$ the 
{\it induced framed local system on $S'$}. 
\vskip 3mm
Similarly, let $S$ be an oriented surface with boundary. Let ${\cal L}$ be positive 
$G(\R)$-local system on $S$. 
Then there is a {\it canonical framing}  on   ${\cal L}$
defined as follows. For each boundary component $\gamma$, whose orientation is induced by the one on $S$, 
take the unique attracting flag invariant by the monodromy along $\gamma$. 
The existence of such a flag is clear when the monodromy is positive hyperbolic, 
and in general follows from the fact that the monodromy around a boundary 
component is conjugate to an element of $B(\R_{>0})$. 
This way we get a canonical embedding 
\begin{equation}
{\cal L}^+_{G,S} \hra {\cal X}^+_{G,S} 
\end{equation}
\vskip 3mm
Recall the complement  $H^0$ to the hyperplanes ${\rm Ker}\chi_{\alpha}$ in $H$, where 
$\chi_\alpha$ runs through all root characters of $H$. It is the part of $H$ where the 
Weyl group $W$ acts freely. 

\begin{definition} \label{11.06.05.5} The submanifold 
$
{\cal X}^+_{G,S'}(\gamma_+, \gamma_-) \subset {\cal X}_{G,S'}(\R_{>0})
$ 
is given  by the following conditions:

(i) The monodromies along the oriented loops $\gamma_+$ and  $\gamma_-$ 
are mutually inverse, and lie  in $H^0(\R_{>0})/W$. 

(ii) The framing along the boundary component $\gamma_+$ (resp. $\gamma_-$) 
is given by the attracting flag for the monodromy $\mu_{\gamma_+}$ (resp. $\mu_{\gamma_-}$). 
\end{definition}

\begin{theorem} \label{10.30.05.2q} Let $S$ be a surface with  or without boundary, $\chi(S) <0$. 
 Let $\gamma$ be a non-trivial loop on $S$, non-homotopic to a boundary component, 
cutting $S$ into one or two 
surfaces. Denote by $S'$ their union. 
Then  

(i) The induced framed  $G(\R)$-local system 
$({\cal L}', \beta')$ on $S'$ is positive, i.e. lies in 
${\cal X}_{G,S'}(\R_{>0})$. There is a well-defined map, called the cutting map:
\begin{equation} \label{11.06.05.15}
C_{\gamma}: {\cal X}_{G,S}(\R_{>0}) \lra {\cal X}^+_{G,S'}(\gamma_+, \gamma_-); 
\qquad ({\cal L}, \beta)\lms ({\cal L}', \beta').
\end{equation} 

(ii) The cutting map  (\ref{11.06.05.15})
is a principal $H(\R_{>0})$-bundle. 

(iii) The space ${\cal X}^+_{G,S'}(\gamma_+, \gamma_-)$ is diffeomorphic to a ball of dimension 
$-\chi(S){\rm dim G}- {\rm dim}H$. 
\end{theorem} 

The part (i) of Theorem \ref{10.30.05.2q} means that, given a positive 
framed local system on $S'$,  cutting along $\gamma$ provides 
a positive 
framed local system on $S'$.

The part (ii) of Theorem \ref{10.30.05.2q} means that, given a 
positive framed local system on $S$, we can glue  it to a positive 
framed local system on $S$ if and only if the gluing conditions 
of Definition \ref{11.06.05.5} are 
 satisfied. Further, the family of positive local systems on  $S$ 
with given restrictions to $S'$ 
is a principal homogeneous space for the group $H(\R_{>0})$.

\vskip 3mm 
The proof of Theorem \ref{10.30.05.2q} will be given in Sections 7.7 - 7.8.

\vskip 3mm
{\bf 7. The cutting and gluing on the level of boundary at infinity of $\pi_1$.} 
Let $S_1$ and $S_2$ be oriented surfaces 
with boundary, and  with negative Euler characteristic. 
Gluing them along boundary circles $\gamma_1$ and $\gamma_2$, 
where $\gamma_i \in \partial S_i$, we get an oriented surface $S$ with $\chi(S) = \chi(S_1)+\chi(S_2)<0$. 
Let $\gamma$ be the loop on $S$ obtained by gluing $\gamma_1$ and $\gamma_2$. 
Cutting $S$ along $\gamma$ we recover  
$S_1$ and $S_2$.

\begin{figure}[ht]
\centerline{\epsfbox{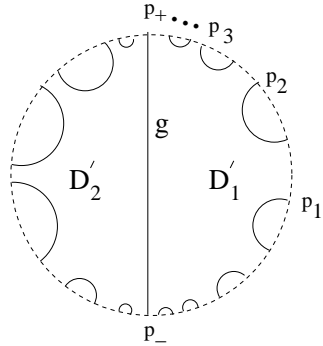}}
\caption{Gluing $D$ out of $D_1$ and $D_2$.}
\label{fg81}
\end{figure}

Let us choose a hyperbolic metric on $S$ such that the boundary of $S$ is geodesic.  
We may  assume that $\gamma$ is a geodesic for this structure. 
Then $S_1$ and $S_2$ inherit hyperbolic metrics with geodesic boundaries. Denote by $D_i$ (resp. $D$)
the universal cover of $S_i$ (resp. $S$). Then $D_i$ is an infinite ``polygon'' whose boundary consists of 
geodesics and segments of the absolute. It is obtained by cutting out of the hyperbolic plane ${\cal H}$ 
discs bounded by geodesics. We say that two boundary geodesics on $D_i$ are 
equivalent if they project to the same boundary component of $S_i$. Both $D_1$ and $D_2$ 
have distinguished equivalence classes of boundary geodesics: the ones which project to $\gamma$.
We call them marked boundary geodesics. 

One can glue $D$ from infinitely many copies of $D_1$ and $D_2$ as follows, see Figure \ref{fg81}. Observe that 
given a geodesic $g$ on ${\cal H}$ we can move by an element of $PGL_2(\R)$ 
the ``polygons'' $D_1$ and $D_2$ so that 
$D_1 \cap D_2 =g$. Now let us take the domain $D_1$, and glue to each 
marked boundary geodesic on it a copy of the domain $D_2$ along a marked geodesic on the latter. 
Then for each 
of the obtained domains we glue a copy of $D_1$ along  every marked boundary 
geodesic on it, and so on infinitely many times. 
\vskip 3mm  
\begin{lemma} \label{11.5.05.1}
The domain glued in this way is identified with the fundamental domain $D$. 
So the group $\pi_1(S)$ acts by  automorphisms of $D$. 
\end{lemma}

{\bf Proof}. Let 
$p_1: D_1 \to S_1$ and $p_2: D_2 \to S_2$ be
 the universal covering maps. Then there is a unique map $p: D \to S$ such that 
its restriction to a copy of the domain $D_i$ in $D$ is provided by the map $p_i$. 
  Since by the very construction the domain $D$ is connected and simply connected, and the map 
$p$ is a covering, $D$ must be the universal cover. The lemma is proved.

\vskip 3mm  
Recall the cyclic $\pi_1(S)$-sets ${\cal G}_{\infty}(S)$ and $\partial_{\infty}\pi_1(S)$. 
Let us address the following 

{\bf Problem}. (i) How to construct the 
cyclic $\pi_1(S)$-set $\partial_{\infty}\pi_1(S)$
out of the ones $\partial_{\infty}\pi_1(S_1)$  and $\partial_{\infty}\pi_1(S_2)$?

(ii) The same question about the 
cyclic $\pi_1(S)$-set ${\cal G}_{\infty}(S)$.

\vskip 3mm  
{\bf Construction}. (i) 
Since $S_i$ has a boundary,  $\partial_{\infty}\pi_1(S_i)$ 
is a Cantor set given by the union of the 
segments of the absolute belonging to the fundamental domain $D_i$. 
The description of $D$ as  the domain obtained by 
gluing infinitely many copies of domains $D_i$, see Lemma \ref{11.5.05.1}, provides us with 
an answer: gluing the absolute parts of the boundaries of the domains $D_i$, we get 
a $\pi_1(S)$-set with a cyclic structure, which is isomorphic to  $\partial_{\infty}\pi_1(S)$.

 (ii) Clearly an answer to (i) implies an answer to (ii). 

\vskip 3mm
Now suppose that  we get $S$ by gluing 
two boundary circles $\gamma_-$ and $\gamma_+$ on a surface $S'$. 
The boundary circles are glued to a loop $\gamma$ on $S$. Choose a hyperbolic metric on 
$S'$ with geodesic boundary. Then $S$ inherits a hyperbolic metric, and $\gamma$ is a geodesic. 
We construct 
the universal cover $D$ of $S$ by gluing copies of the universal cover $D'$ for $S'$ as follows. 
Choose two copies, $D^-$ and $D^{+}$ of $D'$. Each of them 
has a family of marked boundary geodesics: the ones on $D^\pm$ which project to $\gamma_\pm$. 
Then repeat the described above  gluing process: glue copies of $D^{-}$ to a single copy of 
 $D^{+}$ along marked boundary geodesics, and so on. In the end we get 
the fundamental domain $D$. 
This way we also get the boundary at infinity  $\partial_{\infty}\pi_1(S)$, understood as a 
cyclic $\pi_1(S)$-set, glued from infinitely many copies of  $\partial_{\infty}\pi_1(S')$.

\vskip 3mm

{\bf Proof of the part (i) of Theorem \ref{10.30.05.2q}}. It follows immediately 
from the discussion above and Theorem \ref{11.06.05.18}. \vskip 3mm

To prove the part (ii) of Theorem \ref{10.30.05.2q} we need to glue positive 
configurations of flags. 

\vskip 3mm

\begin{figure}[ht]
\centerline{\epsfbox{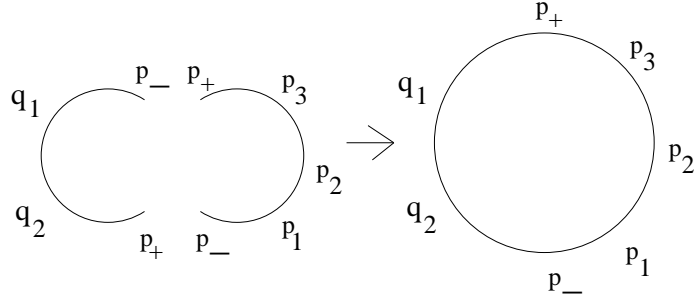}}
\caption{Gluing two cyclic sets along their common two-element subset.}
\label{fg82}
\end{figure}

{\bf 8. Gluing positive configurations of flags}. 
Let 
$$
P = (p_+, p_-, p_1, ..., p_n),\quad Q = (p_-, p_+, q_1, ..., q_m); \qquad P\cap Q = \{p_-, p_+\}
$$ 
be two cyclic sets, 
which may be infinite, i.e. $n =\infty$ or $m =\infty$. We can glue 
$P$ and $Q$ along the subset $\{p_+, p_-\}$, see Figure \ref{fg82}, getting a cyclic set 
$$
P*Q = (p_+, q_1, ..., q_m, p_-, p_1, ..., p_n).
$$
Let 
$$
\beta_P: P \lra {\cal B}(\R), \quad  \beta_Q: Q \lra {\cal B}(\R) 
$$
be two positive configurations of flags. Recall that ``configurations'' means that 
we consider each of them modulo $G(\R)$-action. 
We would like to glue them to get a positive 
configuration of flags 
\begin{equation} \label{11.04.05.2}
\beta_{P*Q}: P*Q \lra {\cal B}(\R) 
\end{equation}
such that its restriction to the subset $P$ (resp. $Q$) is the configuration $\beta_P$ (resp. $\beta_Q$). 
We can always assume that $\beta_{P*Q}(p_+, p_-) = (B^+, B^-)$. The stabilizer 
of this pair of flags is the Cartan subgroup $B^+\cap B^- =H$. 
\begin{theorem} \label{11.04.05.1}
 There is a family of positive configurations of flags (\ref{11.04.05.2}) 
obtained by gluing two given positive configurations $\beta_P$ and $\beta_Q$. 
This family is a principal homogeneous $H(\R_{>0})$-space. 
\end{theorem} 

{\bf Proof}. According to Theorem \ref{11.04.05.4}, we can present the configurations $\beta_P$ and $\beta_Q$ 
as follows:
$$
\beta_P= (B^+, B^-, u_1\cdot B^-, u_1u_2\cdot B^-, ... , u_1 ... u_n\cdot B^-), \quad u_i \in U^+(\R_{>0});
$$
$$
\beta_Q= (B^-, B^+, (v_1 ... v_m)^{-1}\cdot B^-, (v_2 ... v_m)^{-1}\cdot B^-, ... , (v_m)^{-1}\cdot B^-), 
\quad v_i \in U^+(\R_{>0})
$$
Further, we can modify the second configuration by acting by an element 
$h \in H(\R_{>0})$ on it: this changes neither 
the flags $(B^+, B^-)$, nor the positivity of the elements $v_i$. Observe however that the action of any element of $H(\R)-H(\R_{>0})$ would destroy  positivity of the $v_i$'s. Then we can glue 
the two configurations to a single one:
$$
\beta_{P*Q}:= (B^+,      (v_1 ... v_m)^{-1}\cdot B^-, ..., v_m^{-1}\cdot B^-, 
         B^-, u_1\cdot B^-, u_1u_2\cdot B^-, ... , u_1 ... u_n\cdot B^-).
$$
It is positive thanks to Theorem  \ref{11.04.05.4}. The theorem is proved. 
\vskip 3mm

{\bf Proof of the parts (ii) and (iii) of Theorem \ref{10.30.05.2q}}. 
Let us assume first that $S'=S_1 \cup S_2$. 
By Theorems \ref{4.1.04.10} and \ref{11.06.05.18} when $S$ has  boundary, and by Definition 1.10
 when $S$ has no boundary, a positive framed local $G(\R)$-system 
on $S_i$ is the same thing as a $\pi_1(S_i)$-equivariant positive map 
$$
\Psi_i: {\cal G}_{\infty}(S_i) \lra {\cal B}(\R) \qquad \mbox{modulo $G(\R)$-conjugation}.
$$
According to Lemma \ref{11.5.05.1}, ${\cal G}_{\infty}(S)$ 
is glued from infinitely many copies of ${\cal G}_{\infty}(S_i)$ for $i=1,2$. 
Each time when we glue the next domain we glue two cyclic sets along the $2$-element subset given by the endpoints of the unique 
marked boundary geodesic shared by the glued domains.  
 Therefore thanks to Theorem \ref{11.04.05.1} we can 
glue the positive configurations of flags provided by the maps $\Psi_i$ to a positive map $\Psi$. 
The ambiguity of gluing two of these domains, say $D_1'$ and $D_2'$, is given by the action of 
$H(\R_{>0})$. 
As soon as it is fixed, the rest of the gluings are uniquely determined by a  
requirement which we will formulate below. 

\noindent \qquad {\it Gluing requirement}. 
Let $(D_1'', D_2'')$ be another pair of domains glued along a marked geodesic. 
Then there exists an element $g \in \pi_1(S)$ transforming $(D_1', D_2')$ to $(D_1'', D_2'')$. 
It is defined uniquely up to elements of the subgroup of $\pi_1(S)$   
generated by the monodromy around the marked geodesic. 
Denote by $D_1'\ast D_2'$ the domain obtained by gluing of $D_1'$ and $D_2'$. We require that 
{\it the configuration of flags assigned  to (the ${\cal G}_{\infty}$-set on the boundary of) 
$D_1''\ast D_2''$ equals to the one assigned  assigned to $D_1'\ast D_2'$}. 

It follows from Lemma \ref{CCC4} that this requirement, 
 plus the fact that the monodromies along the two boundary components $\gamma_+$ and 
$\gamma_-$ 
coincide, guarantee that there exists a unique representation 
$\rho:\pi_1(S) \to G(\R)$, such that the map $\Psi$, defined by gluing (infinitely many times) of 
the maps $\Psi_1$ and $\Psi_2$, as above,  is $\rho$-equivariant. 
Let us spell out the argument in detail. Let $g \in \pi_1(S,x)$, $x\in \gamma$. 
Let us assume first that $g$ is the homotopy 
class of the loop $\gamma$. Let $D_1'$ and $D_2'$ be two initial fundamental domains glued 
along a geodesic projecting to $\gamma$. Then $g$ provides an automorphism of $D_1'\ast 
D_2'$, and thanks to the condition that the monodromies along $\gamma_+$ and 
$\gamma_-$ 
coincide, there exists  an element $\rho(g) \in G(\R)$ such that the map $D_1'\ast D_2' \to {\cal B}(\R)$ 
intertwines the action of $g$ on $D_1'\ast D_2'$ and the action of $\rho(g)$ on ${\cal B}(\R)$. 
Such an element is unique since the ${\rm Center}(G)$ is trivial.  

Further, for any $g \in \pi_1(S,x)$, the configurations of flags 
assigned to $D_1'\ast D_2'$ and $g(D_1'\ast D_2')$ are isomorphic. 
Moreover, each of these 
configurations   is realised as a part of a bigger configuration of flags, the one assigned to 
${\cal G}_{\infty}(S)$ after the infinite gluing procedure is done. 
Thus there exists a unique element $\rho(g) \in G(\R)$ transforming the first collection
 of flags to the second -- 
the uniqueness follows from the assumption that the ${\rm Center}(G)$ is trivial.  
Clearly $\rho$ is a homomorphism. 

\vskip 3mm
{\bf Remark}. Gluing the positive maps $\Psi_1$ and $\Psi_2$, we must assume that 
the monodromy around loops $\gamma_{\pm}$ preserves the flags $B^+$ and $B^-$, 
and positive configurations, 
and thus is in $H(\R_{>0})$. On the other hand, the monodromy of a  
positive local system on $S_1$ around $\gamma_{+}$  is conjugate to $B(\R_{>0})$. 
One can show that the semisimple elements in $B(\R_{>0})$ are conjugate to $H^0(\R_{>0})$. 
Therefore gluing positive local systems on $S_1$ and $S_2$ we have to assume that their 
monodromies around $\gamma_{\pm}$ are conjugate to an element of $H^0(\R_{>0})$.

To conclude the proof of surjectivity of the cutting map 
it remains to show that the monodromy of a positive local system on $S_1$ 
can be any element of $H^0(\R_{>0})/W$. 
Recall the  canonical map 
\begin{equation} \label{34543ac}
\pi_p: {\cal X}_{G, S} \lra H,
\end{equation} 
provided by the framing and the 
 semi-simple part of the monodromy around a boundary component $p$ of $S$: 
the semi-simple part of the monodromy itself gives a map to $H/W$.   
\begin{lemma} \label{34543a} (i) There exists a coordinate system from the positive atlas on ${\cal X}_{G, S}$ for which 
the map (\ref{34543ac}) is a monomial map, i.e. is given by a homomorphism 
of tori ${\Bbb G}_m^{\chi(S){\rm dim}G} \lra H$. 

(ii) The fibers of the induced map ${\cal X}_{G, S}(\R_{>0}) \to H(\R_{>0})$ are isomorphic to $\R^{\chi(S){\rm dim}G - {\rm dim}H}$. 
\end{lemma} 

{\bf Proof}. (i) We may assume that $S$ is a surface with punctures. Take an ideal  triangulation of 
$S$.  It follows from the construction of the framed local system given in Section 5 that the monodromy around a puncture $p$ is given by the 
product of the $H$-invariants attached to the edges of the triangulation sharing the vertex $p$. 

(ii) The map $z \to z^k$ is an isomorphism on $\R_{>0}$. Thus according to (i) our map is identified with 
a real vector space projection $\R^{\chi(S){\rm dim}G}\to \R^{{\rm dim}H}$.  The Lemma is proved. 

\vskip 3mm

The part i) of the Lemma implies that the monodromy of a positive local system on $S_1$ 
around a given boundary component 
can be any element of $H^0(\R_{>0})/W$. 

The part (iii) of the theorem follows from Lemma \ref{34543a}. Indeed, 
since the monodromy around $\gamma_\pm$ lies in $H^0(\R_{>0})/W$, it is 
semi-simple. Thus, for the surface $S'$, and the boundary component $p$ 
corresponding to $\gamma_+$,  
the fiber of the map $\pi_p$ restricted to ${\cal X}_{G,S}^+$ 
 over an element $m \in H^0(\R_{>0})$ gives 
the space of positive framed local systems on $S'$ 
with given semi-simple monodromy $m$ along $\gamma_+$. 
The case when $S'$ is connected is completely similar. 
Theorem \ref{10.30.05.2q} is proved. 

\vskip 3mm
{\bf Proof of Theorem \ref{4.1.04.10b}}. 
(i). Discreteness is proved by just the same argument as in Theorem \ref{4.1.04.10}. 
Faithfulness follows from Theorem \ref{10.30.05.2q} and the 
similar fact for surfaces with boundary, proved in Theorem \ref{2.1.04.0f}. 
To prove that the monodromy around a non-trivial loop $\gamma$ is positive hyperbolic, 
we cut $S$ along another curve, so that $\gamma$ becomes a non-trivial 
non-boundary loop on the cutted surface, and apply Theorem \ref{10.30.05.2q} and Theorem \ref{2.1.04.0f}. 

(ii). It is sufficient to prove that ${\cal L}_{G, S}^+$ is topologically trivial. 
Since $H(\R_{>0})$ is contractible, this reduces to the similar claim about the moduli space 
${\cal L}_{G, S'}^+(\gamma_+, \gamma_-)$.

(iii). When $S$ is closed the ${\cal X}^+$- and ${\cal L}^+$-Teichm\"uller 
spaces are the same. Hence  the parts (ii) and (iii) of Theorem \ref{10.30.05.2q} show that 
${\cal L}^+_{G,S}$ is a manifold of dimension $-\chi(S'){\rm dim}G$, 
and that it is topologically trivial. 
The theorem is proved. 
\vskip 3mm
\begin{lemma} The canonical map $\pi_p$ 
intertwines the natural action of the Weyl group $W$ on ${\cal X}_{G, S}$ with the one on $H$.
\end{lemma} 

{\bf Proof}. Clear from the definition. 
\vskip 3mm

The canonical map (\ref{34543ac}) provides a map $\pi_p^+: {\cal X}^+_{G, S} \lra H(\R_{>0})$. 

\begin{proposition}
The rational action of the Weyl group $W$ on 
${\cal X}_{G, S}$ corresponding 
to a boundary component $p$ of $S$ preserves $(\pi_p^+)^{-1}H^0(\R_{>0})$.
\end{proposition} 

{\bf Proof}. Let $D$ be the universal cover of a hyperbolic surface $S$. 
Let $g$ be a geodesic lifting the oriented 
boundary geodesic $\gamma_+$ corresponding to the chosen component  of $S$, and 
let $p_-, p_+$ be its ends. Changing the orientation of $\gamma_+$ we get 
an oriented geodesic $\gamma_-$. 
The framing provides a positive map 
$\Psi': {\cal G}_\infty(S) \to {\cal B}(\R)$. Let us alter this map by changing 
 $\Psi'(p_\pm)$ as follows. We set  $\Psi(p_\pm)$ to be the attracting flag 
of the monodromy around the oriented loop $\gamma_\pm$. The attracting flag is well defined since 
the monodromy is in $H^0(\R_{>0})$. The same argument as in the proof of 
Theorem \ref{11.06.05.18} shows that 
the obtained map $\Psi$ is positive. 

We may assume that $\Psi(p_-, p_+) = (B^-, B^+)$, and
 the cyclic order on ${\cal G}_\infty\pi_1(S)$ is given by 
$(p_+, \ldots,  p_-)$. The action of the group $W$ on the frame boils down to replacing $B^-$ by $wB^-w^{-1}$. 
It follows from  Proposition 8.13 in \cite{L1} that  the triple of 
flags $(B_-, wB^-w^{-1}, B^+)$ is non-negative, i.e. in the closure of the triples of positive flags. 
The proposition is reduced to the following claim: Let $(B'', B', B^-, B_t, B^+)$, where $t>0$, 
be a curve in the configuration space of positive $5$-tuples of flags, and $(B'', B', B^-, B_0, B^+)$ 
the limiting configuration. Then the configuration $(B'', B', B_0)$ is positive. 
Conjugating the configuration 
$(B'', B', B^-, B_t, B^+)$ we can write it as  $(B^+, B^-, v_1 \cdot B^-, 
u(t)\cdot B^-, v_2\cdot B^-)$, where 
$u(t),v_1, v_2 \in U^+(\R_{>0})$. Then, in the canonical basis of a representation  of $G$, $u(t)$ has non-negative entries and is bounded from above (resp. below) by the image of $v_2$ (resp.  $v_1$). 
Thus the limit $\lim_{t\to \infty}u(t)$ exists and is in $U^+(\R_{>0})$.  
The proposition is proved. 

{\bf Remark}. We needed to work with $(\pi_p^+)^{-1}H^0(\R_{>0})$ in the proof for 
the existence of the attracting flag for the monodromy around $\gamma_\pm$. Observe that 
 the monodromy around $\gamma_\pm$ is conjugate to an element of $B(\R_{>0})$. 
The  unique attracting flag probably exists for every element of $B(\R_{>0})$, 
at least this is obvious for $G=PGL_m(\R)$.  

\vskip 3mm

\vskip 3mm
 {\bf 9. The space ${\cal L}^+_{G,S}$ for a closed $S$ coincides with the Hitchin component: a proof of 
Theorem \ref{H1}}. 
By  Theorem \ref{10.30.05.2q},  the Teichm\"uller space 
${\cal L}^+_{G,S}$ is an open connected domain 
in ${\cal L}_{G,S}(\R)$.  Furthermore, ${\cal L}^+_{G,S}$ 
lies inside of ${\cal L}^{\rm red}_{G,S}(\R)$. Indeed, by Theorem \ref{4.1.04.10b}i), 
for any representation $\rho$ from 
${\cal L}^+_{G,S}$, and for any $\gamma \in \pi_1(S)$, the element 
$\rho(\gamma)$ is semi-simple. Thus the Zariski closure of $\rho$ is reductive: 
otherwise it is a semi-direct product of a reductive group and a unipotent group $N$, 
and there exists a non-trivial element $\gamma$ such that $\rho(\gamma)\in N$, which contradicts 
to the semi-simplicity of $\rho(\gamma)$.

 So it remains to show that ${\cal L}^+_{G,S}$ is closed in ${\cal L}^{\rm red}_{G,S}(\R)$. 
Indeed, then it is a component of ${\cal L}^{\rm red}_{G,S}(\R)$, and 
since by Corollary \ref{CCC30} ${\cal L}^+_{G,S}$ contains the classical Teichm\"uller space, it is the Hitchin component. 

\vskip 3mm
Let $\{x_i\}$ be a family of points of ${\cal L}^+_{G,S}$ converging when $i \to \infty$ 
to a point of 
${\cal L}^{\rm red}_{G,S}(\R)$. Then 
 there is a family of  
positive representations $\rho_i: \pi_1(S) \to G(\R)$, 
 corresponding to $x_i$'s,  which has the 
limit when $i \to \infty$, denoted by $\rho$.  
Since $G$ has trivial center, each $\rho_i$ 
determines uniquely a positive  map $\psi_i: {\cal G}_{\infty}(S) \to {\cal B}(\R)$. 
Moreover, since there exists the limit $\lim_{i \to \infty}\rho_i$, there exists the limit map 
$\psi:= \lim_{i \to \infty}\psi_i: {\cal G}_{\infty}(S) \to {\cal B}(\R)$ 
(here the convergence means the pointwise convergence). 

\vskip 3mm
For any $s\in {\cal G}_{\infty}(S)$ there exists a unique $s'\in {\cal G}_{\infty}(S)$ 
such that $s,s'$ are the endpoints of a geodesic projecting to a geodesic loop  
on $S$. 
We say that  $s'$  is {\it opposite to $s$  in ${\cal G}_{\infty}(S)$}. 

\begin{lemma} \label{h2} Assume that the map $\psi$ is not positive. 
Then for any pair of points $s_1, s_2 \in {\cal G}_{\infty}(S)$ which are 
not  opposite to each other  the  
flags $\psi(s_1)$ and $\psi(s_2)$ are not in generic position. 
\end{lemma}

{\bf Proof}. Assume the opposite. Then there exists a pair $s_1 \not = s_2 \in {\cal G}_{\infty}(S)$ such 
that the flags 
$\psi(s_1)$ and $\psi(s_2)$ are in generic position, and $s_1, s_2$ are not opposite to each other. 
Conjugating representations $\rho_i$ we may assume 
that $\psi_i(s_1) = B^+$ for all $i$. Still $\{\psi_i\}$ is a family of positive maps, 
and has a limit when $i \to \infty$.

Choose a hyperbolic structure on $S$. 
Take the geodesic $\gamma$ 
on the hyperbolic plane covering $S$,  connecting $s_1$ to its opposite  $s$. 
The stabilizor of 
$\gamma$ in $\pi_1(S)$ is isomorphic to $\Z$. Let $m_{\gamma}$ be its generator. 
Since $s \not = s_2$,  
the points $t_n:= m_\gamma^ns_2$, $n \in \Z$, are different, and we may assume that 
they converge to $s_1$ when $n \to \infty$, and to $s$ when $n \to -\infty$.
We are going to show that this implies that $\psi$ is a positive map. 

Conjugating representations $\rho_i$ 
by elements from $B^+$, we may assume without loss of generality
 that $\psi(s_1)= B^+$ and $\psi(s_2)= B^-$. Then, since  $\psi(t_n)$ is in generic position 
to $B^+$, there exist  
elements $u_i \in U^+(\R)$ such that $\psi(t_n) = u_n \cdot B^-$. Let $U^+_{\geq 0}$ be the closure 
of $U^+(\R_{>0})$ in $U^+(\R)$. The image of any $u \in U^+_{\geq 0}$ 
in any finite dimensional representation of $G$, written in the canonical basis there, 
is a matrix with non-negative coefficients. Since $\psi$ is the limit of 
a family of positive maps, $u_n \in U^+_{\geq 0}$.

Let $C$ be a cyclic set and $a,b,c \in C$. We say that $d \in C$ 
belongs to the arc between $b$ and $c$ and outside of $a$ if the order of $(a,b,c,d)$ 
agrees with the cyclic order induced from $C$.

\begin{figure}[ht]
\centerline{\epsfbox{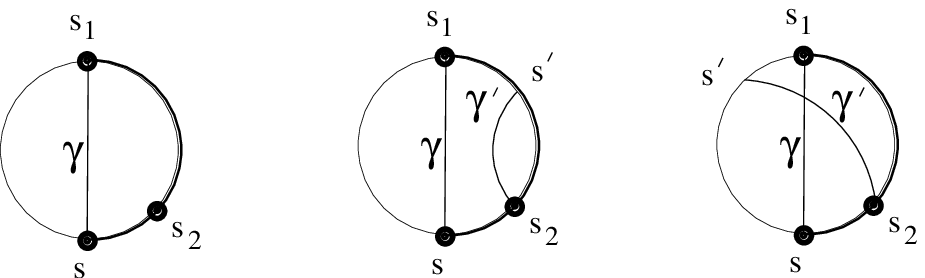}}
\caption{}
\label{fgh1}
\end{figure}

\begin{lemma} \label{h12}
Let $r \in {\cal G}_{\infty}(S)$ be in the arc between $t_{i-1}$ and $t_{i}$ and outside of $s$.  
Then for any $r$ from this arc $\psi(r) = u_r \cdot B^-$, where $u_r \in U^+_{\geq 0}$. 
\end{lemma} 

{\bf Proof}. Choose a section of the projection $U^+_{>0} \to U^+_{>0}/H_{>0}$. 
Writing a flag as $u\cdot B^-$ where $u \in U_{>0}$, we will assume that 
$u$ lies in this section. Therefore if $u^{(n)} \in U^+_{>0}$, the limit when $n \to \infty$ of the flags 
  $u^{(n)} \cdot B^-$ exists, and the limiting flag is in the generic position to $B^+$, then 
the limit $\lim_{n\to \infty}u^{(n)}$ exists, and lies in  $U^+_{\geq 0}$.

Write
 $\psi_n (r):= u_r^{(n)}\cdot B^-$ and $\psi_n (t_i):= u_{t_i}^{(n)}\cdot B^-$,
 where  $u_r^{(n)}, u_{t_i}^{(n)}\in U^+_{> 0}$. 
Then $u^{(n)}_{t_{i}} = u^{(n)}_r \widetilde u^{(n)}_r$ where 
  $\widetilde u^{(n)}_r\in U^+_{> 0}$. 
We know that, as $n \to \infty$, the limit of $u_{t_i}^{(n)}$ exists. Since 
$u^{(n)}_{r}, \widetilde u^{(n)}_r\in U^+_{>0}$, this implies that each of the the limits 
$\lim_{n\to \infty} u^{(n)}_{r}$ and $\lim_{n\to \infty} \widetilde u^{(n)}_r$ 
exists. Indeed, we know that the limit of flags $\psi_n (r)$ exists. Thus 
the first limit will not exist only if one of the matrix coefficients of $u^{(n)}_{r}$  goes to $+\infty$. 
However if this is the case, one of the matrix coefficients of the product 
$u^{(n)}_{r}\widetilde u^{(n)}_r$ also goes to 
$+\infty$, which contradicts to the existence of  $\lim_{n\to \infty} u^{(n)}_{t_i}$. 
If the limits exist, they apparently lie in $U^+_{\geq 0}$. The Lemma is proved. 
\vskip 3mm 
Lemma \ref{h12} implies that  $\psi(r)$ is in generic position to $B^+ = \psi(s_1)$. 
So for any $r$ from 
the arc which is strictly between $s_1$ and $s$ and contains $s_2$, 
$\psi(r)$ is in generic position to $\psi(s_1)$, see Fig. \ref{fgh1}. Take the geodesic $\gamma'$ 
connecting $s_2$ with its opposite $s'$ in ${\cal G}_{\infty}(S)$. 
Since $s' \not = s_1$, $s'$ is either on the right, or on the left of $s_1$, see 
the two pictures on the right of Fig \ref{fgh1}. Assume first that 
$s'$ is on the same arc between $s_1$ and $s_2$ as $s$ (the middle picture on Fig \ref{fgh1}). 
Repeating the above argument we conclude that for any $r\in {\cal G}_{\infty}(S)$ in the arc $A$ 
between 
$s_1$ and $s'$ located outside of $s$ the flag $\psi(r)$ is in generic position to $\psi(s_2)=B^-$. 
Thus for any such $r$, $\psi(r)$ is in generic position to $B^+$ and $B_-$. 
This plus $u_r \in U^+_{\geq 0}$ implies that $u_r \in U^+_{>0}$. 
If $s'$ lies to the left of $s_1$, we apply a similar argument for the arc $A'$ 
between $s_1$ and $s_2$ and outside of $s'$ (the rightmost picture on Fig \ref{fgh1}). 
So the restriction 
of the map $\psi$ to the arc $A$ (respectively $A'$) 
is positive. Since any proper arc of ${\cal G}_{\infty}(S)$ 
can be moved by an element of $\pi_1(S)$ inside of the arcs $A$ or $A'$, 
the map $\psi$ is a positive map. 
Contradiction. The lemma  is proved. 
\vskip 3mm




Let $\Gamma:= \rho(\pi_1(S))$. It is a subgroup of $G(\R)$. 
Let $D_{B^+}$ be the divisor in the flag variety consisting 
of all flags which are not in generic position with 
the flag $B^+$. Observe that if the points $s_1, s \in {\cal G}_{\infty}(S)$ are opposite to each other, 
then $s$ is not on the $\Gamma$-orbit of $s_1$. Thus 
 according to Lemma \ref{h2} one has 
$\Gamma (B^+) \subset D_{B^+}$. (In fact it is easy to avoid using the previous claim). 

A parabolic subgroup in $G$ is a subgroup containing the Borel subgroup $B^+$. 
Let $P_{\alpha}$ be the proper  
maximal parabolic subgroup  corresponding to a simple positive root $\alpha$. 
\begin{lemma} \label{i4}
The subgroup $\Gamma$ is contained in a proper maximal parabolic subgroup of $G$. 
\end{lemma}

{\bf Proof}. Let $\Gamma_{\rm Zar}$ be the Zariski closure 
of  $\Gamma$ in $G$.  Then $\Gamma_{\rm Zar}(B^+) \subset D_{B^+}$. The irreducible components of the 
divisor $D_{B^+}$ are parametrised by the proper maximal parabolic subgroups in $G$. Since
 $\Gamma_{\rm Zar}(B^+) \subset D_{B^+}$, any non-zero vector in the Lie algebra ${\rm Lie}\Gamma_{\rm Zar}$ of 
the algebraic group $\Gamma_{\rm Zar}$ 
lies in a certain proper maximal parabolic Lie algebra of $G$. 

Assume that no proper maximal parabolic subgroup of $G$ contains $\Gamma_{\rm Zar}$. 
Then for every positive simple root $\alpha$ there exists an element $Z_{\alpha} \in {\rm Lie}\Gamma_{\rm Zar}$ 
such that $Z_{\alpha} \not \in {\rm Lie}P_{\alpha}$. Thus for generic numbers $x_\alpha$ the linear combination 
$\sum_\alpha x_\alpha Z_{\alpha}$ does not lie in any maximal parabolic Lie algebra. 
On the other hand it belongs to the Lie algebra of  $\Gamma_{\rm Zar}$, and thus lies in one of the 
proper maximal parabolics. 
Contradiction. The Lemma is proved.

\vskip 3mm
By the definition of the Hitchin component, the algebraic group $\Gamma_{\rm Zar}$ is reductive. 
By Lemma \ref{i4} it is contained in a proper maximal parabolic subgroup 
$P_{\alpha}$, and hence in its Levi radical $M_{\alpha}$. But then its action on the Lie algebra of $G$ 
is not irreducible. On the other hand, $\Gamma_{\rm Zar}$ is irreducible 
by  Lemma 10.1 from 
\cite{Lab} 
(this lemma was formulated in loc. cit. for $SL_n$, but the proof, based on a Higgs bundle argument, 
 works for any $G$). 
This contradiction proves Theorem \ref{H1}.

\vskip 3mm
{\bf 10. Positive curves in ${\cal B}(\R)$ are $G(\R)$-opers on $S^1$}. 
{\it A canonical distribution on ${\cal B}$}. 
The flag variety ${\cal B}$ is equipped with a $G$-invariant non-in\-teg\-ra\-ble 
distribution ${\cal D}$ whose dimension is equal to the rank of the group $G$. 
It is defined as follows. Let $x$ be a point of ${\cal B}$. We may assume that 
the corresponding Borel subgroup is $B^-$. Then the tangent space $T_x{\cal B}$ 
is identified with the Lie algebra ${\rm Lie}(U^+)$. Consider the subspace of 
${\rm Lie}(U^+)$ given by the direct sum of the one dimensional subspaces 
corresponding to the simple roots. Transferring it to the tangent space 
$T_x{\cal B}$, we get 
the fiber ${\cal D}_x$ of the distribution at the point $x$. It is easy 
to see that this definition is independent of the choices involved. 

\begin{proposition} \label{3.23.04.1} 
Let $C$ be a positive differentiable curve  of the flag variety 
${\cal B}(\R)$. Then $C$ is an integral curve of the distribution 
${\cal D}$ on ${\cal B}(\R)$. 
\end{proposition}

{\bf Proof}. Let $c(s) \in {\cal B}(\R)$ be a positive curve. We may assume that 
$c(0) = B^-, c(1) = B^+$. So $c(s) = x(s)B^-$ where $x(s) \in U^+(\R_{>0})$.
 Choose a reduced decomposition $w_0= w_{i_1} ...w_{i_n}$. 
Then the positivity of $x(s)$ implies that we can write it as 
 $x(s) = \prod_{i_k}x_{i_k}(t_k(s))$, where $t_k(0) =0$, $t_k(s) >0$.  Since
$[x_i(s),x_j(t)] = O(st)$, linearizing $x(s)$ at $s=0$ we get  
$\sum_k t'_{i_k}(0) E_{i_k} \in {\rm Lie}(U^+)$, where $x_i(t) = {\rm exp}(tE_{i})$. 
The proposition is proved. 
\vskip 3mm
A smooth map from the circle $S^1$ to ${\cal B}(\R)$ tangent 
to the canonical distribution is nothing else than a {\it $G(\R)$-oper} on $S^1$ in the 
terminology of Beilinson and Drinfeld. So Proposition \ref{3.23.04.1} asserts 
that positive differentiable loops in ${\cal B}(\R)$ are $G(\R)$-opers on $S^1$. 
\vskip 3mm
{\it Describing $PGL_{n+1}(\R)$-opers on real curves}. 
A $C^{n-1}$-smooth (meaning $n-1$ times continuously differentiable) 
curve $K$ in $\R{\Bbb P}^{n}$ gives rise to an 
{\it osculating curve} $\widetilde K$ in the real flag variety
${\cal B}_{n+1}(\R)$ for $PGL_{n+1}$. Namely, $\widetilde K$ 
is formed by the osculating 
flags to the curve at variable points. 
It is an easy and well known fact 
 that any differentiable $PGL_{n+1}(\R)$-oper is  
the osculating curve 
for a certain sufficiently differentiable curve $K$ in $\R{\Bbb P}^{n}$. 
The curve $K$ is recovered by applying the canonical projection $\pi$. 
The $G$-opers for other classical 
groups have a similar description.

 We will show in Section 9.12 that
 positive $PGL_{n+1}(\R)$-opers on $S^1$ come from  convex curves 
in $\R{\Bbb P}^n$. So positive $G(\R)$-opers on real curves  
can be considered as a generalization of convex curves $\R{\Bbb P}^n$ 
to the case of an arbitrary split semi-simple group $G$.

\vskip 3mm
{\bf 11. Remarks on quantum universal Teichm{\"u}ller spaces and $W$-algebras}. 
The space of $G(\R)$-opers  on the parametrised circle $S^1$ 
has a natural Poisson structure, the 
Gelfand-Dikii bracket, defined by Drinfeld and Sokolov in  \cite{DS}. 
On the other hand, since $P^1(\Q)$ is a subset of 
   $S^1$, a positive $G(\R)$-oper provides 
a point of the universal Teichm{\"u}ller space  ${\cal
  X}^+_{G}$.  
 The universal Teichm{\"u}ller space ${\cal X}^+_{G}$ 
has a natural Poisson structure.
For instance if $G = PGL_{m}$ it is  defined 
in the canonical coordinates by the function $\varepsilon_{pq}$ from Definition 
\ref{11.25.02.1} applied to the $m$-triangulation of the Farey triangulation.  
One can show that  
  the Gelfand-Dikii bracket is compatible with the Poisson
bracket on 
${\cal X}^+_{G}$.

The Poisson structure ${\cal X}^+_{G}$ 
admits a natural quantization
equivariant under the action of the Thompson group.
 In the case $G=PSL_m$ the quantization is especially explicit, and 
follows from the results of Chapter 9 and \cite{FG2}.
The case $G = SL_3$ is discussed in \cite{FG3}. 
Therefore we suggest that {\it the quantum universal Teichm\"uller space can be
considered as a combinatorial version of the $W$-algebra corresponding to 
$G$, divided by $G$}.  

\vskip 3mm

{\bf 11. Kleinian $G(\C)$-local systems on hyperbolic threefolds: conjectures.}
Below $S$ is an oriented surface, and $G$ is a semi-simple 
split algebraic group over $\Q$ with trivial center. 
\begin{definition}
The moduli space $Q_{G(\C), S}$ of quasifuchsian $G(\C)$-local systems  on $S$ 
 is the interier part of 
the connected component of the moduli space of $G(\C)$-local systems  on $S$ 
with discrete faithful monodromy 
representations $\pi_1(S) \to G(\C)$, which 
contain the moduli space ${\cal L}_{G, S}^+$ of positive $G(\R)$-local systems  on $S$. 
\end{definition}
Further, let  $Q^{\rm un}_{G(\C),S}$ be the subspace of {\it quasifuchsian unipotent} 
$G(\C)$-local systems  on $S$. It is the interier part of the connected subset of 
${\cal U}_{G,S}(\C)$ which  consists of local systems with 
faithful discrete monodromies and 
 contains ${\cal U}_{G,S}(\R_{>0})$. 

Let ${\cal D}_S$ be the double of the surface $S$, defined as follows. 
Let $\overline S$ be the surface $S$ equipped with the opposite orientation. 
If $S$ has boundary, then ${\cal D}_S$ is a connected surface obtained by gluing 
$S$ and $\overline S$ along the correponding boundary components. 
If $S$ is closed, then ${\cal D}_S:= S \cup \overline S$. 

\begin{conjecture} \label{10.28.05.11}
{\it There are canonical isomorphisms}
\begin{equation} \label{10.28.05.1}
Q_{G(\C), S} \stackrel{\sim}{\lra} {\cal L}_{G, {\cal D}_S}^+, \qquad
Q^{\rm un}_{G(\C), S} \stackrel{\sim}{\lra}
 {\cal L}_{G, S}^{\rm un, +} \times {\cal L}_{G, \overline S}^{\rm un, +}.
\end{equation}
\end{conjecture}

In particular, if $S$ has no boundary, there should be  an isomorphism 
\begin{equation} \label{**}
Q_{G(\C), S} \stackrel{\sim}{\lra} {\cal L}_{G, S}^+ \times {\cal L}_{G, \overline S}^+.
\end{equation}
If $G = PGL_2$, this is equivalent to the Bers 
double uniformization theorem. 

The second isomorphism  in (\ref{10.28.05.1}) should provide $Q^{\rm un}_{G(\C),S}$ 
with a hyperk\"ahler structure: 
One of the complex structures is induced  from the fact that 
$Q^{\rm un}_{G(\C),S}$ is an open domain in ${\cal U}_{G,S}(\C)$. On the other hand, 
${\cal U}_{G,S}(\R_{>0})$ is supposed to be a complex manifold. Finally, the Weil-Petersson 
form on ${\cal U}_{G,S}(\C)$ provides a holomorphic symplectic structure. 

\vskip 3mm
Now let $\Gamma$ be a Kleinian group, that is a torsion-free discrete subgroup of $PGL_2(\C)$. 
So the quotient ${\cal M}_{\Gamma}:= 
\Gamma\backslash {\cal H}^3$, where ${\cal H}^3$ is the hyperbolic space, 
 is a complete hyperbolic threefold. 
We will assume furthermore that ${\Gamma}$ is not {\it elementary}, i.e. the limit set of its action 
on the absolute is infinite, and that ${\cal M}_{\Gamma}$ is 
{\it geometrically finite}, i.e. 
the  convex core of ${\cal M}_{\Gamma}$, defined as the geodesic convex hull of the limit set, 
 is of finite volume. 
The open manifold ${\cal M}_{\Gamma}$ can be compactified by adding a 
surface at each end. We denote the resulting manifold by $\overline {{\cal M}_{\Gamma}}$.  
 We will also assume that ${\cal M}_{\Gamma}$ is {\it incompressible}, that is 
the map $\pi_1(\partial \overline {{\cal M}_{\Gamma}}) \to 
\pi_1(\overline {{\cal M}_{\Gamma}})$ is injective. 

Let $\Omega_\Gamma \subset \partial {\cal H}^3 = \C{\Bbb P}^1$ be the 
discontinuity set for $\Gamma$. It is a union of a finite number of connected, simply-connected components. 
The boundary of $\overline {{\cal M}_{\Gamma}}$ is identified with  
 $S_\Gamma:= \Omega_\Gamma / \Gamma$, which is a union of finite number of surfaces. It inherits 
a complex structure from the absolute, and thus 
provides a point of the Teichm\"uller space ${\cal T}_{S_\Gamma}$. 
Let ${\cal K}_{{\cal M}_\Gamma}$ be the interier part of the deformation space 
of the subgroup $\Gamma$ inside of $PSL_2(\C)$ modulo conjugations. Then the above construction gives rise 
to a map 
\begin{equation} \label{10.31.05.1}
{\cal K}_{{\cal M}_\Gamma} \lra {\cal T}_{S_\Gamma}.
\end{equation}
By the well-known theorem, it is an isomorphism (see \cite{Kr}, \cite{MM} and references therein). 
Below we suggest a conjectural generalization of this picture. 

\vskip 3mm
Pick a principal embedding $\rho: PGL_2 \hra G$. 
Its restriction to $\Gamma$ is a faithful discrete representation $\rho_{\Gamma}: 
\Gamma \to G(\C)$. 
Let ${\cal K}_{G(\C), {\cal M}_\Gamma}$ be the interier part of the 
deformation space of $\rho_\Gamma$ in the class of 
faithful discrete representations $\Gamma \to G(\C)$ 
modulo conjugations.  It is the space of 
$G(\C)$-local systems on  ${\cal M}_{\Gamma}$ with discrete 
faithful monodromies, deforming the one assigned to  $\rho_{\Gamma}$. We call them 
{\it Kleinian local systems} on ${\cal M}_{\Gamma}$.

\begin{conjecture} \label{11.06.05.100}
{\it Assume that $\Gamma$ is a non-elementary Kleinian group, and ${\cal M}_\Gamma$ is geometrically 
finite and incompressible. 
Then there is an isomorphism} 
\begin{equation} \label{10.31.05.2}
{\cal K}_{G(\C), {\cal M}_\Gamma} \stackrel{\sim}{\lra} {\cal L}_{G, S_\Gamma}^+. 
 \end{equation}
\end{conjecture}
Conjecture  \ref{10.28.05.11} is a special case of this conjecture. 
If $G = PGL_2$, the conjectural isomorphism (\ref{10.31.05.2}) reduces to the isomorphism (\ref{10.31.05.1}). In particular, when ${\cal M}_{\Gamma}$ is of finite volume, i.e. is closed or has cusps 
at the ends, this is the Mostow rigidity theorem. Our conjecture claims that, for an arbitrary $G$, 
$\rho_\Gamma$ in  this case has no non-trivial deformations. 

The isomorphism (\ref{10.31.05.2}) should be compatible 
with the one (\ref{10.31.05.1}): the principal embedding $\rho$ should provide 
commutative diagram, where the vertical arrows are the embeddings induced by $\rho$: 
$$
\begin{array}{ccc}
{\cal K}_{{\cal M}_\Gamma} & \stackrel{\sim}{\lra} & {\cal T}_{S_\Gamma}\\
\downarrow &&\downarrow \\
{\cal K}_{G(\C), {\cal M}_\Gamma} &\stackrel{\sim}{\lra}& {\cal L}_{G, S_\Gamma}^+ 
\end{array}
$$

Here is a version of Conjecture \ref{11.06.05.100}. 

\begin{conjecture} \label{11.06.05.1000} {\it Assume that ${\cal M}_\Gamma$ 
is as in Conjecture \ref{11.06.05.100}. Then the restriction of a local $G(\C)$-system 
on ${\cal M}_\Gamma$ to the ends of ${\cal M}_\Gamma$  provides an injective  map 
\begin{equation} \label{11.06.05.21}
{\cal K}_{G(\C), {\cal M}_\Gamma} \hra Q_{G(\C), S_\Gamma}.
\end{equation}
Its image is a Lagrangian submanifold in $Q_{G(\C), S_\Gamma}$. 
The map (\ref{10.31.05.2}) is the composition of the map (\ref{11.06.05.21}) 
followed by the projection $Q_{G(\C), S_\Gamma} \lra {\cal L}^+_{G, S_\Gamma}$
 from Conjecture \ref{10.28.05.11}, see (\ref{**}). }
\end {conjecture}

It is easy to show that the image of the map (\ref{11.06.05.21}) is isotropic. 
Thus the fact that it is Lagrangian should follow from Conjecture \ref{11.06.05.100}. 

\vskip 3mm
{\bf A hyperbolic field theory.} Here is an interpretation of 
Conjecture \ref{11.06.05.100}. 
We assign to each closed surface $S$ a symplectic manifold 
$Q_{G(\C), S}$. Further, let us assign to 
each topological threefold ${\cal M}$ isomprphic to ${\cal M}_\Gamma$ as in 
Conjecture \ref{11.06.05.100} the 
manifold ${\cal K}_{G(\C), {\cal M}} $. By 
 Conjecture \ref{11.06.05.1000} it is a Lagrangian submanifold 
in $Q_{G(\C), \partial{\cal M}}$. 
These Lagrangian submanifolds should satisfy the 
classical field theory axioms. In particular, 
if ${\cal M}$ is glued from threefolds ${\cal M}_1$  and 
${\cal M}_2$ along a boundary component, then 
${\cal K}_{G(\C), {\cal M}}$ is identified with the subset of 
${\cal K}_{G(\C), {\cal M}_1}\times {\cal K}_{G(\C), {\cal M}_2}$ consisting of the pairs 
whose retrictions to the boundary components  coincide. 

\section{A  positive structure on the moduli space ${\cal A}_{G, \widehat S}$ }
\label{posA}

In this section  
  $G$ is  simply connected and  connected. 
We define a positive structure on the moduli space 
 ${\cal A}_{G, \widehat S}$. In particular if $\widehat S$ 
is a disc with $n$ marked points at the boundary we get  
a positive structure on the moduli space  ${\rm Conf}_n({\cal A})$ 
of configurations of $n$  affine 
flags in $G$. We start with a definition of 
 a positive structure on the moduli spaces  ${\rm Conf}_n({\cal A})$ 
for $n=2$ and $n=3$. 
The definition in general is given by using an ideal triangulation 
of  $\widehat S$. We show that it does 
not depend on the choice of an ideal triangulation. 

\vskip 3mm
{\bf 1. Generic configurations of pairs of affine flags}. 
A pair of affine flags is in {\it generic position} 
if the underlying pair of flags is in generic position. 
An $n$-tuple of affine flags is in generic position if every 
two of them are in generic position. 
Let ${\rm Conf}^*_n({\cal A})$ be the variety of 
configurations of $n$  affine 
flags in generic position in $G$. 
Similarly one defines a generic pair $(A_1, B_2)$ where $A_1$ is an affine flag 
and $B_2$ is a flag. 
The  generic pairs $(A_1, B_2)$ 
form a principal homogeneous $G$--space.

If we care only about the birational type of the variety ${\rm Conf}^*_n({\cal A})$, we will usually 
work with the coresponding moduli space ${\rm Conf}_n({\cal A})$.

Recall the element $\overline w_0 \in N(H)$. 
Recall that $[g]_0:= h$ if $g = u_-hu_+$, where $u_{\pm} \in U^{\pm}, h \in H$. 
Let $(A_1, A_2) = (g_1 U^-, g_2\overline w_0 U^-) \in {\rm Conf}^*_2({\cal A})$. We set 
$$
\alpha_2(g_1 U^-, g_2\overline w_0 U^-):= [g_1^{-1}g_2]_0 \in H.
$$
We get  a map $
\alpha_2: 
{\rm Conf}^*_2({\cal A}) {\to} H$, $(A_1, A_2)
\lms \alpha_2(A_1, A_2)$. 

Recall that the Cartan group $H$ acts naturally from the right on the affine 
flag variety: $A \lms A \cdot h$.  Since $H$ is commutative we can also
consider it as a
left action. 

\begin{proposition} \label{4.21.03.1}
The map $\alpha_2$ provides an isomorphism
\begin{equation} \label{9.14.03.2}
\alpha_2: 
{\rm Conf}^*_2({\cal A}) \stackrel{\sim}{\lra} H, \quad  (A_1, A_2)
\lms \alpha_2(A_1, A_2).
\end{equation}
It has the following properties: 
\begin{equation} \label{9.14.03.2er}
\alpha_2(A_1\cdot h_1, A_2\cdot h_2) = h^{-1}_1 w_0(h_2) \cdot \alpha_2(A_1, A_2),
\quad \alpha_2(A_2, A_1) = s_G \cdot w_0\left(\alpha_2(A_1, A_2)^{-1}\right).
\end{equation}
\end{proposition}

{\bf Proof}.
Bruhat decomposition shows that 
 an arbitrary generic pair of affine flags $(A_1, A_2)$ can be 
written uniquely as $(h\cdot U^-, \overline w_0 U^-)$. 
This provides the first claim of the proposition. 
 To check the first of the properties (\ref{9.14.03.2er}) 
notice that $\alpha_2(h_1U^-, \overline w_0 h_2 U^-) = h_1^{-1}w_0(h_2)$. 
The second of the properties (\ref{9.14.03.2er}) follows immediately from this one. 
The Proposition is proved. 
\vskip 3mm
Observe that 
$(A_1 \cdot h,  A_2 \cdot w_0(h)) 
= h(A_1,  A_2)$. 

Recall the twisted cyclic shift map from Definition \ref{4.21.03.157}. 

\begin{lemma} \label{4.21.03.156}
The twisted cyclic shift map is a positive automorphism of 
${\rm Conf}_2({\cal A})$. 
\end{lemma}

{\bf Proof}.   Follows immediately from 
(\ref{9.14.03.2er}). 
\vskip 3mm
If $s_G \not = e$ then the cyclic
shift is not a positive automorphism of ${\rm Conf}_2({\cal A})$. 
\vskip 3mm
{\bf 2. A positive structure on ${\rm Conf}_3({\cal A})$}. 
Let us define a map
\begin{equation} \label{11.10.03.65}
\gamma: U_*^+ \lra H, \quad u^+ \lms [u^+ \overline w_0]_0 \in H.
\end{equation}

\begin{lemma} \label{ytyt}
The map $\gamma$ is a positive map.
\end{lemma}
 
{\bf Proof}. 
Using the isomorphism ${\Bbb G}_m^I  = H$ 
provided by $\omega_1, ..., \omega_r$ we have 
$\gamma(u_+) := \{\Delta_{\omega_i, w_0(\omega_i)(u_+)}\}_{i \in I} 
\in {\Bbb G}_m^I$. It follows from Theorem 
\ref{9.19.03.41} that 
it is a positive map.  The lemma is proved. 

Consider the regular open embeddings
$$
\beta_1^-: B^- \hookrightarrow G/U^-, \quad b_- \lms b_- \overline {\overline w}_0 U^-,
$$
$$
\beta_2^-: B^+ \hookrightarrow G/U^-, \quad b_+ \lms b_+ U^-.
$$
The Borel subgroups $B^-$ and $B^+$ have the standard positive structures. 
They provide $G/U^-$ with two positive structures. 
Similarly the regular open embeddings
$$
\beta_1^+: B^+ \hookrightarrow G/U^+, \quad b_+ \lms b_+ {\overline w}_0 U^+,
$$
$$
\beta_2^+: B^- \hookrightarrow G/U^+, \quad b_- \lms b_- U^+.
$$
provide $G/U^+$ with two positive structures. 
Let $\Phi:B^- \to B^+$ be a rational map such that $\beta_2^+ = \beta_1^+\Phi$
on their common domain of definition. 
\begin{lemma} \label{7.8}

a) The two positive structures on $G/U^-$ 
provided by the maps $\beta_1^-$, $\beta_2^-$ are compatible. 

b) The two positive structures on $G/U^+$ 
provided by the maps $\beta_1^+$ and $\beta_2^+$ are compatible. 

c) The maps $\Phi$ and $\Phi^{-1}$ are positive rational maps. 
\end{lemma}

{\bf Proof}. a) $<=>$ b). Indeed, 
consider the involutive automorphism $\theta:G
\lms G$ which is uniquely determined by  
$$
\theta(h):= h^{-1}, \quad 
\theta(x_i(t)):= y_i(t), \quad \theta(y_i(t)):= x_i(t).
$$
Observe that $\theta(\overline {\overline w}_0 ) = {\overline w}_0 $ and 
$\theta({\overline w}_0 ) = \overline {\overline w}_0$. 
Applying $\theta$ we deduce b) from a) since $\theta(B^{+}) =
B^{-}, \theta(B^{-}) =B^{+}, \theta(\overline {\overline w}_0) =\overline 
w_0$. Similarly b) implies a).

Let us prove b). Recall that $b_+ \overline w_0 U^+ = b_- U^+$. So $b_+ = \Phi(b_-)$. 
We have to show that the maps $b_+ \to b_-$ and 
$b_- \to b_+$ are positive. 
Consider the first map. 
Write $b_+ = h_+u_+, b_- = h_-u_- $. We may assume that $h_+ =1$. Then $h_- = \gamma(u_+)$, 
so it is a positive map by Lemma \ref{ytyt}. The map  $u_+ \lms u_-$ is 
the map $\phi^{-1}$, and hence it is positive. The claim is proved. 

c) Follows from b). The lemma is proved. 
\vskip 3mm

{\bf 3. Positive configurations of affine flags}. 
Let $(A_1, A_2, A_3)$ be a gen\-eric triple of affine flags in $G$. 
Denote by $B_i$ the flag corresponding to $A_i$. 
So every two of the flags $B_1, B_2, B_3$ 
 are in generic position. 
Let ${\rm Conf}^*(A_1, B_2, B_3)$ be the configuration space of generic
triples $(A_1, B_2, B_3)$. Then there is an isomorphism
\begin{equation} \label{4.18.03.4h}
{\rm Conf}^*(A_1, B_2, B_3) \stackrel{\sim}{\lra} U^+_*.
\end{equation}
Indeed, let $(U^-, w_0B^-)$ be a standard 
pair. 
Then there exists unique  $u \in U^+_*$ 
 such that 
\begin{equation} \label{4.18.03.4}
(A_1, B_2, B_3) \sim (U^-, w_0B^-, u_+B^-); \qquad u_+ \in U^+_*.
\end{equation}

\begin{lemma} \label{4.21.03.2}
There is a canonical isomorphism
\begin{equation} \label{9.14.03.1a}
\alpha_3: {\rm Conf}^*_3({\cal A}) \stackrel{\sim}{\lra} {\rm Conf}^*_2({\cal A})
\times {\rm Conf}^*_2({\cal A}) \times {\rm Conf}^*(A_1, B_2, B_3) 
 \end{equation}
$$
(A_1, A_2, A_3) \lms (A_1, A_2) 
\times (A_2, A_3) \times (A_1, B_2, B_3).
$$
\end{lemma}

{\bf Proof}. The 
affine flags $A_2$ and $A_3$ are determined by 
the elements $\alpha_2(A_1, A_2)$ and $\alpha_2(A_1, A_3)$ in $H$. 
The lemma is proved. 
\vskip 3mm
The isomorphism $\alpha_3$ combined with the isomorphisms  
(\ref{9.14.03.2}) and (\ref{4.18.03.4h}) leads to an isomorphism
\begin{equation} \label{9.14.03.1}
\alpha_3: {\rm Conf}^*_3({\cal A}) \stackrel{\sim}{\lra}. 
H \times H \times U^+_*,  
 \end{equation}

\begin{lemma} \label{4.21.03.235}
The inverse to the isomorphism $\alpha_3$ is defined as follows:
$$
\alpha^{'}_3: H \times H \times U^+_*  \lra {\rm Conf}^*_3({\cal A}),
$$
$$
(h_2, h_3, u_+) \lms (U^-, h_2\overline w_0 U^-, h_2 w_0(h_3)s_G u_+ U^-).
$$
\end{lemma}

{\bf Proof.} Consider the  natural projection 
$$
{\bf e}_{ij}: {\rm Conf}^*_3({\cal A}) \lra H; \qquad 
(A_1, A_2, A_3) \lms (A_i, A_j), \quad 1 \leq i<j\leq 3.
$$ It follows from 
the definition of the map $\alpha_2$ that ${\bf e}_{12} \alpha^{'}_3$ 
is the projection on the first $H$-factor. Let us 
check that ${\bf e}_{23} \alpha^{'}_3$ 
is the projection on the second $H$-factor. One has 
$$
(h_2\overline w_0 U^-, h_2 w_0(h_3)s_G u_+ U^-) \sim (U^-, h_3 s_G 
\overline w^{-1}_0 u_+ \overline w_0 \overline w^{-1}_0U^-) = $$ $$ =(U^-, h_3 
\overline w^{-1}_0 (u_+)\overline w_0U^-).
$$
Since $\overline w^{-1}_0 (u_+):= \overline w^{-1}_0 u_+ \overline w_0 
\in U^-$, the statement follows. 
The similar claim about the $U^+_*$-factor is obvious. The lemma is proved. 
\vskip 3mm
\begin{definition} \label{9.14.03.11} 
A positive regular structure on the variety ${\rm Conf}^*_3({\cal A})$ 
is given by the isomorphism (\ref{9.14.03.1}) and the standard positive regular 
structures on $U_*^+$ and $H$. 
\end{definition}

\begin{lemma} \label{9.14.03.15}  
The map 
${\bf e}_{13}: {\rm Conf}_3({\cal A}) \lra H$ 
is a positive rational map. 
\end{lemma}

{\bf Proof}.  Recall the map $\gamma$, see (\ref{11.10.03.65}). 
Using the isomorphism (\ref{9.14.03.1}) we have
$$
{\bf e}_{13}\alpha_3'(h_2, h_3, u_+) = (U^-, h_2 w_0(h_3)s_G u_+ U^-) = 
(U^-, h_2 w_0(h_3) u_+ \overline w_0 \overline w_0U^-). 
$$
Comparing this with the Bruhat decomposition,  we observe that this is equal to 
$h_2 w_0(h_3)\gamma(u_+)$. The lemma is proved.   
\vskip 3mm

\begin{proposition} \label{4.03.6.10}  
The twisted cyclic shift is a positive rational map on the moduli space ${\rm Conf}_3({\cal A})$. 
\end{proposition}

{\bf Proof}. 
To simplify the calculation we may assume without loss of generality that
$h_2 = h_3 = e$. So we start from the configuration $$(U^-,~ \overline w_0U^-,~
u_+s_GU^-).$$
Multiplying from the left by $u_+^{-1}$ 
and then using $\overline w_0^{-1}u_+^{-1}\overline
w_0 \in U^-$ we get 
$$
(U^-,~ \overline w_0 U^-,~ u_+s_GU^-) \sim  
(u_+^{-1}U^-,~  \overline w_0 U^-,~ s_GU^-).
$$
Applying the inversion map $g \lms g^{-1}$ we get 
$$
(U^-u_+,~  U^-~\overline w^{-1}_0,~ U^-s_G).
$$
Applying the antiautomorphism $\Psi$   
we get 
$$
(u_-U^+,~ \Psi ( \overline w^{-1}_0) U^+,~ s_GU^+), \quad u_-:= \Psi(u_+).
$$
According to Lemma \ref{7.8}b) one has $u_-U^+ = v_+h\overline w_0
U^+$. Here $v_+h =  \Phi(u_-)$ and $v_+ = \phi(u_-)$.
The map $\Phi$ is a positive map. Thus, since $v_+^{-1}U^+ = U^+$,
  we can write our configuration as 
$$
(v_+ h \overline w_0 U^+,~ \Psi ( \overline w^{-1}_0) U^+,~ s_GU^+) \sim 
(h \overline w_0 U^+,~ v_+^{-1}\Psi ( \overline w^{-1}_0) U^+,~ s_GU^+). 
$$
Applying the inversion map we write it as 
$$
(U^+\overline w_0^{-1}h^{-1},~  U^+ \Psi ( \overline w_0)v_+ ,~ U^+s_G). 
$$
Then applying $\Psi$ we get 
$$
(h^{-1}\Psi(\overline w_0^{-1}) U^-,~ \Psi (v_+)  \overline w_0 U^-,~ s_G U^-). 
$$
Applying Lemma \ref{7.8}a) we can write $\Psi (v_+)\overline {\overline
  w}_0U^- = h'\phi\Psi(v_+)U^-$, where the map $\Psi(v_+) \lms h'$ is
positive. So we get 
$$
(h^{-1} \overline w_0 U^-,~ s_Gh'(\phi\Psi)^2 (u_+)  U^-,~ s_GU^-). 
$$ 
We conclude that 
\begin{equation}\label{10.16.03.1}
(U^-,~ \overline w_0U^-,~ u_+s_GU^-) 
\sim (h^{-1} \overline w_0 U^-,~ s_Gh'(\phi\Psi)^2 (u_+)  U^-,~ s_GU^-). 
\end{equation}
Since the maps $\phi\Psi$, $u_+\lms h^{-1}$ and $u_+\lms h'$ 
are positive, formula (\ref{10.16.03.1}) implies
that the twisted cyclic shift is a positive map on ${\rm Conf}_3({\cal A})$.
The Proposition is proved.
\vskip 3mm

Let us consider the following regular  open embedding 
$$
H^3 \times U^+_* \times U^-_* \hra {\rm Conf}'_4({\cal A})
$$
\begin{equation}\label{10.16.03.3}
(h_2, h_3, h_4, u_+, v_-) \lms (U^-,~ u_+^{-1}h_2 U^-,~ h_3 \overline w_0U^-,~
v_- h_4 \overline w_0U^-).
\end{equation}
Here ${\rm Conf}'_4({\cal A})$ is the moduli space of configurations of 
four flags $(B_1, ..., B_4)$ such that all flags except the pair 
$(B_2, B_4)$ are in generic position. 

We are going to define 
a positive structure on the moduli space ${\rm Conf}_4({\cal A})$  
by using the birational map 
\begin{equation}\label{10.16.03.13}
H^3 \times U^+_* \times U^-_* \lra {\rm Conf}_4({\cal A}).
 \end{equation}
and the standard positive structures on $H$, 
$U^+_*$ and $U^-_*$. We have to check however that this map is defined on a
complement to a positive rational divisor.

\begin{proposition} \label{4.03.21.5} a) The map (\ref{10.16.03.13}) 
is defined on a
complement to a positive rational divisor. 

b) The twisted cyclic shift provides
a positive automorphism of the moduli space 
${\rm Conf}_4({\cal A})$. 
\end{proposition}

{\bf Proof}. 
Let us define a rational map 
\begin{equation} \label{10.16.03.51q4}
U^+\times U^- \times H\times H  \lra U^-
\times U^+ \times H, \quad (u_+, v_-, h_2, h_4) 
\lms (a_-, a_+, t)
\end{equation} 
 by solving the equation
$
h_2^{-1}u_+v_-h_4 = a_- a_+ t
$. 
It is a positive rational map, since the two Gauss decompositions in $G$ are related by a positive map, 
see Proposition \ref{8.28.03.1}. Let us 
set $$
g:= a_-^{-1}h_2^{-1}u_+ = a_+ t h_4^{-1}v_-^{-1}.$$
 Then multiplying the configuration (\ref{10.16.03.3}) 
by $g$ from the left we get
$$
(a_+th_4^{-1} U^-, ~ U^-, ~a_-^{-1}h_2^{-1}h_3 \overline w_0 U^-, ~ t \overline w_0U^-).
$$ 
Applying to it the twisted cyclic shift we get
\begin{equation} \label{10.16.03.51}
(U^-, ~a_-^{-1}h_2^{-1}h_3 \overline w_0 U^-, ~ t \overline w_0U^-, ~a_+th_4^{-1}s_G U^- ).
\end{equation}
We can write the configuration (\ref{10.16.03.51}) in the form of (\ref{10.16.03.3}):
$$
(\ref{10.16.03.51}) = (U^-,~ \widetilde u_+^{-1}\widetilde h_2 U^-,~ 
\widetilde h_3 \overline w_0U^-,~
\widetilde v_- \widetilde h_4 \overline w_0U^-).
$$
In other words,  one must have 
\begin{equation} \label{10.16.03.514}
i) ~~ a_-^{-1}h_2^{-1}h_3\overline w_0 U^- = \widetilde u_+^{-1}
\widetilde h_2U^-, \quad ii) ~~ t
= \widetilde h_3, \quad iii)~~ a_+th_4^{-1}s_GU^- = 
\widetilde v_- \widetilde h_4 \overline w_0U^-.
\end{equation}

\begin{lemma} \label{10.18.03.1} 
The rational map 
$$(h_2, h_3, h_4, u_+, v_-) \lms (\widetilde h_2, \widetilde h_3, \widetilde
h_4, \widetilde u_+, \widetilde v_-)$$ provided by equations 
(\ref{10.16.03.514}) is a
positive rational map. 
\end{lemma}

{\bf Proof}. 
For $\widetilde h_3$ and the equation ii) this 
follows from the positivity of the map (\ref{10.16.03.51q4}). 
For the equation iii) and $ \widetilde v_-,  \widetilde h_4$ this follows
immediately from Lemma \ref{7.8}a) since
$s_G \overline w_0 = \overline {\overline w}_0$. 
To handle i) means to show that the rational map 
$B^- \to B^+$, $x_- \lms x_+$ 
 determined from the equation 
$x_-^{-1}\overline w_0U^- = 
x_+^{-1}U^-$ is a positive rational map. 
This equation is equivalent to $x_+ = \Psi \Phi \Psi(x_-)$.  
Indeed, applying the inversion antiautomorphism to this equation we get 
$U^- \overline {\overline w}_0 x_- = U^-x_+$. Then applying $\Psi$ we get $\Psi(x_-) {\overline w}_0 U^+ =
 \Psi(x_+)U^+$. So $\Psi(x_+) = \Phi^{-1}\Psi(x_-)$, and hence  $x_+ =
 \Psi\Phi^{-1}\Psi(x_-)$. 
By  Lemma \ref{7.8} the rational map 
$\Psi \Phi^{-1} \Psi: B^+ \to B^-$ is positive. So the lemma follows. 
Both statements of the proposition follow immediately from this lemma. 
The proposition is proved. 
\vskip 3mm
Let us consider the canonical projection (the edge projection) 
\begin{equation} \label{11.6.03.77}
{\bf e} = ({\bf e}_{12}, {\bf e}_{23}, {\bf e}_{13}): 
{\rm Conf}^*_3({\cal A})  \lra H \times H \times H. 
\end{equation}
It follows from Lemma  
\ref{4.21.03.235} and proof of Lemma \ref{9.14.03.15} that one has 
\begin{equation} \label{11.6.03.78}
{\bf e}\alpha'_3(h_2, h_3, u_+) = (h_2, h_3, h_2 w_0(h_3)\gamma(u_+)). 
\end{equation}
Therefore the map {\bf e} is surjective. 
Let us define the variety ${\cal V}_G$ as its fiber over $(e,e,e)$: 
\begin{equation} \label{9.19.03.6}
{\cal V}_G:= {\bf e}^{-1}(e,e,e) \subset {\rm Conf}^*_3({\cal A}).
\end{equation}

\begin{proposition} \label{9.19.03.5}
a) The variety ${\cal V}_G$ has a natural positive structure.

b) The positive structure on the variety ${\cal V}_G$ 
is invariant under the cyclic shift map. 
\end{proposition}

{\bf Proof}. a) We deduce this from Theorem \ref{9.19.03.41} applied to the double Bruhat cell $G^{e, w_0}$. 

Let us identify $G^{e, w_0}$ with the configuration space ${\rm Conf}(A_1, B_2, A_3)$.  
The configuration space of triples $(A_1, B_2, A_3)$ where $(A_1, B_2)$ are in generic position 
is identified with the principal affine space $G/U^-$ by setting 
$(A_1, B_2, A_3)\newline \to (U^-, B^+, g\cdot U^-)$. The inclusion $B^+ \hra G/U^-$ gives rise to a subspace of the configuration space consisting of triples $(A_1, B_2, A_3)$ where in addition $(A_1, A_3)$ are in generic position. 
Then 
the condition that $(B_2, A_3)$ are in generic position singles out the subvariety $B^+ \cap B^-w_0B^- \subset B^+$,  
that is the double Bruhat cell $G^{e, w_0}$.

Given a reduced decomposition 
${\bf i}$ of $w_0$, the collection of the generalized minors 
$F({\bf i}) = \{\Delta(k; {\bf i}), k \in -[1,r] \cup [1, l(w_0)]\}$ 
provides a positive coordinate system 
on $G^{e, w_0}$. The minors $\Delta(k; {\bf i}), k \in -[1,r]$ are 
the ones  $\Delta_{\omega_i, w_0(\omega_i)}$, so they give the map $\gamma$. 
On the other hand among these minors there are $\Delta_{\omega_i, \omega_i}$ for
every $i \in I = \{1, ..., r\}$. To see this write the sequence 
${\bf i} = (i_1, ..., i_{l(w_0)})$. Then for each $p \in \{1, ..., r\}$  take 
the rightmost among all $i_a$'s such that $i_a = p$. Denote it by 
$i_{a(p)}$. Then $\Delta(i_{a(p)}; {\bf i}) = \Delta_{\omega_p, \omega_p}$. 
Indeed, look at the formulas (\ref{9.19.03.1}) and (\ref{9.19.03.2}). For $i_{a(p)}$, 
the elements of the Weyl group defined by these formulas are both equal to $e$. Indeed, 
in the formula (\ref{9.19.03.1}) the product is over the empty set since $\varepsilon(i_l) =+1$ 
for all $l$ involved since this is always so for the element $(e, w)$. 
In the second formula the set of possible $l$'s is empty. 
Thus according to the definition (\ref{bfz03}), we get $\Delta_{\omega_p, \omega_p}$. 
The rest of the generalized minors, that is the minors 
$$
F'({\bf i}) = \{\Delta(k; {\bf i}), k \in  [1, l(w_0)] - \{a(1), ..., a(r)\}\}
$$
give a positive coordinate system on ${\cal V}_G$. The part a) is proved. 

b) Follows easily 
from a) and Proposition \ref{4.03.6.10}. The \!Proposition is proved. 
\vskip 3mm

{\bf 4. A positive structure on the moduli space ${\cal A}_{G, \widehat S}$}.
 Let $\widehat S$ be a  marked hyperbolic surface.
Choose an ideal  triangulation $T$ of $\widehat S$.  
Denote by ${\cal A}_{G, \widehat t}$ the moduli space of the twisted 
decorated $G$-local
systems on the triangle $t$, considered as a disc with three marked points on
the boundary. These marked points are located on the sides of the triangle, one
 point per each side. 
The punctured tangent space to a triangle $t$ of the triangulation $T$ 
 sits in the punctured tangent space $T'S$. So restricting 
the local system ${\cal L}$ on $T'S$ representing an element of 
${\cal A}_{G, \widehat S}$ 
to $T't$ 
we get an element of ${\cal A}_{G, \widehat t}$. 
So we get a projection $ q'_t: 
{\cal A}_{G, \widehat S} \lra {\cal A}_{G, \widehat t}$. 
 Similarly,  each  edge 
${e}$ of the triangulation $T$ can be thickened a bit to became a disc with
 two distinguished points on its boundary provided by the vertices of $e$. 
They cut the boundary of the disc into two arcs. 
Let $\widehat e$ be this disc with two marked points on the boundary, one on
 each of these arcs.  
 

\begin{figure}[ht]
\centerline{\epsfbox{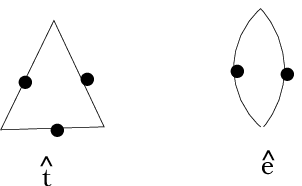}}
\caption{}
\label{fg-14}
\end{figure}

The restriction provides a map 
$q_{e}: {\cal A}_{G, \widehat S} \lra {\cal A}_{G, \widehat e}$. 
All together, they provide a map 
\begin{equation} \label{10.26.03.1}
\varphi_{T}: {\cal A}_{G, \widehat S} \lra 
\prod_{t \in {\rm tr}(T)}{\cal A}_{G, \widehat t} 
\times \prod_{{e} \in {ed}(T)}{\cal A}_{G, \widehat e}. 
\end{equation}
However the target space of this map is too big. Indeed, let $e$ be an edge of
a triangle $t$. Then there is a natural projection $q_{\widehat t,\widehat e}: {\cal A}_{G,
  \widehat t}
\to {\cal A}_{G, \widehat e}$, and $q'_{\widehat t} q_{\widehat
  t,\widehat e} = q_{\widehat e}$. The subvariety 
in the right hand side of (\ref{10.26.03.1}) defined by these
conditions is the
right target space. It contains the following open subset. 
Let us denote by 
 $\widetilde {\rm Conf}^*_n{\cal A}$ the 
variety of twisted cyclic configurations of $n$ affine flags in generic position. 
It is isomorphic to ${\rm Conf}^*_n{\cal A}$, but not canonically. 
\begin{definition}\label{10.26.03.1fa} Let $T$ be an ideal triangulation of 
  $\widehat S$. The subvariety 
$$
{\cal A}_{G, T} \subset \prod_{t \in {\rm tr}(T)}\widetilde {\rm Conf}^*_3{\cal A} 
\times \prod_{{e} \in {ed}(T)}\widetilde {\rm Conf}^*_2{\cal A}
$$
is defined by the conditions $q'_{\widehat t} q_{\widehat
  t,\widehat e} = q_{\widehat e}$ for every  pair $(e, t)$ where $e$ is an
edge of the triangle $t$ of $T$. 
\end{definition}

\begin{theorem} \label{4.03.06.1}  Let $G$ be a
 split simply-connected semi-simple algebraic group. 
 Let $\widehat S$ be a  marked hyperbolic
 surface and ${T}$ an ideal triangulation of $\widehat S$. Then 
there exists a regular open embedding
\begin{equation} \label{9.28.03.11y}
\nu_{{T}}: {\cal A}_{G,  T} \hra {\cal A}_{G, \widehat S}
\end{equation}
such that $\varphi_{{T}} \nu_{{T}} $ is the identity. 
\end{theorem}

\begin{figure}[ht]
\centerline{\epsfbox{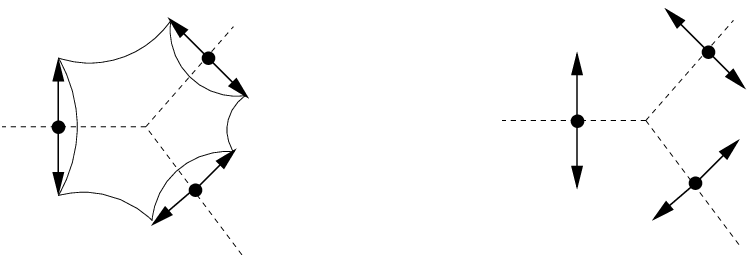}}
\caption{}
\label{fg-10}
\end{figure}

{\bf Proof}. 
Let us introduce a convenient way 
to think of local systems on the punctured tangent bundle $T'S$ 
having the  monodromy $s_G$ around the circle $T'_xS$. 
Let $\Gamma$ 
be a graph on $S$ dual to an ideal triangulation $T$. 
Let us define a collection of points in $T'S$ as follows. 
For every edge $e$ of the graph $\Gamma$ let us choose an internal point 
$x(e)$ 
on this edge, and a pair of 
non zero  tangent vectors $v_1(e)$ and 
$v_2(e)$ at this point looking in the opposite directions from $e$, 
as on the right in Figure \ref{fg-10}. 


We are going to define a groupoid ${\cal G}_{\Gamma}$ with the objects given 
by the points $\{v_i(e)\}$, where $e$ runs through all the edges of 
$\Gamma$ and $i \in \{1,2\}$, which is equivalent to the 
fundamental groupoid of  
$T'S$.  Therefore local $G(F)$-systems on $T'S$ can be understood as 
functors ${\cal G}_{\Gamma} \to G(F)$.  
Here is a system of generators of the groupoid ${\cal G}_{\Gamma}$. 

1). In the torsor of paths in 
$T'_{x(e)}S$ connecting the vectors $v_1(e)$ and 
$v_2(e)$ there are two distinguished elements, $p_+$ and $p_-$, where  $p_+$ 
moves the vector $v_1(e)$ clockwise, and $p_-$ counterclockwise 
towards the vector $v_2(e)$. By definition  $p_{\pm} \in 
{\rm Mor}_{{\cal G}_{\Gamma}}(v_1(e), v_2(e))$. 

2). Let $v$ be a vertex of $\Gamma$. Choose a sufficiently small 
neighborhood ${\cal D}_v$ of $v$ containing the three edges 
sharing $v$. Then the space of nonzero tangent vectors at the points 
of $\Gamma \cap {\cal D}_v$ which are not tangent to $\Gamma$ is a union of three 
connected simply connected components.  Therefore if we denote by $e_1, e_2, e_3$ 
the three 
edges sharing $v$, each connected component has exactly two 
of the six vectors $v_i(e_j)$. By definition the unique path 
connecting each pair belongs to the morphisms in 
${\cal G}_{\Gamma}$. 
 
\begin{figure}[ht]
\centerline{\epsfbox{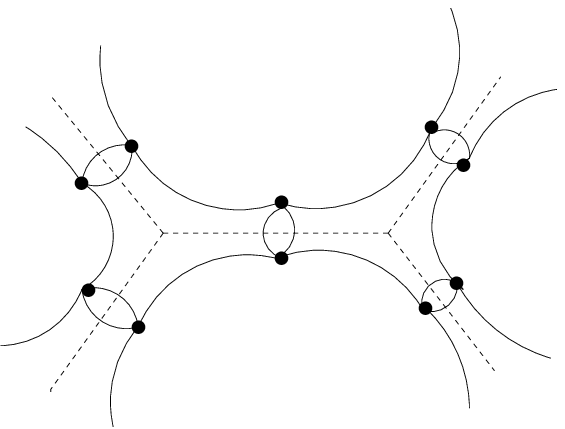}}
\caption{}
\label{fg-11}
\end{figure}

The arcs on the left in Figure \ref{fg-10} illustrate  the 
six generators of the groupoid near a vertex of $\Gamma$.  
They form a hexagon $H_v$. The arcs illustrate how the tangent vectors move. 
So for instance the left arc on Figure \ref{fg-10} intersecting $\Gamma$ 
moves the tangent vector looking up clockwise to the one looking down. 
Putting  the hexagons corresponding to all vertices of $\Gamma$ 
together we get the set of the 
generators of the groupoid ${\cal G}_{\Gamma}$. On Figure
\ref{fg-11} the 
 objects 
of the groupoid   ${\cal G}_{\Gamma}$ are pictured
 by fat points, and the generators by arcs. 


Here are two simple but crucial observations about this picture. 

\begin{lemma} \label{10.21.03.1}
a) The composition of the six path forming the hexagon $H_v$ on the left picture in 
(\ref{fg-10})  amounts to rotation of a tangent vector by $4\pi$. 

b) The composition $p_+p_+$ as well as $p_-p_-$ amounts to rotation of a tangent vector 
by $2\pi$. 
\end{lemma}

\begin{figure}[ht]
\centerline{\epsfbox{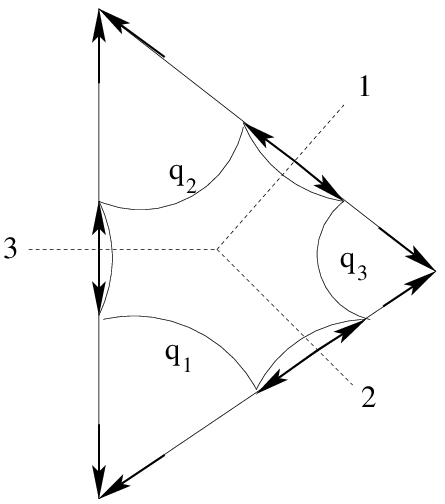}}
\caption{}
\label{fg-13}
\end{figure}

{\bf Proof}. A simple exercise. 
\vskip 3mm
Let us consider a triangle $t$ of the triangulation $T$ dual to a 
vertex $v$ of $\Gamma$.  
We may assume it intersects the edges $e_j$ of  $\Gamma$ 
sharing the vertex $v$ 
at the chosen points 
$x(e_j)$. 
Then there is a canonical homotopy class of 
path in $T'S$ going from the tangent vectors $v_i(e_j)$ along the 
corresponding side of the triangle $t$ to proximity of the corresponding 
vertex of the triangle $t$, as in Figure \ref{fg-13}.


Using this path and the decoration 
at this vertex of $t$ we picking up an affine flag 
over the tangent vector $v_i(e_j)$. 
Let $q$ be a side of the hexagon $H_v$ which 
does not intersect $\Gamma$ on Figure \ref{fg-11} (i.e. sits in one of the 
connected components discussed in 2)). Then the affine flags sitting 
at the endpoints of $q$ 
are obtained one 
from the other by parallel transport along  $q$.

\begin{figure}[ht]
\centerline{\epsfbox{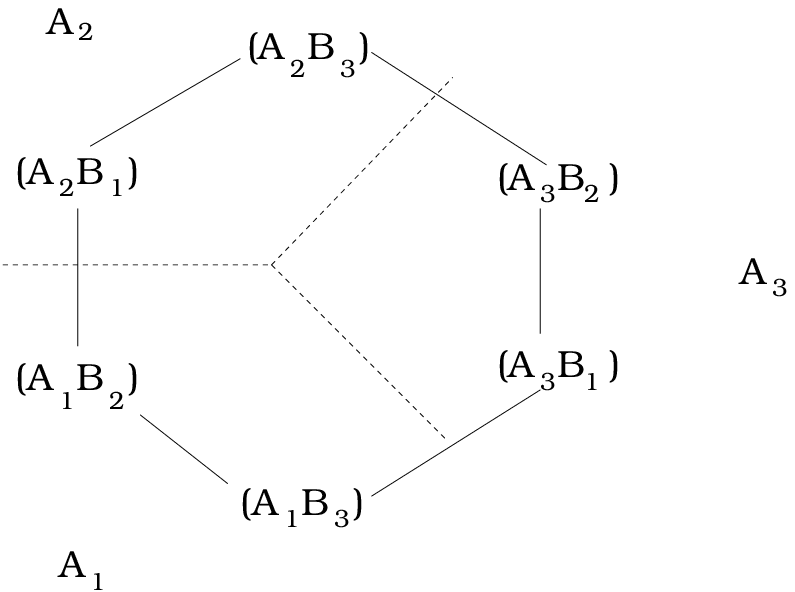}}
\caption{}
\label{fg-12}
\end{figure}

Lemma \ref{10.21.03.1}a) implies that the monodromy 
around the hexagon $H_v$ near a vertex $v$ is the identity. 
So the restriction of the local system ${\cal L}$ to the hexagon $H_v$ is trivial. 
Therefore we can transform all the affine flags into a fiber over one 
point, getting a configuration of three affine flags $(A_1, A_2, A_3)$ in $G$ 
corresponding to the vertex $v$, as on Figure \ref{fg-12}. 
The affine flags are parametrised by the vertices of the
 triangle $t$. 
Denote by $B_1, B_2, B_3$ the corresponding flags. Then we can attach 
the pairs $(A_i, B_j)$ to the vertices of the hexagon $H_v$ 
as in Figure \ref{fg-12}. 


Corresponding to each oriented edge ${\bf a}$ of the 
hexagon $H_v$  there is an element $g_{\bf a}$ 
of $G$ defined as follows. 
Consider a triple of flags $(A_1, A_2, A_3)$ 
representing our configuration.  Let  $(A_i, B_j)$ and $(A_{i'}, B_{j'})$ be the 
pairs assigned to the
first and second vertices of the oriented edge ${\bf a}$. There exists a unique element 
$g \in G$ transforming $(A_i, B_j)$ to the standard pair $(U^-, B^+)$.  
Set $(A'_1, A'_2, A'_3):= (gA_1, gA_2, gA_3)$. 
We define $g_{\bf a}$ as the unique element of $G$ such that 
$(A'_{i'}, B'_{j'}) = g_{\bf a}(U^-, B^+) $. 

There are two types of edges of the graph on Figure \ref{fg-11}: 
the ones intersecting the graph $\Gamma$ on Figure \ref{fg-11},
and the ones which do not intersect it. Each edge
 of the first type determines an edge of the triangulation $T$,
 and each edge of the second type determines a vertex 
of $T$. If ${\bf p}$ (resp. ${\bf q}$) is an edge of the first (resp. second) type,
then 
\begin{equation} \label{10.27.03.1}
g_{\bf p}  \in
\overline w_0H, \quad g_{\bf q}  \in
U^-.
\end{equation}

Each edge of the graph on Figure \ref{fg-11} is assigned 
to a unique hexagon $H_v$. \!Applying the above construction to\! the oriented
edges ${\bf a}$
of these hexagons we get a collection of elements $\{g_{\bf a}\}$ of $G$. It
has the following properties:

i) Reversing orientation of an edge amounts to replacing the corresponding
element 
by its inverse. 

ii) The product of the elements corresponding to a path around the hexagon is 
the identity.

iii) The product of the elements corresponding to a simple 
path around a two-gon is $s_G$.

iv) Let ${\bf q}$ (resp. ${\bf p}$) be an oriented edge corresponding to a
vertex  
(resp. edge) of the
triangulation $T$.  Then one has (\ref{10.27.03.1}).

By Lemma \ref{10.21.03.1}b) 
the monodromy around each of the two-gons on Figure \ref{fg-11} must be $s_G$. 
This  agrees with iii).

We can invert this construction. 

\begin{lemma} \label{10.20.03.20} Let  $\{g_{\bf a}\}$ be an arbitrary collection of 
elements of $G(F)$ assigned to 
oriented edges of the graph on Figure (\ref{fg-11}), which 
satisfy the conditions 
i) - iv). Then  there exists a twisted decorated $G(F)$-local system 
on $\widehat S$ such that the collection of elements  $g_{\bf a}$ constructed for it 
coincide with the original collection.
 \end{lemma}

{\bf Proof}. Thanks to the properties  i) - iii) there exists a unique 
functor ${\cal G}_{\Gamma} \to G(F)$
determined by the condition 
${\bf a} \lms g_{\bf a}$. Here ${\bf a}_1{\bf a}_2$ means that the path ${\bf a}_1$
follows ${\bf a}_2$. It provides a twisted $G(F)$-local system on 
$\widehat S$. The condition  iv) is used to define a decoration. 
The lemma is proved. 
\vskip 3mm
  
Recall that a point of the variety ${\cal A}_{G,T}$ is given by a collection 
of twisted cyclic configurations $(A^t_1, A^t_2, A^t_3)$ of affine flags in
generic position, one for each triangle $t$
of $T$,  compatible over the internal edges of $T$. Given such a configuration 
$(A^t_1, A^t_2, A^t_3)$ attached  to a vertex $v$ of $\Gamma$ we assign 
to the oriented edges of the  hexagon $H_v$ the elements 
$g_{\bf a}$ as above. They satisfy the conditions i), ii), iv). 
Finally, the compatibility over the edges implies the condition iii). It
remains to apply the construction from Lemma \ref{10.20.03.20}. The theorem is
proved. 

It remains to define a positive structure on the variety ${\cal A}_{G,T}$ and prove that
the positive structures corresponding to different ideal triangulations $T$
are compatible.

Consider $\widetilde {\rm Conf}_3{\cal A}$ as a subspace of 
${\cal A}_{G, \widehat t}$. Then choosing 
a vertex  of the triangle $t$, 
 we have 
$$
\widetilde {\rm Conf}^*_3{\cal A} =  {\cal V}_{G} \times 
(\widetilde {\rm Conf}^*_2{\cal A})^3.
$$
Recall that a positive structure on ${\cal V}_{G}$ is given by 
Proposition \ref{9.19.03.5}. The positive coordinate systems are parametrized by
reduced decompositions ${\bf i}$ of the element $w_0$. 
A positive structure on $\widetilde {\rm Conf}_2{\cal
  A}$ is given by Proposition \ref{4.21.03.1} and Lemma \ref{4.21.03.156}.  
This suggests the following definition. 

\begin{definition} \label{9.28.03.7q} Denote by ${\bf T}$ the following data
  on $\widehat S$: 

i) An ideal  triangulation 
$T$ of a marked hyperbolic surface $\widehat S$.

ii) For each triangle $t$ of the triangulation 
$T$ a choice of a vertex of this triangle. 

iii) A reduced decomposition ${\bf   i}$ 
of the maximal length element $w_0
\in W$. 
\end{definition}
So given such a data $\bf T$ and using  Theorem \ref{4.03.06.1} we have an open embedding
\begin{equation} \label{9.28.03.11yx}
\nu_{{\bf T}}: {\cal V}_{G}^{\{\mbox{triangles of $T$}\}} \times 
H^{\{\mbox{edges of $T$}\}} \hra {\cal A}_{G, \widehat S}.
\end{equation}

\begin{theorem} \label{4.03.06.1bn}  Let $G$ be a
 split simply-connected semi-simple algebraic group. 
 Let $\widehat S$ be a  marked hyperbolic
 surface. Then 
the collection of regular open embeddings $\{\nu_{{\bf T}}\}$, when 
${\bf T}$ runs through all the data from Definition \ref{9.28.03.7q} 
assigned to $\widehat S$, provides
a $\Gamma_S$--equivariant 
positive structure on ${\cal A}_{G, \widehat S}$. 
\end{theorem}

\begin{figure}[ht]
\centerline{\epsfbox{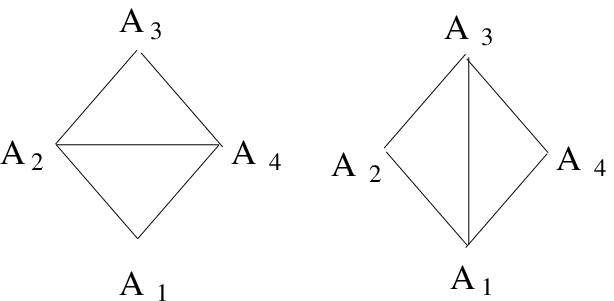}}
\caption{Flip.}
\label{fgo1a}
\end{figure}

{\bf Proof}. We have to check that we get a compatible positive structure by 
changing any of the three 
components in the data  ${\bf T}$ from the definition
\ref{9.28.03.7q}. For iii) this follows from Lusztig's theory \cite{L1}. 
For ii) this follows from Proposition \ref{4.03.6.10}. To prove the statement 
for i) we need to check that a flip as on Figure \ref{fgo1a}
provides a positive transformation. This follows from 
Proposition \ref{4.03.21.5}. The theorem is proved. 
\vskip 3mm

\begin{definition} \label{4.02.03.3wq} 
Let $\widehat S$ be a hyperbolic marked surface and $G$  a 
split simply-connected semi-simple algebraic group.  

a) The decorated Teichm{\"u}ller space  ${\cal A}^+_{G, \widehat S}$ is 
the set of all $\R_{>0}$-points of the positive space  ${\cal A}_{G, \widehat S}$.

b) 
Let ${\Bbb A}^t$ be one of the tropical semifields 
$\Z^t, \Q^t, \R^t$.  
The set  ${\cal A}_{G, \widehat S}({\Bbb A}^t)$ 
of points of the positive variety ${\cal A}_{G, \widehat S}$
with values in the tropical semifield ${\Bbb A}^t$ is called the set 
of ${\cal A}^{{\Bbb A}}$-laminations on $\widehat S$ corresponding to  $G$. 
\end{definition}
The space of real ${\cal A}$-laminations 
serves as  the Thurston boundary of the Teichm{\"u}ller space 
${\cal A}^+_{G, \widehat S}$.

The decorated Teichm{\"u}ller space ${\cal A}^+_{G, \widehat D_n}$ 
is identified with the space of positive twisted cyclic configurations of $n$   affine 
flags in $G(\R)$ on  a disc  $\widehat D_n$ with $n$ marked points on the
boundary.
\begin{corollary}  \label{4.03.21.10} 
The space ${\cal A}^+_{G, \widehat D_n}$ is invariant under the twisted cyclic shift. 
\end{corollary}

\vskip 3mm
{\bf 5. The universal Teichm{\"u}ller ${\cal A}$-space for $G$}. 
Let ${\Bbb S}^1(\R)$ be the set of all rays in $\R^2-\{0,0\}$. It is  a $2:1$ covering of
${\Bbb P}^1(\R)$. Let  ${\Bbb S}^1(\Q)$ be the subset of its rational
points, i.e. arrows with rational slopes. 
Let $s$ be  
the antipodal involution on ${\Bbb S}^1(\R)$. It is the nontrivial 
automorphism of the covering.

Let 
$\{p_1, ..., p_n\}$ be a cyclically ordered subset 
of ${\Bbb P}^1(\Q)$, i.e.  its order 
is compatible with one of the orientations of ${\Bbb P}^1(\R)$. Then a
choice of an initial point $p_1$ 
plus  its lift $\widetilde p_1$ 
determines a lift 
$\{\widetilde p_1, ..., \widetilde p_n\}$ of the set $\{p_{1}, ..., p_n\}$
to  ${\Bbb S}^1(\Q)$: we lift an oriented loop started at $p_1$,
whose orientation agrees with the cyclic order of the points $p_i$, 
to a path starting  at 
$\widetilde p_{1}$, and lift the points $p_{1}, ..., p_n$ using this 
path.  We say $\{\widetilde p_1, ..., \widetilde p_n\}$
is a {\it coherent lift} of
the cyclically ordered set $\{p_1, ..., p_n\}$.

\begin{definition} \label{1.29.04.3f} 
Let $G$ be a split semi-simple simply-connected algebraic group. 
A map 
\begin{equation} \label{1.29.04.1n}
\alpha: {\Bbb S}^1(\Q) \lra {\cal A}(\R)\qquad \mbox{such that
  $\alpha (s p) = s_G \alpha (p)$}
\end{equation} 
is {\em positive} if for any cyclically ordered $n$-tuple of points 
 $\{p_1, ..., p_n\}$ on ${\Bbb P}^1(\Q)$, and a coherent lift $\{\widetilde p_1, ...,
 \widetilde p_n\}$ of these
 points,  
the configuration of affine flags  $ 
(\alpha(\widetilde p_1), ..., \alpha (\widetilde p_n))$ is positive. 

The {\em universal  decorated Teichm{\"u}ller space} ${\cal A}^+_G$ 
     is the quotient of the space of positive maps (\ref{1.29.04.1n})
by the action of the group $G(\R)$. 
\end{definition}

Observe that if the configuration of affine flags  $ 
(\alpha(\widetilde p_1), ..., \alpha (\widetilde p_n))$ is positive for
a certain coherent lift of the cyclically ordered set $\{p_1, ..., p_n\}$, it 
is positive for any such coherent lift. So the positivity is a property of 
 $\{p_1, ..., p_n\}$. Further, it follows from the main results of this
 Section that it is sufficient to check the positivity 
 for $n=2, 3$ only. 

Recall the positive space ${\cal V}_{G}$, see (\ref{9.19.03.6}). Let ${\cal
  V}^+_{G}$ be its positive part.  
Identifying ${\Bbb P}^1(\Q)$ with vertices of the Farey triangulation, we
have a decomposition theorem for the universal  Teichm{\"u}ller space
${\cal A}^+_G$:

\begin{theorem} \label{5.3.04.1} There exists a canonical isomorphism 
$$
{\cal A}^+_G = {{\cal V}^+_{G}}^{\{\mbox{Farey triangles}\}} \times
H(\R_{>0})^{\{\mbox{Farey diagonals}\}}.
$$ 
\end{theorem}
To  fix an isomorphism in this theorem, 
we have to choose a vertex for each of the Farey
triangles. 

The higher Teichm\"uller spaces ${\cal A}^+_{G, S}$ are embedded  
into the universal one just in the same fashion as for the $\cal X$-space. 
Here is a more invariant way to think about this. 

\vskip 3mm
{\bf 6. The space ${\cal A}_{G, \widehat S}$ as a configuration space.} 
Let $\widetilde {\cal F}_{\infty}({\widehat S})$ be the cyclic set obtained as the $2:1$ cover  of the cyclic 
set ${\cal F}_{\infty}({\widehat S})$. Given a hyperbolic metric on $S$, it is induced 
by the $2:1$ cover of the absolute $\partial {\cal H}$. Let $\sigma$ be the non-trivial automorphism 
of this cover. 
Recall the central extension $\overline \pi_1(S)$ of $\pi_1(S)$ by the group $\Z/2\Z$ defined in Section 2.4. 
The generator of the subgroup $\Z/2\Z$ is denoted by $\overline \sigma_S$. The group $\overline \pi_1(S)$ acts by automorphisms of the cyclic set 
$\widetilde {\cal F}_{\infty}({\widehat S})$; in particular the element 
$\overline \sigma_S$ acts by the automorphism $\sigma$.

The moduli space ${\cal A}_{G, \widehat S}(\C)$ parametrises 
the pairs $(\psi, \rho)$ modulo $G(\C)$-conjugation, 
where $\rho: \overline \pi_1(S)\to G(\C)$ is a representation,  
$\rho(\overline \sigma_S) = s_G$, and 
$$
\psi: \widetilde {\cal F}_{\infty}({\widehat S}) \lra {\cal A}(\C), \quad \psi(\sigma) = s_G.
$$
is a $\rho$-equivariant map. 
The Teichm\"uller space ${\cal A}^+_{G, \widehat S}$ parametrises the pairs $(\psi, \rho)$ 
where $\psi$ is a positive map (and hence both $\psi$ and $\rho$ are real). 

\section{Special coordinates on the ${\cal A}$ and ${\cal X}$ spaces when   $G$ is of type $A_m$}
\label{SLN}

Let $\Gamma$ be a marked trivalent ribbon   graph. 
We assume it is not a special graph in the sense of Section 3.8: 
the special graphs are treated as in Section 10.7. 
Then there are varieties ${\cal X}_{G, \Gamma}$ and $
{\cal A}_{G, \Gamma}$ 
parametrising
framed and decorated unipotent $G$-local systems on $\Gamma$. 
If $\Gamma$ is  of type $\widehat S$, 
they  are 
isomorphic to  ${\cal X}_{G, \widehat S}$ and ${\cal A}_{G, \widehat S}$. 
The corresponding  isomorphisms  depend on a choice of an 
isotopy class of  embedding of $\Gamma$ to $\widehat S$. 
They form a principal homogeneous space over the mapping class group. 

In this chapter we introduce some natural functions 
$\{\Delta_{i}\}$ and $\{X_{j}\}$ on the varieties 
${\cal A}_{\rm SL_m, \Gamma}$ and 
${\cal X}_{\rm PGL_m, \Gamma}$. 
The combinatorial 
data parametrising these functions has been 
described already in the  Section 1.13, and will be briefly recalled below.

\vskip 3mm
{\bf 1. The sets ${{\rm I}}_m^{\Gamma}$ and 
${{\rm J}}_m^{\Gamma}$ parametrising the canonical coordinates}. 
Let $\Gamma$ be a marked 
trivalent graph  on $ S$ of type  
$\widehat S$.  
Then there is an ideal  triangulation $T_{\Gamma}$ of 
$\widehat S$ dual to $\Gamma$. It is determined up to an 
 isomorphism 
by the ribbon structure of $\Gamma$. Precisely, let 
 $\Gamma_1$, $\Gamma_2$ be 
two graphs as above such that the corresponding 
marked ribbon graphs $\overline \Gamma_1$ and $\overline \Gamma_2$ are
isomorphic. 
Then an isomorphism 
$\overline \Gamma_1 \to \overline \Gamma_2$ determines uniquely an 
element $g$ of the mapping class group of $S$ such that $g\Gamma_1$ 
is isotopic to $\Gamma_2$, and hence  $gT_{\Gamma_1}$ is isotopic to 
$T_{\Gamma_2}$. 

Summarizing,  a marked ribbon trivalent graph  $\Gamma$ of type $\widehat S$ 
determines an ideal triangulation $T = T_{\Gamma}$ of a marked surface 
$\widehat S_{\Gamma}$ isomorphic to 
$\widehat S$. The isotopy class of an isomorphism $\widehat S_{\Gamma} \to
\widehat S$ is determined 
by  an isotopy class of an embedded graph $\Gamma$.  
In this chapter we work with the marked ribbon trivalent graph $\Gamma$, the 
corresponding triangulation $T= T_{\Gamma}$, and surface $\widehat S = 
\widehat S_{\Gamma}$.  
Recall the $m$-triangulation of an ideal triangulation $T_{\Gamma}$ of the surface 
$\widehat S_{\Gamma}$ defined in the  Section 1.13.  
It provides the oriented graph 
described there and denoted  
$T_m$. Recall the two sets, ${\rm I}_m^{T}$ and ${\rm J}_m^{T}$,
attached to the triangulation. Since the triangulation $T$ and the 
corresponding graph $\Gamma$ determine each other, 
we can use the upper script $\Gamma$ in the notations.  So
\begin{equation} \label{11.21.02.1a}
\begin{array}{rcl}{\rm I}_m^{\Gamma} &:=& \mbox{$\{$vertices of the $m$--triangulation of 
$T_{\Gamma}$$\}$} - \\ & &-  \mbox{$\{$vertices at the  punctures of $S$$\}$},\end{array}
\end{equation}
$$
{\rm J}_m^{\Gamma}:= {\rm I}_m^{\Gamma} - \mbox{$\{$ the 
vertices at the boundary 
of $S$ $\}$}.
$$ 
\vskip 3mm
{\bf Remark}. If $\Gamma \lms \Gamma'$ is a flip then there is canonical bijection 
${\rm I}_m^{\Gamma} \lra {\rm I}_m^{\Gamma'}$. However the oriented graphs $T_m$ 
and $T'_m$ are different. 
\vskip 3mm
{\bf Examples}. 1. The set ${\rm J}_2^{\Gamma}$ is identified with 
the set of all internal (= not ends) edges of $\Gamma$. 

2. $\widehat S = S$ if and only if ${\rm I}_m^{\Gamma} = {\rm J}_m^{\Gamma}$ 
is the same set.

\begin{figure}[ht]
\centerline{\epsfbox{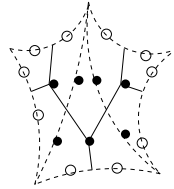}}
\caption{A $3$-triangulation of the disc with $5$ marked points on the boundary.}
\label{fg0}
\end{figure}

3. Figure \ref{fg0} illustrates 
the $3$-triangulation of a disk with $5$ marked points on the boundary. 
The disk is given by a pentagon.  
The ideal triangulation of the disk is given by triangulation 
of the pentagon. 
The boldface points illustrate the elements of the 
set ${\rm J}^{\Gamma}_3$. Adding to them the points 
visualized  by little circles we get the set ${\rm I}^{\Gamma}_3$. 
So $|{\rm J}^{\Gamma}_3|=7$ and $|{\rm I}^{\Gamma}_3|=17$.

{\it A parametrisation of the sets ${{\rm I}}_m^{\Gamma}$ and 
${{\rm J}}_m^{\Gamma}$}.  
Choose a vertex $v$ of $\Gamma$. It is dual to a triangle of the triangulation $T$. The inner  vertices of 
the $m$--triangulation of this triangle 
are described by  triples of nonnegative integers with sum $m$
\begin{equation} \label{11.25.02.10}
(a,b,c), \qquad a+b+c =m, \quad a,b,c \in \Z_{\geq 0}
\end{equation}
such that at least two of the numbers $a,b,c$ are different from zero. 

We say that a triple $(a,b,c)$ is of {\it vertex type} if  $a,b,c >0$, and picture 
 the integers  $a,b,c$ at the pieces of faces of $\Gamma$ sharing 
the vertex $v$:

\begin{figure}[ht]
\centerline{\epsfbox{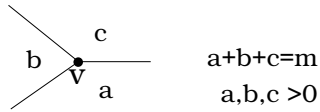}}
\caption{A triple $(a,b,c)$ of vertex type.}
\label{fg5}
\end{figure}
If one (and hence the only one) of the integers $a,b,c$ is zero, the other two 
determine an edge $e$ containing $v$. In this case we say that 
 the triple $(a,b,c)$ is of {\it edge type}:

\begin{figure}[ht]
\centerline{\epsfbox{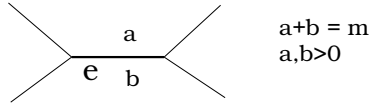}}
\caption{A pair $(a,b)$ of the edge type.}
\label{fg6}
\end{figure}

The vertex/edge type elements correspond to the vertices of 
the $m$--triangulation located inside of the triangles/edges of the dual 
triangulation $\widehat\Gamma $.

\vskip 3mm
{\bf 2. The regular functions $\Delta_{i}$ on 
the space ${\cal A}_{\rm SL_m, \Gamma}$}. They are parameterised by $i \in {\rm I}_m^{\Gamma}$. 

We use the following notation. Let $V$ be a vector space of dimension $n$ 
and $\omega \in {\rm det}V^*$  a volume form in $V$. 
Let $v_1, ..., v_n$ be vectors of $V$. Set
$$
\Delta(v_1, ... , v_n):= \Delta_{\omega}(v_1, ... , v_n):= <v_1 \wedge ... \wedge v_n, \omega>.
$$

Let $i \in {\rm I}_m^{\Gamma}$. Then there is a triangle 
of the triangulation $T$ containing $i$, 
 dual to a vertex $v$ of $\Gamma$. 
This triangle is unique if $i$ is of vertex type. 
Otherwise there are two such triangles. 

\vskip 3mm
{\it The vertex functions}. Let  ${\cal L}_v$ be the fiber of the local system ${\cal L}$ at a vertex $v$. Then, by definition, it carries a monodromy  invariant volume form 
$\omega$. There are three monodromy operators $M_{v, \alpha}$, $M_{v, \beta}$ and 
$M_{v, \gamma}$ 
acting on ${\cal L}_v$. They  correspond to the three face paths $\alpha, \beta, \gamma$ 
starting at  $v$. Given a 
 decorated  ${\rm SL}_m$--local system ${\cal L}$ 
on $\Gamma$ we get   three affine flags  
$$
X = (x_1, ..., x_m); \quad Y = (y_1, ..., y_m); \quad Z = (z_1, ..., z_m)
$$ 
at ${\cal L}_v$. Each of the flags $X$, $Y$ and $Z$ is 
 invariant under the corresponding monodromy operator:
$$
M_{v, \alpha}(X) = X, \quad M_{v, \beta}(Y) = Y, \quad M_{v, \gamma}(Z) = Z.
$$ 
Recall that $i$ is parametrised by a solution 
$(a, b, c)$ of (\ref{11.25.02.10}), and  that the numbers $(a,b,c)$ match the flags $X, Y, Z$. 
Put $x_{(a)}:= x_1 \wedge ... \wedge x_a$. 
The vertex function $\Delta_{i} = \Delta^v_{a,b,c}$ is defined by 
$$
\Delta_{i}({\cal L}): = \Delta^v_{a,b,c}(X, Y, Z) := 
\Delta_{\omega}(x_{(a)}\wedge  y_{(b)} \wedge 
 z_{(c)}) = $$ $$ =
\Delta_{\omega}(x_1 \wedge  ... \wedge x_{a}\wedge  y_1 
\wedge ... \wedge y_{b} \wedge 
 z_{1} \wedge  ... \wedge z_{c}).
$$

Recall that the element $s_G$, which we use to define the twisted cyclic shift, equals 
$(-1)^{m-1}$ for $G=SL_m$.  
The vertex function is invariant under the twisted cyclic shift 
of affine flags:
\begin{equation} \label{9.23.05.1}
\Delta^v_{a,b,c}(X, Y, Z) = \Delta^v_{c,a,b}((-1)^{m-1}Z, X, Y)
\end{equation}
Indeed, the cyclic shift 
multiplies $\Delta^v_{a,b,c}$ by $(-1)^{c(a+b)}$. If $m=a+b+c$ is odd, then $c(a+b)$ is even. 
If $a+b+c$ is even, we get $(-1)^{c^2}= (-1)^{c}$. So 
$\Delta^v_{a,b,c}(X, Y, Z) = (-1)^c\Delta^v_{c,a,b}(Z, X, Y)$ for even $m$. On the other hand \newline
$\Delta^v_{c,a,b}(-Z, X, Y) = (-1)^c\Delta^v_{c,a,b}(Z, X, Y)$. Thus we get (\ref{9.23.05.1}).

\begin{figure}[ht]
\centerline{\epsfbox{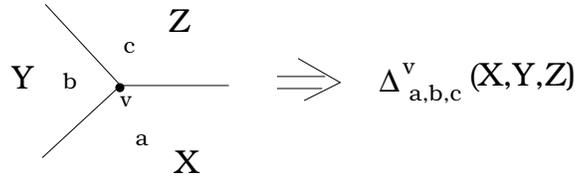}}
\caption{The data giving rise to a vertex function.}
\label{fg1}
\end{figure}

\vskip 3mm
{\it The edge functions}. Let  ${\cal L}_e$ be the fiber of the local system ${\cal L}$ at an edge $e$. The $SL_m$--structure provides a volume form $\omega$ in ${\cal L}_e$. 
There are 
two affine flags 
$$
X = (x_1, ..., x_m); \quad Y = (y_1, ..., y_m)
$$
at  ${\cal L}_e$ invariant under the monodromies around the two face paths $\alpha$ and $\beta$
  sharing the edge $e$:
$M_{v, \alpha}(X) = X, M_{v, \beta}(Y) = Y$. 
 We picture each of the flags located at the corresponding face. 
We define the edge functions by 
$$
\Delta_{i}({\cal L}) = \Delta^e_{a,b}(X,Y) := 
\Delta_{\omega}(x_{(a)}\wedge  y_{(b)} ) = 
\Delta_{\omega}(x_1 \wedge  ... \wedge x_{a}\wedge  y_1 \wedge ... \wedge y_{b}).
$$
The same argument as for the vertex function shows that the edge function is invariant under the twisted cyclic shit of affine flags $(X,Y) \lms ((-1)^{m-1}Y, X)$.

\begin{figure}[ht]
\centerline{\epsfbox{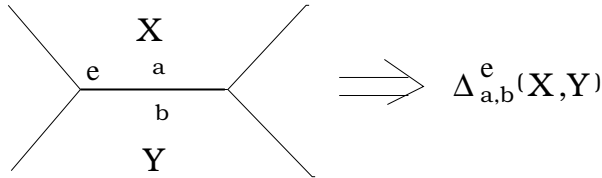}}
\caption{The data giving rise to an edge function.}
\label{fg2}
\end{figure}

\vskip 3mm
{\bf 3. The rational functions  $X_{j}$ on the space 
${\cal X}_{PGL_m, \Gamma}$}. 
Let $i \in {\rm J}_m^{\Gamma}$. Take all the vertices $j\in {\rm I}_m^{\Gamma}$ 
connected with $i$ by an edge in the graph $T_m$. 
Then there is one (if $i$  is of vertex type) or two (if $i$ is of edge type) 
triangles of $T$ containing all these $j$'s. 
Consider the invariant flags attached to the vertices of these triangles. 
There are three such flags if $i$ is vertex type, and four if it is 
of the edge type, see the Figures \ref{fg1} and \ref{fg2}. We can arrange these flags 
to sit in the same vector space 
${\cal L}_v$ (respectively ${\cal L}_e$). 
For each of these flags choose an affine flag dominating it. We denote by $\widetilde F$ the 
affine flag dominating a flag $F$. 
Then for each $j$ as above we have the 
$\Delta_{j}$--coordinate calculated 
for these affine flags. We denote it $\widetilde \Delta_{j}$. 
  We set 
\begin{equation} \label{11.25.02.71}
X_{i}  = \prod_{j \in {\rm I}_m^{\Gamma}}
(\widetilde \Delta_{j})^{\varepsilon_{i j}}.
\end{equation}
The right hand side evidently does not depend on the choice of 
dominating affine flags. 
Below we work out this definition for the vertex and edge type $i$'s.

\vskip 3mm

{\it The vertex functions $X^v_{a,b,c}$}. They are given by 
$$
X^v_{a,b,c}(X,Y,Z):= 
\frac{\Delta_{a-1, b+1, c}(\widetilde X,\widetilde Y,\widetilde Z)\Delta_{a, b-1, c+1}(\widetilde X, 
\widetilde Y, \widetilde Z)
\Delta_{a+1, b, c-1}(\widetilde X, \widetilde Y, \widetilde Z)}{\Delta_{a+1, b-1, c}(\widetilde X, \widetilde Y, \widetilde Z)\Delta_{a, b+1, c-1}(\widetilde X, \widetilde Y, \widetilde Z)\Delta_{a-1, b, c+1}(\widetilde X, \widetilde Y, \widetilde Z)}.
$$
The structure of this formula is illustrated on Figure \ref{fg3}:

\begin{figure}[ht]
\centerline{\epsfbox{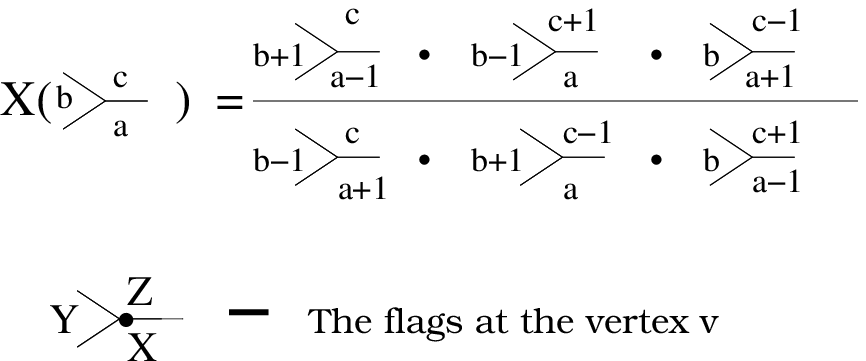}}
\caption{The vertex function.}
\label{fg3}
\end{figure}

\vskip 3mm

{\it The edge functions $X^e_{a,b}$}. They are given by 
$$
X^e_{a,b}(\widetilde X, \widetilde Y, \widetilde Z, \widetilde T):= 
\frac{\Delta_{a, b-1, 1}(\widetilde X, \widetilde Z, \widetilde T) 
\Delta_{a-1, 1, b}(\widetilde X, \widetilde Y, \widetilde Z)}
{\Delta_{a-1, b, 1}(\widetilde X, \widetilde Z, \widetilde T)
\Delta_{a, 1, b-1}(\widetilde X,\widetilde  Y, \widetilde Z)}.
$$
The structure of this formula is illustrated on Figure \ref{fg4}:

\begin{figure}[ht]
\centerline{\epsfbox{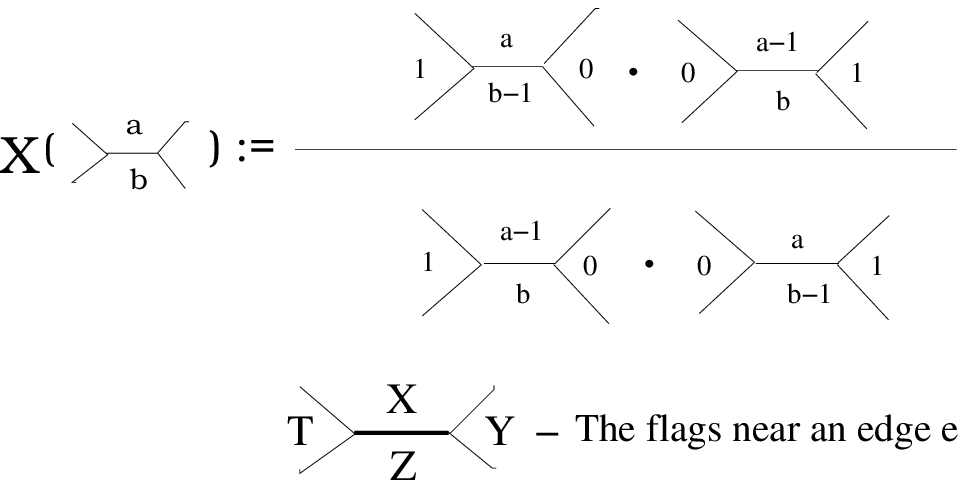}}
\caption{The edge function.}
\label{fg4}
\end{figure}

Recall the canonical 
projection 
\begin{equation} \label{11.5.03.10}
p: {\cal A}_{SL_n, \widehat S} \lra {\cal X}_{PGL_n, \widehat S}.
\end{equation}
\begin{proposition} \label{11.5.03.1} One has 
\begin{equation} \label{11.5.03.2}
p^*X_i = \prod_{j \in {\rm I}^{\Gamma}_m}(\Delta_j)^{\varepsilon_{ij}}.
\end{equation}
\end{proposition}

{\bf Proof}. This is a reformulation of the formula (\ref{11.25.02.71}). 
The proposition is proved. 
\vskip 3mm

It is easy to see using Lemma \ref{9.30.03.1} that 
the total number of the $X$--functions 
equals to the dimension of the 
space of all  $PGL_m$-local systems 
on $S$.

\vskip 3mm
{\bf 4. The cross-ratio and the triple ratio}. 
The space of orbits of the group $PGL_2$ acting on the 
$4$-tuples of distinct points on $P^1$ is one dimensional. 
There is a classical coordinate on the orbit space, called the {\it cross-ratio} 
of four points $x_1, ..., x_4$. To define it we present 
$P^1$ as the projectivisation of a two dimensional vector space $V$, 
and choose four vectors $v_1, ..., v_4$ in $V$ projecting to the given points 
$x_1, ..., x_4$ in $P(V)$. Then
\begin{equation}\label{11.6.02.1}
r^+(x_1, ..., x_4):= \frac{\Delta(v_1, v_2)\Delta(v_3, v_4) }
{\Delta(v_1, v_4)\Delta(v_2, v_3)}.
\end{equation}
One has $r^+(\infty, -1, 0, t) =t$, and 
$$
r^+(x_1, x_2, x_3, x_4) = r^+(x_2, x_3, x_4, x_1)^{-1};\,
r^+(x_1, x_2, x_3, x_4) = - 1- r^+(x_1, x_3, x_2, x_4). 
$$
To check the second equality we use the Pl{\"u}cker relation
$$
\Delta(v_1, v_2)\Delta(v_3, v_4)  - \Delta(v_1, v_3)\Delta(v_2, v_4) 
+ \Delta(v_1, v_4)\Delta(v_2, v_3)=0.
$$

\vskip 3mm
We  also need a {\it triple ratio} of three flags in $P^2$ as defined 
in s. 3.5 of \cite{G2}. 
The space of the orbits of the group $PGL_3$ acting on the 
triples of flags $F_1, F_2, F_3$ in $P^2$ is one dimensional. 
The triple ratio $r^+_3(F_1, F_2, F_3)$ 
is a natural coordinate on this space. To define it, 
 let us present $P^2$ as a projectivisation of 
a three dimensional vector space $V$. Choose three affine flags $A_1, A_2, A_3$ 
dominating the flags $F_1, F_2, F_3$ in $P(V)$. We represent the 
flags $A_i$ as the pairs $(v_i, f_i)$ where $v_i \in V, f_i \in V^*$ and 
$<f_i, v_i>=0$. Then set
$$
r^+_3(F_1, F_2, F_3):= \frac{<f_1, v_2><f_2, v_3><f_3, v_1> }
{<f_1, v_3><f_2, v_1><f_3, v_2>}.
$$
This definition obviously does not depend on the choice of the pair $(v,f)$
 representing a flag $F$. 
It is useful to rewrite  this definition using the presentation 
of the affine flags given by  
$$
A_1 = (x_1, x_1\wedge x_2), \quad A_2 = (y_1, y_1\wedge y_2), \quad 
A_3 = (z_1, z_1\wedge z_2).
$$
Then 
\begin{equation} \label{3.28.04.10}
r^+_3(F_1, F_2, F_3):= \frac{\Delta(x_1, x_2, y_1)\Delta(y_1, y_2, z_1)
\Delta(z_1, z_2, x_1) }
{\Delta(x_1, x_2, z_1)\Delta(y_1, y_2, x_1)\Delta(z_1, z_2, y_1)}.
\end{equation}

\begin{figure}[htbp]
\centerline{\epsfbox{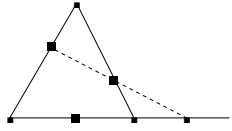}}
\caption{A positive triple of flags.}
\label{fg7}
\end{figure}

\begin{lemma} \label{8.2} A triple $(F_1, F_2, F_3)$ of real flags in $\R P^2$ 
 is positive if and only if the points $x_i$ are
 at the boundary of a convex domain bounded by the lines of the flags, 
as on Figure \ref{fg7}. 
\end{lemma}

{\bf Proof}. A generic triple of flags can be presented as follows: 
$v_1 = (0,1,1); v_2 = (1,0,1); v_3 = (1,t,0); f_1 = x_1, f_2 = x_2, f_3 = x_3$.  
Its triple ratio is $t$. Such a triple of flags satisfies the condition of the 
lemma if and only if $t>0$. The lemma is proved. 
\vskip 3mm
Let $F_i = (x_i, L_i)$ be a flag in $P^2$, so $x_i$ is a point on a line $L_i$, 
and $i=1,2,3$.  
Consider the triangle formed by the lines $L_1, L_2, L_3$. Denote by $y_1, y_2, y_3$ its vertices, so $y_i$ is opposite to the line $L_i$. Let $\widehat x$ be the point of intersection of the line $L_3$ and the line through the points $x_1$ and $x_2$. Then the proof of Lemma \ref{8.2} shows that 
$$
r_3^+(F_1, F_2, F_3) = r^+(y_1, \widehat x, y_2, x_3).
$$
\vskip 3mm
{\bf Remark}. One can find a yet different definition of the triple ratio 
in Section 4.2 of \cite{G5}. 
\vskip 3mm

It follows that the cross- and triple ratios provide isomorphisms
$$
r^+: \{\mbox{$4$--tuples of distinct points of ${\Bbb P}^1$}\}/PGL_2 \stackrel{\sim}{\lra} 
{\Bbb P}^1 - 
\{0, -1, \infty\},
$$
$$
r_3^+: \{\mbox{generic triples of flags of ${\Bbb P}^2$}\}/PGL_3 \stackrel{\sim}{\lra} 
{\Bbb P}^1 - 
\{0, -1, \infty\}.
$$
\vskip 3mm
{\bf Remark}. It is often useful to change the sign in the definition of  
the cross- and triple ratios: $r:= - r^+$ and $r_3:= -r_3^+$. 
\vskip 3mm

{\bf 5. Expressing the $X$-functions via the  cross- and triple ratios}. 
{\it The vertex functions}.  Recall that $F_k$ is the $k$-dimensional subspace of 
a flag $F$. For generic configuration of flags $(X, Y, Z)$ the quotient
\begin{equation} \label{4.18.02.1}
\frac{L_v}{X_{a-1} \oplus Y_{b-1} \oplus Z_{c-1}}
\end{equation}
is a three dimensional vector space. It inherits a configuration of three 
flags  $(\overline X, \overline Y, \overline Z)$. Namely, the flag $\overline X$ is 
given by the projection of the flag \newline $(X_a/X_{a-1} \subset X_{a+1}/X_{a-1} \subset \cdots )$
onto (\ref{4.18.02.1}), and so on. 
 The vertex $X$--function $X_{a,b,c}$ equals to  their triple ratio: 
$$
X_{a,b,c}(X,Y,Z) = 
r^+_3(\overline X, \overline Y, \overline Z).
$$
\vskip 3mm
{\it The edge $X$--functions}. 
For a generic configuration of flags $(X, Y, Z, T)$ in $L_e$ then the quotient 
$$
\frac{L_e}{X_{a-1} \oplus Y_{b-1}}
$$
is a two dimensional vector space. It inherits a configuration $(\overline X_a, \overline Y_1, \overline Z_b, 
\overline T_1)$ of 
four distinct 
one dimensional subspaces:  $\overline X_a := X_a/X_{a-1}$ and so on. 
Their cross--ratio is the edge $X$-function $X_{a,b}$: 
$$
X_{a,b}(X,Y,Z, T) = r^+(\overline T_1, \overline X_a, \overline Y_1, \overline Z_b).
$$
\vskip 3mm
{\bf 6. The main result}. 
Let us choose an isotopy class $\Gamma$ of a marked 
trivalent graph on $\widehat S$ 
isomorphic to a marked ribbon graph $\Gamma$. It provides an 
isomorphism of the moduli spaces 
$f_{\Gamma}: {\cal X}_{PGL_m, \Gamma} \lra {\cal X}_{PGL_m, \widehat S}$. 
Thus the rational function $X_j$ gives rise to a rational function $X^{\Gamma}_{j}$ 
on the moduli spaces 
${\cal X}_{PGL_m, \widehat S}$. 
Similarly there are regular functions $\Delta^{\Gamma}_{i}$.

\begin{theorem} \label{5.6.02.3} 
a) The collections of rational functions $\{X^{\Gamma}_{j}\}$, 
 $j \in {\rm J}^{\Gamma}_m$, where 
$\Gamma$ runs through the set 
of all isotopy classes of marked trivalent 
graphs on $\widehat S$ of type $\widehat S$, 
 provide a positive regular atlas on
 ${\cal X}_{PGL_m, \widehat S}$. 

b) The collections of regular functions $\{\Delta^{\Gamma}_{i}\}$,
 $i \in {\rm I}^{\Gamma}_m$, 
provide a positive regular atlas on
 ${\cal A}_{SL_m, \widehat S}$. 
\end{theorem}

{\bf Plan of the proof}. 
In Sections 9.7 and 9.8 we give a yet another definition of the
$X$--coordinates on the configuration space of triples and quadruples of
  flags in $V_m$. We show in Lemmas \ref{5.20.02.1}, \ref{3.27.04.10}, 
and \ref{3.30.04.1} that it gives the same
  $X$-coordinates as in Section 9.2. 
The results of Section 9.8 make obvious that 
they are related to the general
  definition for a standard reduced decomposition of $w_0$, see Lemma \ref{3.30.04.1}. 
Similarly our $\Delta$-coordinates on the space ${\cal
  A}_{SL_m, \Gamma}$ are evidently related to the 
$A$-coordinates from Section 8 
for a standard reduced decomposition of $w_0$. 
However our special coordinates have some  
  remarkable features which are absent in the general case even for $SL_m$ and
  nonstandard reduced decompositions of $w_0$.  
  For example the coordinates defined in Section 9.2-9.3  
on the configuration space of triples of
  (affine) flags in $V_m$ are manifestly invariant under the cyclic (twisted cyclic) 
shift.

So 
the positivity statements in Theorem \ref{5.6.02.3} are   special cases of  the ones in 
Theorem \ref{7.8.03.1}. On the other hand these statements follow immediately from
the results of Section 10. This way we get simpler proofs of the
positivity results, which are independent from the proof of the
 Theorem \ref{7.8.03.1}. 

We show that 
our $X$-functions provide 
 a regular atlas in Section 9.9. 
A similar statement about 
the $\Delta$-coordinates is reduced to the one 
about the $X$-coordinates using the fact that the  projection 
${\cal A}_{SL_m, \widehat S} \lra {\cal X}_{PGL_m, \widehat S}$ 
is obtained by factorization along the action of a power of the Cartan group.

Let us start the implementation of this plan. 
\vskip 3mm
{\bf 7. Projective bases determined by three flags}. 
A {\it projective basis} in  a vector space $V$  is a basis in $V$  
up to a multiplication by a common scalar. Projective bases in $V$ 
form a principal homogeneous space for $PGL(V)$. 
Let $A, B, C$ be three flags at generic position in an $m$-dimensional vector space 
$V$. Below we introduce several different projective bases in  $V$ 
related to these flags,  
and compute the elements of $PGL_m$ transforming 
one of these bases to the other. 

\begin{figure}[htbp]
\centerline{\epsfbox{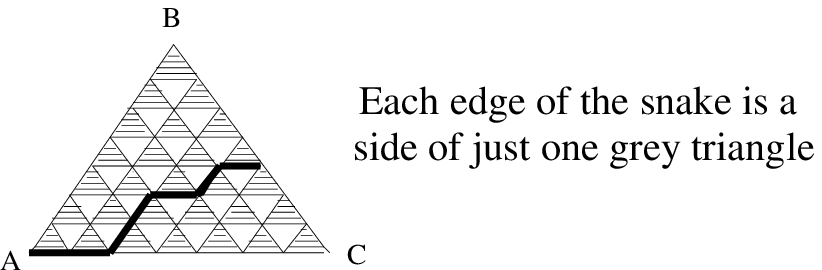}}
\caption{A snake.}
\label{fg15}
\end{figure}

Consider the $(m-1)$--triangulation of a triangle shown on the Figure \ref{fg15}. 
Each side of the triangle carries  $m$ vertices of the  $(m-1)$--triangulation. 
There are two types of the triangles, the triangles  looking down, 
and the triangles looking up. We picture them as the white and grey triangles.    
Let us assign the flags $A, B, C$ 
to the vertices of the big triangle. 
The side  across the $A$-vertex is called the $A$--side, etc.. 
The   vertices of  the $(m-1)$--triangulation  
are parametrized by non negative integer solutions of the equation
$
a+ b+ c =m-1
$. 
The $(a,b,c)$--coordinates of a vertex show
 the distance from the vertex to the $A$, $B$ and $C$ sides of the triangle. 
A vertex $A$ of the big triangle provides two arrows which start at $A$ and 
look along one of  the two sides  sharing  $A$. 
A {\it snake} is a  path  
in the one-skeleton of the $(m-1)$--triangulation from one of the vertices 
of the big triangle  to the opposite side,
 which each time goes in the direction of 
one of the two arrows assigned to the vertex.  

A snake provides 
a projective basis in  $V$. Indeed, 
let  $F^a:= F_{m-a}$ be
 the codimension $a$ subspace of a flag $F$. We attach to each vertex $(a,b,c)$ 
a one--dimensional subspace 
$
V^{a,b,c}:= A^a \cap B^b \cap C^c
$. 
The one--dimensional subspaces $V_1, V_2, V_3$ 
attached to the vertices of any grey triangle $t$
span a two--dimensional subspace $V(t)$. 
The subspaces 
$V_1, V_2, V_3$  provide six projective bases in $V(t)$, assigned to the
  oriented sides of the triangle $t$. 
\begin{figure}[ht]
\centerline{\epsfbox{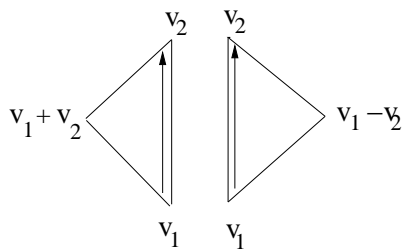}}
\caption{Defining two projective bases in a two dimensional space.}
\label{fgo-17}
\end{figure}
Namely, we assume that the triangle $t$ is counterclockwise oriented. 
Then if the side $V_1V_2$ is oriented according
to the orientation of $t$, we choose vectors  $v_1 \in V_1$ and 
$v_2\in V_2$ so that $v_1-v_2\in V_3$. 
If not, we choose vectors  $v_1 \in V_1$ and 
$v_2\in V_2$ so that $v_1+v_2\in V_3$,  as shown on the picture. 
\vskip 3mm
{\bf Remark}. This definition  of the six projective bases  
agrees with the one given in Section 6.2. 
\vskip 3mm

Each edge of the triangulation 
is a boundary of just one grey triangle.
Using this observation,  let us construct  a projective 
basis corresponding to a snake. 
Let $e_1, ..., e_{m-1}$ be the edges of the snake in their natural order, so 
the $A$--vertex belongs to $e_1$. Denote by $t_1, ..., t_{m-1}$ the 
grey triangles containing these edges. 
Each edge $e_i$ determines a projective basis in the subspace $V({t_i})$. 
A nonzero vector $v_1 \in V^{m-1, 0, 0}$ 
determines a natural basis $v_1, ..., v_m$ parametrized by  
vertices of the snake. Indeed, $v_1$ plus 
a projective basis in  $V(t_1)$ corresponding to  $e_1$  determines a 
vector $v_2$ in $V({t_1})$. The vector $v_2$ 
plus a projective basis in $V({t_2})$ corresponding to  $e_2$ 
provides  $v_3$, and so on.  All together they provide a projective basis.

\vskip 3mm
{\bf 8. Transformations between different projective bases related 
to three flags}. Let us 
 compute the element of $PGL_m$   transforming  the projective basis 
corresponding to a snake to the one  
corresponding to another snake. We assume that the group $PGL_m$ acts on 
projective bases from the left.  
First, let us do it for the snakes sharing the same 
vertex of the big triangle, say the vertex $A$. 
An elementary move of a snake is one of the following two local 
transformations on Figure \ref{fg16}: 

\begin{figure}[ht]
\centerline{\epsfbox{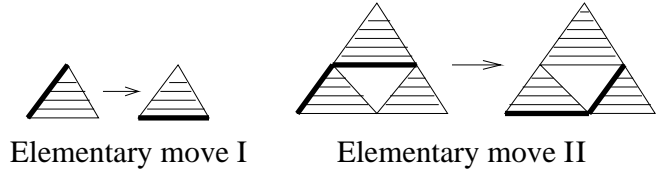}}
\caption{Elementary moves of the snake.}
\label{fg16}
\end{figure}

\begin{figure}[ht]
\centerline{\epsfbox{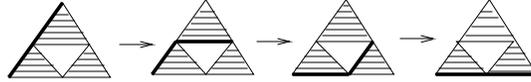}}
\caption{Moving the snake in the case $m=3$.}
\label{fg18}
\end{figure}

We can transform one snake  
to another using a sequence of elementary moves. 
An example for $m=3$ see on Figure \ref{fg18}. 
In general such a  sequence of elementary moves is not unique. 

Let us calculate the elements of $PGL_2$ and $PGL_3$ 
corresponding to the elementary moves ${\rm I}$ and ${\rm II}$. 
For the elementary move ${\rm I}$ we get 
$$
\left( \begin{matrix}v_1\\ v_2 \end{matrix}\right ) \lms  
  \left( \begin{matrix}1&0\\1 &1\end{matrix} \right ),\qquad  
\left( \begin{matrix}v_1\\ v_2 \end{matrix}\right )  = \left( \begin{matrix}v_1\\ v_1+v_2\end{matrix}\right ) .  
$$ 
Here the vectors  $v_1, v_2, v_1+v_2$ 
in $V(t)$ correspond to the vertices of a grey triangle $t$, so  
$(v_1, v_2)$ is the projective basis related to the left edge of $t$, and 
$(v_1, v_1 + v_2)$ is the one  related to the bottom edge of $t$.

Observe that the white triangles, which has been defined using the $(m-1)$-triangulation, 
 are in one
to one correspondance to the inner vertices of the $m$-triangulation, 
which have coordinates $a+b+c =m$ with $a,b,c >0$.  
The vertex $X$--coordinates are assigned to the white triangles. 
Namely, let $\alpha, \beta, \gamma$ 
be the distances from a  white triangle to the sides of the triangle, 
$\alpha + \beta + \gamma = m-3$. Then  the flags $A, B, C$ 
induce the flags $\overline A, \overline B, \overline C$ 
in the three 
 dimensional 
quotient 
\begin{equation} \label{5.20.02.2}
\frac{V}{A_{\alpha} \oplus B_{\beta} \oplus C_{\gamma}}.
\end{equation} 
The triple ratio of these three flags is the 
$X$--coordinate corresponding to the white triangle. 
Thus 
$X_{\alpha , \beta , \gamma} = X_{a,b,c}$ with $a=\alpha +1$,  $b=\beta +1$, 
 $c=\gamma +1$. 

Now consider an elementary move of type $II$ corresponding to 
a white triangle. Let us calculate  the corresponding
element of $PGL_m$, and show that it 
 depends only on the $X$-coordinate assigned to 
this white triangle. 
Denote by 
$a$ a vector spanning the one dimensional 
subspace $\overline A_1$ of the flag $\overline A$. 
There are two  choices for the vectors spanning 
$\overline B_2 \cap \overline C_2$: the vector $r_1$ is defined via the 
first  snake, and 
the vector $r_2$ defined via the second one on the picture 
describing the elementary move II. 
Then 
\begin{equation} \label{5.20.02.3}
r_1 = X r_2. 
\end{equation}

\begin{lemma} \label{5.20.02.1} The number  $X$ in (\ref{5.20.02.3}) 
is  the $X$--coordinate of the flags $\overline A, \overline B, \overline C$. 
\end{lemma}

{\bf Proof}. A choice of vector $a$ provides the  vectors  
$b$, $c$, $p$, $q$ spanning the subspaces 
$\overline B_1$, $\overline C_1$, $\overline A_2 \cap \overline B_2$, 
$\overline A_2 \cap \overline C_2$
respectively, and defined 
using the left (for $b$ and $p$) and  the  right (for $c$ and $q$) 
snakes. The three grey triangles provide the relations 
\begin{equation} \label{5.20.02.4}
r_1  = p + b, \quad q = a + p, \quad c = q + r_2.  
\end{equation}

\begin{figure}[ht]
\centerline{\epsfbox{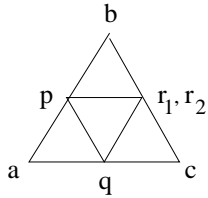}}
\caption{The vectors $p,q,r_1,r_2$.}
\label{fg19}
\end{figure}

Therefore, as easily 
follows from  (\ref{5.20.02.3}) and (\ref{5.20.02.4}), one has 
$$
X(\overline A, \overline B, \overline C)  = 
\frac{\Delta(b,p,c)}{\Delta(c,q,b)} = 
\frac{\Delta(r_1, p,c)}{\Delta(r_2, p, c)}.
$$
The first equality follows from the definition of triple ratio, applied to the 
three affine flags $(a, p), (b, r_1), (c, q)$ dominating the flags $\overline A, \overline B, \overline C$, using 
(\ref{5.20.02.4}). 
The second is obtained using (\ref{5.20.02.4}). Indeed, 
  $\Delta(b,p,c) = - \Delta(r_1,p,c)$, and $\Delta(c,q,b)= \Delta(r_2,q,b) = \Delta(r_2,q,p) = \Delta(r_2,c,p) = 
- \Delta(r_2,p, c)$. and  The lemma is proved.

\vskip 3mm
{\bf Example}. Let us compute the transformation between  the projective bases 
corresponding to the two snakes on Figure \ref{fg18}. 
We present it  as a composition: 
$
(a, p, b) \lms (a, p, r_1) \lms (a,q,r_1) \lms (a, q, r_2) \lms (a, q, c)  
$. 
So thanks to (\ref{5.20.02.4}) and Lemma \ref{5.20.02.1} 
the corresponding matrix is 
\begin{equation} \label{3.27.04.1}
\left ( \begin{matrix}1&0& 0\\ 0&1&0 \\ 0 & 1& 1\end{matrix}\right ) 
\left ( \begin{matrix}1&0& 0\\ 0&1&0 \\ 0 &0& X^{-1}\end{matrix}\right ) 
\left ( \begin{matrix}1&0& 0\\ 1&1&0 \\ 0 &0& 1\end{matrix}\right ) 
\left ( \begin{matrix}1&0& 0\\ 0&1&0 \\ 0 &1& 1\end{matrix}\right ) = 
\left ( \begin{matrix}1&0& 0\\ 1&1&0 \\ 1 & 1+X^{-1}& X^{-1}\end{matrix}\right ) \in PGL_3.
\end{equation}
\vskip 3mm
Let us generalize this computation. 
Let $\varphi_i: SL_2 \hra GL_m$ be the canonical embedding corresponding 
to the $i$-th root $\lambda_i - \lambda_{i+1}$. Set 
$$
F_i = \varphi_i\begin{matrix}1&0\\ 1&1\end{matrix} , \qquad H_i(x):= 
{\rm diag}(\underbrace{x, ..., x}_{i}, 1, ..., 1).
$$
So $H_i: {\Bbb G}_m \lra H$ is the coroot corresponding
to the root $\lambda_i - \lambda_{i+1}$. So $H_i(x)$ commutes with $F_j$ and
$E_j$ if $i \not = j$.  
So the left hand side of (\ref{3.27.04.1}) can be  written as 
$F_2H_2(X)F_1F_2 = F_2F_1H_2(X)F_2 $. 

Let $X_{\alpha, \beta, \gamma}$, where $\alpha +\beta+\gamma=m-3$, $\alpha, \beta, \gamma \geq 0$, 
be the vertex coordinates  
of the configuration of flags $(A, B, C)$. (They are the former coordinates $X_{abc}$,  
where $a=\alpha +1$, $b=\beta+1$, $c=\gamma+1$.) Let $M_{AB\to AC}$  
be the element transforming the projective basis corresponding to the snake
$AB$ to the one for the snake $AC$. 

\begin{proposition} \label{9.30.03.10aw}
a) One has 
\begin{equation} \label{5.18.03.11c}
M_{AB\to AC} = \prod_{j=m-1}^{1}\left(\Bigl(\prod_{i=j-1}^{m-2}H_{i+1}(X_{m-i-2,
      i-j, j-1})F_i\Bigl)F_{m-1}\right).
\end{equation}

 b) The element $M_{AB\to BA} $ is  the transformation 
$S: e_i \lms (-1)^{i-1}e_{m+1-i}$.

c) One has 
\begin{equation} \label{3.28.04.1}
M_{CA\to BA} = SM^{-1}_{AB\to AC}S^{-1}.
\end{equation} 
\end{proposition}

For example for $m=2$ we get $M_{AB\to AC} = F_1$, for $m=3$ we get
$$
M_{AB\to AC} = F_2H_2(X_{000})F_1F_2
$$
and for $m=4$
$$
M_{AB\to AC} = F_3H_3(X_{001})F_2F_3H_2(X_{100})F_1H_3(X_{010})F_2F_3.
$$
\vskip 3mm
{\bf Remark}. In the formula (\ref{5.18.03.11c}) between every two 
elements $F_i$,  $i>1$, which have no $F_i$'s in between,
there is unique element $H_i(X)$ between them. Its precise location 
between the two $F_i$'s is not essential since 
$H_i(X)$ commutes with all the $F$'s located between the two $F_i$'s. 
\vskip 3mm
 {\bf Proof}. 
a) Let us define a sequence of elementary moves transforming the snake 
$AB$ to the snake $AC$. The crucial intermediate 
steps of this sequence  are shown on the picture below.  
We show the snakes which go $j$ steps towards $C$, 
and then $m-j-1$ steps in the direction $AB$, for $i=0, ..., m-1$. 
In between we use the obvious  elementary moves 
transforming  the $(j-1)$-st snake to the the $j$-th. 
\begin{figure}[ht]
\centerline{\epsfbox{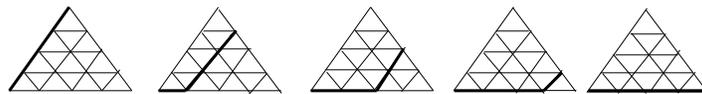}}
\caption{Moving the snake.}
\label{fg200}
\end{figure}

Now, just like in the above example, 
 our formulas and Lemma \ref{5.20.02.1} give the 
formula (\ref{5.18.03.11c}). 

b) Let $(A, B, C)$ be a cyclically ordered triple of 
distinct one dimensional subspaces in a two dimensional vector space. 
Let $[v_1, v_2]$ (respectively $[u_1, u_2]$) be the projective basis 
determined by the snake $AB$ (respectively $BA$). 
This means that 
$v_1 + v_2 \in C$ (respectively $u_1 - u_2 \in C$). It follows that 
$u_1 = v_2, u_2 = -v_1$. This proves the claim in the two dimensional
case. The general case immediately reduces to this. 

c) This is obvious. 
The proposition is proved. 
\vskip 3mm
{\bf Remark}.  Since $SF^{-1}_iS^{-1} = E_{m-i}$, there is a formula 
for $M_{CA\to BA}$ similar to (\ref{5.18.03.11c}). In this formula 
 we use the other standard reduced decomposition
of $w_0$ encoding the sequence of $E_i$'s then the one encoding 
the sequence of $F_i$'s in 
(\ref{5.18.03.11c}).
\vskip 3mm
Consider a configuration of four flags $(A, B, C, D)$. We assign them to
vertices of a quadrilateral triangulated by the diagonal $AC$.  
There are two
projective bases related to the snake $AC$: the one, $P_-(AC)$, comes
from the triangle $ABC$, and the other, 
$P_+(AC)$, comes from the triangle $ACD$. 

\begin{lemma} \label{3.27.04.10} Let  
 $x_a = X_{a, m-2-a}$  
 be the coordinates attached to the edge $AC$. 
Then  
\begin{equation} \label{5.18.03.101}
M_{P_-(AC)\to P_+(AC)} = {\rm diag}(x_{0} x_{1} ... x_{m-2}, ..., x_{m-3}x_{m-2}, x_{m-2}, 1). 
 \end{equation}
\end{lemma}

{\bf Proof}. Follows  from the $SL_2$ case, which was considered in the Example 
in the end of Section 6.6.

\vskip 3mm
{\bf 9. Comparison with the definitions of Sections $5$ and $6$}. 
Recall the subsets $\overline U^-_*({\bf i}) \subset U^-$ and $\overline U^+_*({\bf j}) \subset
U^+$ 
defined by 
reduced 
decompositions ${\bf i}$ and ${\bf j}$ of $w_0$, see (\ref{U^-_*}) and (\ref{U^+_*}). 
Elements $\overline u_+ \in  \overline U^+_*({\bf i})$ and 
$\overline u_- \in \in  \overline U^+_*({\bf i})$ have canonical decompositions:
\begin{equation} \label{3.31.04.13}
\overline u_+ = \prod_{j\in {\bf j}}x_j(s_j), \quad 
\overline u_- = \prod_{i\in {\bf i}}y_i(t_i), \qquad s_j, t_i \in F^*.
\end{equation}
The elements of the sets ${\bf i}$ and ${\bf j}$ are labeled by 
the simple roots $\alpha_1, ..., \alpha_{m-1}$, providing the decompositions
$$
{\bf i} = {\bf i}(\alpha_1) \cup ... \cup {\bf i}(\alpha_{m-1})\qquad 
{\bf j} = {\bf j}(\alpha_1) \cup ... \cup {\bf j}(\alpha_{m-1}).
$$
Our convention was the following: let $\alpha$ be a simple root, 
$j^+_{\alpha}$ is the right element of  ${\bf j}(\alpha)$, 
$i^-_{\alpha}$ is the left element of  ${\bf i}(\alpha)$. Then  
for each simple root $\alpha$ one has 
\begin{equation} \label{3.31.04.12}
x_{j^+_{\alpha}}(s_{j^+_{\alpha}}) =1, \qquad y_{i^-_{\alpha}}(t_{i^-_{\alpha}}) =1. 
\end{equation}

Let $r=m-1$, $\alpha_k = \lambda_k - \lambda_{k+1}$, $k=1, ..., m-1$. 
Consider the following two reduced decomposition of the maximal length element 
in $S_m$: 
$$
{\bf i} = (\alpha_{m-1}, \alpha_{m-2}, \alpha_{m-1}, ..., \alpha_{2}, ..., \alpha_{m-1}, \alpha_{1}, \alpha_{2},..., \alpha_{m-1}),
$$
$${\bf j} = (\alpha_{m-1}, ... , \alpha_1, \alpha_{m-1}, ...,
\alpha_2, ..., \alpha_{m-1}, \alpha_{m-2}, \alpha_{m-1}). 
$$

Consider the $(m-1)$-triangulation of the triangle $ABC$, as on Figure
\ref{fg15}. Its grey triangles 
are presented as a union of $m-1$ strips going parallel to the
side $BC$. We enumerate the strips 
so that the  $k$-th strip contains  $k$ the grey triangles. 
These numbers match the simple roots: $k$ corresponds to $\lambda_k - \lambda_{k+1}$. 
Then the string ${\bf i}$ provides an enumeration of the grey triangles 
from $C$ to $AB$ going in each layer parallel $AB$ from
$A$ to $B$. Let $t$ be a white
triangle and $t'$ the grey triangle sharing a side with $t$, obtained
from $t$ by going towards $B$. 
Recall that our special $X$-coordinates are naturally attached to the white triangles. 
We assign the product of the $X$-coordinates attached to all the white
triangles belonging to the layer parallel to $BC$ from $t$ to $AC$, to the corresponding grey triangle $t'$. 
Finally, we assign $1$ to the grey triangles neighboring the side $AC$. 

\begin{lemma} \label{3.30.04.1q} 
One has 
$$
M_{AB \to AC} = \overline u^-_{AB \to AC} h^-_{AB \to AC}, \qquad 
\overline u^-_{AB \to AC} \in U^-_*({\bf i}), \quad h^-_{AB \to AC}\in H.
$$

The canonical coordinates (\ref{3.31.04.13}) of the element $\overline u^-_{AB \to AC}$ 
are given by the $X$-coordinates assigned as
above to the grey triangles, provided the grey triangles are enumerated 
by the elements of the set ${\bf i}$ as was explained above. 

The element $h^-_{AB \to AC}$ is the product of the elements 
$H_{i_t}(X_{t'})$ assigned to the grey triangles $t'$. 
\end{lemma}

{\bf Proof}. Follows from the very definition and formulas
(\ref{3.31.04.13}), 
(\ref{5.18.03.11c}) 
and  (\ref{3.31.04.12}). 
\vskip 3mm
A projective basis in $V_m$ is the same thing as a pinning for
$PGL_m$. 

\begin{lemma} \label{3.30.04.1} 
The six projective bases  corresponding to the six oriented
edges of the triangle $ABC$ coincide with the six canonical pinnings 
attached to the triple of flags $(A, B, C)$ in $V_m$ in Section 5.2. 
\end{lemma}

{\bf Proof}. 
Lemma \ref{3.30.04.1q} implies that the projective basis assigned to 
the snake $AB$ coincides with the one assigned to the side $AB$ of
the triangle whose vertices are labeled by the configuration of flags
$$
(A, B, C) = (B^-, B^+, u_-B^+u_-^{-1}).
$$
The case of the snakes $BC$ and $CA$ is similar. To treat the case of the
snake $BA$ we use the last remark in Section 9.8. The lemma is proved. 
\vskip 3mm
{\bf 10. Constructing a framed local system with arbitrary 
given non-zero coordinates}. 
Let $\Gamma$ be a trivalent graph on $\widehat S$, of the type of $\widehat
S$.  In Section 9.3 we defined the $X$-functions on the moduli space ${\cal
  X}_{PGL_m, \widehat S}$ corresponding to $\Gamma$. 
Recall that they are 
parametrized by the set ${\rm J}_m^{\Gamma}$. 
Below we are going to show how to invert this construction. Precisely, 
let $F$ be a field. Then given  $X_{i} \in F^*$, 
$ i \in {\rm J}_m^{\Gamma}$ we will construct a framed 
$PGL_m$-local system on $\widehat S$ with the given 
$X$-coordinates $\{X_{i}\}$. 

\begin{figure}[ht]
\centerline{\epsfbox{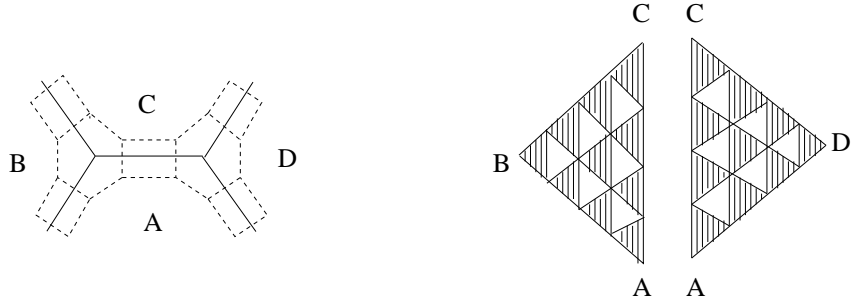}}
\caption{The graph $\Delta$.}
\label{fg21}
\end{figure}

Consider a decomposition of  $\widehat S$ into hexagons  and rectangles 
obtained as follows. Replace the edges of $\Gamma$ by  rectangles and the 
  vertices by hexagons, as on  the left 
picture (Figure \ref{fg21}). Let $\Delta$ be 
the $1$--skeleton of this decomposition. 

Let $G$  a split reductive algebraic group. Then,  given a
framed $G$-local system on $\widehat S$, every vertex $v$ of $\Gamma$ 
provides a configuration of three flags in the vicinity of $v$. 
Therefore there are 
 six canonical pinnings assigned to $v$. These pinnings match the vertices
 of the hexagon of $\Delta$ corresponding to $v$. 

To construct a framed $G(F)$-local system on $\widehat S$ 
we  assign to every oriented edge ${\bf e}$ of $\Delta$ 
an element $M({\bf e}) \in G(F)$. Let 
${\bf p}:= {\bf e}_1{\bf e}_2{\bf e}_3 ... $ be a path on $ \Delta$,  
where  
${\bf e}{\bf f}$ means ${\bf e}$ follows ${\bf f}$. Then
${\bf p}$ 
gives rise to  
$
M({\bf p}):= M({\bf e}_1)M({\bf e}_2)M({\bf e}_3)... $.

\begin{lemma} \label{7.9.03.1}
Suppose that elements $M({\bf e}) \in G(F)$ satisfy
 the following conditions:

i) If ${\bf e}$ does not intersect the graph $\Gamma$ then $M({\bf e}) \in B(F)$.

ii) Let $\overline {\bf e}$ be the edge 
${\bf e}$ equipped with the opposite orientation. 
Then $M({\bf e})M(\overline{\bf  e}) = {\rm Id}$. 

iii) The element assigned to a path around 
a rectangle or hexagon is 
the identity. 

\noindent
Then the elements 
$M({\bf p})$ depend only on the homotopy class of the path ${\bf p}$, 
providing a $G(F)$-local system on $S$. Moreover it has a 
natural framing   on $\widehat S$. 
\end{lemma}

{\bf Proof}. The properties ii) and  iii) guaranty that we get 
a $G(F)$-local system on $S$,  
and  i) provides  a natural framing. The lemma is proved. 
\vskip 3mm

Let us construct such elements $M({\bf e})$. 
The graph $\Delta$ has three types of edges: the 
edges intersecting $\Gamma$, the edges 
of the rectangles which do not intersect $\Gamma$, and the edges 
of the hexagons which do not intersect $\Gamma$. We denote them, respectively,
 by ${\bf a}$,  
  ${\bf b}$ and ${\bf c}$. We will assume that the ${\bf c}$-edges of the 
  hexagons are oriented counterclockwise, i.e. in the 
opposite way to the cyclic order of the edges and hence
the flags in the
vicinity of $v$. The following definition is suggested by Proposition \ref{9.30.03.10aw}
and Lemma \ref{3.27.04.10}. 

a) Set $M({\bf a}):= S$.  

b) Let ${\bf b}$ be an edge of a rectangle 
which goes along an oriented edge ${\bf e}$ of $\Gamma$. Let  
 $x_a = X^{\bf e}_{a, m-2-a}$  
 be the coordinates attached to the edge ${\bf e}$. 
Then  $M({\bf b})$ is given by the right hand side of  (\ref{5.18.03.101}).
  
c) 
Let $\{X_{a,b,c}\}$, where $a+b+c=m-3$ and $a,b,c \geq 0$, be the coordinates 
assigned to a vertex $v$ of $\Gamma$. Then  integers $a,b,c$ match the 
edges of the hexagon surrounding  $v$ which 
do not intersect $\Gamma$. 
Let ${\bf c}$ be such an edge assigned to $a$, counterclockwise oriented. 
Then $M({\bf c})$ is given by the right hand side of (\ref{5.18.03.11c}).

\begin{theorem} \label{5.18.03.10} 
Let $F$ be a field, and we are given  $X_{i} \in F^*$, where 
$ i \in {\rm J}_m^{\Gamma}$. Then 

 a) The elements $M({\bf a}), M({\bf b}), M({\bf c}) \in PGL_m(F)$ 
satisfy the conditions of Lemma \ref{7.9.03.1}. 

\noindent
So they give rise to a framed local system on $\widehat S$. 

b) Let ${\cal L}$ be a framed $PGL_m(F)$--local system on $\widehat S$. 
Then the framed $PGL_m(F)$--local system ${\cal L}(X)$ provided 
by the $X$--coordinates of ${\cal L}$ is isomorphic to ${\cal L}$.
\end{theorem}

{\bf Proof}. Follows immediately from
 Proposition \ref{9.30.03.10aw} and Lemma \ref{3.27.04.10}. 
The theorem is proved. 
\vskip 3mm

 Let $H_m$ (respectively ${\cal V}_m$) 
be the  algebraic torus 
parametrising the internal edge (respectively vertex) $X$--coordinates.  
Then $H_m$ is isomorphic to the Cartan group of $PGL_m$, 
and  
${\cal V}_m = {\Bbb G}_m^{(m-2)(m-1)/2}$. The part a) of Theorem
\ref{5.18.03.10} provides 
 a regular 
 map 
$$
\psi_{\Gamma}: H_m^{\{\mbox{internal edges of $\Gamma$}\}}\times {\cal V}_m^{\{\mbox{vertices of $\Gamma$}\}} = 
{\Bbb G}_m^{{\rm J}^{\Gamma}_m}\lra 
{\cal X}_{PGL_m, \widehat S}. 
$$

\begin{corollary} \label{5.18.03.2}
The map $\psi_{\Gamma}$ is a 
 birational isomorphism, and is an injective regular map. 
\end{corollary}

{\bf Proof}. The part b) of Theorem \ref{5.18.03.10} plus a dimension count 
imply that all  conditions of Lemma \ref{5.18.03.1} are valid for the map 
$\psi_{\Gamma}$.  So it is a 
 birational isomorphism. 
Since $\psi_{\Gamma}$ is regular, it is injective. The corollary is proved.

  \begin{definition} \label{3.27.04.12}
The atlas given by the collection of embeddings $\{\psi_{\Gamma}\}$, 
when $\Gamma$ run through all isotopy classes of marked 
trivalent graphs on $\widehat S$ of type 
$\widehat S$, is called  {\bf special 
atlas} on ${\cal X}_{PGL_m, \widehat S}$. 
\end{definition}

By Lemma \ref{3.30.04.1q} the special atlas is a subatlas of the 
atlas defined in Section 6. This 
 implies that the special atlas is a positive atlas. 
An elementary proof of this claim
is contained in Section 10. 
This positivity claim  plus Theorem  \ref{5.18.03.10} 
gives a complete proof of Theorem \ref{5.6.02.3}.

In the next subsection we investigate the monodromy properties of the
universal $PGL_m$-local system on $S$ with respect to the special positive atlas on 
${\cal X}_{PGL_m, \widehat S}$. 
 \vskip 3mm
{\bf 11. Positive  Laurent property of the monodromy 
of universal $PGL_m$-local system}. 
A {\it special good positive Laurent polynomial on ${\cal X}_{PGL_m, S}$} 
is a rational function which in every coordinate system $\{T_i\}$ 
from the special positive atlas on ${\cal X}_{PGL_m, S}$ 
is a Laurent polynomial with positive integral coefficients. 
We say that a matrix with entries in the field ${\Bbb F}_{PGL_m, S}$ 
is a {\it special good totally positive (integral) Laurent matrix} if 
in every special positive coordinate system, every minor of this matrix is 
a good positive integral Laurent polynomial. There is a similar definition
for the upper/lower triangular matrices. 

\begin{theorem} \label{3.18.04.1q} 
The monodromy of the universal $PGL_m$-local system  ${\cal L}_m$ around a  
non-boundary   loop   on $S$ is conjugated to 
a special good totally positive Laurent matrix. The monodromy around 
a boundary loop is conjugated to an upper/lower triangular  
special good  totally positive Laurent matrix.
\end{theorem}

{\bf Proof}. It is similar to the proof of Theorem \ref{2.1.04.0}. Starting from the
connection on the graph $\Delta$ given by the elements 
$M({\bf a}), M({\bf b}), M({\bf c})$ we define a connection on the
corresponding graph $\Gamma'_T$, shown on Figure \ref{fgo3-11}, as follows. 
Observe that the oriented ${\bf t}$-edges of  $\Gamma'_T$ match the
oriented ${\bf c}$-edges of  $\Delta$. 
Further, 
the ${\bf b}$-edges of $\Delta$ have canonical orientations: the ones 
 compatible with the
clockwise orientation of the holes. So oriented ${\bf e}$-edges of $\Gamma'_T$ match  
the  ${\bf b}$-edges of $\Delta$. 
So, having in mind $M_{AB \to CA} = M_{AC \to CA}M_{AB \to AC}$, we set
$$
M({\bf t}):= SM({\bf c}) = M({\bf a})M({\bf c}); \qquad 
M({\bf e}):= M({\bf b})S = M({\bf b})M({\bf a})
$$
assuming that the ${\bf t}$ and ${\bf e}$-edges of $\Gamma'_T$ 
match the corresponding edges of
$\Delta$. 
Then $M({\bf e})M({\bf t})$ is a lower triangular special totally 
positive integral Laurent
matrix. Indeed, it is obtained by a product of 
the matrices $F_j$, $H_i(x)$ and $M({\bf b})$. Similarly, if ${\bf e}$
follows ${\bf \overline t}$, then  $M({\bf e})M({\bf \overline t}) = 
M({\bf e})M({\bf t})^{-1} = M({\bf b})S M({\bf c})S^{-1}$ 
is an upper triangular special totally 
positive integral Laurent
matrix: see Remark after  Proposition 
\ref{9.30.03.10aw}. The rest  
is as in the proof of Theorem \ref{2.1.04.0}. The theorem is
proved. 

Let  $\rho: PGL_m \to {\rm Aut}(V)$ be a finite dimensional
  representation. Consider the associated local system 
${\cal L}_{\rho}:= {\cal L}_{m}\times_{PGL_m}V$ 
on $S \times {\cal X}_{PGL_m}$. 

\begin{corollary} \label{3.18.04.1} 
a) Let $M_{\alpha}$ be 
the monodromy of the universal local system  around a  loop  $\alpha$. 
Then,  
for any integer $n>0$,  
${\rm tr}(M^n_{\alpha})$ is a special good positive Laurent polynomial 
on  ${\cal X}_{PGL_m, S}$. 

b) For any finite dimensional
  representation $\rho$ of $PGL_m$ the monodromy of the universal 
local system  ${\cal L}_{\rho}$ around any 
loop on $S$ is conjugated to a matrix whose entries are special 
good positive rational Laurent polynomials on 
${\cal X}_{PGL_m, \widehat S}$.

c) The trace of the monodromy of ${\cal L}_{\rho}$ around any 
loop is a special good positive rational Laurent polynomial on 
${\cal X}_{PGL_m, \widehat S}$.
\end{corollary}

{\bf Proof}.  a) Follows immediately from Theorem \ref{3.18.04.1q}. 

b) Take the canonical basis in $V$. Then the elements $E_i$ and $F_i$ 
of $PGL_m$ act by matrices with non-negative rational coefficients, and the
elements $H_i(X)$ act by diagonal matrices whose coefficients are 
monomials in $X^{\pm 1}$.  c) follows from b). The corollary is proved. 

\vskip 3mm
\begin{conjecture} \label{3.18.04.2} 
${\rm tr}(M^n_{\alpha})$ can not be presented as a sum of two non-zero 
elements of ${\Bbb L}_+({\cal X}_{PGL_m, S})$:
$$
{\rm tr}(M^n_{\alpha}) \in {\bf E}({\cal X}_{PGL_m, S}).
$$
\end{conjecture}

\vskip 3mm
{\bf 12. Convex curves in 
$\R{\Bbb P}^{n}$ and positive curves in the flag variety for $PGL_{n+1}(\R)$}. 
A curve $K$ in $\R{\Bbb P}^{n}$ is  {\it convex} if any hyperplane
contains no more then $n$ points of $K$.  
This definition goes back  to Schoenberg \cite{Sch}. It is interesting 
that every convex algebraic curve  is projectively equivalent to the
Veronese curve given by $(x_0:x_1) \lms (x_0^n: x_0^{n-1}x_1: ... : x_1^n)$. 

We start from a useful inductive criteria of positivity of configurations of flags 
in $\R{\Bbb P}^n$. 

\begin{lemma} \label{3.26.04.1}
A configurations of $m$ flags in $\R{\Bbb P}^{2}$ is positive if and only if 
the lines of the flags bound a convex $m$-gon, and the points 
of the flags are at the boundary of this $m$-gon, but not at the vertices, see
Figure \ref{fgo3-15}. 
\end{lemma}

{\bf Proof}. For $m=3$ this is Lemma \ref{8.2}.  
The general case is deduced from it, see 
 Section 2 of \cite{FG3}.

Denote by ${\cal B}_{n+1}(\R)$ the flag variety for  
$PGL_{n+1}(\R)$. There are two natural projections 
 $$
\pi: {\cal B}_{n+1}(\R) \lra \R{\Bbb P}^n, \quad 
\widehat \pi: {\cal B}_{n+1}(\R) \lra \widehat {\R{\Bbb P}}^n.
$$ 

Let $(F_1, ..., F_{m})$  be a generic  configuration of flags 
in $\R{\Bbb P}^n$. Set $x_{k}:= \pi(F_{k})$. We define 
a configuration of flags $(x_k|F_1, ...,  F_{m})$  
in $\R{\Bbb P}^{n-1}$  by projection with the center at 
the point $x_{k}$. Namely, let $F$ be a flag in $\R{\Bbb P}^n$ whose 
  $(n-1)$-dimensional subspace does not contain  $x_{k}$. Then $F$ provides 
 a flag in the projective space of all lines passing 
through $x_{k}$: its $p$-dimensional subspace consists of 
the lines through $x_{k}$ passing through the similar subspace of $F$. This 
way the flags $F_i$, $i \not = k$,  give rise to the flags
$\overline F_i$. Finally, the $p$-dimensional subspace of the 
the flag  $ F_{k}$ provides the $(p-1)$-dimensional subspace of the  flag 
$\overline F_{k}$. Now we set $(x_k|F_1, ...,  F_{m}):= (\overline F_1, ...,
\overline F_{m})$.

\begin{figure}[htb]
\centerline{\epsfbox{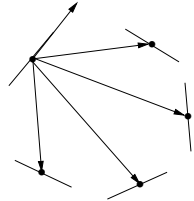}}
\caption{Projecting a configuration of five flags.}
\label{fgo3-16}
\end{figure}

\begin{lemma} \label{3.24.04.3} Let $n \geq 3$. A configuration
  of flags $(F_1, ..., F_m)$ in $\R{\Bbb P}^n$ is positive if and only if 
  all configurations 
$(x_1|F_1, ..., F_m), ..., (x_{m}|F_1, ..., F_m)$ are positive. 
\end{lemma}

{\bf Proof}. Let us prove the implication
 $<=$. Take a convex $m$-gon $P_m$ with vertices labeled by the flags 
$F_1, ..., F_{m}$. Observe that, assuming $n \geq 3$, given a
triangulation of the polygon $P_m$, 
the corresponding canonical coordinates of the configuration $(F_1, ..., F_m)$ 
are defined as follows: we take a triple or quadruple of flags assigned to a
triangle or quadrilateral of the triangulation, project them 
to the flags in $\R{\Bbb P}^2$, and then take the canonical coordinates of the
obtained configuration. It follows from this that 
positivity of the configurations $(x_1|F_1, ..., F_m), ..., (x_{m}|F_1, ...,
F_m)$ 
implies positivity of all the coordinates of the original 
configuration $(F_1, ..., F_m)$ with respect to the given triangulation. 

 Let us prove that if  the configuration $(F_1, ..., F_m)$ is positive, then 
 the one  $(x_{m}|F_1, ..., F_m)$ is also positive. Since the property of 
positivity of
 configurations of flags is cyclically invariant, 
this would imply the implication
 $=>$. 
Consider the triangulation of the polygon $P_m$ by the diagonals from the vertex 
$F_{m}$. Consider the triple 
$(F_{k},  F_{k+1}, F_{m})$ assigned to the vertices 
of one of the triangles. Then, assuming $n \geq 3$,
   the set of the 
coordinates corresponding to the triple $(x_m|F_{k},  F_{k+1}, F_{m})$ is a subset of
  the set of coordinates of the triple 
$(F_{k}, F_{k+1}, F_{m})$. A similar statement is true for the coordinates
assigned to a diagonal 
$(F_{k},  F_{m})$ of the polygon. 
The claim follows immediately from this. 
The
lemma is proved. 
\vskip 3mm
{\bf Remark}. The implication $<=$  of  Lemma \ref{3.24.04.3} is false for $n=2$.
\vskip 3mm

Let $K$ be an oriented continuous convex curve in $\R{\Bbb P}^n$. Let $x_1,
x_2$ be two points of $K$. We say that $x_1 < x_2$ if $x_2$ follows $x_1$ 
according to the
orientation of the curve $K$. Given a collection of points $X = x_1 < x_2 < ... <
x_n$ on $K$, we associate to it a flag 
$$
F_X:= (x_1, [x_1x_2], ..., [x_1 ... x_n])
$$
 where $[x_1 ... x_k]$ denotes the subspace spanned by the points $x_1, ...,
 x_k$. Convexity of $K$ implies that  ${\rm dim}[x_1
 ... x_k] = k-1$, so $F_X$ is indeed a flag. 
More generally, given $n$-tuples of points
\begin{equation} \label{3.24.04.2}
X(1):= \{x_1(1) < ... < x_n(1)\}, \quad ... \quad , X(m):= \{x_1(m) < ... < x_n(m)\}
\end{equation}
 we get the
 flags $F_{X(1)}, ..., F_{X(m)}$. We say that $X(a) < X(b)$ if  $x_i(a) <
 x_j(b)$
for all $1 \leq i, j\leq n$.

\begin{proposition} \label{3.24.04.1} Let $K$ be an oriented continuous convex 
curve in $\R{\Bbb P}^n$. Assume that $X(1) < X(2) < ... < X(m)$ in 
 (\ref{3.24.04.2}). Then the configuration of 
flags $(F_{X(1)}, ..., F_{X(m)})$ is positive. 
\end{proposition}

{\bf Proof}.  Lemma \ref{3.26.04.1} implies  
Proposition \ref{3.24.04.1} for $n=2$. So we may assume that
$n>2$. The projection of a convex curve from a point on
this curve is a convex curve. 
So Lemma \ref{3.24.04.3} implies  Proposition \ref{3.24.04.1} 
by descending induction on $n$. The proposition is proved. 
\vskip 3mm

{\it The osculating curves}. Let $K$ be an oriented continuous convex 
curve in $\R{\Bbb P}^n$. Choose a point $x \in K$, and consider an $n$-tuple of
points $x< x_2 < ... < x_n$ on $K$. It provides a flag. 
We claim that when the points $x_i$ approach the
point $x$ keeping their order, the corresponding family of flags has a limit. 
We prove this claim below. 
The limit is denoted $F_r(x)$ and called the {\it right osculating flag} at
$x$. The map $x \lms F_r(x)$ provides a map $\beta_r: K \to {\cal B}_{n+1}(\R)$,
 which, as one easily sees using the proof of Proposition \ref{3.24.04.1},
 is positive. The curve 
$\widetilde K_r:= \beta_r(K)$ is called the {\it right osculating curve} to
  $K$.
Similarly starting from the points  $x> x_2 > ... > x_n$ we define the {\it left osculating curve} $\widetilde K_l$ to $K$.

To prove the claim, consider a sequence of $(n-1)$-tuples 
of points  $X(k):= x_2(k) < ... < x_n(k)$, $k=1,2,... $, such that 
$x < x_i(k)$ and $x_i(k)$ converges to $x$ as $k\to \infty$. Let us assume that 
$X(a) < X(b)$. Projecting $K$ from  $x$ we get a convex curve
 in $\R{\Bbb P}^{n-1}$, and a sequence of flags 
$F_{X(a)}$ on it. This sequence is positive by Proposition \ref{3.24.04.1}. 
Thus  by Theorem \ref{5.9.03.100} it 
has a limit $\overline F_r(x)$. The limit evidently does not depend
on the choice of sequence $X(k)$. Then  $F_r(x)$ is the unique
flag containing $x$ and projecting to the flag $\overline F_r(x)$.

\begin{theorem} \label{1.30.04.1dwas} 
a) Let $C$ be a positive curve in the flag variety 
${\cal B}_{n+1}(\R)$. Then $\pi(C)$ and $\widehat \pi(C)$ are convex 
curves in $\R{\Bbb P}^{n}$ and $\widehat {\R{\Bbb P}}^n$.

b) Let $K$ be a continuous  convex curve in $\R{\Bbb P}^{n}$. 
 Then  the osculating   curves 
$\widetilde K_l$ and $\widetilde K_r$ are positive 
curves in ${\cal B}_{n+1}(\R)$. 
So if $K$ is a $C^{n-1}$-smooth, then 
$\widetilde K_l = \widetilde K_r =\widetilde K$ is a positive 
curve in ${\cal B}_{n+1}(\R)$.

c) The rule  $K \lms \widetilde K$ gives rise to 
a bijective correspondence 
\begin{equation} \label{3.23.04.10}
\mbox{$C^{n}$-smooth convex curves  
in $\R{\Bbb P}^{n}$} \leftrightarrow \mbox{$C^1$-smooth positive curves in ${\cal B}_{n+1}(\R)$}.
\end{equation}
\end{theorem}

 {\bf Proof}. a) It follows immediately from the following crucial proposition:

\begin{proposition} \label{1.30.04.1da}
Let $(F_1, ..., F_{n+1})$ be a positive configuration of flags 
in  $\R{\Bbb P}^n$. Let $(x_1, ..., x_{n+1})$ be the corresponding
configuration of points in $\R{\Bbb P}^n$, i.e. $\pi(F_i) = x_i$. 
Then these points are not contained in a
hyperplane. 
\end{proposition}

{\bf Proof}. 
Observe that a the flags of a positive configuration of flags are in generic position. 
Lemma \ref{3.24.04.3}  
implies that 
$(x_{n+1}|F_1, ..., F_{n})$ is a positive configuration of flags. 
Further,  the claim of 
Proposition \ref{1.30.04.1da} for the configuration of flags $(x_{n+1}|F_1, ..., F_{n})$
in $\R{\Bbb P}^{n-1}$ evidently implies the one
for the initial flags $F_1, ..., F_{n+1}$. So arguing by the induction on $n\geq 2$ 
we reduce the Proposition to the case $n=2$. Let us prove the Proposition when $n=2$.
\vskip 3mm

Let $V_3$ be a three dimensional vector space. Choose  a 
volume form $\omega$ in $V_3$. For any two vectors $a,b$ 
we define a cross-product $a \times b \in V_3^*$ by
$<a \times b, c>: = \Delta(a,b,c)$. The volume form $\omega$ defines the 
dual volume form in $V_3^*$, so we can define $\Delta(x,y,z)$ for 
any three vectors in $V_3^*$. The following fact is Lemma 5.1 from \cite{G4}. 
\begin{lemma} \label{gz2}
For any $6$ vectors in generic position $a_1,a_2,a_3,b_1,b_2,b_3$ in $V_3$ one has 
$$
\Delta(a_1,a_2,b_2) \cdot \Delta(a_2,a_3,b_3) \cdot \Delta(a_3,a_1,b_1) - 
\Delta(a_1,a_2,b_1) \cdot \Delta(a_2,a_3,b_2) \cdot \Delta(a_3,a_1,b_3)  =
$$
$$
\Delta(a_1,a_2,a_3) \cdot \Delta(a_1 \times b_1,a_2 \times b_2,a_3 \times b_3).
$$
\end{lemma}

{\bf Proof}. Any configuration of six generic vectors in $V_3$ 
is equivalent to the following  one:
$$
\begin{array} {cccccc}
a_1&a_2&a_3&b_1&b_2&b_3\\
-&-&-&-&-&-\\
1&0&0&x_1&y_1&z_1\\
0&1&0&x_2&y_2&z_2\\
0&0&1&x_3&y_3&z_3
\end{array}
$$
Then both the left and the right hand sides are equal to 
$y_3 z_1 x_2 - x_3 y_1 z_2$. The lemma is proved.

Let $(A, B, C)$ 
be a positive configuration of three flags in $\R{\Bbb
  P}^{2}$. The flags $A, B, C$ can be defined by 
pairs of  vectors $(a_1, a_2), (b_1, b_2), (c_1, c_2)$ in $V_3$, so that the flag 
$A$ is determined by the affine flag $(a_1, a_1\wedge a_2)$, and so on. 
Recall the triple ratio (\ref{3.28.04.10}). 
Using Lemma \ref{gz2} we have 
$$
1+r_3^+(A, B, C) =  1 + 
\frac{\Delta(a_1, a_2, b_1)\Delta(b_1, b_2, c_1)
\Delta(c_1, c_2, a_1)}
{\Delta(a_1, a_2, c_1)\Delta(b_1, b_2, a_1)
\Delta(c_1, c_2, b_1)} = 
$$
 $$
\frac{\Delta(a_1, b_1, c_1)\Delta(a_1 \times a_2, b_1 \times b_2, c_1 \times c_2)}
{\Delta(a_1, a_2, c_1)\Delta(b_1, b_2, a_1)
\Delta(c_1, c_2, b_1)}. 
$$
Therefore $\Delta(a_1, b_1, c_1)\not = 0$ if $r_3^+(A, B, C)>0$. 
Proposition \ref{1.30.04.1da}, and hence the part a) of Theorem \ref{1.30.04.1dwas} 
are proved.

b) Observe that projection with the center at a point $x \in K$ maps the
osculating flags on $K$ to the ones on the projected curve. 
So the part b) follows immediately from Lemma \ref{3.24.04.3}. 

c) By Proposition \ref{3.23.04.1} a differentiable positive curve $C$
in ${\cal B}_{n+1}(\R)$  is tangent to the canonical distribution, and hence 
is the osculating curve for the  $n$ times differentiable curve $\pi(C)$ 
in $\R {\Bbb P}^n$. This curve is convex by the part a). 
The converse statement is the part b) of the theorem. The theorem is proved. 



\section{The pair $({\cal X}_{PGL_m, \widehat S}, {\cal A}_{SL_m, \widehat S})$ corresponds to  an orbi-cluster ensemble}
\label{orbi}

In this Section we formulate precisely and prove Theorems \ref{10.29.03.100}
and \ref{10.29.03.100a}. 
 
\vskip 3mm
{\bf 1. Cluster ensembles}. 
We briefly recall some details 
of the cluster ensembles as defined in \cite{FG2}. 
 A cluster ensemble is determined by essentially 
the same combinatorial data, called a {\it seed},  
as the cluster algebras \cite{FZI}. 
A seed $\mathbf i = ({\rm I}, {\rm J}, \varepsilon_{ij}, d_i)$ 
consists of 
a finite set ${\rm I}$,  an integral valued function 
$(\varepsilon_{ij})$ 
on ${\rm I}\times 
{\rm I}$, called a {\it cluster function}, a $\Q_{>0}$-valued {\it symmetrizer function} $d_i$ on ${\rm I}$, such that 
$\widetilde \varepsilon_{ij}:= d_i \varepsilon_{ij}$ is skew-symmetric, 
and a subset ${\rm J} \subset {\rm I}$. The complement ${\rm I} - {\rm J}$ is called the {\it frozen subset} of ${\rm I}$. If $\varepsilon_{ij}$ is skew symmetric, 
like in the case considered below, we set $d_i=1$. 

{\it Mutations}. Given a seed ${\mathbf i}=({\rm I}, {\rm J},\varepsilon, d)$, every non-frozen 
element $k \in {\rm J}$ provides  a mutated in the 
direction $k$ seed $\mu_k({\mathbf i}) = {\mathbf i'} =({\rm I'}, {\rm J'},\varepsilon', d')$ as follows: one has ${\rm I'}:= {\rm I}, {\rm J'}:= {\rm J}, d':=d$ and
\begin{equation} \label{5.11.03.6f}
\varepsilon'_{ij}
= \left\{ \begin{array}{ll} - \varepsilon_{ij} & \mbox{ if $k
      \in \{i,j\}$}
\\ \varepsilon_{ij} + \frac{|\varepsilon_{ik}|\varepsilon_{kj} +
  \varepsilon_{ik}
|\varepsilon_{kj}|}{2} & \mbox{ if $k \not 
      \in \{i,j\}$}
\end{array}\right.
\end{equation}
This procedure is involutive: 
the mutation of  $\varepsilon'_{ij}$ in the direction $k$ is 
the original function $\varepsilon_{ij}$.

Let us mutate 
a given seed in all possible directions, and repeat this infinitely many times. 
Let ${\rm Tr}_{\rm I}$ be a tree 
 defined as follows. 
Its vertices pa\-ra\-me\-te\-rize all seeds obtained from the original one by mutations. 
Two vertices are connected by an edge if the corresponding seeds are related by a mutation. 
This edge inherits a decoration $k$ if the corresponding mutation is in the direction $k$.

A seed ${\mathbf i}$ gives rise to two split algebraic 
tori: 
$$
{\mathcal X}_{\mathbf i} := ({\Bbb G}_m)^{\rm J}, \quad 
{\mathcal A}_{\mathbf i} := ({\Bbb G}_m)^{\rm I}
$$ 
Let $\{X_j\}$ be the natural coordinates on the first torus, and 
$\{A_i\}$ on the second. 
The torus ${\mathcal X}_{\mathbf i}$ is called the {\em seed ${\cal
X}$-torus}, and the torus ${\mathcal A}_{\mathbf i}$ is called the {\em seed ${\cal
A}$-torus}. 
The ${\cal A}$-coordinates parametrized by the complement 
${\rm I} - {\rm J}$ are called the frozen ${\cal A}$-coordinates.
This definition differs slightly  from the one in \cite{FG2}, where 
${\mathcal X}_{\mathbf i} := ({\Bbb G}_m)^I$, since in this paper we do not need 
the frozen ${\cal X}$-coordinates. 

Mutations induce positive 
 rational maps between the corresponding 
 seed ${\cal X}$- and ${\cal A}$-tori, which are  denoted by the same
 symbols $\mu_k$ and defined by the formulas 
\begin{equation} \label{5.11.03.1x}
\mu_k^*X_{i} = \left\{\begin{array}{lll} X_k^{-1}& \mbox{ if } & i=k \\
    X_i(1+X_k)^{-\varepsilon_{ik}} & \mbox{ if } & \varepsilon_{ik}\leq 0 \mbox{ and } i\neq k \\
 X_i(1+X_k^{-1})^{-\varepsilon_{ik}} & \mbox{ if } & \varepsilon_{ik}\geq 0 \mbox{ and } i\neq k
\end{array} \right.
\end{equation}
\begin{equation} \label{5.11.03.1a}
  \mu_k^*A_{k} =  
A_{k}, \qquad \mu_k^*A_{i} = \quad \prod_{j| \varepsilon_{kj} >0} 
A_{j}^{\varepsilon_{kj}} + \prod_{j| \varepsilon_{kj} <0} 
A_{j}^{-\varepsilon_{kj}}; \quad i \not = k
\end{equation}

\vskip 3mm

{\bf 2. A seed for the cluster ensemble 
describing the pair \newline $({\cal A}_{SL_m, \widehat S}, {\cal X}_{PGL_m, \widehat S})$.} 
Choose a marked trivalent graph $\Gamma$ on $S$,
 of type $\widehat S$, and take the corresponding $m$-triangulation 
$T_m = T_m(\Gamma)$ of $S$. Then $T_m$  is  
an oriented graph. It determines a seed ${\mathbf i}$
as follows. Recall that an oriented graph $T$ determines a skew-symmetric  
function $\varepsilon_{pq}(T)$ 
on the set of pairs of its vertices, see  (\ref{11.25.02.1we}). 

\vskip 3mm
i) Let ${\rm I}:= {\rm I}^{\Gamma}_m$ be the set 
of vertices of the graph $T_m(\Gamma)$ minus the punctures of $S$, 
and let ${\rm J}:= {\rm J}^{\Gamma}_m$. 

ii) The cluster function is given by 
$\varepsilon_{pq} = \varepsilon_{pq}(T_m)$.

iii) The coordinates $\Delta_i^{\Gamma}$ are the cluster ensemble coordinates 
$A_i^{\mathbf i}$.

iv)  The coordinates $X_j^{\Gamma}$ are the cluster ensemble coordinates $X_j^{\mathbf i}$.

v) Formula (\ref{11.5.03.2}) shows that the function 
$X_j^{\Gamma}$ and $\Delta_i^{\Gamma}$ are related under the projection $p$ 
the same way as the cluster ensemble coordinates $X_j^{\mathbf i}$ and $A_i^{\mathbf i}$, 
see \cite{FG2}. 
\vskip 3mm

An interpretation of a flip in the cluster language is a non 
trivial problem:   
we are going to prove that a flip can be presented as a composition of mutations. 

\begin{figure}[ht]
\centerline{\epsfbox{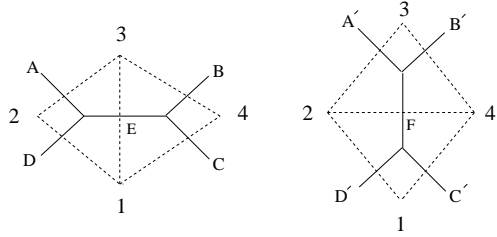}}
\caption{Describing a flip.}
\label{fg14a}
\end{figure}

\vskip 3mm
{\it Describing a flip and its action on the cluster function}.  
Let $\Gamma'$ be a graph obtained from  $\Gamma$ by a
 flip at an edge $E$. Denote by $F$ the new edge of $\Gamma'$. 
Let $\{A, B, C, D\}$ (resp. $\{A', B', C', D'\}$) be the set of edges of  $\Gamma$ 
(resp. $\Gamma'$) adjacent to $E$ (resp. $F$). Then there is a canonical bijection 
$\{A, B, C, D\} \lra \{A', B', C', D'\}$ such that $X \lms X'$. 
Let us number by $1,2,3,4$ the pieces of faces in the vicinity of 
the edge $E$. They match the  ones in the vicinity of  $F$. 
We  assume that their order is compatible with the 
cyclic order provided by the ribbon structure of $\Gamma$, see Figure  
(\ref{fg14a}). Denote by $T_m'$ the $m$-triangulation corresponding to $\Gamma'$.

The claim that  a flip can be obtained as 
 a composition of mutations means 
the following: 

i) The composition of mutations transforms the function $\varepsilon_{pq}(T_m)$
 to  $\varepsilon_{pq}(T'_m)$

ii) The transformation rules for the  ${\cal A}$- and ${\cal X}$-coordinates under a flip coincide 
with the ones provided by the sequence of mutations. 

\vskip 3mm
{\bf 3.  A sequence of mutations 
representing a flip at an edge $E$}. We present it as a composition of $m-1$ steps. 
Consider the diamond with the vertices $1,2,3,4$ on Figure \ref{fg14a}. 
One can inscribe into the $m$-triangulation of this  diamond 
rectangles of the size $i \times (m-i)$, where $i=1, ..., m-1$, centered at 
the center of the diamond.   
Here $i$ (resp. $(m-i)$) is the length in direction of the edge $F$ (resp. $E$). 
We subdivide this rectangle into $i  (m-i)$ equal squares.

{\it Step $i$}.  Do mutations in the direction of the elements of the set
${\rm I}
={\rm I}^{\Gamma}_m$ 
corresponding to the  centers of  little squares of the $i \times (m-i)$ rectangle. 

For example at the step $1$ (resp. $m-1$) we do mutations 
corresponding to the $m-1$ vertices located inside of the vertical
(resp. horizontal) diagonal 
of the diamond.

\begin{proposition} \label{12.8.02.20} The described sequence of mutations 
transforms the cluster function $\varepsilon_{pq}(T_m)$ to the one 
$\varepsilon_{pq}(T'_m)$.  
\end{proposition}

\begin{figure}[ht]
\centerline{\epsfbox{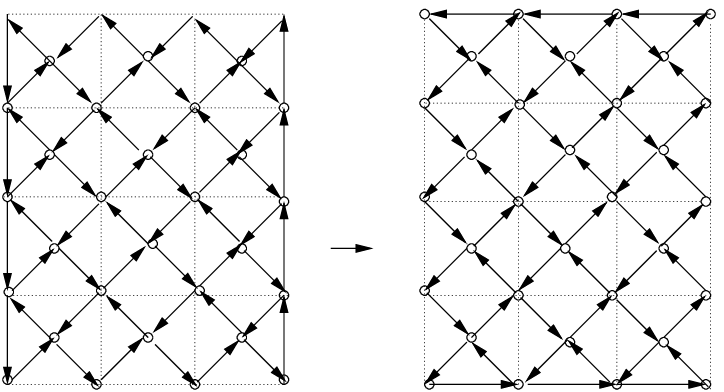}}
\caption{}
\label{fg44}
\end{figure}

{\bf Proof}. Consider an $m\times n$ rectangle where $m, n \geq 1$ 
are integers. Subdivide it into $mn$ equal squares,  
 mark the vertices and centers of the squares. 
Connect the marked points by arrows as on Figure \ref{fg44}. 
Denote by $R_{m,n}$ the obtained oriented graph. 
Observe that the right picture on Figure \ref{fg44} 
is obtained from the left one by inverting directions 
of the arrows at each square center, and replacing the arrows at the 
sides of the rectangle by the ones  at the top and the bottom, 
keeping the counterclockwise orientation.

\begin{lemma} \label{11.4.03.1} 
Consider the cluster function $\varepsilon_{pq}$ corresponding to the 
oriented graph on the left of Figure \ref{fg44}. Let us perform a sequence 
of mutations at all centers of the $mn$ squares, in any order. Then the 
obtained cluster function $\varepsilon'_{pq}$ is described by the picture 
on the right of Figure \ref{fg44}.
\end{lemma}

{\bf Proof}. Follows immediately from the definitions. 

\begin{figure}[ht]
\center
\hspace{8mm}\begin{minipage}{5.5cm}
\centerline{\epsfbox{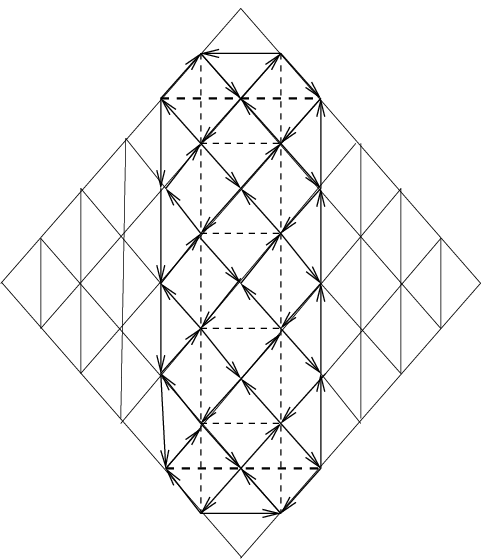}}
\end{minipage}
\begin{minipage}{5.5cm}
\centerline{\epsfbox{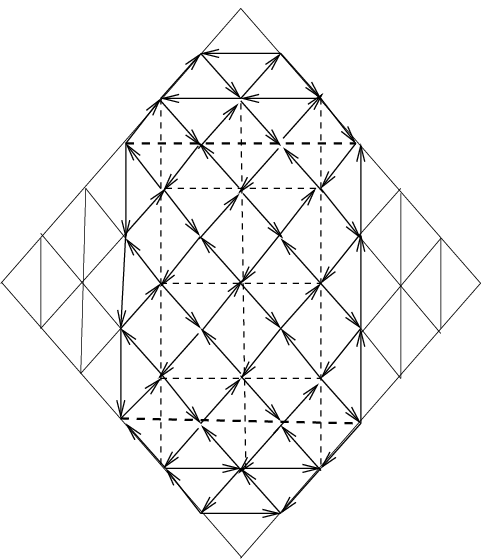}}
\end{minipage}
\vspace{3mm}
\caption{}
\label{fg49}
\end{figure}

\vskip 3mm
{\bf Example}. The two pictures on Figure \ref{fg49}  
illustrate the case $m=6$. We show a part of the mutated graph after the steps 1 and 2. 
On each picture we show the distinguished 
rectangles corresponding to this and previous steps.
\vskip 3mm
Contemplation will convince the reader that after $m-1$ steps we get 
the graph of the function $\varepsilon_{pq}(T_m')$. The proposition is proved. 

\vskip 3mm
{\bf 4. Hypersimplices and a decomposition of a simplex}. Recall \![GGL] that 
a hypersimplex $\Delta^{k,l}$, where $k+l=n-1$, is a section of 
the $n+1$-dimensional cube $0 \leq x_i \leq 1$ by the hyperplane 
$\sum x_i = l+1$. The hypersimplex $\Delta^{k,l}$ can be 
obtained as the convex hull of the centers of $k$-dimensional faces 
of an $(k+l+1)$--dimensional simplex. 
\vskip 3mm
{\bf Example}.  $\Delta^{1,1}$ is an octahedron. 
It is the convex hull of the centers of the edges of a tetrahedron. 
\vskip 3mm
Let $m$ be a positive integer. 
Consider the simplex
\begin{equation} \label{11.16.02.1}
\Delta^n_{(m)}:= \{(x_0, ..., x_n) \in \R^n\quad |\quad  \sum_{i=0}^n x_i =m, \quad x_i \geq 0\}
\end{equation}

\begin{lemma} \label{5.10.02.2} 
The hyperplanes $x_i =k$ where $k$ is an integer cut the simplex (\ref{11.16.02.1}) 
into a union of hypersimplices. 
\end{lemma} 

{\bf Proof}. Consider a decomposition of $\R^n$ into the unit cubes 
with the faces $x_i =k$ where $k$ is an integer. The hyperplane 
$\sum x_i =m$ intersect each of these cubes either by an empty set, or by a 
hypersimplex. The lemma is proved. 
\vskip 3mm
We call it the hypersimplicial decomposition of the simplex 
$\Delta^n_{(m)}$. 
\vskip 3mm
{\bf Example}. When $n=2$ we get an $m$--triangulation of a triangle defined in 
Section 1.13.
\vskip 3mm

\begin{figure}[ht]
\centerline{\epsfbox{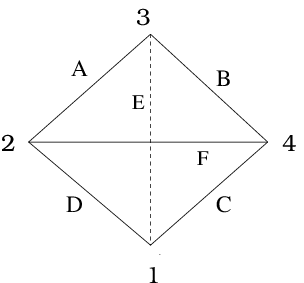}}
\caption{}
\label{fg14b}
\end{figure}

Consider the three dimensional simplex $\Delta_{(m)}^3(E)$ given by 
\begin{equation} \label{11.5.03.31}
x_1+x_2+x_3+x_4 =m, \quad x_i\geq 0, 
\end{equation} 
whose vertices are labeled by $1, 2, 3, 4$, as on Figure \ref{fg14b}. 
Its  edges are dual to the $A, B, C, D, E, F$ edges, see Figure \ref{fg14a}, and 
get the same labels. 
The faces match  the vertices of the 
edges $E$ and $F$. The two faces matching the vertices of $E$ are 
called the $e$--faces, and the other two are called the $f$--faces.

The hypersimplicial decomposition of $\Delta_{(m)}^3(E)$ 
induces the $m$--tri\-an\-gu\-la\-tion of its faces. 
The faces of $\Delta_{(m)}^3(E)$ are identified with the corresponding triangles 
of $T$ or $T'$, so that the 
$m$--triangulations of the faces match  the ones of the 
triangles. Thus to make a flip we glue the simplex $\Delta_{(m)}^3(E)$ 
so that its $e$--faces match the corresponding two triangles of $T$. 
Then the  $f$--faces of $\Delta_{(m)}^3(E)$ provide the two 
triangles of $T'$. 

\vskip 3mm
 {\bf 5. The action of a flip on the $A$-coordinates}.  
Let us attach to the integral points inside of the simplex (\ref{11.5.03.31}) 
functions on the configurations space of four affine flags in an 
$m$-dimensional vector space $V_m$ with a volume form $\omega$. 
Given an integral solution 
$(a,b,c,d)$ of (\ref{11.5.03.31}) and a $4$-tuple of affine flags $(X, Y, Z, T)$
we set
$$
\Delta_{a,b,c,d}(X,Y,Z,T):= \Delta_{\omega}(x_{(a)} \wedge y_{(b)} \wedge 
z_{(c)} \wedge t_{(d)})
$$
where $x_{(a)} = x_1 \wedge ... \wedge x_{a}$ and so on. 
If one of the numbers $a,b,c,d$ is zero we get the 
functions defined in Section 9.2. The Pl{\"u}cker relations allow 
to derive some identities between these functions which we are going 
to describe now. 

\noindent \qquad The hypersimplicial decomposition of the three dimensional 
simplex \!(\ref{11.5.03.31}) is a decomposition into tetrahedrons of type $\Delta^{2,0}$, octahedrons $\Delta^{1,1}$, and tetrahedrons 
of type $\Delta^{0,2}$. The octahedrons play an important role in our story. 
They are parametrized by the solutions 
$(a_1, a_2, 
a_3, a_4)$ 
of the equation
\begin{equation} \label{11.22.02.31}
a_1 + a_2 +  
a_3 + a_4 =m-2, \quad a_i \in \Z, a_i \geq 0.
\end{equation} 
The collection of vertices of these octahedrons coincides with the set of 
the non negative integral solutions 
of the equation (\ref{11.5.03.31}). The vertices of the $(a_1, a_2, 
a_3, a_4)$--octahedron have the  coordinates
\begin{equation} \label{uks}
(a_1, a_2, 
a_3, a_4) + (\varepsilon_1, \varepsilon_2, 
\varepsilon_3, \varepsilon_4), \quad \varepsilon_i \in \{0,1\}, \quad 
\sum_{i=1}^4 \varepsilon_i = 2.
\end{equation} 
We will use a shorthand 
$
\Delta_{b_1, b_2, b_3, b_4}:= \Delta_{b_1, b_2, b_3, b_4}(X, Y, Z, T)
$. 
\begin{lemma} \label{11.5.03.32}
For any  solution $\overline a = (a_1, a_2, a_3, a_4)$ of 
(\ref{11.22.02.31})
 one has 
\begin{equation} \label{11.5.03.35}
\Delta_{\overline a + (1,1,0,0)}\Delta_{\overline a + (0,0, 1,1)} 
+ \Delta_{\overline a + (1,0,0,1)}\Delta_{\overline a + (0, 1,1, 0)} =
\Delta_{\overline a + (1,0,1,0)}\Delta_{\overline a + (0,1,0,1)}. 
\end{equation} 
\end{lemma}

{\bf Proof}. Let $X = (x_1, x_1 \wedge x_2, x_1\wedge x_2 \wedge x_3, ... )$ 
be an affine flag in $V_m$. Let 
$X_1 \subset X_2 \subset ... \subset X_m = V_m$ be the corresponding flag. 
The affine flag $X$ induces an affine flag in each of the quotients $V_m/X_i$. 
The four affine flags $(X, Y, Z, T)$  in generic position 
in  $V_m$ induce four affine flags 
in the quotient 
\begin{equation} \label{11.22.02.35}
\frac{V_m}{X_{a_1} \oplus  Y_{a_2} \oplus  Z_{a_3} \oplus  T_{a_4}}.
\end{equation} 
Since $a_1+a_2+a_3+a_4 = m-2$, this quotient is two dimensional, and 
the induced affine flags provide a  configuration of vectors 
$(\overline x_{a_1+1}, ... , \overline t_{a_4+1})$. 
We define a volume form $\omega'$ in the space (\ref{11.22.02.35}) 
by 
\begin{equation} \label{11.5.03.36}
\omega'(v_1 \wedge v_2):= \Delta(x_{(a_1)}
\wedge y_{(a_2)}\wedge z_{(a_3)}\wedge t_{(a_4)}
\wedge v_1 \wedge v_2).
\end{equation}
Set $(v_1, v_2, v_3, v_4):= 
(\overline x_{a_1+1}, ... , \overline t_{a_4+1})$ 
in the two dimensional space (\ref{11.22.02.35}).  One easily checks
that multiplying  (\ref{11.5.03.35}) by $(-1)^{a_2+a_4}$
we get the Pl{\"u}cker relation
$$
\omega'(v_1\wedge v_2) \omega'(v_3\wedge v_4) +
 \omega'(v_1\wedge v_4) \omega'(v_2\wedge v_3) = 
\omega'(v_1\wedge v_3) \omega'(v_2\wedge v_4). 
$$
The lemma is proved. 
\vskip 3mm

\begin{figure}[ht]
\centerline{\epsfbox{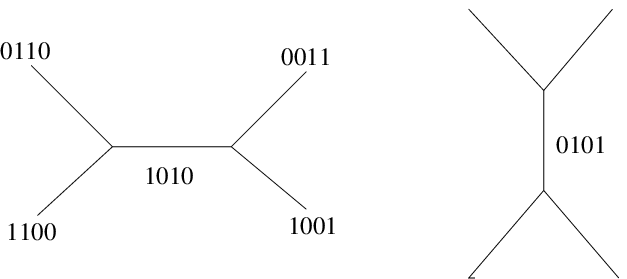}}
\caption{}
\label{fg31}
\end{figure}

One can interpret it as follows. 
The six vectors $(\varepsilon_1, \varepsilon_2, 
\varepsilon_3, \varepsilon_4)$ from (\ref{uks}) 
match the edges of the graphs $\Gamma$ and $\Gamma'$  as shown on 
Figure \ref{fg31}. 
For instance $(1100)$ matches the edge $D$ on the 
Figure (\ref{fg14a}) and so on. 
Set $\Delta_{\overline a, A}:= \Delta_{\overline a + (1,1,0,0)}$ and so on.
Then we have 
\begin{equation} \label{11.5.03.45}
\Delta_{\overline a, A}\Delta_{\overline a, C} 
+ \Delta_{\overline a, B}\Delta_{\overline a, D} =
\Delta_{\overline a, E}\Delta_{\overline a, F}. 
\end{equation}

The sequence of mutations defined in Section 10.2 can be described 
geometrically as follows. The octahedrons of the hypersimplicial decomposition
of $\Delta^3_{(m)}(E)$ 
form $m-1$ horizontal 
layers, if we assume that the edges $E$ and $F$ parallel to 
a horizontal plane, and on the $i$-th step we add the octahedrons on the $i$-th
level. 

Recall that a mutation in the direction $k$ 
changes just one $A$-coordinate, the one parametrised by $k$, 
via the exchange relation given by (\ref{5.11.03.1a}). 
In our situation the mutation corresponding to an octahedron 
replaces $\Delta_{\overline a, E}$ by $\Delta_{\overline a, F}$, so the
exchange relation is given by (\ref{11.5.03.45}). Presenting a flip as a composition of mutations, we 
describe how the ${\cal A}$-coordinates change under a flip.

\vskip 3mm
{\bf 6. The action of a flip on the $X$-coordinates}.  Recall the canonical 
projection $p: {\cal A}_{SL_n, \widehat S} \lra {\cal X}_{PGL_n, \widehat S}$. 
Observe that in the special case when $\widehat S$ is a disc with $n$ marked points 
the canonical projection $p: {\rm Conf}_n({\cal A}) \to {\rm Conf}_n({\cal B})$ is 
a map onto. Thus to prove any identity between the functions $X^{\Gamma}_i$
on ${\rm Conf}_n({\cal B})$ it is sufficient to prove them 
for the pull backs $p^*X^{\Gamma}_i$. 
We already proved that the 
$A$-coordinates and the cluster functions transform the same way by a flip. We have 
the cluster type formula (\ref{11.5.03.2}) for $p^*X^{\Gamma}_i$. So the 
$X$-coordinates also transform by a flip the same way as 
by the corresponding composition of mutations. 

{\it The functions corresponding to the octahedrons of 
hypersimplicial decomposition}. The cross--ratio of the vectors $
\overline x_{a_1+1}, ... , \overline t_{a_4+1}$ provides a function on configurations of $4$--tuples of flags in $V_m$:
\begin{equation} \label{11.22.02.33}
R_{a_1, a_2, a_3, a_4}(X, Y, Z, T) := r^+(\overline x_{a_1+1}, \overline y_{a_2+1}, \overline z_{a_3+1}, \overline t_{a_4+1}).
\end{equation} 
These functions correspond to the octahedrons of the hypersimplicial decomposition of the simplex $\Delta^3_{(m)}(E)$. They are 
determined by the part of the graph $\Gamma$ in the vicinity of $E$. 
A flip $\Gamma \to \Gamma'$ at the edge $E$ changes these functions to their inverses. 

\vskip 3mm
{\bf Proof of Theorem \ref{10.29.03.100}}. 
Since the Farey triangulation is not special 
in the sense of Section 3.8, the results obtained above 
imply the proof of Theorem
\ref{10.29.03.100}a). To get Theorem \ref{10.29.03.100}b) we use in addition 
Lemma \ref{5.16.04.1}. 
\vskip 3mm
To formulate precisely and prove Theorem 
\ref{10.29.03.100a}  it remains to consider the 
special triangulations. We do it below.

\vskip 3mm

 {\bf 7. Treatment of the special graphs}. Recall the special trivalent ribbon 
graphs defined
 in Section 3.8. By definition, they contain a part given by 
 either an eye or a virus. 
There are  well
defined ${\cal X}$- and ${\cal A}$-coordinates assigned to edges and 
triangles of an eye. For the
${\cal X}$-coordinates assigned to an edge $E$ of an eye this deserves a comment.   
Namely, let $Q_E$ be the rectangle having $E$ as an edge, see the middle graph 
on Figure \ref{fgo-100}. The two vertices
opposing the edge $E$ are glued. However the four flags assigned to the
vertices of $Q_E$ are generic: the flags assigned to the two glued 
vertices differ by
a monodromy transformation. Thus the corresponding 
${\cal X}$-coordinates on the edge
$E$ are well
defined. Contrary to this,  to define the ${\cal X}$-coordinates 
corresponding to the loop of a
virus we have to go to a $2:1$ covering of $S$ pictured on 
Figure \ref{fgoeye} and described below. So it is no 
surprise that we face problems comparing the effects at this loop of cluster and geometric
transformations, and that these problems are resolved by going
to a cover. 

\begin{figure}[ht]
\centerline{\epsfbox{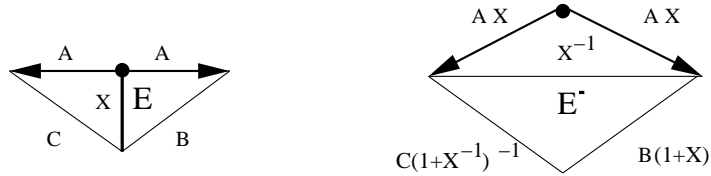}}
\caption{A virus  is obtained from an eye by the flip at the edge $E$.}
\label{fgo-102}
\end{figure}

\begin{figure}[ht]
\centerline{\epsfbox{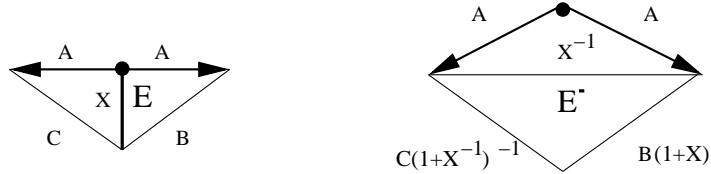}}
\caption{Cluster transformation of coordinates for the flip at $E$ making 
a virus  from an eye.}
\label{fgo-103}
\end{figure}

Let us consider first the case $G = PGL_2$. Figure \ref{fgo-102} 
shows how we obtain a virus by 
making a flip on an eye at the edge $E$, and 
tells us how the coordinates on the moduli space 
${\cal X}_{PGL_2, \widehat S}$ are transformed by this flip. 
Figure  \ref{fgo-103} tells us the effect of the 
cluster transformation, mutation,  for the cluster ensemble corresponding to 
the same moduli space, under the same flip. 
Observe that the only difference is for the 
${\cal X}$-coordinate $A$ assigned to the edge 
obtained by gluing the two edges shown by arrows on the Figure.  
In the geometric case the transformation is $A \lms AX$, while in 
the  cluster case the coordinate $A$ does not change. Since ${\cal
  X}_{PGL_2, \widehat S}$ 
is embedded into ${\cal X}_{PGL_m, \widehat S}$, the 
cluster rule for the transformation of the ${\cal X}$-coordinates 
for $G = PGL_2$ also 
does not agree with the geometric one. There is a similar problem  
for the ${\cal A}$-coordinates. 
Below we are going to show how to 
resolve this problem. 

Let $v$ be the vertex of triangulation shown by the 
fat point on Figure \ref{fgo-102}. 
Let us construct a $2:1$ cover $\widetilde S \to S$ ramified at the 
single point $v$ of the compactification of $S$. Then there is an embedding 
$
i: {\cal X}_{G, S} \hra {\cal X}_{G, \widetilde S} 
$. 
Its image is identified with the set of the stable points of the 
involution $\sigma$ of the covering. 

\begin{figure}[ht]
\centerline{\epsfbox{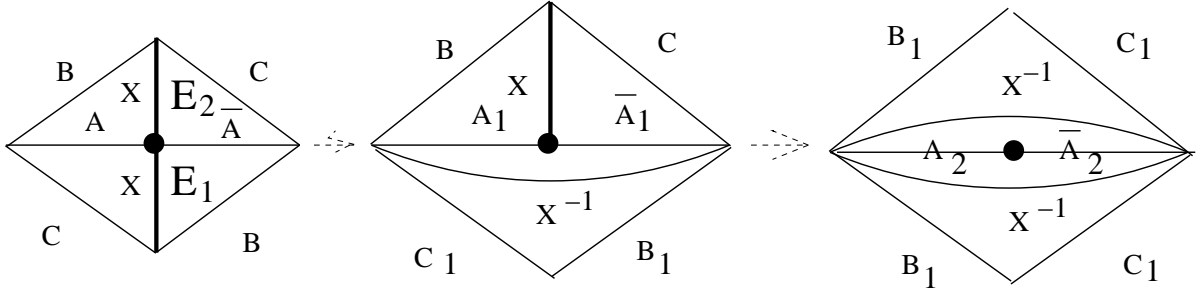}}
\caption{A flip on an eye is a composition of two mutations on a $2:1$ cover.}
\label{fgo-104}
\end{figure}

The coordinate system on either 
${\cal A}$- or ${\cal X}$-space 
corresponding to a trivalent graph $\Gamma$ on $S$ is 
described by a torus $T^*_{S, \Gamma}$, where $*$ 
stays either for ${\cal A}$ or for ${\cal X}$.  
A flip $\Gamma \to \Gamma'$ at an edge $E$ provides a birational automorphism 
$\varphi^*_E: T^*_{S, \Gamma} \to T^*_{S, \Gamma'}$,  
The coordinate system on $\widetilde S$ corresponding to the lifted graph 
$\widetilde \Gamma$ on $\widetilde S$ is described by the torus 
$T^*_{\widetilde S, \widetilde \Gamma}$. One has 
\begin{equation} \label{5.14.04.2} 
T^*_{S, \Gamma} = (T^*_{\widetilde S, \widetilde \Gamma})^{\sigma = {\rm Id}}.
\end{equation} 
Lifting the trivalent graph on $S$ to the cover $\widetilde S$, we get the 
graph on the left of Figure \ref{fgo-104}. The involution $\sigma$ 
interchanges the edges marked by the same letter, and acts as the 
central symmetry with respect to the $4$-valent vertex. We will make the flip at 
the edge $E_1$, 
followed by the flip at the edge $E_2$. These two flips are 
illustrated on Figure \ref{fgo-104}. 

We claim that, using the identification (\ref{5.14.04.2}),
the  geometric coordinate transformation corresponding to the flip at 
the edge $E$ on Figure \ref{fgo-102} equals to  
the  composition of two flips on the cover $\widetilde S$. 
Each of these flips is a product of mutations for the cluster ensemble 
$({\cal X}_{PSL_m, \widetilde S}, {\cal A}_{PSL_m, \widetilde S})$.

\begin{proposition} \label{5.11.04.1}
i) The composition of the flips at the edges $E_1$ and $E_2$, restricted to the 
subtorus (\ref{5.14.04.2}), is equivalent to the flip at the edge $E$ on $S$,
i.e. the following diagram is commutative:
$$
\begin{array}{ccc}
T^*_{S, \Gamma} & \stackrel{\varphi^*_E}{\lra} & T^*_{S, \Gamma'}\\
\downarrow = & & \downarrow = \\
(T^*_{\widetilde S, \widetilde \Gamma})^{\sigma = 
{\rm Id}}&\stackrel{\varphi^*_{E_2}\circ 
\varphi^*_{E_1}}{\lra} 
& (T^*_{\widetilde S, \widetilde \Gamma'})^{\sigma = {\rm Id}}
\end{array}
$$
The flips at $E_1$ and $E_2$ commute. 
 
 ii) Each of the flips at the edges $E_1$ and $E_2$ is a product of 
mutations in the cluster ensemble defined using the trivalent graph
on $\widetilde S$ (partly shown on 
 Figure \ref{fgo-104}) and the function $\varepsilon_{ij}$ related to it. 

\end{proposition} 

This proposition allows to deduce many 
 properties of the coordinate transformations of the 
flips from  an eye to a virus and back 
from the similar properties of 
mutations on the cover $\widetilde S$.

{\bf Proof}. i) This boils down to the very definitions. 

ii)  
It is instructive to consider first the 
$PGL_2$ case. The transformation rule of ${\cal X}$-coordinates under a flip in this case is 
given in (\ref{10.30.03.23}) below. So 
the coordinate 
transformations for the ${\cal X}$-coordinates  are as follows: 
$$
B_1 = B(1+X), \quad C_1 = C(1+X^{-1})^{-1}, \quad \overline A_1 = A(1+X^{-1})^{-1}, \quad 
A_1 = A(1+X)
$$
$$
\overline A_2 = A(1+X)(1+X^{-1})^{-1} = AX, \quad 
A_2 = A(1+X^{-1})^{-1}(1+X)= AX.
$$
Thus the composition of two mutations on Figure \ref{fgo-104} 
commutes with the automorphism $\sigma$ of 
the covering, and thus corresponds to a transformation on $S$. 
The obtained transformation 
is the geometric transformation of coordinates corresponding 
to the flip at the edge $E$, see Figure \ref{fgo-102}. 

Now let us address the general case, when $G=PGL_m/SL_m$. 
We say that an edge $E$ is {\it generic} if for the quadrilateral $Q_E$ 
of the triangulation having $E$ as the diagonal both  maps 
${\cal A}_{G, \widehat S} \to {\cal A}_{G, Q_E}$ and 
${\cal X}_{G, \widehat S} \to {\cal X}_{G, Q_E}$ provided 
by the restriction to $Q_E$ are surjective at the generic point. 
For the ${\cal A}$-space a 
 flip at a generic edge is a composition of mutations. 
Indeed, the ${\cal A}$-coordinates on a quadrilateral are defined using the
 four affine flags assigned to the vertices of the quadrilateral. 
The quadrilaterals corresponding to the edges $E_1$ 
on the left picture and $E_2$ on the middle picture of Figure \ref{fgo-104} are generic.  
This gives the claim ii) for the ${\cal A}$-coordinates. 

\begin{figure}[ht]
\centerline{\epsfbox{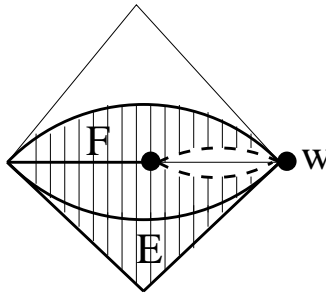}}
\caption{The two vertices of the shaded pentagon $P_{E,F}$ are identified at
  the vertex $w$.}
\label{fgo-105}
\end{figure}

The ${\cal X}$-coordinates in the internal part of the quadrilateral $Q_E$ 
also depend only on the four flags assigned to the vertices of $Q_E$. 
However the ${\cal X}$-coordinates at 
an external edge $F$ depend also on the fifth flag, and so we have to take 
into account the pentagon $P_{E, F}$ containing as the internal diagonals the edge $E$ and the 
external edge $F$. The natural projection  ${\cal X}_{G, \widehat S} 
\to {\cal X}_{G, P_{E, F}}$ provided 
by the restriction to the pentagon $P_{E, F}$ is
 surjective at the generic point. This is clear for the edge $E_1$.  
To check this for the edge $E_2$ observe that  although the pentagon $P_{E,F}$ 
has two vertices identified (at the vertex $w$ on Figure \ref{fgo-105}), 
the flags at these vertices differ
 by the  monodromy around the puncture $v$, and thus they are generic. 
The claim ii) for the ${\cal X}$-coordinates 
follows from this. The proposition, and hence 
Theorem \ref{10.29.03.100a} are proved. 

\vskip 3mm
{\bf 8. Integral points of higher Teichm\"uller spaces}. The following pro\-position
is due to Grothendieck. We will prove it in a slightly different way. 
\begin{proposition} \label{4.2.04.1} 
There is canonical bijection between the ribbon graphs 
and finite index torsion free subgroups of $SL_2(\Z)$ modulo conjugation. 
\end{proposition}

{\bf Construction}. Recall that ribbon graphs encode
pairs $(S, T)$, where $S$ is a surface with punctures, and $T$ an ideal 
triangulation of $S$,
considered modulo diffeomorphisms of surfaces. A torsion free finite
index subgroup $\Delta \subset SL_2(\Z)$ provides a surface 
$S_\Delta := {\cal H}/\Delta$ with triangulation given by the image of the
Farey triangulation of ${\cal H}$. Vice versa, 
given a pair $(S, T)$ we consider the point of the Teichm\"uller space ${\cal X}_{PGL_2, S}^+$ 
whose coordinates with respect to the triangulation $T$ are equal to $1$ 
at every edge of $T$. According to, say, Section 9.11, the monodromy representation 
of such a local system provides a subgroup $\Delta_{(S, T)}$ of $PSL_2(\Z)$,
well defined up to a conjugation. One checks that 
the   triangulated surface corresponding to $\Delta_{(S, T)}$ 
is isomorphic to the pair $(S, T)$. Indeed,  to compute the 
coordinate assigned to an edge $E$ of the triangulation $T$, we take 
the two triangles sharing $E$ to the universal cover, which is the 
hyperbolic plane with the Farey triangulation.  There we get 
the two triangles $(\infty, -1, 0)$ and $(\infty, 0, 1)$. So 
the coordinate of $E$ is $r^+(\infty, -1, 0, 1) =1$.  
The proposition is proved.

\begin{definition} \label{integ}
A point of the Teichm\"uller space ${\cal X}^+_{G, S}$ is integral if the 
corresponding monodromy representation can be represented by a homomorphism 
$\pi_1(S) \hra G(\Z)$. 
\end{definition}

So the pairs $(S, T)$ considered modulo ${\rm Diff}_0(S)$ 
describe the integral points of the Teichm\"uller space ${\cal X}^+_{PGL_2, S}$. 
Here is a hypothetical generalization. 

Consider a cluster coordinate system on ${\cal X}_{PGL_m, S}$, and take 
a faithful representation $\pi_1(S) \hra PGL_m(\Q)$ 
corresponding to the point with 
all cluster ${\cal X}$-coordinates equal to $1$. 

\begin{question} \label{4.2.04.1f}
Is it true that the above construction 
provides  an inclusion $\pi_1(S) \hra PGL_m(\Z)$, and every 
integral point of the Teichm\"uller space ${\cal X}^+_{PGL_m, S}$ appears this way? 
\end{question}
\vskip 3mm
{\bf Example}. The 
cluster coordinate system assigned 
to an ideal triangulation of $S$ provides a map
given as a composition $\pi_1(S) \hra PSL_2(\Z) \hra PSL_m(\Z)$, where the 
first map comes from Proposition \ref{4.2.04.1}, and the 
second is given by the irreducible
$m$-dimensional representation of $SL_2$. 
\vskip 3mm

\section{The classical Teichm{\"u}ller spaces}
\label{classical}

\vskip 3mm
{\bf 1. Proof of Theorem \ref{3.27.02.1}a)}. 
Let $p$ be a point of  ${\cal X}^+_{PGL_2, S}$. Let us assign
 to $p$ the corresponding representation $\rho_p: \pi_1(S) \to PSL_2(\R)$.
By Theorem \ref{7.8.03.2} we get a discrete faithful representation 
of $\pi_1(S)$, and hence a point of the classical Teichm{\"u}ller space. 
A framed structure for $\rho_p$ is equivalent to a choice of orientations
 of those holes in $S$ where 
 the monodromy $M_{\rm hole}$ around the hole is not unipotent. Indeed, a framed structure near such 
a hole is given by a choice of a 
point on $\R P^1$ invariant under $M_{\rm hole}$. 
To see this, we identify $\R P^1$ with the absolute of ${\cal H}$. Then the 
$M_{\rm hole}$-invariant points are the ends of the unique $M_{\rm hole}$-invariant
 geodesic covering the boundary geodesic on the surface ${\cal H}/
\rho_p(\pi_1(S))$. (An element $M_{\rm hole} \in PSL_2(\R)$ is defined up to a 
conjugation by the image of $\pi_1(S)$). So choosing one of them we 
orient the geodesic towards the chosen point. So we get an inclusion 
$$
A: {\cal X}^+_{PGL_2, S} \hra {\cal T}^+_S.
$$
Similarly 
there is a tautological inclusion  $B$ left inverse to  $A$:
$$
B: {\cal T}^+_S \hra {\cal X}_{PGL_2, S}(\R), \qquad B \circ A = {\rm Id}.
$$
It remains to show that $A$ is surjective or, equivalently, 
the image of the map $B$ lies in 
${\cal X}_{PGL_2, S}(\R_{>0})$. So we have to show that 
our coordinates on  ${\cal X}_{PGL_2, S}(\R)$ take positive values 
on the image of  $B$.

\begin{figure}[ht]
\centerline{\epsfxsize29\baselineskip\epsfbox{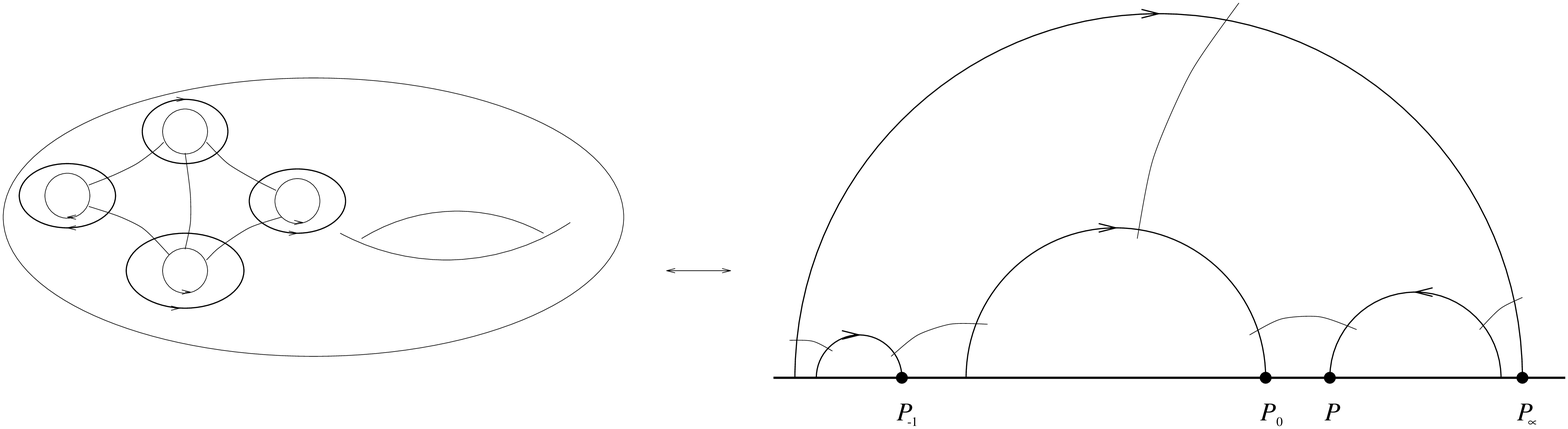}}
\caption{A geodesic octagon corresponding to a rectangle of the triangulation.}
\label{fg-20}
\end{figure}

Let $S$ be a surface equipped with a hyperbolic metric 
with geodesic boundary.  Choose a trivalent graph $\Gamma$ on $S$ and cut $S$ along the 
 sides of the dual graph $\Gamma^{\vee}$ into hexagons, as on 
Figure \ref{fg-20}.  
The two hexagons
 sharing an edge form an octagon. It has four geodesic sides surrounding
 holes. Let us lift the octagon to the hyperbolic plane and continue the 
geodesic sides till the absolute $\R P^1$.

The orientations of the holes 
induce the orientations of the geodesics on ${\cal H}$, and hence the
endpoints of these four geodesics.  Therefore we get four
points $p_1, ..., p_4$ on the absolute $\R P^1$.  The orientation of $S$ induces 
a cyclic order of geodesic
sides of the octagon on $S$, and hence a cyclic  order of the points $p_1,
..., p_4$. 
The key fact is that this cyclic  order is compatible 
with the orientation of the absolute. Therefore the 
 cross-ratio $r^+(p_1, ...., p_4)$ is
positive. Finally, if we assume that $p_1$ is the point incident to the
diagonal $E$ of the octagon, then $r^+(p_1, ...., p_4)$ 
is the restriction of our coordinate function $X^{\Gamma}_{E}$ to ${\cal
  X}^+_S$.  The part a) of the theorem is proved.

Before we start the proof of  the part b) of Theorem \ref{3.27.02.1}, it is
instructive to consider the following 
toy example.
\vskip 3mm
{\bf 2. Describing the decorated Teichm\"uller space 
for a disc with $n$ marked points}. 
Let $V_2$ be a real two dimensional 
vector space with a symplectic form $\omega$. 
There is a canonical map
$
s: V_2 \lra S^2V_2, \quad v \lms v\cdot v
$. 
A symplectic structure in $V_2$ determines a quadratic form $(\cdot, \cdot)$ 
in $S^2V_2$:
$$
(v_1\cdot v_2, v_3 \cdot v_4) := 
\frac{1}{2}\Bigl( \omega(v_1, v_3)\omega(v_2, v_4) + 
\omega(v_1, v_4)\omega(v_2, v_3)\Bigr).
$$
In particular 
\begin{equation} \label{5.7.02.3}
(s(v_1), s(v_2)) = \omega(v_1, v_2)^2.
\end{equation}
It is of signature $(2,1)$. 
The vectors $s(v)$ are isotropic and lie on one half of the isotropic 
cone in $S^2V_2$, denoted $Q_+$.  We identify 
$Q_+$ with the quotient of $V_2 -0$ 
by the reflection $v \lms -v$.

Observe that the projectivisation $PQ_+$ of the cone $Q_+$ is identified 
with the projective line $PV_2 = \R P^1$. So an orientation of the vector
space $V_2$ induces an orientation of $PQ_+$. A point $q$ on the cone $Q_+$
projects to the  point $\overline q \in PQ_+$. 
A configuration $(q_1, ..., q_n)$ of $n$ points on the cone
$Q_+$ is {\it positive} if the corresponding projective configuration of points 
$(\overline q_1, ..., \overline q_n)$ in $PQ_+ = PV_2$ is positive, i.e. its
cyclic order coincides with the one induced by the orientation of the
projective line $PV_2$. 

Denote by ${\rm Conf}^+_n(Q_+)$ the space of positive configurations  
of $n$ points $(q_1, ..., q_n)$ on the cone $Q_+$, and by ${\rm CConf}^+_n(Q_+)$ 
the corresponding
space of cyclic positive configurations. 

Recall    the space $\widetilde {\rm Conf}_n(V_2)$ of twisted cyclic
configurations of $n$ vectors in athe vector space $V_2$. 
Let $\widetilde {\rm Conf}^+_n(V_2)$ 
be the set of $\R_{>0}$-points for the positive atlas defined in Section 8.

\begin{proposition} \label{1.7.04.1} There is a canonical isomorphism
\begin{equation} \label{1.7.04.2}
{\rm CConf}^+_n(Q_+) \stackrel{\sim}{\lra} \widetilde {\rm Conf}^+_n(V_2).
\end{equation}
\end{proposition}

{\bf Proof}. Let us interpret the space ${\rm CConf}^+_n(Q_+)$ 
in terms of the real vector space $V_2$. Since $Q_+ = (V_2 -
\{0\})/\pm 1$, a configuration of $n$ points in $Q_+$ is the same thing as a 
$PSL_2(\R)$-configuration $(\pm  v_1, ..., \pm v_n)$ of $n$ non-zero vectors in $V_2$, each
considered up to a sign. Now starting from such a cyclic 
$PSL_2(\R)$-configuration $(\pm
v_1, ..., \pm v_n)$ we make a twisted cyclic $SL_2(\R)$-configuration of $n$
vectors in $V_2$ as follows. Consider a one dimensional subspace $L$ in $V_2$
which does not contain any of the vectors $\pm v_i$. Let $L'$ be one of the two
connected components of 
$V_2 - L$. Let $v_i'$ be the one of the two vectors $\pm v_i$ located in $L'$. 
We may assume that  
$(v'_1, ..., v'_n)$ is the order of the vectors $v_i'$ in $L'$: otherwise we
change cyclically the numerations of the vectors $v_i$. 
 Then the correspondence 
$$
(\pm v_1, ..., \pm v_n) \lra (v'_1, ..., v'_n)
$$
 provides  map (\ref{1.7.04.2}). To check that it is well defined notice that 
if we move the line $L$ so that it crosses the line spanned by the vector
 $v_1$, but not $v_2$, then we will get a configuration $(v'_2, ..., v'_n, -v'_1)$ which is
 twisted cyclic equivalent to $(v'_1, ..., v'_n)$. The proposition is
 proved. 
\vskip 3mm
The isomorphism (\ref{1.7.04.2}) transforms the real 
coordinates $|\Delta(\pm v_i,
\pm v_j)|$, $i<j$,  
on the left to the positive coordinates $\Delta(v'_i,
v'_j)$, $i<j$, on the right.

\begin{figure}[ht]
\centerline{\epsfbox{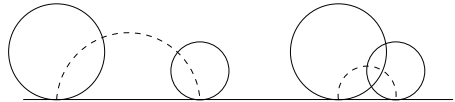}}
\caption{The distance between the horocycles can be any real number.}
\label{fg12}
\end{figure}

\vskip 3mm
{\bf 3. Proof of Theorem \ref{3.27.02.1}b)}. 
Let us realize the hyperbolic plane 
as the sheet of the hyperboloid $$
-x_1^2 - x_2^2 + x_3^2 = 1, \quad x_3>0.
$$ 
Denote by 
$<x, x'>$ the symmetric bilinear form 
corresponding to the qua\-dra\-tic form $-x_1^2 - x_2^2 + x_3^2$. 
 The horocycles on the hyperbolic plane are described by 
vectors $x$ on the cone $Q_+$, defined by the equation 
$<x,x>=0$ in the upper half space $x_3>0$. Namely, such a vector $x$ 
provides a horocycle 
$$
h_x:= \{\xi \quad | \quad <\xi, x>=1, \quad <\xi,\xi>=1, \quad \xi_3>0\}.
$$  
The distance $\rho(h_{x}, h_{x'})$ between the two horocycles $h_{x}$ and $h_{x'}$ 
is defined as follows. Take the unique   geodesic perpendicular to both horocycles. 
The horocycles cut out  on it a geodesic segment of finite length. Then  $\rho(h_{x}, h_{x'})$ 
is the length of this segment taken with the sign plus if the segment is outside 
of the discs inscribed into the horocycles, and with the sign minus otherwise, see Figure \ref{fg12}. 

Observe that $<x, x'>$ is positive unless  $x$ is proportional to $x'$. 
The distance between the two horocycles can be 
computed by the formula
\begin{equation} \label{5.7.02.2}
\rho(h_{x}, h_{x'})= \frac{1}{2} \log <x, x'>.
\end{equation}

Now take a hyperbolic topological surface $S$ with $n$ punctures 
equipped with  a complete hyperbolic metric $g$. 
Its universal cover  is the hyperbolic plane. 
The punctures are described by  $\pi_1(S)$--orbits 
of some  arrows on the cone $<x, x>=0$, $x_3>0$. 

According to Penner \cite{P1}, a decoration of $(S,g)$ is given by 
a choice of a horocycle $h_i$ for each of the punctures $s_i$. 
Let ${\cal A}^+_S$ be the moduli space of decorated  complete hyperbolic surfaces of  
given topological type $S$. 
Penner's  coordinates $a^{\Gamma}_E$ on it  are 
defined using a  trivalent graph $\Gamma$ on $S$ and parametrised by the set $\{E\}$ of 
edges of $\Gamma$. Namely, 
consider the ideal triangulation of $S$ dual to the graph $\Gamma$. 
For each side of the triangulation 
there is a unique geodesic between the punctures 
isotopic to the side. 
The coordinate $a_E^{\Gamma}$ is the length of the segment of the 
geodesic $\gamma_E$ contained between the two horocycles, for the triangulation dual to $E$.  

We can interpret this definition as follows. 
Choose a lift $\widetilde \gamma_E$ of 
the geodesic $\gamma_E$ on the hyperbolic plane. It is defined up to the action of 
$\pi_1(S)$. Let $l_1$ and $l_2$ be the arrows on the cone $<x, x> =0$ 
corresponding to the boundary of $\widetilde \gamma_E$. 
So  the plane spanned by $l_1$ and $l_2$ intersects 
the hyperboloid by a geodesic $\widetilde \gamma_E$. 
The arrows $l_1, l_2$  project to the punctures $s_1, s_2$ 
bounding $\gamma_E$.  A choice of horocycle(s) at 
the puncture $s_i$, $i=1,2$,  is equivalent to a choice of 
vector(s) $x_i$ on the arrow $l_i$.

We need to show that the 
coordinates $a_E^{\Gamma}$ 
coincide with our $\Delta$--coor\-di\-nates. We can spell the definition of 
the co\-or\-di\-na\-tes $a_E^{\Gamma}$ as follows. 
Choose an edge $E$ of the ribbon graph $\Gamma$. Let $F_1$ and $F_2$ be the 
two face paths sharing the edge $E$. Then a decoration on a twisted local 
$SL_2(\R)$--system ${\cal L}$ on $S$ 
provides  vectors $\pm v_1$ and $\pm v_2$,  
each defined only up to a sign, in the fiber ${\cal L}_E$ of the local 
system  ${\cal L}$ 
over the edge $E$. The vectors $\pm v_i$ are invariant up to a sign under 
the monodromy along the face path $F_i$. Then we assign to the edge $E$ the 
coordinate
$$
\Delta^{\Gamma}_E:= |\Delta(v_1, v_2)|.
$$
Observe that  (\ref{5.7.02.3}) implies 
$
\Delta(v_1, v_2)^2 = (s(v_1), s(v_2))
$. 
So thanks to (\ref{5.7.02.2})
$$
a_E^{\Gamma} = \log |\Delta(v_1, v_2)| = \log |\Delta^{\Gamma}_E|.
$$
It is well known that the monodromies along the conjugacy 
classes 
determined by  punctures are  unipotent. The theorem is proved.

\section{Laminations and canonical pairings}
\label{lam}

In this Section we give geometric definitions of two types of laminations,  
${\cal A}$- and ${\cal X}$-laminations, with integral, rational and real 
coefficients, on a punctured surface 
$S$. Our rational ${\cal A}$-laminations are essentially the same as 
Thurston's laminations \cite{Th}, which are usually considered for a 
surface $S$ without boundary. Our definition of  ${\cal A}$-laminations differs 
 in the 
treatment of holes, and is more convenient for our purposes. 

Our main results about  laminations are the following: 

i) The rational and real lamination spaces 
are canonically isomorphic to the points of the positive moduli spaces 
${\cal A}_{SL_2, S}$ and ${\cal X}_{PSL_2, S}$ with values in the tropical 
semifields $\Q^t, \R^t$. 
The same is true for the integral ${\cal X}$-laminations, while 
the set of integral ${\cal A}$-laminations contains ${\cal A}_{SL_2, S}(\Z^t)$  
as an 
explicitly described  subset. 

ii) The real ${\cal A}$-lamination space is  
canonically isomorphic to the space of Thurston's transversely 
measured laminations \cite{Th}, \cite{PH}. 

iii) We define a canonical pairing between the real ${\cal A}$- (resp. ${\cal X}$-) 
laminations and the ${\cal X}$- (resp. ${\cal A}$-)  Teichm\"uller spaces.

iv)  We define a canonical map from the integral ${\cal A}$- 
(resp. ${\cal X}$-) 
laminations to functions on the ${\cal X}$- (resp. ${\cal A}$-)  
moduli spaces. Its image 
consists of functions on the moduli space which are Laurent polynomials 
with positive integral coefficients in each of our coordinate systems. 
We conjecture that the image of this map provides a canonical basis 
for the space of such functions,  
and that the canonical  map 
admits a natural $q$-deformation. 

In \cite{FG4} we generalized results of this Section to the case of surfaces 
with marked points on the boundary.

The definitions of 
laminations are similar in spirit to the definitions of the singular
homology groups. There are two different ways to define the notion of
 laminations for surfaces with boundary, which are
similar to the definition of homology group with compact and
closed support respectively.

\vskip 3mm
 Throughout this section $S$ is a surface with $n>0$ holes, $\Gamma$ is
a trivalent graph embedded into $S$ in such a way that $S$ is
retractable onto it, and $\varepsilon_{ij}$, where $i$ and $j$ run
through the set of edges of $\Gamma$, is a skew-symmetric matrix
given by
$$
\begin{array}{rcl}
\varepsilon_{ij} &=& \mbox{ (number of vertices where the edge $i$
is
to the left of the edge $j$) }-\\
 & -& \mbox{ (number of vertices where the edge $i$ is to the right
of the edge $j$) }.
\end{array}
$$
The entries of $\varepsilon_{ij}$ take five possible values: $\pm
1,\pm 2$ or $0$.
\vskip 3mm
{\bf 1. ${\cal A}$- and ${\cal X}$-laminations and tropicalisations of the 
${\cal A}$- and ${\cal X}$- moduli spaces.}

\begin{definition} \label{9.8.03.10}
A {\em rational ${\cal A}$-lamination} on a surface $S$ is the 
homotopy class of a finite collection of disjoint simple unoriented closed curves  with rational
weights, subject to the following conditions and  equivalence
relations.

\begin{enumerate}
\item Weights of all curves are positive, unless a curve surrounds
a hole.

\item A lamination containing a curve of zero weight is equivalent
to that with this curve removed.

\item A lamination containing two homotopy equivalent curves of
weights $a$ and $b$ is equivalent to a lamination with one of
these curves removed and with the weight $a+b$ on the other.

\end{enumerate}

\end{definition}

The set of all rational ${\cal A}$-laminations on a  surface $S$ is
denoted by ${\cal A}_L(S,{\mathbb Q})$. It contains  a 
subset ${\cal A}_L(S,{\mathbb Z})$ of {\it integral ${\cal
A}$-laminations}, defined as  rational ${\cal
A}$-lamination with integral
weights.

Any rational ${\cal A}$-lamination can be represented by
a collection of $3g-3+n$ curves. Any integral ${\cal
A}$-lamination can be represented by a finite collection of curves
with unit weights.

\begin{definition} \label{9.8.03.11} A 
{\em Rational ${\cal X}$-lamination} on a  surface $S$  is
a pair consisting of 

a) A homotopy class of a finite collection of nonselfintersecting 
and pairwise non-intersecting curves either
closed or connecting two boundary components (possibly coinciding)
with positive rational weights assigned to each curve and subject
to the following equivalence relations:

\begin{enumerate}
\item A lamination containing a curve retractable to a boundary
component is equivalent to the one with this curve removed.

\item A lamination containing a curve of zero weight is equivalent
to the one with this curve removed.

\item A lamination containing two homotopy equivalent curves of
weights $a$ and $b$ is equivalent to the lamination with one of
these curves removed and with the weight $a+b$ on the other.
\end{enumerate}

b) A choice of orientations of all boundary
components but those that do not intersect curves of the
lamination. 

\end{definition}

Denote the space of rational ${\cal X}$-laminations on $S$ by 
${\cal X}_L(S,{\mathbb Q})$. It contains  a subset ${\cal X}_L(S,{\mathbb
Z})$ of { integral ${\cal X}$-laminations}, defined as 
rational ${\cal X}$-laminations representable by collections of curves
with integral weights.

 Any rational ${\cal X}$-lamination can be represented
by a collection of no more than $6g-6+2n$ curves (for Euler
characteristic reasons).  Any integral ${\cal X}$-lamination can
be represented by a finite collection of curves with unit weights.

\begin{theorem} \label{10.07.03.1} There are canonical isomorphisms
$$
{\cal A}_{SL_2, S}({\Bbb Q}^t) \stackrel{\sim}{=} {\cal A}_L(S, {\Bbb Q}), \qquad 
{\cal X}_{PSL_2, S}({\Bbb Q}^t) \stackrel{\sim}{=} {\cal X}_L(S, {\Bbb Q}).
$$
\end{theorem} 

The proof of this theorem will occupy the rest of this subsection. 
The idea is this. We are going to introduce 
natural coordinates on the space of ${\cal A}$- and ${\cal X}$-laminations 
corresponding to a given isotopy class  of a trivalent graph $\Gamma$  on $S$, 
which is homotopy equivalent to $S$. 
Our construction is a
modification of Thurston's ``train track'' (\cite{Th}, section 9
and \cite{PH}), and borrows from  \cite{F} 
with some modifications. 
Then we observe that the transformation 
rules for this coordinates under the flips 
are precisely the tropicalisations of the corresponding  
transformation rules for the canonical 
coordinates on the ${\cal A}$- and ${\cal X}$-moduli spaces. 
This proves the theorem. 
Let us start the implementation of this plan. 

\vskip 3mm
\paragraph{\bf Construction of coordinates for ${\cal A}$-laminations} 
We are going to assign to  a given rational ${\cal A}$-lamination rational
numbers on the edges of the graph $\Gamma$ and show that these
numbers are coordinates on the space of laminations.

 Retract the lamination to the graph in such a way that each curve
retracts to a path without folds on edges of the graph, and no
curve goes along an edge and then, without visiting another edge,
back. Assign to each edge $i$ the sum of weights of curves going
through it (Picture  (\ref{blam})), divided by two. The collection of these numbers,
one for each edge of $\Gamma$, is the desired set of coordinates.

\begin{equation}
{\epsfxsize5\baselineskip\epsfbox{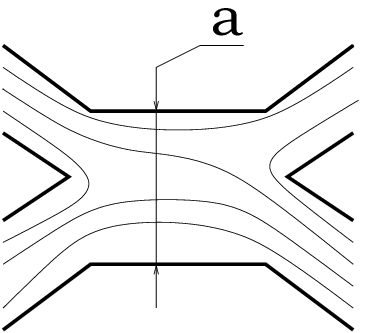}} \label{blam}
\end{equation}
\vskip 3mm
{\bf Remark}. We divided the natural coordinates of ${\cal A}$-laminations 
by two to simplify some formulas 
for the multiplicative canonical pairing, see below. 
\vskip 3mm
\paragraph{\bf Reconstruction of ${\cal A}$-laminations}
We need to prove that these numbers  are coordinates. For 
this purpose we describe an inverse construction which gives a
lamination starting from the numbers on edges.

 First of all note that if we are able to reconstruct a
lamination corresponding to a set of numbers $\{a_i\}$, we can do
this as well for the set $\{ra_i\}$ and $\{a_i+t\}$ for any
rational $r \ge 0$ and $t$. Indeed, multiplication of all numbers
by $r$ can be achieved by multiplication of all weights by $r$ and
adding $t$ is obtained by adding loops with weight $t/2$ around
each hole. Therefore we can reduce our problem to the case when
$\{a_i\}$ are positive integers and any three numbers $a_1, a_2,
a_3$ on the edges incident to each given vertex satisfy the
following triangle and parity conditions
\begin{eqnarray}
\vert a_1-a_2 \vert \le a_3 \le  a_1+a_2  \\  a_1+a_2+a_3  \hbox{
is an integer}
\end{eqnarray}
 Now the reconstruction of the lamination is almost obvious. Draw
$a_i$ lines on the $i$-th edge and connect these lines at vertices
in an non-intersecting way (Picture (\ref{vertex})), which can be done
unambiguously.

\begin{equation}
{\epsfxsize9\baselineskip\epsfbox{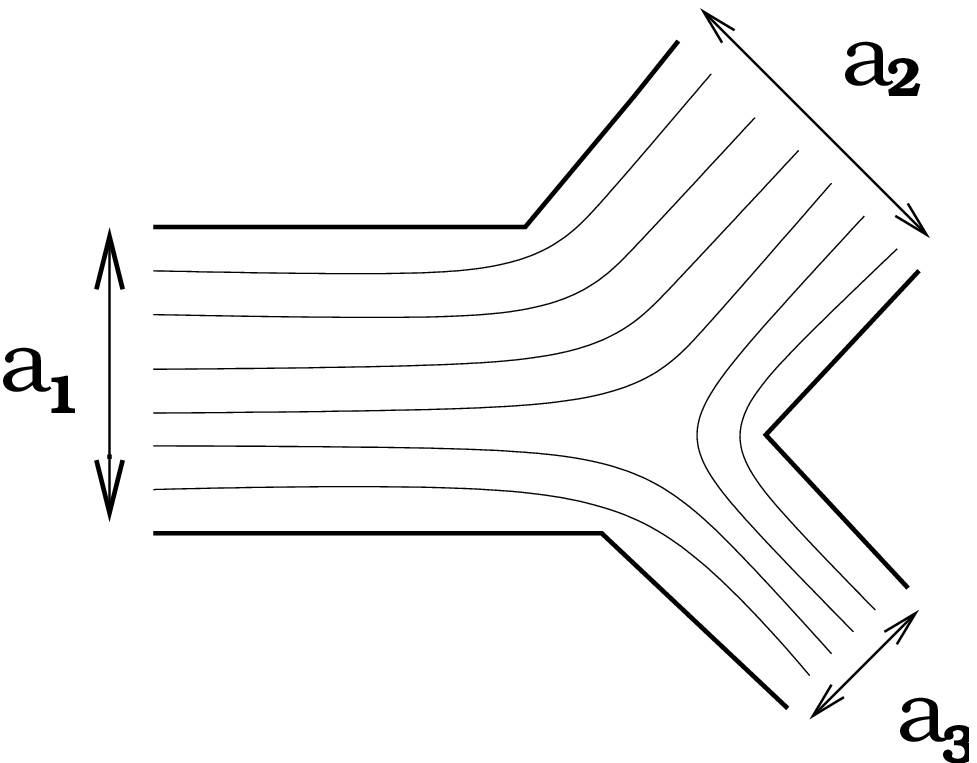}} \label{vertex}
\end{equation}
\vskip 3mm
\paragraph{Graph change for ${\cal A}$-laminations.}
The flip of the edge $k$ changes the coordinates according to the
formula:
$$
a_i \rightarrow \left\{ \begin{array}{lll} a_i & \mbox{ if
}&
k \neq i,\\
-a_k + \max(\sum\limits_{j | \varepsilon_{k j}>0} \varepsilon_{k j}
a_j,-\sum\limits_{j | \varepsilon_{k j}<0} \varepsilon_{k j} a_j)
&\mbox{ if }& k = i.
\end{array}\right.
$$
 Or in a graphical form:
\begin{equation}
\setlength{\unitlength}{0.15mm}%
\begin{picture}(500,275)(20,485)
\thicklines \put(280,700){\line( 1,-2){ 40}} \put(320,620){\line(
1, 0){280}} \put(600,620){\line( 1, 2){ 40}} \put(600,620){\line(
1,-2){ 40}} \put(320,620){\line(-1,-2){ 40}} \thinlines
\put(180,620){\vector(-1, 0){  0}} \put(180,620){\vector( 1, 0){
55}} \thicklines \put(100,540){\line( 2,-1){ 80}}
\put(100,700){\line( 0,-1){160}} \put(100,540){\line(-2,-1){ 80}}
\put( 20,740){\line( 2,-1){ 80}} \put(100,700){\line( 2, 1){ 80}}
\put( 40,745){\makebox(0,0)[lb]{$a_1$}}
\put(160,745){\makebox(0,0)[rb]{$a_2$}}
\put(120,620){\makebox(0,0)[lc]{$a_0$}}
\put(160,485){\makebox(0,0)[rt]{$a_3$}}
\put(40,485){\makebox(0,0)[lt]{$a_4$}}
\put(300,540){\makebox(0,0)[lt]{$a_4$}}
\put(620,540){\makebox(0,0)[rt]{$a_3$}}
\put(620,690){\makebox(0,0)[rb]{$a_2$}}
\put(300,690){\makebox(0,0)[lb]{$a_1$}}
\put(460,640){\makebox(0,0)[cb]{$\max(a_1+a_3,a_2+a_4)-a_0$}}
\end{picture}
\label{dflipl}
\end{equation}
(Only part of the graph is shown here, the numbers on the other
edges remain unchanged.)

\vskip 3mm
\paragraph{\bf Construction of coordinates for ${\cal X}$-laminations} 
We are going to assign for a given point of the space ${\cal
X}_L(S,{\mathbb Q})$ a set of rational numbers on edges of the
graph $\Gamma$, and show that these numbers are global coordinates
on this space.

 Straightforward retraction of an  ${\cal X}$-lamination onto $\Gamma$
is not good because some curves may shrink to points or finite
segments. To avoid this problem, let us first rotate each oriented
boundary component infinitely many times in the direction
prescribed by the orientation as shown on Picture (\ref{twist}).
\begin{equation}
{\epsfxsize10\baselineskip\epsfbox{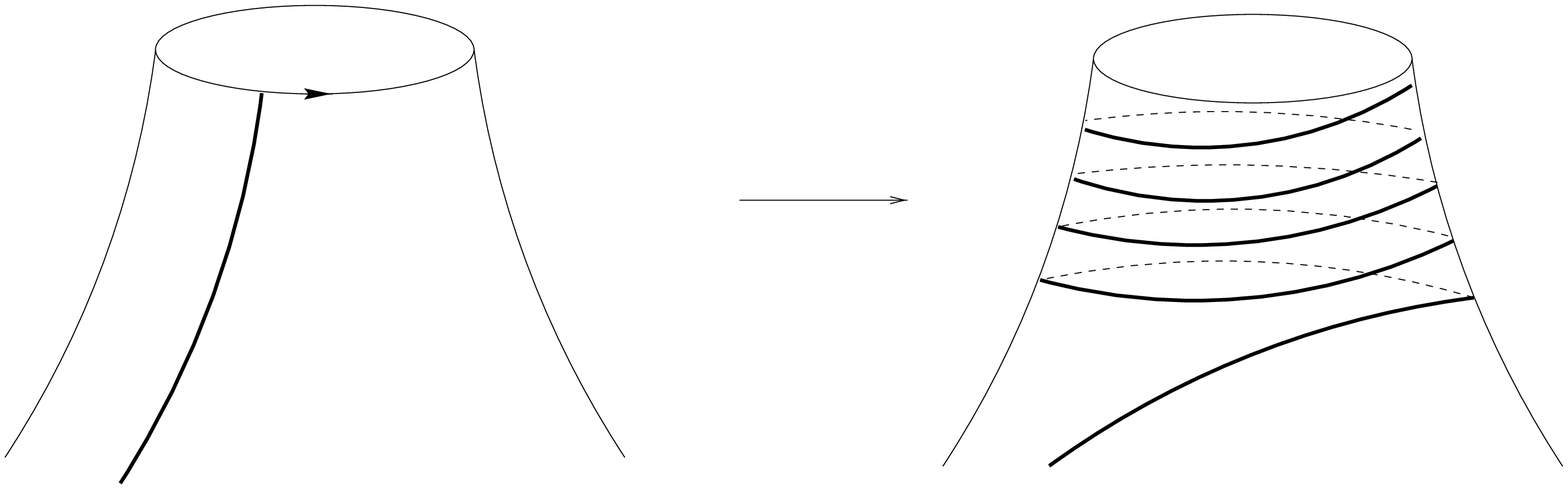}} \label{twist}
\end{equation}
 The resulting lamination can be retracted on $\Gamma$ without
folds. However we can get infinitely many curves going through an
edge. Now assign to the edge the sum with signs of weights of
curves that go diagonally, i.e., such that being oriented turn
to the left at one end of the edge and to the right at another one
(in this case we choose the plus sign), or first to the right and
then to the left (in this case we choose the minus sign), as shown
on Picture (\ref{ublam}).
\begin{equation}
{\epsfxsize5\baselineskip\epsfbox{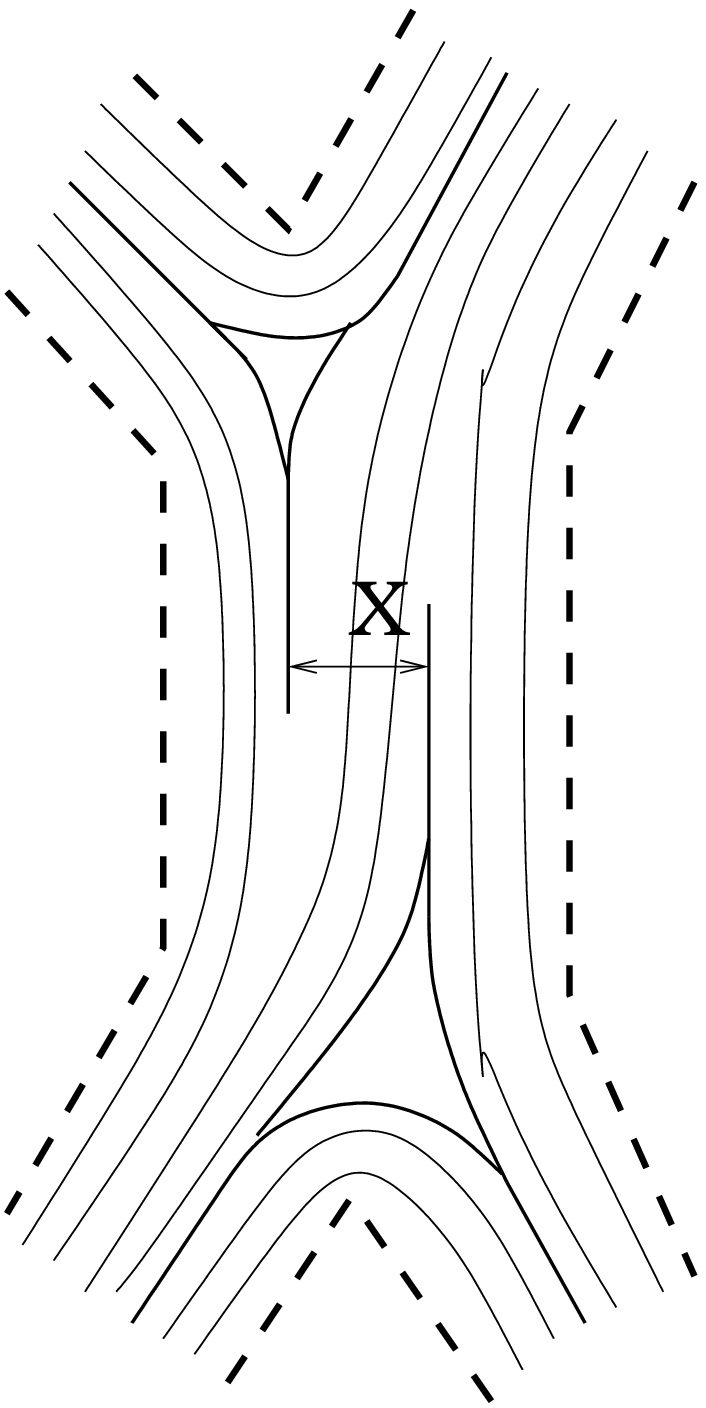}} \label{ublam}
\end{equation}
  The collection of such rational numbers, one for
each edge, is the desired coordinate system on ${\cal
X}_L(S,{\mathbb Q})$.
  Note that the number of curves going diagonally is always finite and
therefore the numbers assigned to an edge are well defined.
Indeed, consider the curves retracted on the graph. We can mark a
finite segment of each non-closed curve in such a way that each of
two unmarked semi-infinite rays goes only around a single face and
therefore never goes diagonally along the edges. Therefore only
the finite marked parts of curves contribute to the numbers on the
edges.
\vskip 3mm
\paragraph{\bf Reconstruction of ${\cal X}$-laminations} 
We need to prove that these numbers are coordinates indeed. We
will do it by describing an inverse construction.  Note that if we
are able to construct a lamination corresponding to a set of
numbers $\{x_i\}$, we can equally do it for the set $\{r x_i\}$
for any rational $r > 0$.  Therefore we can reduce our task to the
case when all numbers on edges are integral.  Now draw ${\mathbb
Z}$-infinitely many lines along each edge.  In order to connect
these lines at vertices we need to split them at each of the two
ends into two ${\mathbb N}$-infinite bunches to connect them with
the corresponding bunches of the other edges.  Let us do so at
the $i$-th edge, such that $x_i \leq 0$ (resp.  $x_i
\geq 0$), in such a way that the intersection of the right
(resp. left) bunches at the both ends of the edge consist of $x_i$
lines (resp. $-x_i$ lines).  Here the left and the right side are
considered from the center of the edge toward the corresponding
end. The whole procedure is illustrated on Picture  (\ref{ublam-}). The
resulting collection of curves may contain infinite number of
curves surrounding holes, which should be removed in accordance
with the definition of an unbounded lamination.

\begin{equation}
{\epsfxsize13\baselineskip\epsfbox{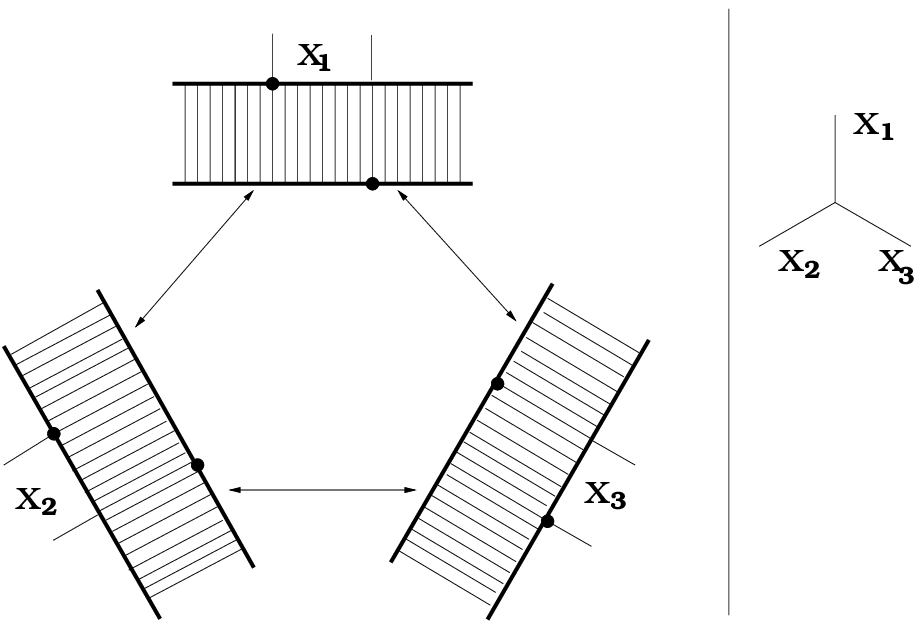}} \label{ublam-}
\end{equation}

 Note that although we have started with infinite bunches of curves the
resulting lamination is finite. All these curves glue together
into a finite number of connected components and possibly infinite
number of closed curves surrounding punctures. Indeed,  any curve of
the lamination is either closed or goes diagonally along at least
one edge. Since the total number of pieces of curves going
diagonally $I=\sum_i\vert x_i\vert$ is finite,  the resulting
lamination contains no more than this number of connected
components. (In fact the number of connected components equals $I$
provided the numbers $x_i$ are all non-positive or all
nonnegative.)
\vskip 3mm
\paragraph{\bf Graph change for ${\cal X}$-laminations} 
The flip of the edge $k$ changes the
coordinates according to the formula
$$
x(i) \rightarrow \left\{ \begin{array}{lll}
-x_i & \mbox{ if } & i = k,\\
x_i & \mbox{ if } & \varepsilon_{k i}=0 \mbox{ and }k \neq i, \\
x_i + \max(0,x_k) & \mbox{ if } & \varepsilon_{k i}>0,\\
x_i - \max(0,-x_k) & \mbox{ if } & \varepsilon_{k i}<0,
\end{array}\right.
$$
or, in the graphical form:
\begin{equation}
\setlength{\unitlength}{0.15mm}%
\begin{picture}(500,275)(20,485)
\thicklines \put(280,700){\line( 1,-2){ 40}} \put(320,620){\line(
1, 0){280}} \put(600,620){\line( 1, 2){ 40}} \put(600,620){\line(
1,-2){ 40}} \put(320,620){\line(-1,-2){ 40}} \thinlines
\put(180,620){\vector(-1, 0){  0}} \put(180,620){\vector( 1, 0){
55}} \thicklines \put(100,540){\line( 2,-1){ 80}}
\put(100,700){\line( 0,-1){160}} \put(100,540){\line(-2,-1){ 80}}
\put( 20,740){\line( 2,-1){ 80}} \put(100,700){\line( 2, 1){ 80}}
\put( 40,745){\makebox(0,0)[lb]{$x_1$}}
\put(160,745){\makebox(0,0)[rb]{$x_2$}}
\put(120,620){\makebox(0,0)[lc]{$x_0$}}
\put(160,485){\makebox(0,0)[rt]{$x_3$}} \put(
40,485){\makebox(0,0)[lt]{$x_4$}}
\put(275,545){\makebox(0,0)[lt]{$x_4 - \min(x_0,0)$}}
\put(645,545){\makebox(0,0)[rt]{$x_3 + \max(x_0,0)$}}
\put(645,700){\makebox(0,0)[rb]{$x_2 - \min(x_0,0)$}}
\put(275,700){\makebox(0,0)[lb]{$x_1 + \max(x_0,0)$}}
\put(470,640){\makebox(0,0)[cb]{$-x_0$}}
\end{picture}
\label{Hflipl}
\end{equation}
(Only part of the graph is shown here, the numbers on the other
edges remain unchanged.)

\vskip 3mm
{\bf End of the proof of Theorem \ref{10.07.03.1}}. 
Recall that positive structures on the moduli spaces 
${\cal A}_{SL_2, S}$ and ${\cal X}_{PSL_2, S}$ have been defined 
using the coordinate systems corresponding 
to trivalent ribbon graphs $\Gamma$ embedded to $S$. 
The coordinate systems on the lamination spaces are defined 
using the same graphs. 
The transformation rules (\ref{dflipl}) and (\ref{Hflipl})
are tropical limits of the ones for the respective moduli spaces 
${\cal A}_{SL_2, S}$ and ${\cal X}_{PSL_2, S}$, as one sees comparing the pictures (\ref{dflipl}) and (\ref{Hflipl}) with the 
pictures (\ref{10.30.03.22}) and (\ref{10.30.03.23}) 
where these rules are exhibited. The theorem is proved. 

\begin{equation} \label{10.30.03.22}
\setlength{\unitlength}{0.15mm}%
\begin{picture}(500,275)(20,485)
\thicklines \put(280,700){\line( 1,-2){ 40}} \put(320,620){\line(
1, 0){280}} \put(600,620){\line( 1, 2){ 40}} \put(600,620){\line(
1,-2){ 40}} \put(320,620){\line(-1,-2){ 40}} \thinlines
\put(180,620){\vector(-1, 0){  0}} \put(180,620){\vector( 1, 0){
55}} \thicklines \put(100,540){\line( 2,-1){ 80}}
\put(100,700){\line( 0,-1){160}} \put(100,540){\line(-2,-1){ 80}}
\put( 20,740){\line( 2,-1){ 80}} \put(100,700){\line( 2, 1){ 80}}
\put( 40,745){\makebox(0,0)[lb]{$A_1$}}
\put(160,745){\makebox(0,0)[rb]{$A_2$}}
\put(120,620){\makebox(0,0)[lc]{$A_0$}}
\put(160,485){\makebox(0,0)[rt]{$A_3$}}
\put(40,485){\makebox(0,0)[lt]{$A_4$}}
\put(300,540){\makebox(0,0)[lt]{$A_4$}}
\put(620,540){\makebox(0,0)[rt]{$A_3$}}
\put(620,690){\makebox(0,0)[rb]{$A_2$}}
\put(300,690){\makebox(0,0)[lb]{$A_1$}}
\put(460,640){\makebox(0,0)[cb]{$\frac{A_1A_3+A_2A_4}{A_0}$}}
\end{picture}
\end{equation}

\begin{equation} \label{10.30.03.23}
\setlength{\unitlength}{0.15mm}%
\begin{picture}(500,275)(20,485)
\thicklines \put(280,700){\line( 1,-2){ 40}} \put(320,620){\line(
1, 0){280}} \put(600,620){\line( 1, 2){ 40}} \put(600,620){\line(
1,-2){ 40}} \put(320,620){\line(-1,-2){ 40}} \thinlines
\put(180,620){\vector(-1, 0){  0}} \put(180,620){\vector( 1, 0){
55}} \thicklines \put(100,540){\line( 2,-1){ 80}}
\put(100,700){\line( 0,-1){160}} \put(100,540){\line(-2,-1){ 80}}
\put( 20,740){\line( 2,-1){ 80}} \put(100,700){\line( 2, 1){ 80}}
\put( 40,745){\makebox(0,0)[lb]{$X_1$}}
\put(160,745){\makebox(0,0)[rb]{$X_2$}}
\put(120,620){\makebox(0,0)[lc]{$X_0$}}
\put(160,485){\makebox(0,0)[rt]{$X_3$}} 
\put(40,485){\makebox(0,0)[lt]{$X_4$}}
\put(245,545){\makebox(0,0)[lt]{$X_4(1+(X_0)^{-1})^{-1}$}}
\put(645,545){\makebox(0,0)[rt]{$X_3(1+X_0)$}}
\put(685,700){\makebox(0,0)[rb]{$X_2(1+(X_0)^{-1})^{-1}$}}
\put(275,700){\makebox(0,0)[lb]{$X_1(1+X_0)$}}
\put(470,640){\makebox(0,0)[cb]{$(X_0)^{-1}$}}
\end{picture}
\end{equation}

\vskip 3mm
\paragraph{Relations and common properties of
${\cal X}_L(S,{\mathbb Q})$ and ${\cal A}_L(S,{\mathbb Q})$}

{\bf 1.}  Since the transformation rules for coordinates
(\ref{dflipl}) and (\ref{Hflipl}) are continuous with respect to  the
standard topology of ${\mathbb Q}^n$, the coordinates define a
natural topology on the lamination spaces.  One can now define the
spaces of {\em real laminations} as completions of the
corresponding spaces of rational laminations. These spaces are
denoted as ${\cal A}_L(S,{\mathbb R})$ and ${\cal X}_L(S,{\mathbb
R})$, respectively. Of course, we have the coordinate systems on
these spaces automatically. 
Evidently there are canonical isomorphisms
$$
{\cal A}_{SL_2, S}({\Bbb R}^t) \stackrel{\sim}{=} {\cal A}_L(S, {\Bbb R}), \qquad 
{\cal X}_{PSL_2, S}({\Bbb R}^t) \stackrel{\sim}{=} {\cal X}_L(S, {\Bbb R}).
$$

The constructed spaces of real laminations coincide  with the
spaces of transversely measured laminations introduced by
W.Thurston \cite{Th}. We'll prove this statement below.


{\bf 2.} An ${\cal X}$-lamination is integral if and only if it
has integral coordinates.  
However an ${\cal A}$-lamination is integral if
and only if it has half-integral coordinates and their
sum at every vertex is an integer. So we have 
$$
{\cal X}_{PSL_2, S}({\Bbb Z}^t) \stackrel{\sim}{=} {\cal X}_L(S, {\Bbb Z}), 
\qquad 
{\cal A}_{SL_2, S}({\Bbb Z}^t) \subset {\cal A}_L(S, {\Bbb Z})\subset 
{\cal A}_{SL_2, S}(\frac{1}{2}{\Bbb Z}^t).
$$ 

{\bf 4.} There is a  canonical  map 
$$
p: {\cal A}_L(S,{\mathbb Q}) \rightarrow {\cal X}_L(S,{\mathbb Q})
$$
which removes all curves retractable to  boundary components. 
In coordinates it is given by
\begin{equation}
x_i = \sum_j \varepsilon_{ij}a_j.
\end{equation}

{\bf 5.} For any given ${\cal X}$-lamination and any given puncture
one can calculate the sum
\begin{equation}
\sum_{i\in \gamma}x_i,
\end{equation}
where $\gamma$ is the set of edges surrounding the puncture. The
absolute value of this sum gives the total weight of the curves
entering the puncture, and its sign tells the difference between 
the orientation of the puncture and the one induced by the orientation of the surface.

{\bf 6.}
 For any given ${\cal A}$-lamination and any given puncture consider the
expression
\begin{equation} \label{signA}
\max_{i\in \gamma} (a_{i''}-a_i-a_{i'}),
\end{equation}
where $i'$ is the edge next to $i$ in counterclockwise direction
and $i''$ is the edge between $i$ and $i'$. This expression  gives
the weight of the loop surrounding the puncture.

\vskip 3mm
{\bf 2. Additive canonical pairings and the intersection pairing}.  
An additive canonical pairing  is a
function of two arguments. One argument is a rational (or, later on, real) 
lamination, which can be 
either ${\cal A}$- or ${\cal X}$-lamination, and the other is a 
point of the opposite type  Teichm\"uller space, ${\cal X}^+:= {\cal X}^+_{PSL_2, S}$ or 
${\cal A}^+:= {\cal A}^+_{SL_2, S}$:
$$
I:{\cal A}_L(\Q)\times{\cal X}^+\rightarrow {\mathbb R}, \quad 
I:{\cal X}_L(\Q)\times {\cal A}^+\rightarrow {\mathbb R}.
$$
We will also define the {\it intersection pairing}
$$
{\cal I}: {\cal A}_L(\Q)\times{\cal X}_L(\Q)\rightarrow {\mathbb R}
$$
 which should be thought of as a degeneration of an additive pairing. 
Abusing notation,  we shall denote the first two 
of them by a single letter $I$. We shall also denote the canonical maps
$p: {\cal A}^+\rightarrow {\cal X}^+$ 
and $p: {\cal A}_L\rightarrow {\cal X}_L$ by
the same letter $p$. 

Let us define these pairings. Below by the length of a curve on 
a surface with metric of curvature $-1$ we always mean the length of the unique 
geodesic isotopic to this curve. 
\begin{definition}
\begin{enumerate}
\item Let  $l \in {\cal X}_L $ be a single closed curve. Then $I(l,m)$ is 
its length w.r.t the hyperbolic
metric on $S$ defined by $m$.

\item Let $l \in {\cal A}_L $ be a single closed curve. 
Then $I(l,m)$ is equal to $\pm$ its length
w.r.t the hyperbolic metric on $S$ defined by $m$. 
The sign is $+$ unless $l$ surrounds a negatively
oriented hole.

\item Let $l$ be a curve connecting two boundary components and $m
\in {\cal A}^+$. Let us assume that orientations of these boundary components 
are induced by the ones of the surface. 
Then $I(l, m)$ is the signed length of the curve $l$ between
the horocycles. 

\item Let $l_1$ and $l_2$ be non-intersecting collections
of curves. Then 
$I(\alpha
l_1+\beta l_2,m) = \alpha I(l_1,m)+\beta I(l_2,m)$, and similarly for ${\cal I}$. 

\item Let $l_1$ and $l_2$ be
two curves. Then ${\cal I}(l_1,l_2)$ is the minimal number of
intersection points of $l_1$ and $l_2$.  

\end{enumerate}
\end{definition}

In the proposition below by coordinates we always mean 
the coordinates corresponding to a certain given
trivalent graph $\Gamma$ such that 
$S$ is retractable on $\Gamma$. All the statements are 
assumed to be valid for every such graph $\Gamma$. 
Recall the logarithms $x_i:= \log X_i$ and $a_i:= \log A_i$ 
of the canonical positive coordinates 
on the ${\cal X}$- and ${\cal A}$- Teichm\"uller spaces. 
The coordinates on the lamination spaces are also denoted by $x_i$ and $a_i$. 

\begin{proposition}~\label{properties-of-XL}
\begin{enumerate}

\item\label{continuity} The pairings  $I$ and ${\cal I}$ are
 continuous.

\item\label{symmetry-add} 
$I(l, p(m))=I(p(l), m)$, where
$l \in {\cal A}_L$ and $m\in {\cal A}^+$.

\item\label{explicit} If a point $l$ of the space ${\cal X}_L$ 
has positive coordinates $x_1,\ldots,x_n$
and a point $m$ of the Teichm\"uller 
${\cal A}$-space has coordinates $a_1,\ldots,a_n$,
then $I(l, m)= \sum_\alpha a_\alpha x_\alpha $.

\item\label{limit-add} Let $m$ be a point of a Teichm\"uller space 
with coordinates $u_1,\ldots,u_n$ with respect to a  
graph $\Gamma$. Let $C$ be a positive real number. Denote by 
$C\cdot m$ a point with coordinates $C u_1,\ldots, Cu_n$.
Let $m_L$ be a lamination with the coordinates $u_1,\ldots, u_n$. Then
$\lim\limits_{C \rightarrow \infty}I(l, C\cdot m)/C= {\cal I}(l, m_L)$.
\end{enumerate}
\end{proposition}
Observe that the last property implies the properties 
\ref{symmetry-add} and \ref{explicit} for ${\cal I}$.

 {\bf Proof.} The property \ref{symmetry-add} 
follows from the definition
and the property \ref{continuity}. The property \ref{explicit} follows from
the observation that if all coordinates of an ${\cal X}$-laminations 
are positive w.r.t. a
given graph $\Gamma$, all its leaves can be deformed to paths
across the edges of $\Gamma$. Indeed, if all coordinates are $\geq 0$, a path in $\Gamma$ is of the form $...l...lr...r...$ where $lr$ corresponds to a positive coordinate, then such a path is homotopic to a curve from one puncture to another. The proof of  property \ref{limit-add} will be
postponed to our consideration of the multiplicative pairing. 


Let us first prove the continuity of the pairing 
between ${\cal A}$-laminations and the ${\cal X}$-Teichm\"uller
space. To do this it suffices to prove it for laminations without
curves with negative weights. Indeed, if we add such curve to a
lamination, the length obviously changes continuously. We are
going to show that the length of an integral lamination is a convex
function of its coordinates, i.e. that

\begin{equation} \label{ineq}
I(l_1, m) + I(l_2, m) \geq I(l_1 +_\Gamma l_2, m) 
\end{equation}
were by $l_1 +_\Gamma l_2$ we mean a lamination with coordinates
being sums of the respective coordinates of $l_1$ and $l_2$.
Taking into account the homogeneity of $I$, one sees that
the inequality (\ref{ineq}) holds for all rational laminations and
therefore can be extended by continuity to all real laminations.

Let us prove the inequality (\ref{ineq}).  Draw both laminations
$l_1$ and $l_2$ on the surface and deform them to be geodesic.
These laminations in general intersect each other in finite number
of points.  Then retract the whole picture to the ribbon graph in
such a way that no more intersection points appear and the
existing ones are moved to the edges.  Now it becomes obvious
that at each intersection point we can rearrange our lamination
cutting both intersecting curves at the intersection point and
gluing them back in another order in a way to make the resulting
set of curves homotopically equivalent to a non-intersecting
collection.  We can do it at each intersection point in two
different ways, and we use the retraction to the ribbon graph to
choose one of them.

Indeed, connect them as shown on Picture (\ref{convex}).
\begin{equation}
{\epsfxsize13\baselineskip\epsfbox{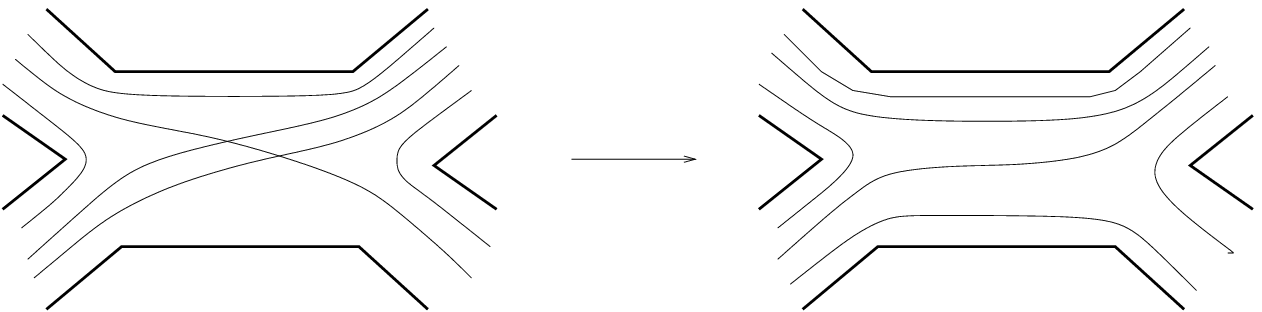}} \label{convex}
\end{equation}

One can easily see that the numbers on edges, corresponding to the
new lamination $l$ are exactly the sums of the numbers
corresponding to $l_1$ and $l_2$. On the other hand, the
lamination on the original surface is no longer geodesic, because
its curves may have breaks. But its length is exactly $I(l_1,m)
+ I(l_2,m)$. When we deform $l$ to a geodesic lamination, its
length can only decrease, which proves the inequality (\ref{ineq}).

The proof of the continuity when the first argument is an ${\cal
X}$-lamination and the second is a point of the ${\cal
A}$-Teichm\"uller space is  similar.  Call a pair of laminations
$l_1$ and $l_2$ {\em similar} if at each hole both laminations $l_1$
and $l_2$ have the same orientation or at least one of them does
not touch the hole. Once the inequality (\ref{ineq}) 
 is proved for
similar laminations $l_1$ and $l_2$, the continuity follows since the set of
laminations similar to a given one is a closed subset in the space
${\cal X}_L$. Conjecturally the inequality (\ref{ineq}) holds for 
any pair $l_1$ and $l_2$. But for a pair of similar laminations the arguments
used above works: one can first forget about orientations and move
the leaves in order to minimize the number of intersection points.
Then one can first spiralize the leaves according to the
orientations, retract the resulting collection of curves to the
graph and then resolve the intersection points just as above (Picture 
(\ref{convex})). The inequality (\ref{ineq}) follows.

Using 4  we deduce from this the convexity 
of the pairing ${\cal I}$. 

The continuity implies the following 
\begin{corollary} Our real 
${\cal A}$-laminations are canonically identified with Thurston's transversally 
measured laminations. 
\end{corollary}

{\bf Proof}. Thurston proved that rational laminations
are dense in the space of all transversally 
measured laminations w.r.t. the topology given by the length
function. Our rational ${\cal A}$-laminations coincide
 with Thurston's. Our canonical pairing of ${\cal A}$-laminations with the 
Teichm\"uller space is
nothing but Thurston's length function. Therefore the  continuity
shows that our topology coincides with Thurston's. The corollary is proved.

\begin{corollary} The additive canonical pairings and the intersection pairing 
 provide continuous 
pairings between the real lamination and Teichmuller/lamination spaces:
$$
I: {\cal A}_L(\R)\times {\cal X}\rightarrow {\mathbb R}, \quad 
I: {\cal X}_L(\R)\times {\cal A}\rightarrow {\mathbb R}, \quad 
{\cal I}: {\cal A}_L(\R)\times{\cal X}_L(\R)\rightarrow {\mathbb R}.
$$
\end{corollary}

We say that an element of ${\cal X}_{PSL_2, S}(\C)$ is quasifuchsian 
if the corresponding representation of $\pi_1(S)$ is quasi-Fuchsian. 
The set of all quasifuchsian elements forms an open domain 
of  ${\cal X}_{PSL_2, S}(\C)$ which we call the quasifuchsian domain. 

\begin{conjecture} {\it Let $l$ be a real ${\cal A}$-lamination. The
 function $I(l, *)$ is 
analytically continuable exactly to the quasifuchsian domain of ${\cal X}_{PSL_2, S}(\C)$.}
\end{conjecture}
\vskip 3mm
{\bf 3. Multiplicative canonical pairings.}
A multiplicative canonical pairing ${\Bbb I}$ is a pairing 
between integral ${\cal A}_{SL_2, S}$-laminations and ${\cal X}_{PSL_2, S}$-moduli space
or between integral ${\cal X}_{PSL_2, S}$-laminations and 
${\cal A}_{SL_2, S}$-moduli 
space. These pairings are better understood as  maps
$$
{\Bbb I}_A: {\cal A}_{SL_2, S}(\Z^t) \lra \Q({\cal X}_{PSL_2, S}), \qquad 
{\Bbb I}_X: {\cal X}_{PSL_2, S}(\Z^t) \lra \Q({\cal A}_{SL_2, S}) 
$$
 equivariant with respect to  the action of the mapping class group of $S$. 
The group $(\Z/2\Z)^n$ acts naturally on the 
space ${\cal X}_{PSL_2, S}$. We will define in Section 12.6 an action of this group 
on the space ${\cal A}_{PSL_2, S}$, so that the above maps intertwine these actions. 

Before defining the map ${\Bbb I}_A$ we will define an auxiliary map $$\widetilde {\Bbb
I}_A: {\cal A}_{SL_2, S}((\Z^t)/2) \lra \widehat{\Q({\cal X}_{PSL_2,
S})},$$ where the hat denotes certain covering of
 the space ${\cal X}_{PSL_2,S}$. 
 Afterwords we
shall check that the restriction of this pairing to ${\cal A}_{SL_2,
  S}(\Z^t)$ 
comes from ${\cal X}_{PSL_2,S}$, thus defining the desired pairing.

Here is our strategy. Say that a Laurent polynomial 
is {\it sign definite} if all its non zero coefficients are of the same sign. 
Given a sign definite Laurent polynomial we can multiply it by $\pm 1$ to make 
a positive Laurent polynomial out of it. 
We are going to define first the related maps 
${\Bbb I}'_{\cal A}$ and ${\Bbb I}'_{\cal X}$. Then we will show that 
 ${\Bbb I}'_{\cal A}(l)$ and ${\Bbb I}'_{\cal X}(l)$ are sign definite 
Laurent polynomials. By definition the elements ${\Bbb I}_{\cal A}$ and ${\Bbb I}_{\cal X}$  are the corresponding positive  
Laurent polynomials. We will produce our sign definite Laurent polynomials 
as products of  
sign definite Laurent polynomials, so we allow them to be 
well defined only up to a sign.   For example 
the trace of a matrix from $PSL_2(F)$ is defined up to a sign function. 

Recall that a curve $l$ connecting two punctures on $S$ 
provides a well defined regular function $\Delta(l)$ on the space 
${\cal A}_{SL_2, S}$. 
Further, let $l$ be a closed curve surrounding a 
puncture on $S$. Then the monodromy operator for a framed local system 
${\cal L}$ on $S$ 
around $l$ has a distinguished eigenvalue $\lambda_l$ determined by the framed structure: the eigenvector with this eigenvalue $\lambda_l$ generates 
the invariant flag assigned to this puncture. 

 Let $l$ be a loop on $S$. Let $\widetilde l$ be a section 
of the punctured tangent bundle $T'S$ of $S$ restricted to $l$. 
 The homotopy class of $\widetilde l$ is well defined.  
 By monodromy of a twisted $SL_2$-local system ${\cal L}$ 
on $S$ we understand the monodromy 
 of ${\cal L}$ around $\widetilde l$. 

\begin{definition} \label{17.27} 
\begin{enumerate}
\item Let $l$ be a loop with the weight $k$ which 
does not surround a puncture. Then 

i) the value of the function 
$\pm \widetilde {\Bbb I}'_{\cal A}(l)$ 
on the  framed $PSL_2$-local system ${\cal L}$ 
is the trace of the $k$-th power of the monodromy of ${\cal L}$ along $l$.  

ii) the value of the function 
${\Bbb I}'_{\cal X}(l)$ 
on the twisted decorated $SL_2$-local system ${\cal L}$ 
is the trace of the $k$-th power of the monodromy of ${\cal L}$ along 
the lifted loop $\widetilde l$ in $T'S$.

\item Let $l$ be a closed curve surrounding a puncture on $S$, and 
$k \cdot l$ the ${\cal A}$-lamination 
determined by this curve and a weight $k$. Then  
$\widetilde {\Bbb I}'_{\cal A}(k \cdot l) = \lambda_l^k$.

\item Let $l$ be a curve with
weight $k$ connecting two punctures. Let us assume that orientations of these  punctures are induced 
by the orientation of the surface. Then ${\Bbb I}'_{\cal X}(l) = \Delta(l)^k$. 

\item Let $l_1$
and $l_2$ be non-intersecting collections of curves, such that no 
isotopy class of a 
curve enters both $l_1$ and $l_2$. Then $
{\Bbb I}'_*(l_1+ l_2)= {\Bbb I}'_*(l_1)
{\Bbb I}'_*(l_2)$, where 
$*$ stands for ${\cal A}$ or ${\cal X}$, and $\widetilde {\Bbb I}'$ 
is needed in the ${\cal A}$-case.  
\end{enumerate}
\end{definition}
Here we assume that we have reduced the lamination 
to the form without equivalent curves. 

In this definition the map $\widetilde {\Bbb I}'_{\cal A}(l)$ is defined for all 
integral laminations. In the definition of the map ${\Bbb I}'_{\cal X}(l)$ 
we assume so far that if a curve of the lamination connects two punctures, then orientations  
of 
both punctures are induced by the orientation of the surface $S$. 
We will extend the definition to all ${\cal X}$-laminations in the end of the Section 12.6. Namely, 
recall that the group $(\Z/2\Z)^n$ acts on the ${\cal X}$-space by changing the orientations of the boundary components. We 
will  define an action of the group $(\Z/2\Z)^n$ on the ${\cal A}$-space, and 
extend the map ${\Bbb I}_{\cal X}$ to a map equivariant with respect to the 
action of the group $(\Z/2\Z)^n$. 


\vskip 3mm
In the theorem below by coordinates we always mean 
the coordinates corresponding to a certain fixed 
trivalent graph $\Gamma$ such that 
$S$ is retractable on $\Gamma$. All statements are 
assumed to be valid for every such graph $\Gamma$. The 
function ${\Bbb I}'_{\cal X}(l)$ is restricted to those laminations $l$ where
it has been defined.

\begin{theorem}~ \label{10.10.03.2}
\begin{enumerate}
\item\label{x-positivity} Let $a_1,\ldots,a_n \in \Z$ be coordinates 
of an integral 
${\cal
A}$-lamination $l$. Then the function $\widetilde {\Bbb I}'_{\cal A}(l)$ is a
sign definite Laurent polynomial with integral coefficients in the
coordinates 
 $X_1,\ldots,X_n$  on ${\cal X}$. Its highest term
is equal to $\prod_i X_i^{a_i}$ and the lowest one is
$\prod_i X_i^{-a_i}$. 
\item\label{a-positivity} The function ${\Bbb 
I}'_{\cal X}(l)$ is a sign definite Laurent polynomial with integral
coefficients in the coordinates $A_1,\ldots,A_n$ on ${\cal A}$.

\item\label{i-compatibility} Let $C \in \Z$. Then 
$\lim_{C \rightarrow \infty}\log{\Bbb I}_*(l, C\cdot m)/C= {\cal I}(l,m_L)$, where the
definitions of $C\cdot m$ and $m_L$ are the same as for the property 
4 in the additive case. 

\item\label{symmetry-mult} ${\Bbb I}_{\cal A}(l, p(m))={\Bbb I}_{\cal X}(p(l),m)$ 
where $l \in {\cal A}_{SL_2, S}(\Z^t)$ and $m \in {\cal A}_{SL_2, S}$. 

\item\label{multiplication} ${\Bbb I}_*(l_1){\Bbb I}_*(l_2)=
\sum_l c_*(l_1, l_2; l) {\Bbb I}_*(l)$, 
where $*$ stands for either ${\cal X}$ or 
${\cal A}$, and $c_*(l_1, l_2; l)$ are non-negative integers, and the
sum is finite.
\end{enumerate}
\end{theorem}

Observe that the property
\ref{x-positivity} implies that if an ${\cal
 A}$-lamination is integral, the corresponding function is lifted from 
${\cal X}_{PSL_2,S}$, and the map ${\Bbb I}_{\cal A}(l)$ is defined.
\vskip 3mm
{\bf Proof.}
Statement \ref{i-compatibility} immediately follows from 
the analogous additive statement
and the fact that for a curve $l$ connecting punctures
$\exp( I(l,m))={\Bbb I}_{\cal X}(l,m)$, and for a closed curve $l$ we have
$\cosh( I(l,m))={\Bbb I}_*(l,m)$.

As a corollary one deduces property \ref{limit-add} of 
Proposition \ref{properties-of-XL}. Property \ref{x-positivity} of Theorem \ref{10.10.03.2}
implies property \ref{i-compatibility} for integral laminations. This implies that 
the identity $\lim\limits_{C\rightarrow \infty}I(l, C\cdot m)/C= {\cal I}(l,m)$ holds
for an integral lamination $l$. By continuity it can be extended to
any laminations.

Statement \ref{symmetry-mult} is obvious if $l$ is a loop: if this loop
does not surround a puncture, both sides of
the equality are equal to the trace of the monodromy of the corresponding
local system; if $l$ surrounds a puncture, it is  $1$ for both sides. The
general case follows by multiplicativity.  

The properties \ref{x-positivity} and \ref{a-positivity} 
for the ${\cal X}$-moduli space follow essentially from 
 Theorem \ref{3.18.04.1q} and Corollary \ref{3.18.04.1}, valid for $PGL_m$. 
However for the sake of completeness of this Section, we prove them  below. 
 
Before presenting the proof of the 
properties \ref{x-positivity} and \ref{a-positivity}
recall the construction of monodromy matrices starting from coordinates of our 
moduli spaces given in Sections 8 and 9, specified to the case of 
$G=PGL_2$ and $G=SL_2$. 

Let us do it first for the ${\cal X}$-moduli
space. Replace the graph $\Gamma$ by another one $\Gamma'$ 
by gluing small triangles
into each vertex and define a graph connection on it as shown 
on Figure \ref{x-monodromy}. 
(The graph $\Gamma$ is shown by thin lines.)
\begin{equation}\label{x-monodromy}
\setlength{\unitlength}{0.3mm}%
\parbox{150pt}{
\begin{picture}(150,100)(-75,-50)
\thinlines
\put(-40,5){\line( 1,0){ 80}}
\put(-40,5){\line( -1,-2){ 20}}
\put(-40,5){\line( -1,2){ 20}}
\put(40,5){\line( 1,2){ 20}}
\put(40,5){\line( 1,-2){ 20}}
\thicklines
\put(50,20){\vector(0,-1){ 40}}
\put(50,-20){\vector(-2,1){40}}
\put(10,0){\vector(2,1){40}}
\put(-50,-20){\vector(0,1){ 40}}
\put(-50,20){\vector(2,-1){40}}
\put(-10,0){\vector(-2,-1){40}}
\put(-10,0){\line(1,0){20}}
\put(-50,20){\line(-1,2){10}}
\put(-50,-20){\line(-1,-2){10}}
\put(50,20){\line(1,2){10}}
\put(50,-20){\line(1,-2){10}}
\put(0,6){\makebox(0,0)[cb]{$k$}}
\put(0,-1){\makebox(0,0)[ct]{$f_k$}}
\put(-30,13){\makebox(0,0)[cb]{$g$}}
\put(-30,-12){\makebox(0,0)[ct]{$g$}}
\put(-50,0){\makebox(0,0)[rc]{$g$}}
\put(30,13){\makebox(0,0)[cb]{$g$}}
\put(30,-12){\makebox(0,0)[ct]{$g$}}
\put(52,0){\makebox(0,0)[lc]{$g$}}
\end{picture}}
\end{equation}
Here $f_k$ and $g$ are given by 
$$f_k=\left( \begin{array}{cc} 0&-X_k^{1/2}\\
X_k^{-1/2}&0\end{array}\right), \qquad 
g=\left( \begin{array}{cc} 1&1\\
-1&0\end{array}\right).
$$
Here, since we work with $PGL_2$, we employ a bit more symmetric notation for the matrix 
$f_k= \left( \begin{array}{cc} 0&-X_k\\
1&0\end{array}\right)$, using coordinates  $X_k^{\pm 1/2}$ on a covering. 
The monodromy around small triangles is the identity since $g^3= -1$ in $SL_2$, i.e. 
$g^3=1$ in $PSL_2$.
The orientations of the edges with the matrices $f_i$ do not play any role since
$f_i=f_i^{-1}$ in $PSL_2$.

Here is a similar construction for the ${\cal A}$-moduli space.
Replace the graph $\Gamma$ by another one by replacing each vertex by a hexagon and
define a graph connection on it as shown on Picture (\ref{a-monodromy}).
(The graph $\Gamma$ is shown by the thin lines).
\begin{equation}\label{a-monodromy}
\setlength{\unitlength}{0.3mm}
\parbox{150pt}{
\begin{picture}(150,120)(-75,-60)
\thicklines
\qbezier(0,-30)(12,0)(0,30)
\put(0,-30){\vector(-1,-2){0}}
\qbezier(0,-30)(-12,0)(0,30)
\put(0,30){\vector(1,2){0}}
\put(0,-30){\vector(2, -1){40}}
\put(-40,-50){\vector(2, 1){40}}
\put(40,50){\vector( -2, -1){40}}
\put(0,30){\vector(-2,1){ 40}}
\qbezier(40,-50)(53,-30)(80,-30)
\put(80,-30){\vector(1,0){0}}
\qbezier(40,-50)(67,-50)(80,-30)
\qbezier(40,50)(53,30)(80,30)
\put(40,50){\vector(-1,1){0}}
\qbezier(40,50)(67,50)(80,30)
\qbezier(-40,50)(-53,30)(-80,30)
\put(-80,30){\vector(-1,0){0}}
\qbezier(-40,50)(-67,50)(-80,30)
\qbezier(-40,-50)(-53,-30)(-80,-30)
\put(-40,-50){\vector(1,-1){0}}
\qbezier(-40,-50)(-67,-50)(-80,-30)
\put(-80,30){\vector( 0,-1){ 60}}
\put(80,-30){\vector( 0,1){ 60}}
\thinlines
\put(-40,0){\line( 1,0){ 80}}
\put(-40,0){\line( -1,-2){ 30}}
\put(-40,0){\line( -1,2){ 30}}
\put(40,0){\line( 1,2){ 30}}
\put(40,0){\line( 1,-2){ 30}}
\thicklines
\put(40,50){\line(0,1){ 20}}
\put(40,-50){\line(0,-1){ 20}}
\put(-40,50){\line(0,1){ 20}}
\put(-40,-50){\line(0,-1){ 20}}
\put(-80,-30){\line( -2,-1){ 20}}
\put(80,-30){\line( 2,-1){ 20}}
\put(80,30){\line( 2,1){ 20}}
\put(-80,30){\line( -2,1){ 20}}
\put(-20,7){\makebox(0,0)[cc]{$k$}}
\put(75,70){\makebox(0,0)[cc]{$i$}}
\put(75,-70){\makebox(0,0)[cc]{$j$}}
\put(7,6){\makebox(0,0)[lb]{$d_k$}}
\put(62,25){\makebox(0,0)[cc]{$d_i$}}
\put(62,-25){\makebox(0,0)[cc]{$d_j$}}
\put(83,0){\makebox(0,0)[lc]{$u_{ijk}$}}
\put(19,-47){\makebox(0,0)[cc]{$u_{jki}$}}
\put(19,47){\makebox(0,0)[cc]{$u_{kij}$}}
\end{picture}}
\end{equation}

Here the matrices $d_k$ and $u_{ijk}$ are given by:
$$d_k=\left( \begin{array}{rr}
0&-A_k\\ \frac{1}{A_k}&0\end{array}\right), \qquad 
u_{ijk}=\left( \begin{array}{cc} 1&0\\
\frac{A_k}{A_iA_j}&1\end{array}\right).
$$
The monodromy around the hexagons is the unity 
since $d_ju_{jki}d_ku_{kij}d_iu_{ijk}=1$ in
$SL_2$. The monodromy around the $2$-gons is $-1$. 
 The monodromy along a face path fixes the vector 
$\left( \begin{array}{c} 0\\1\end{array}\right)$.

Now proceed with the proof. Due to multiplicativity it suffices to prove 
statements 1 and 2 for a connected curve $l$.  Consider first  statement
\ref{x-positivity}.
In this case $l$ must be a loop. Let $\alpha$ be a loop on  
the graph $\Gamma'$ homotopic, on $S$, to the loop $l$. 
We may assume that the loop $\alpha$ contains no consecutive edges 
of little triangles on $\Gamma'$. Observe that 
\begin{equation} \label{3.19.04.1}
f_kg = \left( \begin{array}{cc} X_k^{1/2} & 0 \\
X_k^{-1/2} & X_k^{-1/2}\end{array}\right), \qquad f_kg^{-1} = 
\left( \begin{array}{cc}X_k^{1/2}  & X_k^{1/2} \\
0 & X_k^{-1/2}\end{array}\right).
\end{equation}
(The second equality is in $PGL_2$.) So calculation of the monodromy around 
such a  loop reduces to taking products of matrices of type (\ref{3.19.04.1}), 
just one of them per each edge of the original graph $\Gamma$ which appears 
in the loop $\alpha$. These matrices 
are positive and Laurent.  
Evidently the trace of such a monodromy operator is a positive Laurent
polynomial, and  the coefficients of the highest and  lowest
monomials are equal to 1.
\vskip 3mm
{\bf Remark}. In general the coordinates of an ${\cal A}$-lamination $l$  
are half-integers $a_1, ..., a_n \in \Z/2$, and the function 
$\widetilde {\Bbb I}'_{\cal A}(l)$ is defined only on 
 a covering of the space ${\cal X}$. 
However if $a_i \in \Z$, the proof shows that it is defined on the space 
${\cal X}$ itself. So we get the function ${\Bbb I}'_{\cal A}(l)$. 
\vskip 3mm
Consider now property \ref{a-positivity}. 
Call a matrix with Laurent-polynomial
coefficients {\em sign definite} if each its entry is a Laurent
polynomials with all positive or all negative coefficients. Let the
{\em type}         of a {\em sign definite} matrix $M$ be a matrix $|M|$ of the same size
as $M$ with entries equal to $+,-$ or $0$ indicating the sign of the
monomials of the corresponding entries of $M$. Define a partial product
on the set of types by the requirement, that $ab=c$ if for any two matrices
$M_1$ and $M_2$, such that $|M_1|=a$ and $|M_2|=b$, we have $|M_1M_2|=c$.

Consider  a path $l$ connecting two holes, retracted to the graph $\Gamma'$. The desired
intersection index is given by the top right element of the
monodromy matrix. The sequence of matrices along the path is
composed  of the matrices $u,u^{-1}$ and $d$. Since it connects two holes, 
we may assume that it starts and
ends with a $d$. The path can be deformed to make the
sequence  avoid the subsequences $uu^{-1},u^{-1}u,udu, d^2$ and
$u^{-1}du^{-1}$. To investigate the sign properties of the  product,
let us  replace every matrix in the sequence
by its type and investigate the type of the product. Let
$$
D:=|d_i|=\left( \begin{array}{cc} 0&-\\+&0\end{array}\right),\quad 
U^+=|u_{ijk}|=\left( \begin{array}{cc} +&0\\+&+\end{array}\right), \quad 
U^-=|(u_{ijk})^{-1}|=\left( \begin{array}{cc} +&0\\-&+\end{array}\right).
$$
Now observe that $(U^+)^n=U^+$ and $(U^-)^n=U^-$ for any $n\in {\mathbb N}$,
and replace the sequences  of $U^+$ and $U^-$ by single letters. The resulting
sequence has one of the following four forms:
$$
\begin{array}{c}DU^-D(U^+DU^-D)^kU^+D, \quad D(U^+DU^-D)^kU^+D, \quad DU^-D(U^+DU^-D)^k,\\
D(U^+DU^-D)^k.\end{array}
$$ 
 The crucial observation is this:
\begin{equation} \label{3.19.04.10}
\mbox{the products $U^+DU^-D$ and   
$DU^-DU^+$  are of the type $-P=\left( \begin{array}{cc} 
-&-\\-&-\end{array}\right)$}.
\end{equation}
 In particular the type of $(U^+DU^-D)^k$ is $(-1)^kP$. Our problem now is reduced
to the verification that in each of the four cases the top right 
entry is sign definite. Since $DP$ and $PD$ are sign definite,  this takes care of  cases one and four. 
Observe that $DU^-D= \left( \begin{array}{cc} 
-&-\\0&-\end{array}\right)$, so $DU^-DP$ is sign definite, and so the claim in case three follows. 
Finally, $DPU^+D$ is sign definite, so the claim in case two follows.

Let us prove the statement for a 
closed loop $l$. In this case ${\Bbb I}'_{\cal X}(l)$ is 
the trace of a cyclic word of type $... DU^{\pm}DU^{\mp}DU^{\pm}DU^{\mp}DU^{\pm}$. 
We will assume that its length is minimal possible for a given loop $l$. 
So the following cases may occur: 

i) The word is $U^+$, $U^-$ or $D$. These are trivial cases.

ii) The cyclic word is $DU^+$ or $DU^-$. One easily checks that 
in this case  
the trace is sign-definite.

iii) The cyclic word has length greater than two. Then the length is even: 
otherwise the length of the cyclic word can be reduced. If the length is bigger then 
$2$ and not divisible by $4$, there must be a segment $...DU^{\pm}DU^{\pm}\!..$. 
Therefore the sequence can be made shorter. If the length is two, the claim follows from ii). 
Finally, if the length is $4k$, we can write it as 
a product of $k$ segments of type $DU^+DU^-$, each of which is sign definite by (\ref{3.19.04.10}).

\begin{figure}[hc]
\centerline{\epsfbox{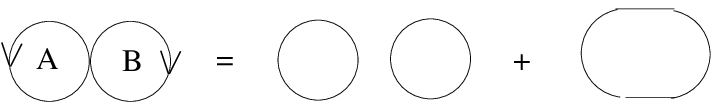} \qquad \qquad \qquad \epsfbox{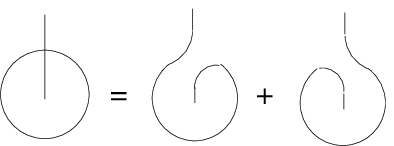}}
\caption{Skin relations.}
\label{fg87}
\end{figure}

Let us prove the property \ref{multiplication}. Let us spell out the proof in the more complicated case 
of ${\cal X}$-laminations. 
For this we  introduce the notion of a {\em quasi-lamination}.\!
A quasi-lamination is \!a collection of curves, connecting 
punctures or closed. Unlike the laminations, 
they may intersect or have self-intersections. 
One can define the canonical multiplicative pairings for single 
curve  
quasi-laminations by the same rules as for laminations, and extend it to the
collections of curves by multiplicativity, i.e., by definition 
${\Bbb I}_*(l_1){\Bbb  I}_*(l_2) = {\Bbb  I}_*(l)$, where 
$l$ is a quasi-lamination given by 
the union of $l_1$ and $l_2$. The functions on the moduli spaces 
corresponding to  quasi-laminations satisfy certain 
linear {\em skein relations}, described below. We will show that using
these relations one can represent any function corresponding to 
 a quasi-lamination 
as a positive integral linear combination of similar functions corresponding to 
 laminations, thus proving the property. Let us describe the 
relations. 
A lamination in the vicinity of an intersection or self-intersection point 
looks like a cross
$
\setlength{\unitlength}{0.15ex}
\begin{picture}(20,20)(-10,-5)
\put(-10,0){\line(1,0){ 20}}
\put(0,-10){\line(0,1){ 20}}
\end{picture}
$
.
Let $l_1$ and $l_2$ be the quasi-laminations obtained by replacing this cross by
$
\setlength{\unitlength}{0.15ex}
\begin{picture}(20,20)(-10,-5)
\put(10,10){\oval(20,20)[bl]}
\put(-10,-10){\oval(20,20)[tr]}
\end{picture}
$
and
$
\setlength{\unitlength}{0.15ex}
\begin{picture}(20,20)(-10,-5)
\put(10,-10){\oval(20,20)[tl]}
\put(-10,10){\oval(20,20)[br]}
\end{picture}
$
respectively. The relation is
$$
{\Bbb I}_*(l, m)\pm {\Bbb I}_*(l_1, m)\pm {\Bbb I}_*(l_2, m)=0,
$$
for some choice of signs depending on $l,l_1,l_2$. To prove this, and determine the signs, 
we consider the following three cases: (i) an intersection of two loops; (ii) an intersection of a path connecting 
two punctures, each puncture equipped with a decoration as a vector $v_i$, and a loop; (iii) 
an intersection of two paths connecting decorated punctures. 
The proof in the cases (i) and (ii) follows from the two pictures on  Figure \ref{fg87}, and the following 
two identities for $2\times2$ matrices with determinant $1$ 
 and a pair of vectors $v_1, v_2$ in a two dimensional vector space:
\begin{equation} 
{\rm tr} A{\rm tr} B ={\rm tr} AB +{\rm tr} AB^{-1}, \qquad 
(v_1\wedge v_2)\cdot {\rm tr}A = v_1 \wedge Av_2 + v_1 \wedge A^{-1}v_2.
\end{equation}
The second identity follows from ${\rm tr}A \cdot {\rm Id} = A + A^{-1}$. 
The proof in the case (iii) reduces to the  Pl\"ucker identity 
$\Delta(v_1, v_3) \Delta(v_2, v_4) = \Delta(v_1, v_2) \Delta(v_3, v_4) + 
\Delta(v_1, v_4) \Delta(v_2, v_3)$.  
\begin{figure}[ht]
\centerline{\epsfbox{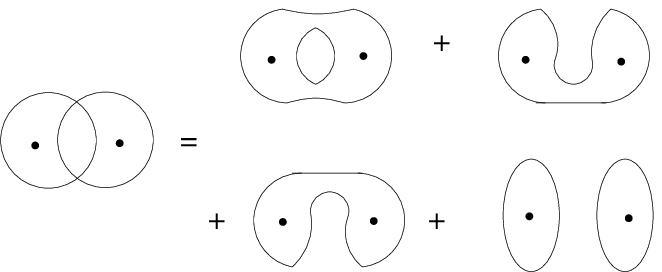}\qquad \qquad \qquad   \epsfbox{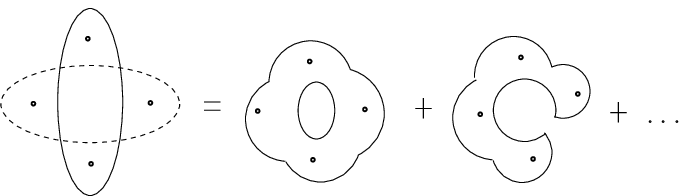}}
\caption{Skin relations producing contractible loops.}
\label{fg95}
\end{figure}
Since $l_1$ as well as $l_2$ have one intersection point less than $l$, repeating this procedure
one can express the pairing with any quasi-la\-mi\-na\-tion by a linear combination of pairings with
quasi-laminations without intersection points. Now to reduce the pairing to 
the pairing with a lamination
we need to replace equivalent curves by a single curve with a multiplicity. 
Observe that ${\rm tr} A^k=P_k({\rm tr} A)$,
where $P_k$ is the $k$-th Tchebychev polynomial. 
Furthermore, in the case of ${\cal X}$-laminations 
each simple contractible curve which might appear in the process 
contributes a factor $-2$ to the pairing. Indeed, lifting such a curve to the punctured tangent bundle 
of $S$ we get a generator of the center of $\pi_1$, 
so the trace of the monodromy of a twisted $SL_2$-local system along it is $-2$. 
Thus the  sum obtained obviously has integral coefficients. 
The contractible curves are  the only source of negative coefficients, so we have to show 
that,  nevertheless, all coefficients are positive. 
For example, consider the skin relation shown on the left of Figure \ref{fg95}. Contracting 
the loop inside of the 
the first lamination on the right of the equality sign, 
we get the outer curve with the coefficient $-2$, which cancels with the 
next two curves, which are isotopic to that curve. To prove the claim in general we may assume, without 
 loss of generality, that our quasi-la\-mi\-na\-tion has two components, and thus 
a simple contractible loop appears as a $2n$-gon, $n \geq 1$, obtained by their intersection. 
Arguing as  above (see the right picture on Figure \ref{fg95}) we see that contracting the loop we get 
a curve  with coefficient 
$2n-2$. The property 5 of Theorem \ref{10.10.03.2} for the ${\cal X}$-laminations 
is proved. The case of ${\cal A}$-laminations is similar but simpler: we do not have to worry about 
contractible loops. Theorem \ref{10.10.03.2} is  proved. 

\vskip 3mm
{\bf 4. Conjectures}. 
 Let $L$ be a set. Denote by $\Z_+\{L\}$ the abelian 
semigroup generated by $L$. Its elements are expressions 
$\sum_i n_i \{l_i\}$, where $l_i \geq 0$, the sum is finite, and 
$\{l_i\}$ is the generator corresponding to an element $l_i \in L$. 
Similarly denote by $\Z\{L\}$ the abelian group generated by $L$.

Let ${\cal X}$ be a positive space. 
Recall from Section 4.3 the ring ${\Bbb L}[{\cal X}]$ of 
all good Laurent polynomials,   
the semiring ${\Bbb L}_+[{\cal X}]$ of positive good Laurent polynomials and 
the set of extremal elements 
${\bf E}({\cal X})$.

It follows from Definition \ref{17.27} and 
Theorem \ref{10.10.03.2} that the 
canonical multiplicative pairings 
provide canonical maps
\begin{equation} \label{10/9/03/1a}
{\Bbb I}^+_{\cal A}: \Z_+\{{\cal A}_{SL_2, S}(\Z^t)\} \lra 
{\Bbb L}_+[{\cal X}_{PSL_2, S}], \quad 
{\Bbb I}^+_{\cal X}: \Z_+\{{\cal X}_{PSL_2, S}(\Z^t)\} \lra {\Bbb L}_+[{\cal A}_{SL_2, S}], 
\end{equation}
and hence 
\begin{equation} \label{10/9/03/2a}
{\Bbb I}_{\cal A}: \Z\{{\cal A}_{SL_2, S}(\Z^t)\} \lra 
{\Bbb L}[{\cal X}_{PSL_2, S}], \quad 
{\Bbb I}_{\cal X}: \Z\{{\cal X}_{PSL_2, S}(\Z^t)\} 
\lra {\Bbb L}[{\cal A}_{SL_2, S}], 
\end{equation}
which are multiplicative in the sense that property 4. of 
Definition \ref{17.27} holds. 
Conjecturally the polynomials provided by laminations in both cases are {\em minimal},
i.e., any difference of such polynomials have at least one
negative coefficient. Moreover, we conjecture the following.

 \begin{conjecture} \label{10.10.03.11}
{\it The maps (\ref{10/9/03/1a}) and (\ref{10/9/03/2a}) are isomorphisms. They provide the isomorphisms 
$$
{\cal A}_{SL_2, S}(\Z^t) = {\bf E}({\cal X}_{PSL_2, S}), 
\qquad {\cal X}_{PSL_2, S}(\Z^t) = {\bf E}({\cal A}_{SL_2, S}).
$$}
\end{conjecture}

\vskip 3mm
In particular 
this means that, in both  ${\cal A}$ and ${\cal X}$ cases,
 any good Laurent polynomial is a linear 
combination with 
integer coefficients of the images of the corresponding type laminations under the canonical map. For positive good Laurent polynomials all the 
coefficients are positive.

Observe that if $G$ is simply-connected then the Langlands dual 
 $^L G$ always has trivial center, so the corresponding moduli space 
${\cal X}_{^L G, S}$ has the positive structure defined in Section 6.

\begin{conjecture} \label{10.10.03.10}
{\it Let $G$ be a connected, simply-connected, split semi-simple algebraic group. Let $S$ be an open surface. 
Then 

a) There exist canonical additive pairings between the 
real lamination and Teichm\"uller spaces
$$
I:{\cal A}_{G,   S}(\R^t)\times {\cal X}^+_{^L G,   S}
\rightarrow {\mathbb R}, \quad 
I: {\cal X}_{^L G,   S}(\R^t) \times {\cal A}^+_{G,   S}
\rightarrow {\mathbb R}, \quad 
$$
as well as the intersection pairing between the real lamination spaces
$$
{\cal I}: {\cal A}_{G,   S}(\R^t)\times {\cal X}_{^L G,   S}(\R^t) 
\rightarrow {\mathbb R}.
$$
All of them are continuous. 

b) There exist canonical isomorphisms 
\begin{equation} \label{10/9/03/3}
{\Bbb I}^+_{\cal A}: \Z_+\{{\cal A}_{G,   S}(\Z^t)\} \lra 
{\Bbb L}_+[{\cal X}_{^L G,   S}], \qquad 
{\Bbb I}^+_{\cal X}: \Z_+\{{\cal X}_{^L G,   S}(\Z^t)\} 
\lra {\Bbb L}_+[{\cal A}_{G,   S}], 
\end{equation}
\begin{equation} \label{10/9/03/2}
{\Bbb I}_{\cal A}: \Z\{{\cal A}_{G,   S}(\Z^t)\} \lra 
{\Bbb L}[{\cal X}_{^L G,   S}], \qquad 
{\Bbb I}_{\cal X}: \Z\{{\cal X}_{^L G,   S}(\Z^t)\} 
\lra {\Bbb L}[{\cal A}_{SL_2,   S}]. 
\end{equation}
They provide the isomorphisms 
$$
{\cal A}_{G,   S}(\Z^t) = {\bf E}({\cal X}_{^L G,   S}), 
\qquad {\cal X}_{^L G,   S}(\Z^t) = {\bf E}({\cal A}_{G,   S}). 
$$}
\end{conjecture}

 Since  $^L SL_2 = PSL_2$, 
Conjecture \ref{10.10.03.10} reduces to 
Conjecture \ref{10.10.03.11}.

\vskip 3mm
{\em Quantum canonical map}. We use notations from \cite{FG2}. 
So  ${\cal X}^q_{PSL_2, S}$ denotes the $q$-deformation 
of the moduli space ${\cal X}_{PSL_2, S}$, defined by gluing 
the quantum tori given by the  generators ${X_i}$ and relations
$$
q^{{- \varepsilon_{ij}}}{X_i}{X_j}=q^{-
{\varepsilon_{ji}}}{X_j}{X_i}.
$$

\begin{conjecture} \label{10.10.03.10fgf}
{\it There exists quantum canonical map 
$$
\widehat {\Bbb I}^q: {\cal A}_{SL_2, S}(\Z^t) \lra {\cal X}^q_{PSL_2, S}.
$$
So for any positive 
integral {\cal A}-lamination $l$ one can associate
a Laurent polynomial $\widehat{\Bbb I}^q(l)$ of the variables 
${X_1}, \ldots, {X_n}$, satisfying
the following properties:
\begin{enumerate}
\item $\widehat{\Bbb I}^1(l)={\Bbb I}_{\cal A}(l)$
\item The highest term of $\widehat{\Bbb I}^q(l)$ is $q^{-\sum_{i<j}
{\varepsilon}_{ij}a_ia_j}\prod_i X_i^{a_i}$.
\item All coefficients of $\widehat{\Bbb I}^q(l)$ are 
positive Laurent polynomials of $q$.
\item Let $\ast$ be the canonical involutive antiautomorphism on the fraction field 
$\Q({\cal X}^q_{PSL_2, S})$. Then 
$\ast \widehat{\Bbb I}^q(l)=\widehat{\Bbb I}^q(l)$.
\item $\widehat{\Bbb I}^q(l_1)\widehat{\Bbb I}^q(l_2)=\sum_l c^q(l_1, l_2;l)
 \widehat{\Bbb I}^q(l)$, where the $c^q(l_1, l_2;l)$
are Laurent polynomials of $q$ with positive integral coefficients and  the sum is finite.
\end{enumerate}}
\end{conjecture}
\vskip 3mm

Of course, we expect a similar conjecture to be valid for an arbitrary $G$. 
The most general and precise conjectures are formulated in the setup of 
cluster ensembles, see Section 4 of \cite{FG2}.  We expect that the positive 
atlases defined in this paper provide the same rings of good positive Laurent polynomials as the ones 
provided by the cluster structure.  Therefore
the corresponding conjectures describe the same  objects.

There exists a generalization of these 
 conjectures 
where $S$ is replaced by $\widehat S$. In the case of $G=SL_2$ the 
canonical pairings in this set up were defined in [FG4].

Finally, there is a (partial) generalisation of these 
 conjectures when $S$ is a closed surface. In this case the 
space of integral laminations for a group $G$ with trivial center was defined in 
Section 6.9. We conjecture that it parametrises 
a canonical basis on the quantum space ${\cal A}_{^L G,S}$: for a closed $S$ the 
 space ${\cal A}_{^L G,S}$ has a symplectic structure, and can be quantised. 
For $G=PGL_2$ we get Turaev's algebra.

\vskip 3mm
{\bf 5. Bases in the space of functions on moduli spaces which are parametrised by laminations}. 
Recall that there are two versions of the spaces of ${\cal A}$-laminations: 
$$
{\cal A}_{SL_2, S}(\Z^t) \subset {\cal A}_L(S, \Z).
$$
Let ${\cal O}(X)$ be the algebra of regular functions on a stack $X$. Below we prove the following result

\begin{theorem} \label{10.6.05.1}
(i) The functions ${\Bbb I}_{\cal A}(l)$, when $l\in {\cal A}_{SL_2, S}(\Z^t)$, provide a basis in 
${\cal O}({\cal X}_{PGL_2, S})$.

(ii) The functions ${\Bbb I}_{\cal A}(l)$, when $l\!\in\! {\cal A}_{L}(S, \Z)$, provide a basis in 
${\cal O}({\cal X}_{SL_2, S})$.
\end{theorem}

{\bf Remark}. The set ${\cal A}_L(S, \Z)$ should be thought of as the space of 
integral $PGL_2$ ${\cal A}$-laminations: Indeed, it  provides a basis on the space 
${\cal X}_{SL_2, S}$. 
\vskip 3mm
Let us describe the stack of $G$-local systems on a graph for a reductive algebraic group $G$.  

Let us start from an arbitrary group  $G$. 
Let $\Gamma$ be any graph. 
Recall the set ${\cal L}(G,\Gamma)$ of $G$-local systems  on $\Gamma$, understood as local systems of right principal 
$G$-bundles. Let $E(\Gamma)$ be the set of all edges, and $V(\Gamma)$ the set of all vertices of $\Gamma$. 

\begin{lemma} \label{10.6.05.1sd} Given an orientation ${\rm Or}$ of the graph $\Gamma$, 
there exists a natural bijection
\begin{equation} \label{12.7.04.1}
\kappa_{\rm Or}: {\cal L}(G,\Gamma) \stackrel{\sim}{\lra} G^{E(\Gamma)}/G^{V(\Gamma)}.
\end{equation} 
\end{lemma}

{\bf Proof}. 
Let us define the right action of the group 
$G^{V(\Gamma)}$ on the set $G^{E(\Gamma)}$, which we use in 
(\ref{12.7.04.1}). It is the only thing in (\ref{12.7.04.1}) depending on 
an orientation of $\Gamma$. 
Let $(v_1(E), v_2(E))$ be vertices of the edge $E$, ordered according to the orientation of $E$. 
Then an element  $\{g_v\} \in G^{V(\Gamma)}$, where $v\in V(\Gamma)$, 
acts on an element $\{g_E\} \in G^{E(\Gamma)}$, where $E\in E(\Gamma)$, 
by the formula
\begin{equation} \label{12.7.04.2}
\{g_v\}\ast \{g_E\} = \{g'_E\}, \quad \mbox{where} \quad g'_E = g^{-1}_{v_2(E)}g_Eg_{v_1(E)}
 \end{equation} 
Let ${\cal L}$ be a $G$-local system on $\Gamma$. Let us trivialize its fibers at every vertex of 
 $\Gamma$. Then the parallel transport along an oriented edge $E$ of $\Gamma$ 
is uniquely described by an element $g_E\in G$. 
The collection $\{g_E\}$ of these 
elements provides a point of $G^{E(\Gamma)}$. 
Changing trivializations at vertices 
amounts to the action (\ref{12.7.04.2}) of the 
group $G^{V(\Gamma)}$. The lemma follows. \vskip 3mm
\vskip 3mm

\begin{corollary-definition} \label{12.7.04.4} Let $G$ be an affine algebraic group. Then 
the space of regular functions on the stack 
${\cal L}(G,\Gamma)$ is given by the space of $G^{V(\Gamma)}$-invariants
\begin{equation} \label{12.7.04.3}
{\cal O}\left({\cal L}(G,\Gamma)\right) = 
\Bigl({\cal O}(G)^{E(\Gamma)}\Bigr)
^{G^{V(\Gamma)}}.
\end{equation} 
\end{corollary-definition}
\vskip 3mm 

Now let  $G$ be a reductive group. Let us work out a convenient description of the space (\ref{12.7.04.3}). 
The Peter-Weyl 
theorem implies an isomorphism
$$
{\cal O}(G) \stackrel{\sim}{\lra}\bigoplus_{\lambda\in \widehat G}V_{\lambda}\bigotimes V_{\lambda}^*,
$$
where the summation is over the set $\widehat G$ of isomorphism classes of 
finite dimensional irreducible algebraic 
$G$-modules, and $V_{\lambda}$ is a representation space 
corresponding to $\lambda$. 
 Combining with (\ref{12.7.04.3}), we get an isomorphism
$$
{\cal O}\left({\cal L}(G,\Gamma)\right) \stackrel{\sim}{\lra} \Bigl(\bigoplus_{\{\lambda: E(\Gamma) \to \widehat G\}}\bigotimes_{E \in E(\Gamma)}
V_{\lambda(E)}\bigotimes V_{\lambda(E)}^*\Bigr)^{G^{V(\Gamma)}},
$$
where the summation is over the set of all functions $\lambda: E(\Gamma) \to \widehat G$. 

Let us rewrite this formula using invariants of a tensor product over the set $F(\Gamma)$ 
of flags of $\Gamma$. Say that a {\it $\widehat G$-coloring} of $\Gamma$ is a map 
$\lambda: F(\Gamma) \to \widehat G$ which has the following property: 
for every edge $E$ of the graph it 
assigns contragradient representations to the two flags assigned to $E$.

We assign the factors of the product $V_{\lambda(E)}\otimes V_{\lambda(E)}^*$  
to the two flags of the edge $E$, so that $V_{\lambda(E)}$ is assigned to the flag $(v_1(E), E)$, and 
$V^*_{\lambda(E)}$ to the flag $(v_2(E), E)$. Then  we get an isomorphism
\begin{equation} \label{12.7.04.7}
{\cal O}\left({\cal L}(G,\Gamma)\right) \stackrel{\sim}{\lra}
 \Bigl(\bigoplus_{\mbox{$\lambda\in  \{\widehat G$-colorings of $\Gamma$}\}}\bigotimes_{F \in F(\Gamma)}V_{\lambda(F)}\Bigr)^{G^{V(\Gamma)}}.
\end{equation} 
Collecting the flags corresponding to a given vertex $v$ of the graph, we get the desired formula
\begin{equation} \label{12.7.04.8}
{\cal O}\left({\cal L}(G,\Gamma)\right) \stackrel{\sim}{=} \bigoplus_{\mbox{$\lambda\in  \{\widehat G$-colorings of 
$\Gamma$}\}}\bigotimes_{v \in V(\Gamma)}\Bigl(\bigotimes_{\{E\to v\}}V_{\lambda(v, E)}\Bigr)^{G}. 
\end{equation} 
Here the set $\{E\to v\}$ consists of all edges of $E$ incident to the given vertex $v$, and $(v, E)$ means the flag 
at the vertex $v$ corresponding to such edge $E$.

\vskip 3mm
Our next step is to describe the space ${\cal O}({\cal X}_{G,S})$ of regular functions on ${\cal X}_{G,S}$ 
as an ${\cal O}({\cal L}_{G,S})$-module. 
Let $\widetilde G$ be Grothendieck's simultaneous resolution of $G$. Then 
there are isomorphism of stacks $\widetilde G/{\rm Ad}G = H$ and $G/{\rm Ad}G = H/W$. 

Given a puncture $p$ and an edge $E$ surrounding this puncture, there is a map $G^{E(\Gamma)} \lra G$ provided by the monodromy around the puncture, starting at $E$. Thus choosing edges 
for each of the $n$ punctures we get a map $G^{E(\Gamma)} \to G^{n}$. 
Let $ {{\cal X}'_{G,S}}$ be the fibered product in 
the following Cartesian diagram of varieties: 
$$
\begin{array}{ccc}
{{\cal X}'_{G,S}} & \lra & {\widetilde  G}^{n}\\
\pi' \downarrow & & \downarrow \\
G^{E(\Gamma)} & \lra & G^{n}
\end{array}
$$
The group  $G^{V(\Gamma)}$ acts on it, as explained in the proof of Lemma \ref{10.6.05.1sd}.  
Taking the quotient by the action of $G^{V(\Gamma)}$ we arrive at a Cartesian square of stacks, 
which is independent of the above choices:
$$
\begin{array}{ccc}
{\cal X}_{G,S} & \lra & H^{n}\\
\pi \downarrow & & \downarrow \\
{\cal L}_{G,S} & \lra & (H/W)^{n}
\end{array}
$$
It follows that 
$$
{\cal O}({\cal X}_{G,S}) =  {\cal O}({\cal L}_{G,S}) \otimes_{{\cal O}((H/W)^{n})} {\cal O}(H^{n})
$$
By Chevalley's theorem $\Q[H]$  is a free $\Q[H]^W$-module with $|W|$ generators. 
Thus the space ${\cal O}({\cal X}_{G,S})$ of regular functions on ${\cal X}_{G,S}$ 
is a free ${\cal O}({\cal L}_{G,S})$-module of rank $|W|^n$. A set of  generators is obtained by 
pull backs of 
generators of the $\Q[(H/W)^{n}]$-module $\Q[H]^n$. 
\vskip 3mm

Now let us assume that $G = SL_2$ and $\Gamma$ is a trivalent graph 
homotopy equivalent to $S$.

Let ${\cal A}^0_L(S, \Z) \subset {\cal A}_L(S, \Z)$ be the subset consisting of all integral 
${\cal A}$-laminations such that the weights of all curves of the lamination, including the boundary curves, 
 are positive integers. 
Given a loop $\alpha$ of a lamination $l$, we assign to it the weight 
$n_{\alpha}$, which equals to the number of loops in $l$ homotopic to $\alpha$. 
So we write a lamination $l$ as $l = \sum_{\alpha}n_\alpha \alpha$. Then $l \in {\cal A}^0_L(S, \Z)$ if and only if $n_\alpha >0$ for all $\alpha$. Given such a lamination $l$, set
\begin{equation} \label{A1/3}
{\Bbb I}^0_{\cal A}(l):= \prod_{\alpha}{\rm Tr}(M^{n_{\alpha}}_{\alpha}) \in {\cal O}({\cal L}_{SL_2,S}).
\end{equation}
Here the product is over all 
isotopy classes of loops in $l$.
Further, set $$
{\cal A}^0_{SL_2, S}(\Z^t):= {\cal A}_{SL_2, S}(\Z^t)\cap {\cal A}^0_L(S, \Z).
$$ 
So we have ${\cal A}^0_{SL_2, S}(\Z^t) \subset {\cal A}^0_L(S, \Z)$. 
\vskip 3mm
\begin{proposition} \label{A1/2} The map ${\Bbb I}^0_{\cal A}$ provides 
isomorphisms
$$
{\Bbb I}^0_{\cal A}: \Q\{{\cal A}^0_L(S, \Z)\} \stackrel{\sim}{\lra} {\cal O}({\cal L}_{SL_2,S}); \qquad 
{\Bbb I}^0_{\cal A}: \Q\{{\cal A}^0_{SL_2, S}(\Z^t)\} \stackrel{\sim}{\lra} {\cal O}({\cal L}_{PGL_2,S}).
$$
\end{proposition}
In other words, the functions $\{{\Bbb I}^0_{\cal A}(l)\}$, where $l \in 
{\cal A}^0_L(S, \Z)$ (respectively $l \in {\cal A}^0_{SL_2, S}(\Z^t)$), 
form a basis of ${\cal O}({\cal L}_{SL_2,S})$ (respectively ${\cal O}({\cal L}_{PGL_2,S})$). 
\vskip 3mm
{\bf Proof}. 
The list of all irreducible 
representations of $SL_2$ is given by the representations $V_{n/2}$ in 
the space of degree $n$ polynomials 
in two variables, $n \geq 0$. One has $V_{n/2} = S^nV_{1/2}$. 
Recall that 
\begin{equation} \label{12.7.04.9}
{\rm dim}(V_a\otimes V_b \otimes V_c)^{SL_2} = 
\left\{ \begin{array}{ll} 1 & \mbox{ if $a+b+c\in \Z$} \\
                            & \mbox{ and $a,b,c$ satisfy the triangle inequalities},\\
0& \mbox{ otherwise}.
\end{array}\right.
\end{equation} 
Recall that  laminations $l \in {\cal A}^0_L(S, \Z)$ are described by 
collections of non-negative half-integers $\{a_E\}$ on the edges of $\Gamma$, such that 
for every vertex 
the sum of the three numbers at the edges incident to the vertex is an integer, and these numbers satisfy the triangle inequalities. Such data  gives rise to the following two different functions on the moduli space 
${\cal L}_{G, S}$:

i) We assign to each edge $E$ of $\Gamma$ the irreducible representation $V_{a_E}$. 
Then thanks to the conditions on $\{a_E\}$ and (\ref{12.7.04.9}), 
there is a unique (up to a scalar) 
invariant $S_{a,b,c}$ at every vertex. Taking the tensor product of these invariants, we get 
a function $S_{\{a_E\}}$. (Here $S$ stands for ``symmetric'': $V_{n/2}= S^nV_{1/2}$).

ii) We assign to each edge $E$ of $\Gamma$ the reducible representation 
$\otimes^{2a_E}V_{1/2}$. Then given a triple of non-negative half-integers 
$(a,b,c)$ satisfying the triangle inequalities, and with an integral sum, 
we define an invariant in the triple tensor product 
\begin{equation} \label{12.7.04.11}
\otimes^{2a}V_{1/2} \otimes \otimes^{2b}V_{1/2} \otimes \otimes^{2c}V_{1/2} 
\end{equation} 
as follows. Let us draw an ${\cal A}$-lamination near the vertex 
$v$ whose coordinates at the edges are $2a, 2b, 2c$. Let us assign a standard 
representation $V_{1/2}$ to both ends of each arc of the lamination. 
Then we assign to each arc of the lamination the invariant in the tensor product 
$V_{1/2} \otimes V_{1/2}$ of the two representations assigned to the ends of this arc. 
Then taking the tensor product we get an $SL_2$-invariant 
$G_{a,b,c}$ in (\ref{12.7.04.11}). Taking the tensor product 
of these invariants over all vertices of $\Gamma$, we get  a function $T_{\{a_E\}}$ on the moduli space. 
(Here $T$ stands for ``tensor'').

The functions $T_{\{a_E\}}$ are linear combinations of the $S_{\{a_E\}}$, and vice versa. Indeed, this follows from the fact that the invariants $\{T_{a,b,c}\}$ 
are linear combinations of the ones $\{S_{a,b,c}\}$, and vice versa. 
(In fact they 
 are related by upper triangular transformations.) Since the 
functions $S_{\{a_E\}}$ form a basis 
in the space of functions on  ${\cal L}_{G, S}$, the 
functions $T_{\{a_E\}}$  also form a basis.  

The basis $T_{\{a_E\}}$
has a more geometric description, which does not 
rely on the graph $\Gamma$. Namely, given 
a set of half-integral coordinates $\{a_E\}$ on the edges of $\Gamma$, 
with an integral sum and satisfying the triangle inequalities,  
let us construct the corresponding  ${\cal A}$-lamination 
$l$ on $S$.   
Then we assign to each loop $\alpha$ of the resulting lamination the trace 
${\rm Tr}M_{\alpha}$ 
of the monodromy $M_{\alpha}$
of the local system around this loop, and take 
the product of obtained functions over all loops of the lamination. 
This way we get a yet another collection of functions $T'_{\{a_E\}}$.

The following lemma follows immediately from the definitions.

\begin{lemma} \label{12.8.04.1}
The functions $T_{\{a_E\}}$ and $T'_{\{a_E\}}$ coincide, up to a non-zero scalar.
\end{lemma}

Let us compare the previous construction with (\ref{A1/3}). It gives instead 
of the product of ${\rm Tr}(M_{\alpha}^{n_{\alpha}})$ the product of $\Bigl({\rm Tr}M_{\alpha}\Bigr)^{n_{\alpha}}$. However 
since both the Laurent polynomials 
$$
(\lambda + \lambda^{-1})^n \qquad \mbox{and} \quad \lambda^n + \lambda^{-n}, \quad n \in \Z_{\geq 0}
$$ 
form a basis in the space of all Laurent polynomials in $\lambda$ 
invariant under the involution $\lambda \lms \lambda^{-1}$, the functions 
${\Bbb I}^0_{\cal A}(l)$ are linear combinations of the  $T_{\{a_E\}}$ 
(and vice versa), and hence also form a basis.  

In the case when $G=PGL_2$ the proof is similar. The only difference is that 
the list of irreducible representations of $PGL_2$ is given by the representations 
$V_a$, where $a \in \Z_{\geq 0}$, and thus in (\ref{12.7.04.9}) we can drop the condition 
that $a+b+c \in \Z_{\geq 0}$. The proposition is proved. 
\vskip 3mm

Now we can complete the proof of Theorem \ref{10.6.05.1}. We start from  claim i). We claim that 
the products of ${\Bbb I}_{\cal A}(\alpha_i)$ over some of the boundary components provide all 
generators of the ${\cal O}({\cal L}_{PGL_2,S})$-module ${\cal O}({\cal X}_{PGL_2,S})$.
Observe that 
${\Bbb I}^0_{\cal A}(l)$ differs from ${\Bbb I}^0_{\cal A}(l)$ only in the treatment 
of boundary curves. Precisely, let $\alpha_i$ be the boundary loop for the $i$-th boundary component. 
Let $\lambda_i$ be the monodromy around $\alpha_i$   
restricted to the invariant one-dimensional subspace provided by the framing. 
Then 
${\Bbb I}_{\cal A}(n \alpha_i) = \lambda_i^n + \lambda_i^{-n}$, while 
${\Bbb I}^0_{\cal A}(n \alpha) = \lambda_i^n$. 
So the theorem is reduced to the following obvious claim: The space $\Q[\lambda, \lambda^{-1}]$ is a free $\Q[\lambda, \lambda^{-1}]^{\Z/2\Z}$-module of rank two, 
with generators $1$ and $\lambda$. Here the generator of $\Z/2\Z$ acts by 
$\lambda \lms \lambda^{-1}$. 

The second claim of the theorem is proved in a completely similar way.
The theorem is proved. 

\vskip 3mm
{\bf 6. The action of the group $(\Z/2\Z)^n$}. 
This group 
acts by birational transformations on the moduli space ${\cal X}_{PGL_2, S}$. 
Its generator assigned to a  boundary component $p$ exchanges 
the two monodromy-invariant flags at $p$. 
\begin{lemma} \label{12.25.05.1} 
The group $(\Z/2\Z)^n$ 
acts on ${\cal X}_{PGL_2, S}$ by positive birational transformations.
\end{lemma}

{\bf Proof}. Let $x_1, ..., x_n$ be the coordinates at the edges entering a puncture $p$. 
The monodromy around $p$ is given by 
$\prod_{i=1}^n \left( \begin{matrix}x_i& 0\\ 1&1\end{matrix}\right ) = 
\left(\begin{matrix}x_1... x_n& 0\\ c&1\end{matrix}\right )$, 
where $c= 1+x_n+x_nx_{n-1} + \ldots + x_n ... x_2$. 
The two eigenvectors of this matrix are $[0,1]$ and $[\alpha, 1]$, where $\alpha = (x_1... x_n - 1)/c$. 
We need to calculate how the ${\cal X}$-coordinates change when we 
alter the framing at $p$ by switching 
the first eigenvector with the second.
 \vskip 3mm
\begin{figure}[ht]
\centerline{\epsfbox{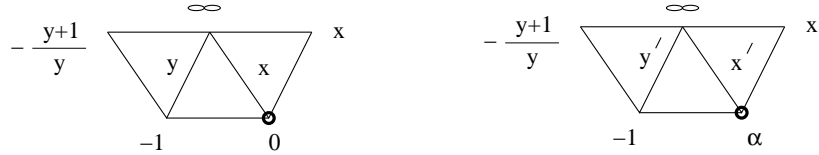}}
\caption{The action of $\Z/2\Z$ at a puncture (located at the fat vertex).}
\label{fg94}
\end{figure}
Observe that only the coordinates assigned to the edges of the triangles 
with a vertex at $p$ may change. On Figure \ref{fg94}  $x:=x_1$ and $y:=y_1$, and 
at vertices shown the points on ${\Bbb P}^1$ corresponding to the flags at the punctures. 
An easy computation shows that 
$$
y_1' = (1+\alpha)y_1 = \frac{(x_n+x_nx_{n-1} + x_n...x_1)y_1}{1+ x_n+x_nx_{n-1} + \ldots + x_n ... x_2} 
=\frac{[x_n, \ldots, x_1]}{[1, x_n, x_{n-1}, \ldots, x_2]}y_1,
$$
$$
x_1'= \frac{\alpha - x_1}{-1-\alpha} = \frac{cx_1+1 - x_1 ... x_n}{c-1+x_1...x_n} = 
$$
$$
=\frac{1+x_1+x_1x_n +  \ldots + x_1x_n...x_3}{x_n+x_nx_{n-1} + \ldots +x_n...x_1} 
=\frac{[1, x_1, x_n, \ldots, x_3]}{[x_n, x_{n-1}, \ldots, x_1]}
$$
where $[a_1, ..., a_m]:= a_1+ a_1a_2+ \ldots + a_1 ... a_m$. 
The lemma is proved.

\vskip 3mm
{\it The rational action of the group $(\Z/2\Z)^n$ on the moduli 
space ${\cal A}_{SL_2, S}$}. 
Denote by $A_{a, b}$ the $A$-coordinate assigned to an edge connecting punctures 
$a$ and $b$. Let $p$ be a puncture. The monodromy $M_p$ around $p$ is given by 
$$
M_p = \prod_{j=1}^q\left ( \begin{matrix}1& 0\\ \frac{A_{j, j+1}}{A_{p, j}A_{p, j+1}}&1\end{matrix}\right ) = 
\left ( \begin{matrix}1& 0\\\alpha(p) &1\end{matrix}\right ). \qquad 
 \alpha(p):= \sum_{j=1}^q \frac{A_{j, j+1}}{A_{p, j}A_{p, j+1}},
$$  
where the product is over the ordered set $\{1, \ldots , j\}$ of all 
edges entering the puncture $p$, so that $\{(p, j, j+1)\}$ is the set of triangles
 of the triangulation sharing $p$. 
It follows that $\alpha(p)$
does not depend on the choice of triangulation. Observe that 
 $\alpha(p) = \langle M_pv_p', v_p \rangle$, where 
$v_p$ is the $M_p$-invariant vector defining 
a decoration at $p$, $v_p'$ is a complementary vector, $\langle *, *\rangle$ the invariant area. 

Let us define an involution $\sigma_p$ 
of the ${\cal A}$-space. It acts by changing the decoration 
$v_p$ at the 
puncture $p$ by $v_p\lms  \alpha(p)v_p$, leaving untouched
 the other decorations and the twisted local system. The map $\sigma_p$ 
changes the invariant $\alpha(p)$ to its inverse:  $\sigma^*_p(\alpha(p)) = \alpha(p)^{-1}$. 
It follows that $\sigma^2_p = {\rm Id}$. The involutions $\sigma_p$ for different punctures commute, 
and thus give rise to an action of the group $(\Z/2\Z)^n$ on the moduli space ${\cal A}_{SL_2, S}$. 
The projection of this action to the ${\cal U}$-space is trivial.

Here is a slightly different way to see this action. Given a puncture $p$, 
the canonical projection $p: {\cal A}_{SL_2, S}\to 
{\cal U}_{SL_2, S}$ has a $1:2$ ``section''. It is given by choosing a multiple $\widetilde v_p$ of $v_p$ such that 
the corresponding invariant $\alpha(p)$ is $1$. Such a vector $\widetilde v_p$ is 
well defined up to a sign. The group ${\Bbb G}_m$ acts on ${\cal A}_{SL_2, S}$ 
by altering the decoration at $p$ by 
$\lambda \widetilde v_p \lms \lambda^{-1} \widetilde v_p$. 

\vskip 3mm
{\bf Remarks}. 1. For a point of ${\cal A}_{SL_2, S}(\R_{>0})$, the number $\alpha(p)$ 
is the area enclosed by the decorating 
horocycle at $p$. 

2. If $S$ has $n>1$ holes, the group $(\Z/2\Z)^n$ acts by 
cluster transformations. Indeed, given a hole $p$, we can find 
a trivalent graph $\Gamma$ which has a virus subgraph whose circle encloses the hole. Then flip at the cricle is the cluster transformation providing the 
generator of $\Z/2\Z$ corresponding to $p$. If $S$ has just one hole, 
changing the frame vector at this hole is not 
a cluster transformation.

\vskip 3mm
{\it A definition of the map ${\Bbb I}_{\cal X}$ for all ${\cal X}$-laminations}. The group $(\Z/2\Z)^n$ acts by positive birational transformations of both moduli spaces 
${\cal A}_{SL_2, S}$ and ${\cal X}_{PGL_2, S}$.  The tropicalization of this action provides 
an action of the group $(\Z/2\Z)^n$ on the lamination spaces. 
In the ${\cal X}$-case (resp. ${\cal A}$-case) 
it coincides with the one given by changing orientations of the holes (changing the signs of the weights
 of 
the loops surrounding the holes - see (\ref{signA})). 
The map ${\Bbb I}_{\cal X}$ was defined only for ${\cal X}$-laminations with the standard orientations of the boundary components. 
We extend it to all ${\cal X}$-laminations by imposing 
the equivariance with respect to the action of the group $(\Z/2\Z)^n$. 
Since $\alpha(p)$ is a positive Laurent polynomial, 
its image consists of good positive Laurent polynomials. It is straightforward that the other 
claims of Theorem \ref{10.10.03.2} also remain valid for the extended map ${\Bbb I}_{\cal X}$. 

\vskip 3mm
{\bf Conclusion}. Combining this with 
Theorems \ref{10.10.03.2} and \ref{10.6.05.1}, we have defined 
the canonical maps ${\Bbb I}_{\cal A}$ and 
${\Bbb I}_{\cal X}$ and proved that they satisfy all 
the desired properties but one: we do not know 
that laminations give rise to indecomposable elements.

$$
$$

\section{Completions of Teichm\"uller spaces and canonical bases}
\label{fgolam1}

{\bf 1. Cyclic sets, laminations and cacti-sets}. 
Let $C$ be a cyclic set, perhaps infinite. 
A {\it simple lamination $l$ in $C$} is a disjoint collection, perhaps infinite,  of pairs 
of elements $(c_i, c_i')$ of $C$
so that for any $i \not = j$  the order of the quadruple  
$(c_i, c_j, c_i', c_j')$ is not  compatible with the cyclic order of $C$. (If $C= S^1$, this just means 
that the chords $(c_i,c_i')$ and $(c_j, c_j')$ do not intersect).

\begin{figure}[ht]
\centerline{\epsfbox{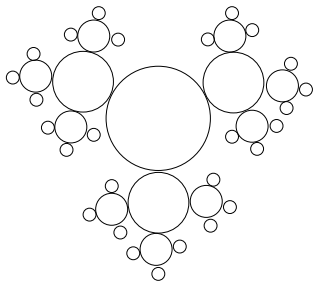}}
\caption{A cacti-set $C/l$. }
\label{fg10}
\end{figure}
Gluing every pair of points $c_i, c'_i$ 
we get a set $C/l$. It no longer has a cyclic structure compatible with the projection $\pi: C \to C/l$. 
A subset $M$ of $C/l$ is {\it cyclic} if it has a cyclic structure such that the projection 
$\pi^{-1}M \to M$ is a map of cyclic sets.  So a cyclic 
subset has a canonical cyclic structure inherited from $C$. 
Cyclic subsets are ordered by inclusion. The maximal elements 
are called the {\it maximal cyclic} subsets of $C/l$. 
The set $C/l$ is a union of its maximal cyclic subsets; any two of them 
are either disjoint, or have a single element intersection. We call such sets {\it cacti-sets}. 
\vskip 3mm
{\bf Example}. Let $C$ be a circle. A simple lamination $l$ in $C$ is given by the endpoints 
of a collection of disjoint 
chords inside the circle. Shrinking all chords to points we get a cacti-set 
$C/l$ -- see Figure \ref{fg10}. Its circles are the maximal cyclic subsets. 
\vskip 3mm

Let $S$ be a surface, with or without boundary. 
A {\it simple lamination} on $S$ is a finite collection of simple, non-contractible, 
mutually non-isotopic, and non-isotopic to boundary components,  
disjoint loops on $S$, considered  modulo isotopy.

A simple lamination $l$ on $S$ gives rise to a simple lamination in the cyclic set ${\cal G}_{\infty}(S)$. 
Indeed, let $\widetilde l$ be 
the preimage of $l$ in the universal cover $\widetilde S$ of $S$. The endpoints of $\widetilde l$ 
provide a lamination in 
the cyclic set ${\cal G}_{\infty}(S)$. 
Denote by 
${\cal G}_{\infty}(S,l)$ the corresponding cacti-set. 
The connected components of 
$\widetilde S-\widetilde l$ 
are in bijection with maximal cyclic subsets of ${\cal G}_{\infty}(S, l)$.
\vskip 3mm
{\bf 2. The classical case}. Let $\overline {\cal M}_{0, X}$ be the Knudsen-Deligne-Mumford moduli space 
of configurations of points on ${\Bbb C}{\Bbb P}^1$ parametrised by a finite set $X$. An embedding of sets 
$X \subset X'$ gives rise to
 a canonical projection  $\overline {\cal M}_{0, X'}\lra \overline {\cal M}_{0, X}$ -- 
forgetting the points parametrized by $X'-X$. 

Let $C$ be an infinite countable set. Set 
$$
\overline {\cal M}_{0, C}:= \lim_{\longleftarrow} \overline {\cal M}_{0, X}
$$ 
where the projective limit is over finite subsets $X$ of $C$, and the maps 
 correspond to inclusions of finite sets. 
(If $C$ is a finite set, we recover $\overline {\cal M}_{0, C}$).

Now let $C$ be a cyclic set. Then ${\cal M}_{0, C}$ is equipped with a positive atlas, so 
there is its real positive part ${\cal M}^+_{0, C}$. Its closure   in $\overline {\cal M}_{0, C}(\R)$
is denoted $\overline {\cal M}^+_{0, C}$. The set $\overline {\cal M}^+_{0, C}$ has a decomposition 
(into ``cells'' of {\em a priori} infinite dimensions) parametrized by laminations in $C$.

\vskip 3mm 
{\bf Example.} 
When $C$ is a finite cyclic set, $\overline {\cal M}^+_{0, C}$ is (the closure of) the Stasheff polytope, 
with its natural cell decomposition parametrised by the laminations in $C$ (\cite{GM}).
\vskip 3mm 

The points of the Teichm\"uller space ${\cal X}^+_{PGL_2, S}$ are given by positive 
$\pi_1(S)$-equivariant maps 
${\cal G}_{\infty}(S) \to {\Bbb P}^1(\R)
$. 
Thus they are points of the moduli space ${\cal M}^+_{0, C}$ for $C = {\cal G}_{\infty}(S)$.

\begin{definition} \label{y12}
The set $\overline {\cal X}^+_{PGL_2, S}$ is the closure of ${\cal X}^+_{PGL_2, S}$ 
in $\overline {\cal M}_{0, C}$ for $C = {\cal G}_{\infty}(S)$. 
\end{definition}

We need the following general definition. Given a simple lamination $l$ on $S$, 
the boundary components of the curve $S-l$ correspond either to the loops of the lamination $l$, 
or to the boundary components of $S$. The former are called $l$-boundary components of $S-l$.

\begin{definition} \label{y22}
Let $l$ be a simple lamination on $S$. The moduli space ${\cal X}_{G, S-l}^{{\rm un}}$ 
parametrizes framed $G$-local systems on $S-l$ 
with unipotent monodromies around the $l$-boundary components  of $S-l$. 
The moduli space  ${\cal A}_{G, S-l}^{{\rm un}}$ parametrizes  twisted unipotent 
$G$-local systems on $S-l$ with decorations 
at those  boundary components which are inherited from $S$. 
\end{definition}

So 
${\cal X}^{{\rm un}}_{G,S-l}$ is the preimage of the identity element under the 
projection $
{\cal X}_{G,S-l} \lra H^k$, 
given by the semi-simple parts of the monodromies around the $l$-boundary 
components  on $S-l$. 
The moduli space 
${\cal A}^{{\rm un}}_{G,S-l}$ is the quotient of ${\cal A}^{}_{G,S-l}$ by the action of the group $H^k$ 
provided by the $l$-boundary 
components  on $S-l$. Therefore each of them carries an induced  positive atlas, and thus there are the
 real positive parts  ${\cal X}_{G, S-l}^{+, {\rm un}}$ and ${\cal A}_{G, S-l}^{+, {\rm un}}$.

\begin{proposition} \label{y13}
$\overline {\cal X}^+_{PGL_2, S}$ 
is a disjoint union of cells parametrised by simple laminations in $S$:
$$
\overline {\cal X}^+_{PGL_2, S} = \coprod_l {\cal X}_{PGL_2, S-l}^{+, {\rm un}}; 
\qquad \mbox{where $l$ are simple laminations on $S$}.
$$
\end{proposition} 
An embedding of the right hand side to the left 
for any $G$ is defined  in Section 13.3.

Apparently the action of the mapping class group $\Gamma_S$ 
extends to the closure $\overline {\cal X}^+_{PGL_2, S}$. 
The quotient $\overline {\cal U}^+_{PGL_2, S}/\Gamma_S$ is identified with the Knudsen-Deligne-Mumford 
moduli space $\overline {\cal M}_S$, 
usually denoted $\overline {\cal M}_{g,n}$, where $g$ is the genus of 
$\overline S$ and $n$ is the number of punctures on $S$. 

\vskip 3mm
{\bf 3. Completions of higher Teichm\"uller spaces.} Recall that a positive 
configuration of flags parametrised by a cyclic set $C$ is the same thing as a positive map 
$C \to {\cal B}(\R)$ modulo $G(\R)$-conjugation. Recall the canonical projection 
$
\pi: {\cal G}_{\infty}(S) \to {\cal G}_{\infty}(S, l).
$ 

\begin{definition} \label{y15} Let $C$ and $C'$ be cyclic sets related by a 
surjective map of cyclic sets $\pi: C \to C'$. 

A sequence $\{\psi_n\}$ of positive configurations of flags 
parametrised by $C$ is convergent to 
a positive configuration $\psi$ parametrised by $C'$ 
if there are maps $\widetilde \psi_n: C \to {\cal B}(\R)$ 
representing configurations $\psi_n$ such that the limit 
$\widetilde \psi:= \lim_{n \to \infty}\widetilde \psi_n$ exists, factors through $C'$, i.e. 
$\widetilde \psi (c_i) = \widetilde \psi (c'_i)$ if  $\pi(c_i) = \pi(c_i')$, and the induced map 
 $\widetilde \psi: C' \to {\cal B}(\R)$ represents the configuration $\psi$. 
\end{definition}

\begin{definition} \label{y14} A sequence of points $\psi_n \in {\cal X}^+_{G, S}$ is convergent 
if there exists a simple lamination $l$ on $S$ such that for every maximal cyclic subset $\mu$ of 
${\cal G}_{\infty}(S, l)$ the sequence of positive configurations 
\begin{equation} \label{y16}
\psi_n^\mu: \pi^{-1}(\mu) \lra {\cal B}(\R),
\end{equation}
obtained by 
restrictions of $\psi_n$'s to 
$\pi^{-1}(\mu)$, is convergent to  a positive configuration  
$
\psi^\mu: \mu \lra {\cal B}(\R).
$
\end{definition}

Since the maps (\ref{y16}) correspond to points of ${\cal X}^+_{G, S}$, they are 
$\pi_1(S)$-equivariant. Therefore the limiting map $\psi^\mu$ is $\pi_1(S_\mu)$-equivariant, where 
$S_{\mu}$ is the  component of $S-l$ corresponding to $\mu$.  

Recall that there exists a unique regular unipotent conjugacy class in $G(\R)$: 
it is given by the maximal Jordan block for $G= PGL_n$. 

\begin{lemma} \label{y117} The monodromy of the local system on $S_\mu$ corresponding to $\psi^\mu$ 
around an $l$-boundary component is regular unipotent. 
\end{lemma}

{\bf Proof}. The monodromy $M$ around a boundary component is conjugate to an element of $B(\R_{>0})$. 
Since the canonical flags corresponding to $M$ and $M^{-1}$ coincide, $M$ is unipotent. 
Any element of $U(\R_{>0})$ is regular. The Lemma is proved. 

\vskip 3mm

Lemma \ref{y117} implies 
that a convergent sequence of points $\psi_n \in {\cal X}^+_{G, S}$ corresponds to a 
point of the space ${\cal X}^{+, {\rm un}}_{G, S-l}$, for a certain simple 
lamination $l$ on $S$, which we view as its limit. 
The gluing from Section 7 shows that we can get any point of 
${\cal X}^{+, {\rm un}}_{G, S-l}$. Similar questions for the ${\cal A}$-spaces 
are reduced to the ones for the ${\cal X}$-spaces. 
So we arrive at the following definition.

\begin{definition} Each of the sets  ${\cal B}({\cal X}^+_{G, S})$ and 
 ${\cal B}({\cal A}^+_{G, S})$ is a  disjoint unions 
of cells, which are parametrised by simple laminations $l$ in $S$:
$$
{\cal B}({\cal X}^+_{G, S}) = \coprod_l {\cal X}_{G, S-l}^{+, {\rm un}}; \qquad 
{\cal B}({\cal A}^+_{G, S}) = \coprod_l {\cal A}_{G, S-l}^{+, {\rm un}}, 
$$
where $l$ are simple laminations on $S$.

The topology on these sets is provided by Definition \ref{y14}. The obtained  topological space 
are called completions of the higher Teichmuller spaces ${\cal X}^+_{G, S}$ and ${\cal A}^+_{G, S}$. 
\end{definition}
We use here a different notation, comparing to Definition \ref{y12}, 
to emphasize that for $G \not = PGL_2$ we do not know 
whether the space ${\cal B}({\cal X}^+_{G, S})$ is the closure of the space ${\cal X}^+_{G, S}$ in any sense. 

\vskip 3mm
Here is how one should define a full completion of the higher Teichm\"uller space for a group $G$ different from $PSL_2$. 
One should be able to define a natural compactification $\overline {\rm Conf}_n({\cal B})$ of 
the space of configurations of $n$ flags in generic position, so that for $G=PSL_2$ we recover the 
moduli space $\overline {\cal M}_{0,n}$. Then we would define a completion of the space ${\cal L}^+_{G,S}$ 
as the closure $\overline {\rm Conf}^+_{{\cal G}_{\infty}(S), \pi_1(S)}({\cal B})$ of the positive $\pi_1(S)$-equivariant 
configuration space 
${\rm Conf}^+_{{\cal G}_{\infty}(S), \pi_1(S)}({\cal B})$ in the above compactification.

\vskip 3mm
{\bf 4. The canonical map for closed surfaces, and its behavior at the boundary}. 
For a closed surface $S$ the moduli space 
${\cal A}_{G,S}$ is isomorphic to the moduli space of $G$-local systems on $S$, 
although this isomorphism is non-canonical if $s_G \not = e$. 
So it has a canonical symplectic structure. When $G$ has trivial center, 
the set ${\cal X}_{G,S}(\Z^t)$ was defined in Section 6.9.

\begin{conjecture} \label{y8} {\it Let $S$ be a surface without  boundary. Assume that $G$ has trivial center.  
Then there exists a canonical map 
\begin{equation} \label{y10}
{\Bbb I}: {\cal X}_{G,S}(\Z^t) \lra {\cal O}({\cal A}_{^L G,S}).
\end{equation}
Its image is a basis of ${\cal O}({\cal A}_{^L G,S})$.
The map (\ref{y10}) has a quantum deformation.} 
\end{conjecture} 
\vskip 3mm

We conjecture that the canonical map (\ref{y10}) behaves nicely at the boundary of 
the moduli space ${\cal A}_{^L G,S}$, and this requirement determines it 
 uniquely. Let us formulate this precisely. 

\vskip 3mm
Recall the embedding 
\begin{equation} \label{y18}
\beta_l: {\cal A}^{+, {\rm un}}_{G,S-l}\hra {\cal B}({\cal A}^+_{G,S})
\end{equation}
Let ${\cal O}_l({\cal A}_{G,S})$ be the 
subspace of functions which have limit at the component (\ref{y18}). 
So by the very definition there is a restriction map 
\begin{equation} \label{y17}
R_l: {\cal O}_l({\cal A}_{^L G,S})\lra {\cal O}({\cal A}^{{\rm un}}_{^L G, S-l})
\end{equation}

The space on the right should carry a canonical basis: 
the Duality Conjecture for the surface $S-l$ implies 
that there should be a canonical map
\begin{equation} \label{y30}
{\cal X}^{{\rm un}}_{G,S-l}(\Z^t) \lra {\cal O}({\cal A}^{{\rm un}}_{^L G, S-l}).
\end{equation}
Indeed, the canonical map ${\Bbb I}_{\cal X}: {\cal X}^{}_{G,S-l}(\Z^t) \lra {\cal O}({\cal A}^{}_{^L G, S-l})$ should have the following property: if $x \in {\cal X}^{}_{G,S-l}(\Z^t)$, 
and $y\in H^k(\Z^t)$ is its image under the map provided by the projection ${\cal X}^{}_{G,S-l}\to H^k$, then the function ${\Bbb I}_{\cal X}(x)$ 
transforms under the action of the dual torus $^L H^k$ by the character $\chi_y$ corresponding to $y$. 
Thus the functions corresponding to points of the subset  ${\cal X}^{{\rm un}}_{G,S-l}(\Z^t)$ 
should be invariant under the action of $^L H^k$. So
 the map ${\Bbb I}_{\cal X}$ should restrict to the map (\ref{y30}).

\vskip 3mm
Let $V$ be a vector space with a given basis ${\Bbb B}$, $V_0$ a subspace of $V$, and $r:V_0\to W$ 
a surjective linear map. We say that the basis ${\Bbb B}$ restricts under the map $r$ to a basis in $W$ if:

(i) The basis elements which lie inside of $V_0$ generate $V_0$, i.e. provide a basis ${\Bbb B}_0$ in $V_0$.

(ii) Let ${\Bbb B}'_0$ be the subset of the basis ${\Bbb B}_0$ 
consisting of the elements which are not killed by $r$. 
Then $r({\Bbb B}'_0)$ is a basis in $W$.

\begin{conjecture} \label{y88} {\it Let $S$ be a surface without boundary. Then, 
for any simple lamination $l$ on $S$, the restriction map $R_l$ is surjective, and 
the canonical basis (\ref{y10}) in ${\cal O}({\cal A}_{^L G,S})$ restricts 
under the map $R_l$ to the canonical basis in ${\cal O}({\cal A}^{{\rm un}}_{^L G, S-l})$.} 
\end{conjecture} 
\vskip 3mm
In other words,  for any simple lamination  $l$ on $S$, we should have a commutative diagram
\begin{equation} \label{y2}
\begin{array}{ccc}
{\cal X}^{\perp l}_{G,S}(\Z^t)& \lra & {\cal O}_l({\cal A}_{^L G,S})\\
&&\\
\downarrow &&\downarrow R_l\\
&&\\
{\cal X}^{{\rm un}}_{G,S-l}(\Z^t)&\lra & {\cal O}({\cal A}^{{\rm un}}_{^L G, S-l})
\end{array}
\end{equation}
Here the bottom arrow is the  map (\ref{y30}). 
The top arrow is given by the restriction of the map ${\Bbb I}$. 

The canonical maps for closed surfaces should be uniquely determined by Conjecture \ref{y88} and 
the canonical maps 
for surfaces with boundary. 
\vskip 3mm
{\bf Example}. Conjectures \ref{y8} and  \ref{y88} are valid for $G=PGL_2$. Indeed, in this case 
 the set ${\cal X}_{G,S}(\Z^t)$ 
coincides with the set of integral laminations on $S$, i.e. collections of simple disjoint curves with positive integral weights on $S$, considered modulo isotopy. 
The canonical map assigns to a lamination $\sum_i n_i\{\alpha_i\}$ the function 
$\prod_i {\rm Tr}M_{\widetilde \alpha_i}^{n_i}$ on the moduli space ${\cal A}_{SL_2,S}$ of twisted 
$SL_2$-local systems on $S$. Here $\widetilde \alpha_i$ is the lift of the loop 
$\alpha_i$ to the punctured tangent bundle of $S$. The quantum version of this map 
is nothing else but the Turaev algebra of $S$. The subset ${\cal X}^{\perp l}_{G,S}(\Z^t)$ 
consist of the laminations with the zero intersection index with $l$. The left vertical map 
in (\ref{y2}) kills the lamination $l$. So in this case the diagram (\ref{y2}) is 
 well-defined and commutative.

\section{Positive coordinate systems for ${\cal X}$-laminations}
\label{poslamX}

By Lemma \ref{12.25.05.1}, 
the group $(\Z/2\Z)^n$ acts by positive birational transformations 
on the moduli space ${\cal X}_{PGL_2, S}$.  
We extend the positive atlas on ${\cal X}_{PGL_2, S}$ by adding coordinate systems 
obtained from the standard ones by the action of $(\Z/2\Z)^n$. 
\vskip 3mm
{\bf 1. Positive coordinate systems: definitions and the main result}.  Let ${\cal C}_S$ be the set of pairs (an ideal  triangulation of
$S$ considered modulo isotopy, a choice of 
orientations of boundaries of all holes of $S$). 
The direct product of the
mapping class group of $S$ and the group $(\Z/2\Z)^n$ acts naturally on the 
set ${\cal C}_S$. The coordinate systems on the space of ${\cal X}$-laminations 
on $S$ are parametrized by the set ${\cal C}_S$. So there are 
$2^n$ coordinate systems corresponding to an ideal  triangulation of $S$. 

\begin{definition} \label{5.9.04.24a}
Let $l$ be an ${\cal X}$-lamination on $S$. 
A coordinate system corresponding to an ideal  triangulation of $S$ is 
{\em positive (non-negative) for the lamination $l$}, if 
all coordinates of $l$ in this coordinate system are $>0$ (respectively $\geq 0$).
\end{definition}

\begin{theorem} \label{5.31.05.1}
a) There exists a dense (resp. dense open)
 subset in the set of all real ${\cal X}$-laminations 
on $S$ such that every lamination from it admits a non-negative
 (resp. positive) coordinate system.

 b) Let $T$ be an ideal triangulation of $S$ corresponding to a
non-negative coordinate system for a finite real ${\cal X}$-lamination $l$. 
Then the isotopy class of the collection of 
edges of $T$ with $>0$ coordinates is isotopic to  $l$, and hence 
uniquely determined by $l$.

c) A positive coordinate system for an ${\cal X}$-lamination is unique if
exists. 
\end{theorem}

To get the part a) of Theorem \ref{5.31.05.1} we 
need the following result: for a generic real ${\cal X}$-lamination all 
curves of the lamination go between the punctures. 
For this we use some results from the theory of Morse-Smale foliations.

\vskip 3mm
{\bf 2. Morse-Smale foliations and real laminations}. 
Let us recall the notion of a Morse-Smale foliation on $S$ and some
results about them. 

Given a foliation on a punctured disc, the {\it winding
number} of the foliation at the puncture is an integer given by the 
 rotation number of the fibers of the foliation
restricted to a small circle near the puncture. For example the radial
foliation has the winding number $+1$, and the foliation tangent to the family
of concentric circles also has the winding number $+1$. 

Let $S$ be a surface with holes. 
A {\it Morse-Smale foliation} on $S$ is a foliation ${\cal F}$ on $S$ minus a
finite number of points, called singular points of the foliation, with 
the following structures/properties: 

i) ${\cal F}$ is equipped with a real transversal measure.

ii) The fibers of ${\cal F}$ are transversal to the boundary of $S$. 

iii) The winding number of the
foliation at each singular point is $<0$.

A fiber of a Morse-Smale foliation on $S$ is {\it singular} if 
it ends at a singular point of the foliation. The condition iii) implies that the number of 
singular fibers of a Morse-Smale foliation ending in a singular point is
finite. 
The condition ii) easily implies that a Morse-Smale foliation does not have fibers
contractible onto the boundary. 

\begin{lemma} \label{MS} The number of 
homotopy types of the leaves of a Morse-Smale foliation going 
between the boundary components is finite. 
\end{lemma} 

{\bf Proof}. The above remarks show that we may assume that the fibers can not
be contracted to the boundary. 
Cutting the surface along a non-singular leaf going between the boundary
components 
we decrease either the number of holes, or the genus. The lemma follows. 
\vskip 3mm

We need the following result.

\begin{theorem} \label{5.31.05.2}
There exists an open dense domain $U_{MS}(S)$ in the space of all Morse-Smale foliations
on $S$ 
such that every non-singular curve of a foliation from this domain 
joins two boundary components. 
\end{theorem}

Now let us replace the homotopic non-singular 
fibers of a Morse-Smale foliation ${\cal F}$  by a
single curve with the weight provided by the transversal measure. 
According to a well known theorem of Thurston, this way we get 
a bijection between the Morse-Smale foliations on $S$ modulo isotopy and
Thurston's $\R$-laminations: every
$\R$-lamination is representable as a Morse-Smale foliation on $S$, and 
this representation is unique up to an isotopy.  It follows
from Lemma \ref{MS} that if ${\cal F}$ was from the set $U_{MS}(S)$, we will get
a Thurston $\R$-lamination $l_{\cal F}$ with a finite number of curves.

\vskip 3mm
{\bf 3. Proof of Theorem \ref{5.31.05.1}}. 
Obviously c) is a particular case of b). Moreover, clearly
in b) not
only the set of the edges of $T$ carrying positive coordinates, 
but also the set of the orientations of the holes serving as the
endpoints of these edges is determined 
by $l$. 

a)  We will assume without loss of generality that an ${\cal X}$-lamination $l$
is defined using the orientations of the holes induced by the orientation of
the surface $S$. The general case is reduced to this using the action of the
group $(\Z/2\Z)^n$ changing  orientations of the holes.

A real ${\cal X}$-lamination is called {\it finite} if it can be presented by a finite 
collection $l$ of mutually and self non-intersecting curves with real weights on them. 
 Recall that to get an ${\cal X}$-lamination out of such a collection of curves with real weights 
$l$ we have to choose in addition 
orientations of all holes. As we already said, we do it by taking the orientations induced by
the orientation of the surface $S$, getting an ${\cal X}$-lamination
$l^+$. To define the ${\cal X}$-coordinates 
of the finite ${\cal X}$-lamination $l^+$ we have 
to 
perform  the
following procedure at every hole of $S$: we replace a segment of the 
lamination going to the hole by a ray winding 
infinitely many times around the puncture in the direction prescribed by the
orientation of $S$, and then 
retract the winded lamination to a trivalent graph on $S$. 

Let us call finite real 
${\cal X}$-laminations such that every leaf of the lamination joins two boundary component 
{\it special real ${\cal X}$-laminations}. 
Theorem \ref{5.31.05.2} and Lemma \ref{MS} imply that the space of 
special finite real ${\cal X}$-laminations 
is an open dense subset in the space of real ${\cal X}$-laminations on $S$.

\begin{figure}[ht]
\centerline{\epsfbox{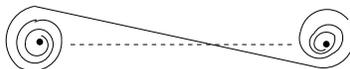}}
\caption{Calculating a coordinate of an ${\cal X}$-lamination assigned to an
  edge of  the dual graph.}
\label{fgo-150}
\end{figure}

Now, given  a  special ${\cal X}$-lamination $l^+$,
 let us exhibit a non-negative triangulation of $S$ for it. 
The curves of the finite lamination $l$ provide a decomposition $D_{l}$ 
of $S$. Adding a finite number of curves to this decomposition we obtain 
a triangulation $T$ of $S$ containing the decomposition $D_{l}$. 
Let $\Gamma$ be the dual graph for this triangulation. 
We claim that it provides a non-negative 
coordinate system for the ${\cal X}$-lamination $l^+$: 
  all coordinates assigned to the edges of
$D_{l}$ are 
$>0$, while the rest of the coordinates are equal to zero. 
Indeed, a curve of the lamination, being retracted onto the graph $\Gamma$,
 will cross only one edge of $\Gamma$, the one dual to the 
edge of the triangulation $T$ given by this curve. Figure \ref{fgo-150} 
illustrates this claim, by showing  
the winding procedure for a curve of the lamination on $S$ traveling between the holes.
The part a) of the theorem is proved.

 b) Given a special real ${\cal X}$-lamination $l$, there is
 the following correspondence
\begin{equation} \label{6.22.04.4}
{\cal I}_l \subset \{\mbox{edges of $\Gamma$}\} \times \{\mbox{connected components of 
the lamination $l$}\}.
\end{equation}
Namely, a pair (an edge $e$ of $\Gamma$, a connected component
 $\alpha$ of $l$) belongs to ${\cal I}_l$
if and only if there exists a part of  $\alpha$ 
which goes diagonally along the edge $e$. 
Let us denote by $\pi_2$ the natural projection 
of the correspondence ${\cal I}_l$ to the second factor in (\ref{6.22.04.4}):
$$
\pi_2: {\cal I}_l \lra \{\mbox{connected components of the lamination $l$}\}.
$$

\begin{lemma} \label{6.22.04.1} Let $l$ be a special real ${\cal X}$-lamination on $S$. 
Then: 
 
a) The map $\pi_2$ is surjective.

b) If $\pi_2^{-1}(\alpha)$ is a single edge $e$, then 
this edge is isotopic to an edge of the dual triangulation $T$. 

c) If all coordinates of $l$ 
are $\geq 0$, then the map $\pi_2$ is injective. 
\end{lemma} 

{\bf Proof}. a) Take a connected component $\beta$ of $l$ which is not in the 
image of $\pi_2$. This means that there is no edge of $\Gamma$ such that $\beta$ goes diagonally along this edge. Thus $\beta$ can be shrunken into the boundary of a hole. Indeed, 
$\beta$ either turns to the  right at every edge of $\Gamma$, or it turns to the  left at every edge 
of $\Gamma$, and thus it is isotopic to a curve which does not intersect $\Gamma$. 

b) If $\pi_2^{-1}(\alpha) = e$, then $\alpha$ is homotopic to a curve which intersects 
$\Gamma$ just once. Indeed, it either goes to the  left before $e$ and to the  right after $e$, 
or to the  right before and to the  left after $e$. The claim follows from this. 

c) Assume the opposite. Then there is a connected component $\alpha$ of $l$ such that 
$|\pi_2^{-1}(\alpha)|>1$. But this is impossible: after the curve $\alpha$ 
crossed diagonally the first edge on its way, it can not cross the other edge since 
all its coordinates are $\geq 0$. 
The lemma is proved.
\vskip 3mm

\begin{corollary} \label{6.22.04.2}
Let $l$ be a special real ${\cal X}$-lamination on $S$. 
Suppose that all coordinates of $l$ in the coordinate system corresponding to 
an ideal triangulation $T$ 
are $\geq 0$. Then the set of the edges of $T$ which carry $>0$ coordinates 
of $l$ is isotopic to the lamination $l$. 
\end{corollary} 

{\bf Proof}. Follows immediately from the parts b) and c) of Lemma \ref{6.22.04.1}. 

The part b) of the theorem is an immediate consequence of Corollary \ref{6.22.04.2}. 
The theorem is proved. 
\vskip 3mm

\begin{corollary} \label{5.31.05.10}
a) Any rational ${\cal X}$-lamination without loops admits a non-negative coordinate system.

b) The set of rational ${\cal X}$-laminations 
without loops is dense in the space ${\cal X}_L(S; \Q)$ of all
rational ${\cal X}$-laminations.  
\end{corollary}

{\bf Proof}. a) It has been proved during the proof of Theorem \ref{5.31.05.1}. 

b) The intersection of the domain described in Theorem
\ref{5.31.05.1} 
with the space of
all rational ${\cal X}$-laminations is described as the set of the ones 
without closed loops. 
The corollary is proved. 

So a
little perturbation destroys loops in a lamination. 

\vskip 3mm 
It is tempting to avoid the use of the Morse-Smale theory in the proof, making
 the proof self-contained. 
Here is perhaps a first step in this direction.  

\begin{lemma} \label{5.31.05.11s} Let $l$ be a real  ${\cal X}$-lamination and $\alpha$ be a closed
leaf of it. Take any 3-valent graph and let now  $x_i$ be the coordinates of
$l$ and $a_i$ be the coordinates of $\alpha$ considered as an ${\cal
A}$-lamination. Then $\sum_i x_i a_i=0$.
 \end{lemma}

{\bf Proof}. The intersection pairing of $\alpha$ with $l$ equals zero. 
On the other hand it is computed by the formula $\sum_i x_i a_i$. The lemma is proved. 
 
\vskip 3mm
This lemma implies that, since all $a_i$'s are integers, if all $x_i$'s are
linearly independent over ${\Q}$, there are no closed curves. 
However this does not prove yet that for a generic real ${\cal X}$-lamination all 
curves of the lamination go between the punctures: we have to 
take care of the non-closed curves winding inside of $S$. 
\vskip 3mm
{\bf 4. An application to the canonical map ${\Bbb I}_{\cal X}$}. 
It follows from the part a) of Corollary \ref{5.31.05.10} 
 that the restriction of the canonical map ${\Bbb I}_{\cal X}$ 
to the integral laminations without loops is especially simple. Namely, 
given such a
lamination $l$,  in every 
non-negative coordinate system for $l$ the element ${\Bbb I}_{\cal X}(l)$ is 
given by a product of the 
coordinate functions in this coordinate system. So ${\Bbb I}_{\cal X}(l)$ 
is in the cluster algebra corresponding to the 
space ${\cal A}_{SL_2, S}$. The image of laminations containing loops 
does not belong to it.

\section{A Weil--Petersson form on ${\cal A}_{G, \widehat S}$ and its motivic avatar}
\label{WPform}

{\bf 1. The $K_2$-invariant of a $G$--local system on a compact surface}. 
Let $S$ be a compact smooth oriented surface of genus $g>1$. 
A representation $\rho: \pi_1(S, x) \to G(F)$ of $\pi_1(S, x)$ 
to the group of 
$F$--points of a reductive group $G$ induces a map of the classifying spaces
$\rho_B: B\pi_1(S, x) \to BG(F)$. Since $S$ is a $K(\pi_1(S, x))$-space, applying the functor $H_2$  we get a map 
$$
\rho_*: H_2(S) = H_2(B\pi_1(S, x)) \lra  H_2BG(F) = H_2G(F).
$$

According to the stabilization theorem \cite{Sus}, the natural map $H_2GL_2(F) \to H_2GL_n(F)$ is an isomorphism for any infinite field $F$ provided $n \geq 2$. There is a natural map $H_2GL_2(F) \to K_2(F)$ (see e.g. loc. cit.). Thus if $n\geq 2$ there is a map $\pi_n: H_2GL_n(F) \to K_2(F)$. A finite dimensional representation $V$ of $G$ provides a map $BG \to BGL(V)$. Composing it with  $\pi_{\rm dim V}$ we get a homomorphism 
$\pi_{G, V}: H_2G(F) \to K_2(F)$. Let 
$ 
\pi_{G}: H_2G(F) \to K_2(F)
$  
be the homomorphism corresponding to the adjoint representation of $G$. 
Composing it with  $\rho_*$ we get a map
$ 
\pi_{G} \circ  \rho_*: H_2(S)  \lra K_2(F)
$. 
Let $[S] \in H_2(S)$ be the fundamental class of $S$. 
\begin{definition} Let $G$ be a split reductive group and $F$ is a field. 
The $K_2$--invariant 
$ 
w_{\rho} \in K_2(F)
$  of 
a representation  $\rho:\pi_1(S) \to G(F)$ is the element 
element $\pi_{G} \circ  \rho_*([S]) \in K_2(F)$. 
\end{definition} 
Since the mapping class group $\Gamma_{S}$ preserves the 
fundamental class of $S$,  $w_{\rho}$ is invariant under the 
action of $\Gamma_{S}$, i.e. $w_g{\rho} = w_{\rho}$ for $g \in \Gamma_{S}$.

Applying this to the field ${\Bbb F}$ of rational functions on the moduli space ${\cal L}_{G, S}$ and the corresponding tautological representation  
$\pi_1(S, x) \to G({\Bbb F})$ we get a 
${\Gamma_{S}}$-invariant class 
$W_{G, S} \in K_2({\Bbb F})$. 
The restriction of this class 
to any point ${\rho} \in {\cal L}_{G, S}$ is the class 
$w_{{\rho}}$. Therefore the class $W_{G, S}$ is regular 
on ${\cal L}_{G, S}$, i.e. provides an element 
$$
W_{G, S} \in H^2({\cal L}_{G, S}, \Q_{\cal M}(2))^{\Gamma_{S}} \subset K_2({\Bbb F})^{\Gamma_{S}}.
$$
Applying to it the map $d\log$ we get a regular $2$-form $\Omega_{G, S}$ on ${\cal L}_{G, S}$.

\vskip 3mm
{\bf 2. The weight two motivic 
complex}. Let $F$ be a field. Recall the Bloch complex:
of $F$ 
$$
\delta: B_2(F) \lra \Lambda^2F^*; \qquad B_2(F) := \frac{\Z[F^*]}{R_2(F)}.
$$
The subgroup $R_2(F)$ is generated by the elements 
$$
\sum_{i=1}^5 (-1)^i\{r(x_1, ..., \widehat x_i, ..., x_5)\}, \qquad x_i \in P^1(F), 
\quad x_i \not = x_j.
$$
where $\{x\}$ is the generator of $\Z[F^*]$ corresponding to $x \in F^*$. 
Let $\{x\}_2$ be  the 
projection of $\{x\}$  to $B_2(F)$. Then  
$\delta \{x\}_2 := (1-x) \wedge x$.  One has $\delta (R_2(F) =0$ (see \cite{G1}). 

Let $X_k$ be the set of all codimension $k$ irreducible subvarieties of $X$. 
The weight two motivic complex $\Gamma(X;2)$ of a regular irreducible 
  variety $X$ with the field of functions $\Q(X)$ is  the following 
cohomological complex
\begin{equation} \label{12.12.02.11}
\Gamma(X;2):= \qquad B_2(\Q(X)) \stackrel{\delta}{\lra} \Lambda^2 \Q(X)^* \stackrel{\rm Res}{\lra}  \prod_{Y \in X_1}\Q(Y)^*
\stackrel{{\rm div}}{\lra}  \prod_{Y \in X_2}\Z
\end{equation}
where the first group is in the degree $1$ and ${\rm Res}$ is the tame symbol map 
(\ref{12.20.02.1es}). If $Y \in X_1$ is normal, 
the last map 
is given by the divisor ${\rm div}(f)$ of $f$. 
In $Y$ is not normal,  we take its normalization $\widetilde Y$, compute ${\rm div}(f)$ 
on $\widetilde Y$, 
and then push it down to $Y$.

\vskip 3mm
{\bf 3. The second motivic Chern class}. 
According to Milnor the universal $G$--bundle over the classifying space 
for an algebraic group $G$ 
can be thought of as  the following (semi)simplicial variety $EG_{\bullet}$:
$$
\stackrel{\lra}{\lra}G^4 \stackrel{\lra}{\lra} G^3 \stackrel{\lra}{\lra} G^2 \stackrel{\lra}{\lra} G.
$$

The group $G$ acts diagonally from the left on $EG_{\bullet}$, 
and the quotient $BG_{\bullet}$ is a simplicial model for 
the classifying space of $G$. 
The simplicial structure includes the maps 
$$
s^n_i: G^n \lra G^{n-1}; \quad (g_1, ..., g_n) \lms (g_1, ..., \widehat g_i, ... , g_n). 
$$

Let ${\cal F}^{\bullet}(X)$ be a complex of sheaves on $X$ 
which depends functorially on $X$ for  smooth projections. 
The hypercohomology  $H^*(BG_{\bullet}, {\cal F}^{\bullet})$ are 
the cohomology of the total complex associated with the bicomplex:
$$
\begin{array}{cccccccc}
&...&... &... &...&... &...&... \\
&\uparrow &&\uparrow &&\uparrow &&\uparrow \\
\stackrel{d}{\longleftarrow} & {\cal F}^3(G^3) &\stackrel{d}{\longleftarrow} 
& {\cal F}^3(G^2) &\stackrel{d}{\longleftarrow} 
 & {\cal F}^3(G) &\stackrel{d}{\longleftarrow} 
& {\cal F}^3(*)\\
&\uparrow &&\uparrow &&\uparrow &&\uparrow \\
\stackrel{d}{\longleftarrow} & {\cal F}^2(G^3) &\stackrel{d}{\longleftarrow}
& {\cal F}^2(G^2) &\stackrel{d}{\longleftarrow} & {\cal F}^2(G) 
&\stackrel{d}{\longleftarrow} & {\cal F}^2(*)\\
&\uparrow &&\uparrow &&\uparrow &&\uparrow \\
\stackrel{d}{\longleftarrow} & {\cal F}^1(G^3) &\stackrel{d}{\longleftarrow} 
& {\cal F}^1(G^2) &\stackrel{d}{\longleftarrow} 
 & {\cal F}^1(G) &\stackrel{d}{\longleftarrow} &{\cal F}^1(*)\\
\end{array}
$$
The differential $d: {\cal F}^*(G \backslash G^{n-1})\lra {\cal F}^*(G \backslash G^n)$ is defined by
$$
d:= \sum_{i=1}^n(-1)^{i-1}(s^n_i)^*.$$ 
So $H^n(BG_{\bullet}, {\cal F}^{\bullet})$ 
is represented by cocycles in $\oplus_{p+q=n}{\cal F}^p(G^q)$. 

It is well known that 
$$
H^*(BG_{\bullet}, \Q_{\cal M}(2)) = 
H^*(BG_{\bullet}, \Gamma(X;2)\otimes \Q)
$$
where the right hand side  is computed via the bicomplex above. 
Moreover, it is known that if $*=4$, this $\Q$--vector space  can be identified 
with the one of $G$--invariant quadratic forms 
on the Lie algebra of $G$. In particular the Killing form provides 
a distinguished cohomology class
\begin{equation} \label{12.12.02.10q}
c_{2, G}^{\cal M} = c_{2}^{\cal M} \in H^4(BG_{\bullet}, \Q_{\cal M}(2)).
\end{equation}
  A cocycle representing the class $c_2^{\cal M}$ is given by the following data:
\begin{equation} \label{12.12.02.111}
C_4 \in B_2(\Q(G^4)^*)^G, \quad C_3 \in \Bigl(\Lambda^2 
\Q(G^3)^*\Bigr)^G, \quad C_2 \in  
\Bigl(\prod_{Y \in (G^2)_1}\Q(Y)^*\Bigr)^G.
\end{equation}

\begin{lemma} \label{12.14.02.1} There exists a cocycle 
$C^U_{\bullet} = (C^U_4, C^U_3, C^U_2)$ 
with $C^U_2 = 0$ 
representing 
the class 
$c_2^{\cal M}$ 
such that its component on $G^k$ is invariant under the 
 right action of $U^k$, and skew symmetric under the 
permutations of the factors of  $G^k$ modulo $2$-torsion. 
\end{lemma}

{\bf Proof}. 
Let $Q(X_*(H))^W$ be the group  
of $W$-invariant quadratic forms on the group of cocharacters of the 
Cartan group $H$. 
Recall (\cite{BrD}, Section 4; \cite{EKLV}, 3.2, 4.7, 4.8) that for a split
semi-simple group $G$ there exists 
canonical isomorphism, functorial in $G$,
\begin{equation} \label{9.20.03.1}
H^2(BG, \Z_{\cal M}(2)) \stackrel{\sim}{\lra} H^2(BH, \Z_{\cal M}(2))^W = Q(X_*(H))^W
\end{equation}
 
\vskip 3mm
{\bf Remark}.  The work \cite{BrD} employs 
$H^0(X, {\bf K}_2)$, where ${\bf K}_2$ is the sheaf of $K_2$
groups in the Zariski topology.  The Gersten resolution shows that it is
isomorphic to our group 
$H^2(X, \Z_{\cal M}(2))$. In  (\cite{BrD}, 4.1) the isomorphism (\ref{9.20.03.1}) is
formulated for simply-connected $G$, but it is valid for all $G$. 
\vskip 3mm

If $G = GL_{n}$, a cocycle  satisfying
all the conditions of the lemma   was constructed 
in \cite{G3}. We will recall its construction below. 
It provides a class 
in $H^2(BGL_n, \Z_{\cal M}(2))$ denoted by $w_{GL_n}$.

 We may assume that $G$ is simple. Let $\rho: G \lra GL_n$ be a nontrivial homomorphism.  
The Killing form 
on $GL_n$ restricts to $\lambda_{\rho}$ times the Killing form on $G$. 
Since  $\lambda_{\rho} \not = 0$,  the
class $\rho^*w_{GL_n}$ is not zero. 
The class $\rho^*w_{GL_n}$ is represented by an element in 
$C_3' \in K_2(\Q(G^3))$. Since $\rho (U)$ lies in a maximal 
unipotent subgroup $U'$ of $GL_n$, it is invariant under the right 
action of $U^3$. Therefore there exists an
$U^3$-invariant   
lift $C^U_3 \in \Lambda^2(\Q(G^3))$ of the class $\rho^*w_{GL_n}$. 
Consider the element 
\begin{equation} \label{11.27.03.1}
\sum_{i=1}^4(-1)^{i-1} (s^4_i)^*C^U_3.
\end{equation}
 It lies in  
$\Lambda^2\Q((G/U)^4)^*$. We claim that its projection to $K_2$ is zero. 
This means that  
there exists an 
$U^4$-invariant 
element $C^U_4 \in B_2(\Q(G^4))$ 
such that 
$\delta(C^U_4) = \sum_{i=1}^4(-1)^{i-1} (s^4_i)^*C^U_3$.  
It remains to skewsymmetrize it.  
To check the  claim observe that the projection of the element
(\ref{11.27.03.1}) 
to $K_2$ lies in the subspace $H^2((G/U)^4, \Z_{\cal M}(2))$ of the elements
of $K_2$ with zero tame symbols. On the other hand we know that its pull back
to $G^4$ is zero. Since $H^2(X, \Z_{\cal M}(2))$ is homotopy invariant, and
the variety $U$ is isomorphic to an affine space, we are done. 
The lemma is proved. 
\vskip 3mm
Observe that $G\backslash G^3/U^3$ is birationally isomorphic
to ${\rm Conf}_3({\cal A})$. So it comes from a class in $K_2$ of the field of rational
functions on ${\rm Conf}_3({\cal A})$. 
Set ${\Bbb F}_{G,k}:= \Q({\rm Conf}_k({\cal A}))$. Lemma \ref{12.14.02.1} 
just means that  there exist elements 
$$
C^U_4 \in B_2({\Bbb F}_{G,4}), \quad C^U_3 \in \Lambda^2{\Bbb
  F}^*_{G,3}
$$
which satisfies the following three conditions:
$$
{\rm i}) \sum_{i=1}^5(-1)^{i-1} (s^5_i)^*C^U_4  = 0; \quad {\rm ii}) 
\sum_{i=1}^4(-1)^{i-1} (s^4_i)^*C^U_3  = 
\delta C^U_4; \quad {\rm iii}) {\rm Res} (C^U_3) = 0.
$$

\vskip 3mm
{\bf 4. The motivic avatar of the Weil--Petersson form}. 
The mapping class group $\Gamma_S$ acts on  ${\cal A}_{G, \widehat S}$. 
Thanks to Proposition \ref{12.02.15.30} we can define the  
second rational 
 {\it $\Gamma_S$--equivariant weight two motivic cohomology} of ${\cal A}_{G,
   \widehat S}$ 
by 
$$
H_{\Gamma_S}^2({\cal A}_{G, \widehat S}, \Q_{\cal M}(2)):= H^2({\Bbb G}_{\bullet}(S) 
\times_{\Gamma_S} {\cal A}_{G, \widehat S}, \Q_{\cal M}(2)):=
$$
$$
:= H^2({\Bbb G}_{\bullet}(S) 
\times_{\Gamma_S} {\cal A}_{G, \widehat S}, \Gamma_{\cal M}(2)).
$$
\vskip 3mm
{\bf Remarks}. 1. ${\cal A}_{G, \widehat S}$ and ${\Bbb G}_{\bullet}(S)$  
are objects of different nature: 
the first is a stack whose generic point is an algebraic variety, while the second 
is a polyhedral complex. Their product is a mixed object. 
To define its hypercohomology with coefficients in a complex of 
sheaves ${\cal F}^{\bullet}$ one proceeds similarly to the definition 
of hypercohomology of a simplicial algebraic variety. 

2. Our cocycles for the classes in  $H_{\Gamma_S}^2$ will sit at the generic 
point of  ${\cal A}_{G, \widehat S}$. Therefore we can forget about the complications
related to the stack nature of  ${\cal A}_{G, \widehat S}$. 
\vskip 3mm
\begin{theorem} \label{12.14.02.4} 
\!A choice of $W$-invariant quadratic form \!$Q$ on \!$X_*(H)$ provides a class
$$
{\Bbb W}_G(S) \in H_{\Gamma_S}^2({\cal A}_{G, \widehat S}, \Q_{\cal M}(2)).
$$
It depends linearly on the form $Q$. 
\end{theorem}

Here is a refined  version of Theorem \ref{12.14.02.4}. 
Let ${\Bbb F}$ be the field of rational functions on ${\cal A}_{G, \widehat S}$. 
So the generic  ${\Bbb F}$--point of the variety 
${\cal A}_{G, \widehat S}$ describes  generic decorated unipotent $G$--local system 
on $\Gamma$. 

\begin{theorem} \label{3.27.02.3} 
A  cocycle $C^U_{\bullet}$ as in Lemma \ref{12.14.02.1} 
representing the class $c_2^{\cal M}$ 
provides a homomorphism of complexes 
$$
\begin{array}{cccccccc}
  G_{5}(S)& \lra & G_{4}(S) & \lra & G_{3}(S)&\lra & 0\\
\downarrow &&\downarrow w_4&& \downarrow w_3&&\downarrow \\
 0 &\lra & B_2({\Bbb F}) &\stackrel{\delta}{\lra} 
&\Lambda^2{\Bbb F}^*& \stackrel{\partial}{\lra} & \prod_{Y \in {\cal A}_m(S)_1}\Q(Y)^*
\end{array}
$$
Different cocycles lead to homotopic homomorphisms.
\end{theorem} 

{\bf Proof of Theorem \ref{3.27.02.3}}. 
A local system on a graph $\Gamma$ is  
the following: 

i) a collection  $\{L_v\}$ of vector spaces
attached to  the vertices $v$ of $\Gamma$. 

ii) For any edge $e$ an  orientation of $e$ gives rise to  an 
operator $M_{\stackrel{\to}{e}}: L_{v_1} \lra L_{v_2}$. Here $v_1$ and $v_2$ are the 
vertices
 of  $e$, and the chosen orientation  is from $v_1$ to $v_2$. 
Reversing the orientation 
of $e$ we get the inverse operator. 

Let $\Gamma$ be a graph embedded to $S$, so that $S$ retracts to $\Gamma$. 
Then 
there is a natural bijection between 
the local systems on  $S$ and 
  $\Gamma$.

We picture a valence $k$ vertex $v$ of a ribbon graph by a little 
oriented circle $S^1_v$
with $k$ marked points on it corresponding to the edges shared by $v$. 
The orientation of the circle is compatible with 
the cyclic order at $v$.  
The marked points cut the circle  into 
$k$ arcs $\alpha_1, ..., \alpha_k$. 
A decorated unipotent $G$--local system on 
a ribbon graph $\Gamma$ is given by a collection of 
affine flags $X_{v, \alpha_i}$ over  
${v}$. The affine  flag $X_{v, \alpha_i}$
is invariant under the monodromy around the face path provided by $\alpha_i$.

\begin{figure}[ht]
\centerline{\epsfbox{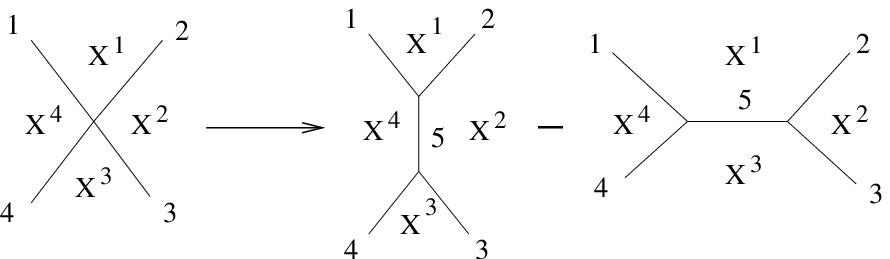}}
\caption{}
\label{fg50}
\end{figure}

Take a generator $(\Gamma, \varepsilon_{\Gamma})$ of $G_3(S)$. We assign to it 
an element of $\Lambda^2{\Bbb F}^*$ as follows. For each vertex $v$ of $\Gamma$ 
there are three cyclically ordered affine flags $X^1_v, X^2_v, X^3_v$ over $v$. 
Evaluating the component $C_3^U$ on them we get an element 
$C^U_3(X^1_v, X^2_v, X^3_v) \in \Lambda^2{\Bbb F}^*$. A trivalent ribbon graph has a 
canonical orientation provided by the 
cyclic structure at the edges. The sign ${\rm sgn}(\varepsilon_{\Gamma})$ as the ratio of 
this orientation to the one 
$\varepsilon_{\Gamma}$. We define the map $w_3$ by 
$$
w_3(\Gamma, \varepsilon_{\Gamma}) := \sum_{v \in V(\Gamma)} {\rm sgn}(\varepsilon_{\Gamma}) C^U_3(X^1_v, X^2_v, X^3_v)\in \Lambda^2{\Bbb F}^*.
$$
Let $(\Gamma, \varepsilon_{\Gamma}) \in G_4(S)$. Then $\Gamma$ is represented 
by a ribbon graph with exactly one $4$--valent vertex $v_0$. Let us enumerate 
the edges at this vertex in a way compatible with the cyclic order. 
Then the graph inherits a canonical orientation, and we will assume that 
$\varepsilon_{\Gamma}$ is this orientation.  One has 
$$
\partial: (\Gamma, \varepsilon_{\Gamma}) \lms (\Gamma_1, \varepsilon_{\Gamma_1}) - 
(\Gamma_2, \varepsilon_{\Gamma_2})
$$
where $\Gamma_1$ and $\Gamma_2$ are the two 
trivalent ribbon graphs shown on Fig. \ref{fg50}, 
and $\varepsilon_{\Gamma_i}$ is the canonical orientation of the trivalent 
ribbon graph $\Gamma_i$.

A decorated unipotent $G$--local system on $\Gamma$ provides us four affine flags 
$X^1, ..., X^4$ in the domains near the vertex $v_0$. 
We enumerate them using the chosen order 
of the edges at $v_0$ and set 
$$
w_4(\Gamma, \varepsilon_{\Gamma}) := C^U_4(X^1, X^2, X^3, X^4)\in B_2({\Bbb F}).
$$
The skew symmetry property of the component $C_4^U$ guarantees that this definition does not 
depend on the ordering of the edges at $v_0$. 
The condition
$
\delta \circ w_4(\Gamma, \varepsilon_{\Gamma}) = w_3 \circ \delta (\Gamma, \varepsilon_{\Gamma}) 
$ 
is equivalent to the cocycle condition ii) for $C^U_{\bullet}$.

\begin{figure}[ht]
\centerline{\epsfbox{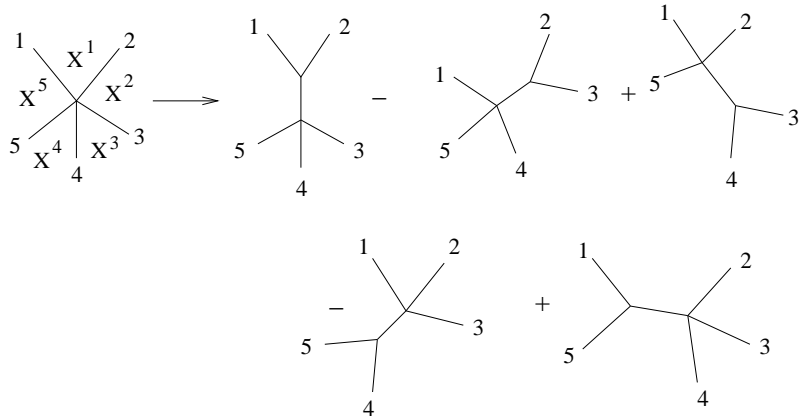}}
\caption{The pentagon relation.}
\label{fg51}
\end{figure}

Let $(\Gamma, \varepsilon_{\Gamma}) \in G_5(S)$. 
We assign to $(\Gamma, \varepsilon_{\Gamma})$ zero. One needs to check 
that
$$
w_4\Bigl(\partial (\Gamma, \varepsilon_{\Gamma})\Bigr) =0 \quad 
\mbox{in $B_2({\Bbb F})$}.
$$
Observe that $\Gamma$ has either one  $5$--valent vertex $v_0$, or 
two $4$--valent vertices, and all the other vertices are of valence three. 
Suppose that  $\Gamma$ has a $5$--valent vertex $v_0$. 
Choose an order of the edges at $v_0$ compatible with the cyclic
structure at the vertex. It provides 
an orientation of $\Gamma$, and we can assume that 
$\varepsilon_{\Gamma}$ is this orientation. 
Then $\partial(\Gamma, \varepsilon_{\Gamma})$ is given by the five generators shown on Fig. \ref{fg51}.  
Let $X^1, ..., X^5$ be the five affine flags over the vertex $v_0$, 
enumerated by using the numeration of the edges, see Fig. \ref{fg51}. Then

$$
w_4\Bigl(\partial (\Gamma, \varepsilon_{\Gamma})\Bigr) = \sum_{i=1}^5(-1)^{i-1}
C^U_4(X^1, ..., \widehat X^i, ..., X^5) = 0.
$$
The right hand side  is zero thanks to the condition i)  
for the cocycle $C^U_{\bullet}$.  

If $\Gamma$ has two $4$--valent vertices then it is easy to check, using  
only  the skew symmetry of $C_4^U$,  that 
$w_4\Bigl(\partial (\Gamma, \varepsilon_{\Gamma})\Bigr) = 0$. 
We conclude that the properties i)--iii) for the 
 cocycle $C_{\bullet}^U$ imply that $w_{\bullet}$ 
is a morphism of complexes. The theorem is proved. 
\vskip 3mm
{\bf Proof of Theorem \ref{12.14.02.4}}. Similar to  the proof 
of Theorem \ref{3.27.02.3}. 
\vskip 3mm
{\bf 5. The $K_2$--avatar of the degenerate symplectic form on 
${\cal A}_{G, \widehat S}$}. 
\begin{corollary} \label{12.13.02.2} Let $G$ be a split semi-simple algebraic 
 group. Then 
a choice of a $W$-invariant quadratic form $Q$ on $X_*(H)$ provides a class 
$$
W_{G, \widehat S} \in H^2({\cal A}_{G, \widehat S}, \Q_{\cal M}(2))^{\Gamma_S}.
$$
The form $\Omega_{G, \widehat S}:= d\log(W_{G, \widehat S})$ 
is a regular ${\Gamma_S}$-invariant $2$-form on the nonsingular part of  ${\cal A}_{G, \widehat S}$. 
\end{corollary}

{\bf Proof}. For a field $F$ one has 
$
K_2(F) = {\rm Coker} (\delta)
$ by Matsumoto's theorem.  
Thus for a regular variety $X$ we have 
$
H^2(X, \Q_{\cal M}(2)) = H^2\Gamma(X,2)\otimes \Q
$  
where the left hand side was defined in (\ref{12.20.02.1}). 

We define the class $W_{G, \widehat S}$ 
as the image of the class ${\Bbb W}_{G, \widehat S}$ 
under the canonical projection 
$$
H^2({\Bbb G}_{\bullet}(S) 
\times_{\Gamma_S} {\cal A}_{G, \widehat S}, \Q_{\cal M}(2)) \lra 
H^2({\cal A}_{G, \widehat S}, \Q_{\cal M}(2))^{\Gamma_S}.
$$
The properties of the map 
$d\log$ provide the second statement. The corollary is proved. 

{\it A direct construction of the class $W_{G,\widehat S}$}. 
Recall the canonical 
orientation of a trivalent ribbon graph $\Gamma$. 
Let $\Gamma \to \Gamma'$ be a flip defined by an edge $e$, and 
 $\overline \Gamma$ the ribbon graph obtained by shrinking 
of the edge $e$. 
Then $\overline \Gamma$  has just one $4$--valent vertex, the shrinked edge $e$. 
According to the proof of 
Theorem \ref{3.27.02.3} we have  
\begin{equation} \label{12.13.02.1}
w_3(\Gamma) - w_3(\Gamma') = \delta w_4(\overline \Gamma).
\end{equation} 
So the 
projection of 
$w_3(\Gamma) \in \Lambda^2{\Bbb F}^*$ to $K_2({\Bbb F})$ does not depend on 
the choice of  $\Gamma$. 
It is the class $W_{G, \widehat S}$. 

 In the next subsection we look at the proof of theorem \ref{3.27.02.3} from a bit 
different perspective, getting a more general result. 
\vskip 3mm
{\bf 6. The graph and the affine flags complexes}. 
Let $G$ be a group and $X$ a $G$--set. 
Let $C_n(X)$ be the coinvariants of the diagonal action of $G$ on 
$\Z[X^n]$. So the generators of $C_n(X)$ are elements $(x_1, ..., x_n)$ 
such that 
$(x_1, ..., x_n) = (gx_1, ..., gx_n)$ for any $g \in G$. 
The groups $C_{n}(X)$ form a complex $C_{\bullet}(X)$, called the {\it configurations complex}, 
  with a differential $\partial$ given by 
\begin{equation} \label{12.12.02.1}
\partial: (x_1, ..., x_n) \lms \sum_{i=1}^n(-1)^{i-1}(x_1, ..., \widehat x_i, ... x_n).
\end{equation}
Applying this construction to the set 
of ${\Bbb F}$-points of the 
affine flag variety ${\cal A}_G$ for $G$ we get a  configuration complex 
$C_{\bullet}({\cal A}_G({\Bbb F}))$. 
We define the related to $\widehat S$ { affine flags complex } 
$$
A^{\widehat S}_{\bullet}(G) := \quad ... \stackrel{\partial}{\lra} 
A^{\widehat S}_{n+2}(G) \stackrel{\partial}{\lra}
 A^{\widehat S}_{n+1}(G) \stackrel{\partial}{\lra}  A^{\widehat S}_{n}(G)  
\stackrel{\partial}{\lra} ... 
$$  
as its subcomplex 
generated by generic configurations of affine flags over ${\Bbb F}$. 

Let us define a homomorphism $f_k: G_k(S) \lra A^{\widehat S}_k(G)$. 
If $\Gamma$ has at least 
two vertices of valence $\geq 4$ we set $f_k(\Gamma, \varepsilon_{\Gamma})=0$. 
Suppose that $\Gamma$ has just one vertex $v_0$ of valence $\geq 4$. 
Then ${\rm val}(v_0) =k$. Choose an order of the edges sharing $v_0$ 
compatible with the cyclic order at $v_0$. This order provides an orientation 
$\widetilde \varepsilon_{\Gamma}$ of $\Gamma$. The order of the edges 
provides  an order of  the domains near $v_0$, and hence the affine flags 
$X^1, ..., X^k$ attached to the domains. 
Set 
$ 
f_k(\Gamma, \varepsilon_{\Gamma}):= \widetilde \varepsilon_{\Gamma}/\varepsilon_{\Gamma}
\cdot (X^1, ..., X^k)
$. 
It does not depend on the chosen order of the edges. 

\begin{proposition} \label{12.15.02.11}The map $f_{\bullet}$ 
provides a map of complexes
$$
\begin{array}{ccccccccc}
  \lra & G_{5}(S)& \lra & G_{4}(S) & \lra & G_{3}(S)& \lra & 0 \\
&\downarrow f_5 &&\downarrow f_4&& \downarrow f_3 && \downarrow \\
\lra & A^{\widehat S}_5(G) &\stackrel{}{\lra} & A^{\widehat S}_4(G) &\stackrel{}{\lra} & A^{\widehat S}_3(G)&\lra &A^{\widehat S}_2(G)
\end{array}
$$
\end{proposition}

\begin{figure}[ht]
\centerline{\epsfbox{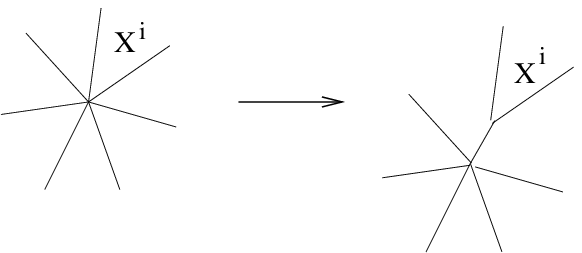}}
\caption{}
\label{fg52}
\end{figure}

{\bf Proof}. Let us check that for a generator 
$(\Gamma, \varepsilon_{\Gamma}) \in G_k(S)$ one has 
\begin{equation} \label{12.15.02.20}
\partial\circ  f_{k-1}(\Gamma, \varepsilon_{\Gamma}) = 
 f_k \circ \partial(\Gamma, \varepsilon_{\Gamma}).
\end{equation}
Suppose first that $k>3$. Then both terms in (\ref{12.15.02.20}) are zero
 if either there are  three vertices in $\Gamma$ 
of valence $\geq 4$, or there are two vertices of valence  $\geq 5$. 
Suppose there are just two vertices,   $v_1$ and $v_2$, of valence $\geq 4$. 
We may assume that one of them, say for $v_1$, is of valence $4$. 
Then $v_1$ provides two terms 
to $\partial(\Gamma)$. They enter  with the opposite signs, and hence 
$f_{k-1}$ kills their sum. 
It remains to consider the case 
when there is just one vertex $v$ of valence $k > 3$. 
Then the only component of  $\partial(\Gamma)$ 
which is not immediately killed by 
$f_{k-1}$ is the one shown on Fig. \ref{fg52}. 
It matches the $i$-th term in the formula for $\partial \circ f_k (\Gamma)$. 
So (\ref{12.15.02.20}) holds. 
If $k=3$ then each edge of $\Gamma$ contributes 
two terms to $f_k \circ \partial(\Gamma, \varepsilon_{\Gamma})$, and they 
appear with the opposite signs,  canceling each other. 
The proposition is proved.

\vskip 3mm
{\bf 7. Calculation of the ${\Bbb W}$-class for $G = SL_m$}.  
Recall that $\Gamma$ is a ribbon graph embedded to $\widehat S$, 
and ${\rm I}_m^{\Gamma}$ is the set
parametrising the corresponding $\Delta$- and $X$-coordinates, see Section 6. Recall the canonical map 
$p:  {\cal A}_{SL_m, \widehat S} \to {\cal X}_{ PSL_m, \widehat S}$.

\begin{theorem} \label{11.25.02.2} 
The $W$--class  on ${\cal A}_{PSL_m, \widehat S}$ is given by 
\begin{equation} \label{9.20.03.2}
\begin{array}{rcl}W_{SL_m, \widehat S} &=& \sum_{i, j \in {\rm I}_m^{\Gamma}}\varepsilon_{i j}  \{\Delta^{\Gamma}_{i}, d\log \Delta^{\Gamma}_{j}\} = \\ &=& 
\sum_{i \in {\rm I}_m^{\Gamma}} \{\Delta^{\Gamma}_{i} , 
p^*(X^{\Gamma}_{\alpha})\} \in K_2({\Bbb F}),\quad {\Bbb F}= \Q(
{\cal A}_{ SL_m, \widehat S}).
\end{array}
\end{equation}
It does not depend on the choice of the graph $\Gamma$.
\end{theorem}

\begin{corollary} \label{11.25.02.2c} 
The Weil-Petersson form on ${\cal A}_{ PSL_m, \widehat S}$ is given by 
$$
\Omega_{ SL_m, \widehat S} = \sum_{i, j \in {\rm I}_m^{\Gamma}}
\varepsilon_{i j}  
d\log \Delta^{\Gamma}_{i} \wedge d\log \Delta^{\Gamma}_{j}.
$$
It does not depend on the choice of the graph $\Gamma$.
\end{corollary}

{\bf Proof}.  We are going to  apply our construction to 
the cocycle $C^U_{\bullet}$  
for the second motivic Chern class $c^{\cal M}_2$  defined in  \cite{G3}. 
This cocycle satisfies the conditions of Lemma \ref{12.14.02.1}. 
We will get a cocycle
representing the ${\Bbb W}$-class whose $\Lambda^2{\Bbb F}^*$-component is given by the
following formula:
\begin{equation} \label{11.25.02.4}
\sum_{i, j \in {\rm I}_m^{\Gamma}} \varepsilon_{i, j} 
\Delta^{\Gamma}_{i} \wedge \Delta^{\Gamma}_{j} = 
\sum_{i \in {\rm I}_m^{\Gamma}} \{\Delta^{\Gamma}_{i} , 
p^*(X^{\Gamma}_{i})\} \in 
\Lambda^2 {\Bbb F}^*.
\end{equation} 
It, of course, depends on the choice of the graph $\Gamma$. 
Theorem \ref{11.25.02.2} follows immediately from this formula. 
Observe that the equality of the two terms in the formula is a straitforward
consequence of the formula (\ref{11.5.03.2}) for $p^*(X^{\Gamma}_{i})$. 
To prove (\ref{11.25.02.4}) we will give an explicit construction of the
$w_{\bullet}$-map.

Denote by $A_{\bullet}(m)$ the affine flag complex for $PSL_m$. 
We will construct the following
 diagram of morphisms of complexes, where the $A^{\widehat S}_{\bullet}(m)$ and 
$BC^{\widehat S}_{\bullet}(2)$ 
are complexes defined using  ${\Bbb F}$-vector spaces, 
as explained  below: 
$$
\begin{array}{ccccccccc}
  \lra & G_{5}(S)& \lra & G_{4}(S) & \lra & G_{3}(S)&&\\
&\downarrow f_5 &&\downarrow f_4&& \downarrow f_3 &&\\
\lra & A^{\widehat S}_5(m) &\stackrel{}{\lra} & A^{\widehat S}_4(m) &\stackrel{}{\lra} & A^{\widehat S}_3(m)&&\\
&\downarrow g_5 &&\downarrow g_4&& \downarrow g_3&&\\
\lra & BC^{\widehat S}_5(2) &\stackrel{}{\lra} & BC^{\widehat S}_4(2) &\stackrel{}{\lra} & BC^{\widehat S}_3(2)&&\\
&\downarrow &&\downarrow h_4&& \downarrow h_3&&\\
 &0 &\lra & B_2({\Bbb F}) &\stackrel{\delta}{\lra} &\Lambda^2 
{\Bbb F}^* & \stackrel{\rm Res}{\lra} & 
\prod_{Y \in {({\cal A}_{SL_m, \widehat S}})_1}\Q(Y)^*
\end{array}
$$
Then $w_{\bullet}$ is  the composition $w_{\bullet}:= 
h_{\bullet} \circ g_{\bullet}\circ f_{\bullet}$.

The complexes in the diagram are 
given by horizontal lines, while morphisms are given by vertical arrows.  
The top complex is the graph complex on $S$. 
The next one is the complex of affine flags in the generic position 
in an $m$-dimensional vector space. Then follows the rank $2$ bigrassmannian complex 
([G1]) and the complex $\Gamma({\cal A}_{SL_m, \widehat S}; 2)$. 
The homomorphisms $g_{\bullet}$ 
and $h_{\bullet}$ were defined in \cite{G3} and \cite{G1} - \cite{G2}. We recall their definition 
below.

{\it The grassmannian complex}. Take $X = V_m$, where $V_m$ is an $F$--vector space.  
The grassmannian complex 
$$
C_{\bullet}(m)_{/F}:= \quad ... \stackrel{\partial}{\lra} C_{m+2}(m) \stackrel{\partial}{\lra}
 C_{m+1}(m) \stackrel{\partial}{\lra}  C_{m}(m)  
$$ 
is the  subcomplex of the $m$--truncated complex $\sigma_{\geq m}C_{\bullet}(V_m)$ given by the configurations of vectors 
in the generic position in $V_m$.

{\it The rank two bigrassmannian complex $BC_{\bullet}(2)$} (\cite{G1}).
It is  the total complex associated with 
the bicomplex
$$
\begin{array}{ccccccc}
...&& ... && ... && ... \\
...&&\downarrow \partial'&&\downarrow \partial'&&\downarrow \partial'\\
...&\stackrel{\partial}{\lra}  &C_6(3)&\stackrel{\partial}{\lra} 
& C_5(3)&\stackrel{\partial}{\lra} & C_4(3)\\
...&&\downarrow \partial' &&\downarrow \partial'&&\downarrow \partial'\\
...&\stackrel{\partial}{\lra} 
& C_5(2) &\stackrel{\partial}{\lra} & C_4(2) & \stackrel{\partial}{\lra} & C_3(2)
\end{array}
$$
The horizontal differential are given by (\ref{12.12.02.1}). The vertical ones are  
$$
\partial': C_n(k) \lra C_{n-1}(k-1); \quad (x_1, ..., x_n) \lms \sum_{i=1}^n(-1)^{i-1}(x_i|x_1, ..., \widehat x_i, ... x_n).
$$
Here $(x_i|x_1, ..., \widehat x_i, ... x_n)$ denotes the configuration of 
vectors obtained by projecting $x_j$, $j \not = i,$ to the quotient 
$V_k/<x_i>$.

\vskip 3mm
 {\it The homomorphism $g_{\bullet}$} (\cite{G3}). Let us define its component
$$
g_n^k: A_n(m) \lra C_n(k); \quad k \leq m.
$$
Choose integers $p_1, ..., p_n\geq 0$ 
such that $p_1+ ... + p_n = m-k$. Let $X_{\bullet}^1, ..., X_{\bullet}^n$ be the 
flags underlying  affine flags 
$\widetilde X^1, ..., \widetilde X^n$. The quotient
\begin{equation} \label{12.02.02.2}
\frac{V_m}{X_{p_1} \oplus ... \oplus X_{p_n}}.
\end{equation} 
is a vector space of dimension $k$. Our affine flags induce  affine flags 
in this quotient. Taking the first vector of each of them, 
 we get a configuration of $n$ vectors in (\ref{12.02.02.2}). The corresponding element 
of the group $C_n(k)$ is, by definition, 
$g_n^k(\widetilde X^1, ..., \widetilde X^n)$. 
According to the Key Lemma from \cite{G3} this way we get a morphism of complexes.

\vskip 3mm
{\it The homomorphism $h_{\bullet}$}.  Consider the following 
diagram: 
$$
\begin{array}{cccccc}
&C_5(2)& \stackrel{\partial}{\lra} &C_4(2)& \stackrel{\partial}{\lra}& C_3(2)\\
&\downarrow l_5 &&\downarrow l_4 &&\downarrow l_3\\
&0& \lra &B_2({{\Bbb F}}) &\stackrel{\delta}{\lra} &\Lambda^2 
{{\Bbb F}}^* 
\end{array}
$$
$$
l_3: (l_1, l_2, l_3) \lms \Delta(l_1, l_2) \wedge \Delta(l_2, l_3) + 
\Delta(l_2, l_3) \wedge \Delta(l_3, l_1) + \Delta(l_3, l_1)\wedge \Delta(l_1, l_2)
$$
$$
l_4: (l_1, l_2, l_3, l_4) \lms -\{r(l_1, l_2, l_3, l_4)\}_2.
$$
Using the Pl\"ucker relation we see that,  modulo $2$--torsion, 
it provides a morphism of complexes,  (see s. 2.3 in \cite{G1}). 
The rank two grassmannian complex is  a subcomplex of the rank two 
grassmannian bicomplex. Extending it by zero 
to the rest of the rank two grassmannian bicomplex we get the map 
 $h_{\bullet}$. 
By s. 2.4  
of \cite{G1} it is  a morphism of complexes. 

Formula \ref{11.25.02.4} follows from the construction of the map
$w_3$. Theorem \ref{11.25.02.2} and hence 
Corollary \ref{11.25.02.2c} are proved. 

The class ${\Bbb W}_{PSL_m, \widehat S}$
is lifted from 
a similar class on the space ${\cal U}_{PSL_m, \widehat S}$: 
\begin{theorem} \label{9.7.03.1}
a) There exists a unique class 
$$
{\Bbb W}^{\cal U}_{PSL_m, \widehat S} \in 
H^2_{\Gamma_S}({\cal U}_{PSL_m, \widehat S}, \Q_{\cal M}(2)) 
\quad \mbox{such that ${\Bbb W}_{SL_m, \widehat S} = 
p^*{\Bbb W}^U_{PSL_m, \widehat S}$}.
$$

b) The form ${\Omega}^U_{SL_m, \widehat S}:= d\log ({W}^U_{PSL_m, \widehat S})$ is a 
symplectic form on  ${\cal U}_{PSL_m, \widehat S}$. 
\end{theorem}

{\bf Proof}. a) The class ${\Bbb W}_{SL_m, \widehat S}$ has two nontrivial
components, 
with values in $B_2$ and $\Lambda^2$. The $B_2$ component comes from 
the ${\cal U}_{PSL_m, \widehat S}$-space by the very definition. 
For the $\Lambda^2$-component this is deduced from the formula
(\ref{11.25.02.4}). 

b) The class $W^U_{SL_m, \widehat S}$ coincides with the one defined in 
the set-up of cluster ensembles, using the results of Section 9. 
Then the claim  is a general property of the $W$-class in $K_2$ of a 
cluster ensemble, see [FG2], Chapter 5. 

\newpage





\begin{thebibliography}{BGSV}

\bibitem[BAG]{BAG} I. \smalltextsc{Biswas}, P. \smalltextsc{Ares-Gastesi}, S. \smalltextsc{Govindarajan}, Parabolic Higgs bundles and Teichmuller spaces for punctured surfaces., \textit{Trans. Amer. Math. Soc.}, \textbf{349,} (1997), no. 4, 1551--1560, \texttt{alg-geom/9510011}.
\bibitem[BD]{BD} A.A. \smalltextsc{Beilinson}, V.G. \smalltextsc{Drinfeld}, {\it Opers}, \texttt{math.AG/0501398}.
\bibitem[BeK]{BeK} A.\smalltextsc{Berenstein}, D. \smalltextsc{Kazhdan}, Geometric and unipotent crystals., GAFA, Special Volume GAFA2000, 188-236. 
\bibitem[BeZ]{BeZ} A. \smalltextsc{Berenstein}, A. \smalltextsc{Zelevinsky}, Tensor product multiplicities, canonical bases and totally positive algebras., \textit{Invent. Math.}, \textbf{143,} (2001), no. 1, 77--128, \texttt{math.RT/9912012}.
\bibitem[BFZ96]{BFZ96} A. \smalltextsc{Berenstein}, S. \smalltextsc{Fomin}, A. \smalltextsc{Zelevinsky}, Parametrizations of canonical bases and totally positive matrices., \textit{Adv. Math.}, \textbf{122,} (1996), no. 1, 49--149, 
\bibitem[BFZ03]{BFZ03} A. \smalltextsc{Berenstein}, S. \smalltextsc{Fomin}, A. \smalltextsc{Zelevinsky}, Cluster algebras. III: Upper bounds and double Bruhat cells., \textit{Duke Math. J.}, \textbf{126,}(2005),  no. 1, 1--52, \texttt{math.RT/0305434}.
\bibitem[Bers]{Bers} L. \smalltextsc{Bers}, 
Universal Teichm"uller space.  Analytic methods in mathematical physics (Sympos., Indiana Univ., Bloomington, Ind., 1968),  pp. 65--83. Gordon and Breach, 1970. 
\bibitem[Bers2]{Bers2} L. \smalltextsc{Bers}, On the boundaries of Teichm\"uller spaces and on Kleinian groups., \textit{Ann. of Math.}, \textbf{91,} (1970), 670--600. 
 \bibitem[Bon]{Bon} F. \smalltextsc{Bonahon}, The Geometry of Teichmuller space via geodesic currents., \textit{Inv. Math.},  \textbf{92,}  (1988),  no. 1, 139--162. 
\bibitem[Bo]{Bo} N. \smalltextsc{Bourbaki}, \textit{Lie groups and Lie algebras. Chapters 4-6.}, 
Springer, (2002). 
\bibitem[Bor]{Bor} M. \textit{Borovoi}, Abelianization of the second nonabelian Galois cohomology, 
\textit{Duke Math. J.}, \textbf{72,} (1993), 217--239. 
\bibitem[BrD]{BrD} J.-J \smalltextsc{Brylinsky}, \smalltextsc{Deligne}, Central extensions of
  reductive groups by $K_2$, \textit{Publ. Math. IHES}, (2001).
\bibitem[CG]{CG} N. \smalltextsc{Chriss}, V. \smalltextsc{Ginzburg}, \textit{Representation theory and complex geometry.},
Birkhauser Boston, Inc., Boston, MA, (1997). 
\bibitem[FCh]{FCh}  L.O. \smalltextsc{Chekhov}, V.V. \smalltextsc{Fock}, Quantum Teichm{\"u}ller spaces.,
 \textit{Teoret. Mat. Fiz.}, \textbf{120,} (1999), no. 3, 511--528, \texttt{math.QA/9908165}.
\bibitem[C]{C} K. \smalltextsc{Corlette}, Flat $G$-bundles with canonical metrics, 
\textit{J. Diff. Geom.}, \textbf{28,} (1988), 361--382. 
\bibitem[Del]{Del} P. \smalltextsc{Deligne}, {\'E}quations diff{\'e}rentielles 
{\`a} points singuliers r{\'e}guliers. Springer LNM, Vol. 163. 1970. 
\bibitem[DS]{DS} V.G. \smalltextsc{Drinfeld}, V.V. \smalltextsc{Sokolov}, Lie algebras and 
equations of Korteweg-de Vries type., \textit{Current problems in mathematics,} 
\textbf{24,} (1984),81--180, In Russian.
\bibitem[D]{D} S. \smalltextsc{Donaldson}, Twisted Harmonic Maps and the Self-Duality Equations., \textit{Proc. London Math. Soc.},  \textbf{55,} (1987), 127--131.
\bibitem[EKLV]{EKLV} H. \smalltextsc{Esnault}, B. \smalltextsc{Kahn}, M. \smalltextsc{Levine}, E. \smalltextsc{Viehweg}, The Arason invariant and mod $2$ algebraic cycles, \textit{JAMS}, \textbf{11,} (1998) no 1, 73-118.
 Lett. Math. Phys. 34 (1995), no. 3, 249--254.
\bibitem[F]{F} V.V. \smalltextsc{Fock}, Dual Teichm\" uller spaces., \texttt{dg-ga/9702018}.
\bibitem[FR]{FR} V.V. \smalltextsc{Fock}, A.A. \smalltextsc{Rosly}, Poisson structure on moduli of flat connections on Riemann 
surfaces and $r$-matrix. \textit{AMS Transl. Ser. 2}, \textbf{191,} (1999), 67--86, \texttt{math.QA/9802054}.
\bibitem[FG2]{FG2} V.V. \smalltextsc{Fock}, A.B. \smalltextsc{Goncharov}, Cluster ensembles, 
quantization and the dilogarithm, \texttt{math.AG/0311245}. 
\bibitem[FG3]{FG3} V.V. \smalltextsc{Fock}, A.B. \smalltextsc{Goncharov}, Moduli spaces of convex
    projective structures on surfaces, 
To appear in \textit{Adv. in Math.}, (2006), \texttt{math.AG/0405348}. 
\bibitem[FG4]{FG4} V.V. \smalltextsc{Fock}, A.B. \smalltextsc{Goncharov}, Dual Teichm\"uller and lamination spaces, To appear in the \textit{Handbook on Teichm\"uller theory}, \texttt{math.AG/0510312}. 
\bibitem[FG5]{FG5} V.V. \smalltextsc{Fock}, A.B. \smalltextsc{Goncharov}, Cluster ${\cal X}$-varieties, amalganations, and Poisson-Lie groups. \textit{Progress in Mathematics}, Birkhauser. The volume dedicated to V.G. Drinfeld, \texttt{math.RT/0508408}.
\bibitem[FG6]{FG6} V.V. \smalltextsc{Fock}, A.B. \smalltextsc{Goncharov}, To appear. 
\bibitem[FZ99]{FZ99} S. \smalltextsc{Fomin}, A. \smalltextsc{Zelevinsky}, Double Bruhat cells and total positivity., \textit{JAMS}, \textbf{12,} (1999), no. 2, 335--380, \texttt{math.RA/9912128}.
\bibitem[FZI]{FZI} S. \smalltextsc{Fomin}, A. \smalltextsc{Zelevinsky}, Cluster algebras. I., 
\textit{JAMS}, \textbf{15,} (2002), no. 2, 497--529, \texttt{math.RT/0104151}.
\bibitem[FZII]{FZII} S. \smalltextsc{Fomin}, A. \smalltextsc{Zelevinsky}, Cluster algebras. II: Finite type classification., \textit{Invent. Math.}, \textbf{154,} (2003), no. 1, 63--121, \texttt{math.RA/0208229}. 
\bibitem[FZ3]{FZ3} S. \smalltextsc{Fomin}, A. \smalltextsc{Zelevinsky}, The Laurent  phenomenon. \textit{Adv. in Appl. Math.}, \textbf{28,} (2002), no. 2, 119--144, \texttt{math.CO/0104241}.
\bibitem[GGL]{GGL} A.M. \smalltextsc{Gabrielov}, I.M. \smalltextsc{Gelfand}, M.V.\smalltextsc{Losik}, Combinatorial computation of characteristic classes. I, II. (Russian), \textit{Funk. Anal. i Pril.}, \textbf{9,} (1975), no. 2, 12--28; no. 3, 5--26.
\bibitem[GKr]{GKr} F.R. \smalltextsc{Gantmacher}, M.G. \smalltextsc{Krein},  \textit{Oscillation matrices and kernels and small vibrations of mechanical systems}. Revised edition of the 1941 Russian original. 
\bibitem[GKr2]{GKr2} F.R. \smalltextsc{Gantmacher}, M.G. \smalltextsc{Krein}, Sur les matrices oscillatores., \textit{CR Acad. Sci. Paris}, \textbf{201,} (1935). AMS Chelsea Publ., Providence, RI, (2002).
\bibitem[GSV1]{GSV1}  M. \smalltextsc{Gekhtman}, M. \smalltextsc{Shapiro}, A. \smalltextsc{Vainshtein}, Cluster algebras and Poisson geometry.,   \textit{Mosc. Math. J.}, \textbf{3,}  (2003),  no. 3, 899--934; \texttt{math.QA/0208033}.  
\bibitem[GSV2]{GSV2} M. \smalltextsc{Gekhtman}, M. \smalltextsc{Shapiro}, A. \smalltextsc{Vainshtein}, Cluster algebras and Weil-Petersson forms., \textit{Duke Math. J.}  \textbf{127,} (2005),  no. 2, 291--311, \texttt{math.QA/0309138}. 
\bibitem[Gui]{Gui} O.\smalltextsc{Guichard}, Sur les r\'epresentations de groupes de surface., Preprint. 
\bibitem[Gol1]{Gol1} W. \smalltextsc{Goldman}, The symplectic nature of fundamental groups of surfaces., \textit{Adv. in Math.}, \textbf{54,} (1984), no. 2, 200--225. 
\bibitem[Gol2]{Gol2} W. \smalltextsc{Goldman}, Convex real projective structures on compact surfaces.,
\textit{J. Diff. Geom.}, \textbf{31,} (1990) 126--159. 
\bibitem[G1]{G1} A.B. \smalltextsc{Goncharov}, Geometry of configurations, polylogarithms, and motivic cohomology., \textit{Adv. Math.}, \textbf{114,} (1995), no. 2, 197--318.
\bibitem[G2]{G2} A.B. \smalltextsc{Goncharov}, Polylogarithms and motivic Galois groups.,  
Motives (Seattle, WA, 1991),\textit{Proc. Sympos. Pure. Math.}, \textbf{55,} Part 2, Amer. Math. Soc., Providence, RI, (1994),43--96,
\bibitem[G3]{G3} A.B. \smalltextsc{Goncharov}, Explicit construction of characteristic classes.
\textit{I. M. Gelfand Seminar,  Adv. Soviet 
Math.}, \textbf{16,} Part 1, Amer. Math. Soc., Providence, RI, (1993),169--210.
\bibitem[G4]{G4}  A.B. \smalltextsc{Goncharov}, Deninger's conjecture of $L$-functions of elliptic curves at $s=3$. Algebraic geometry, 4. \textit{J. Math. Sci.}, \textbf{81,} (1996), no. 3, 2631--2656; \texttt{alg-geom/9512016}.
\bibitem[G5]{G5} A.B. \smalltextsc{Goncharov}, Polylogarithms, regulators and Arakelov motivic complexes., \textit{J. Amer. Math. Soc.},  \textbf{18,}  (2005),  no. 1, 1--6; \texttt{ math.AG/0207036}
\bibitem[GM]{GM}  A.B. \smalltextsc{Goncharov}, Yu.I. \smalltextsc{Manin}, Multiple $\zeta$-motives and moduli spaces ${\cal M}_{0,n}$., \textit{Compos. Math.}, \textbf{140,} (2004), no. 1, 1--14, \texttt{math.AG/0204102}.
\bibitem[Ha]{Ha} J. \smalltextsc{Harer}, The virtual cohomological dimension of the mapping class group of an orientable surface., \textit{Invent. Math.}, \textbf{84,} (1986) no. 1, 157--176. 
\bibitem[H1]{H1} N.J. \smalltextsc{Hitchin}, Lie groups and Teichmuller space., 
\textit{Topology} \textbf{31,} (1992), no. 3, 449--473. 
\bibitem[H2]{H2} N.J. \smalltextsc{Hitchin}, The self-duality equation 
on a Riemann surface., \textit{Proc. London Math. Soc.} \textbf{55,} (1987), 59--126.
\bibitem[K]{K} R.M.\smalltextsc{Kashaev}, Quantization of Teichm{\"u}ller spaces and 
the quantum dilogarithm., \textit{Lett. Math. Phys.}, \textbf{43,} (1998), no. 2,
105--115.
\bibitem[Kr]{Kr} I. \smalltextsc{Kra}, Deformation spaces.  A crash course on Kleinian groups (Lectures at a Special Session, Annual Winter Meeting, Amer. Math. Soc., San Francisco, Calif., 1974), \textit{Lecture Notes in Math.}, \textbf{400,} Springer, Berlin, (1974),48--70.. 
\bibitem[Ko]{Ko} M. \smalltextsc{Kontsevich}, Formal (non)commutative symplectic geometry.,  
\textit{The Gelfand Mathematical Seminars 1990--1992}, Birkh\"auser Boston, Boston, MA, (1993).173--187.
\bibitem[Lab]{Lab} F. \smalltextsc{Labourie},  {\em Anosov flows, surface groups and curves in projective spaces}. Preprint, Dec. 8, (2003). 
\bibitem[L1]{L1} G. \smalltextsc{Lusztig}, Total positivity in reductive groups,
\textit{Lie theory and geometry}, Progr. Math., \textbf{123,} Birkh{\"a}user Boston, 
Boston, MA, (1994),531--568.
\bibitem[L2]{L2} G. \smalltextsc{Lusztig}, Total positivity and canonical bases., 
\textit{Algebraic groups and Lie groups,} \textit{Austral. Math. Soc. Lect. Ser.}, 
\textbf{9,} Cambridge Univ. Press, Cambridge, (1997), 281--295.
\bibitem[MM]{MM} C. \smalltextsc{McMullen}, Iteration on Teichmüller space.,  \textit{Invent. Math.},  \textbf{99,}  (1990),  no. 2, 425--454. 
\bibitem[Mi]{Mi} J. \smalltextsc{Milnor}, \textit{Introduction to algebraic $K$-theory.}, Annals
  of Mathematics Studies, No. 72. Princeton University Press, Princeton, N.J.,
  University of Tokyo Press, Tokyo, (1971). 
\bibitem[NZ]{NZ} I. \smalltextsc{Nikolaev},E.  \smalltextsc{Zhuzhoma}, Flows on 2-dimensional manifolds. LNM1705 1999.
\bibitem[P1]{P1} R.C. \smalltextsc{Penner}, The decorated Teichm{\"u}ller space of punctured surfaces., 
\textit{Comm. Math. Phys.}, \textbf{113,} (1987), no. 2, 299--339. 
\bibitem[P2]{P2} R.C. \smalltextsc{Penner}, Weil-Petersson volumes., \textit{J. Differential Geom.}, \textbf{35,} (1992), no. 3, 559--608.  
\bibitem[P3]{P3} R.C. \smalltextsc{Penner}, Universal constructions in Teichm{\"u}ller theory.,
 \textit{Adv. Math.} \textbf{98,} (1993), no. 2, 143--215. 
\bibitem[P4]{P4} R.C. \smalltextsc{Penner}, The universal Ptolemy group and its completions. Geometric Galois actions, \textbf{2,} 293--312, London Math. Soc. Lect. Note Ser., 243,
   Cambridge Univ. Press,  1997.
\bibitem[PH]{PH} R.C. \smalltextsc{Penner}, J.L. \smalltextsc{Harer}, Combinatorics of train tracks. \textit{Annals of Mathematics Studies,}  \textbf{125,} Princeton University Press, Princeton, NJ, 1992. 
\bibitem[Sch]{Sch} I.J. \smalltextsc{Schoenberg}, Convex domains and linear 
combinations of continuous functions. \textit{BAMS}, \textbf{39,}(1933), 273-280. 
\bibitem[Sch1]{Sch1} I.J. \smalltextsc{Schoenberg}, \"Uber variationsvermindernde lineare Transformationen. 
\textit{Math. Z.} \textbf{32,} (1930) 321-322. 
\bibitem[S]{S} C. \smalltextsc{Simpson}, Constructing Variations of Hodge Structures Using Yang-Mills Theory and Applications to Uniformization, \textit{JAMS} \textbf{1,} (1988) 867--918. 
\bibitem[Se]{Se}J.-P. \smalltextsc{Serre}, Cohomologie Galoisienne, Lect. Notes Math., \textbf{5}. 
\bibitem[Str]{Str} K. \smalltextsc{Str\``bel}, \textit{Quadratic Differentials}, Springer,
Berlin--Heidelberg--New York (1984).
\bibitem[SZ]{SZ} P. \smalltextsc{Sherman}, A. \smalltextsc{Zelevinsky}, 
Positivity and canonical bases in rank 2 cluster algebras of finite and affine types,  \textit{Mosc. Math. J.}  \textbf{4,}  (2004),  no. 4, 947--974, \texttt{math.RT/0307082}.
\bibitem[Sus]{Sus} A.A. \smalltextsc{Suslin}, Homology of ${\rm GL}\sb{n}$, characteristic classes and Milnor $K$-theory. \textit{Algebraic geometry and its applications. Trudy Mat. Inst. Steklov.} \textbf{165,} (1984),  188--204. 
\bibitem[Th]{Th} W. \smalltextsc{Thurston}, The geometry and topology of three-manifolds,
\textit{Princeton University notes,} 
\texttt{http://www.msri.org/publications/books/gt3m}. 
\bibitem[W]{W} A.M. \smalltextsc{Whitney}, A reduction theorem for totally positive matrices,
\textit{J. d'Analyse Math.} \textbf{2,} (1952), 88-92.
\bibitem[Wo]{Wo} S.  \smalltextsc{Wolpert}, Geometry of the 
Weil-Petersson completion of the Teichm\"uller space. 
\textit{Surveys in Differential Geometry.J.} \textbf{VIII,} (2002), 357-393.

\end{thebibliography}
\end{document}